\documentclass[reqno,10pt]{amsart}
\usepackage{amsmath,mathrsfs}
\usepackage{amssymb, amsmath}
\usepackage{graphicx}
\usepackage{pdfsync}
\usepackage{color}
\usepackage[colorlinks,
            linkcolor=red,
            anchorcolor=blue,
            citecolor=green
            ]{hyperref}
\usepackage{colonequals}
\def\cs#1#2{\big(#1\big)^{\frac{1}{2}}\big(#2\big)^{\frac{1}{2}}}

\def\inte#1{
\displaystyle\mathop{#1\kern0pt}^\circ }

\let\pa=\partial

\let\f=\frac

\def\pa{\partial}

\def\virgp{\raise 2pt\hbox{,}}
\def\cdotpv{\raise 2pt\hbox{;}}

\def\C{\mathop{\mathbb C\kern 0pt}\nolimits}
\def\DD{\mathop{\mathbb D\kern 0pt}\nolimits}
\def\EE{\mathop{{\mathbb E \kern 0pt}}\nolimits}
\def\K{\mathop{\mathbb K\kern 0pt}\nolimits}
\def\N{\mathop{\mathbb N\kern 0pt}\nolimits}
\def\Q{\mathop{\mathbb Q\kern 0pt}\nolimits}
\def\R{\mathop{\mathbb R\kern 0pt}\nolimits}
\def\SS{\mathop{\mathbb S\kern 0pt}\nolimits}
\def\ZZ{\mathop{\mathbb Z\kern 0pt}\nolimits}
\def\TT{\mathop{\mathbb T\kern 0pt}\nolimits}
\def\P{\mathop{\mathbb P\kern 0pt}\nolimits}

\DeclareMathOperator*{\esssup}{ess\,sup}

\newcommand{\beq}{\begin{equation}}
\newcommand{\eeq}{\end{equation}}
\newcommand{\ben}{\begin{eqnarray}}
\newcommand{\een}{\end{eqnarray}}
\newcommand{\beno}{\begin{eqnarray*}}
\newcommand{\eeno}{\end{eqnarray*}}

\newtheorem{thm}{Theorem}[section]
\newtheorem{lem}{Lemma}[section]
\newtheorem{col}{Corollary}[section]
\newtheorem{prop}{Proposition}[section]

\theoremstyle{definition}

\newtheorem{rmk}{Remark}[section]

\numberwithin{equation}{section}

\begin{document}
\title[Quantum Boltzmann equation]
{Global well-posedness of the quantum Boltzmann equation for bosons interacting via inverse power law potentials}

\author[Y.-L. Zhou]{Yu-long Zhou}
\address[Y.-L. Zhou]{School of Mathematics, Sun Yat-Sen University, Guangzhou, 510275, P. R.  China.} \email{zhouyulong@mail.sysu.edu.cn}

\begin{abstract}  We consider the spatially inhomogeneous quantum Boltzmann equation for bosons with
a singular collision kernel, the weak-coupling limit of a large system of Bose-Einstein particles interacting through inverse power law. Global well-posedness of the corresponding Cauchy problem is proved in a periodic box near equilibrium for initial data satisfying high temperature condition.
\end{abstract}

%
%
%
%

\maketitle
\markright{Quantum Boltzmann equation}

\setcounter{tocdepth}{1}
\tableofcontents






 %

\section{Introduction}

%
%
%
%

Quantum Boltzmann equations are proposed to describe the time evolution of a dilute system of weakly interacting bosons or fermions. The derivation of such equations dates back to as early as 1920s by Nordheim \cite{nordhiem1928kinetic} and 1933 by Uehling-Uhlenbeck \cite{uehling1933transport}.
As a result, the quantum Boltzmann equations are also called Boltzmann-Nordheim equations or Uehling-Uhlenbeck equations.
Later on, further developments were made by Erd{\H{o}}s-Salmhofer-Yau \cite{erdHos2004quantum}, Benedetto-Castella-Esposito-Pulvirenti \cite{benedetto2004some}, \cite{benedetto2005weak},  \cite{benedetto2006some}, \cite{benedetto2008n} and \cite{benedetto2007short},  Lukkarinen-Spohn \cite{lukkarinen2009not}.
One can refer the classical book \cite{chapman1990mathematical} for physical backgrounds.

In this article, we consider the Cauchy problem of the quantum Boltzmann equation for bosons
\ben \label{quantum-Boltzmann-UU}
\partial _t F +  v \cdot \nabla_{x} F=Q_{\Phi, \hbar}(F,F), ~~t > 0, x \in \mathbb{T}^{3}, v \in \R^3 ; \quad
F|_{t=0}(x,v) = F_{0}(x,v).
\een
Here $F(t,x,v)\geq 0$ is the density
function of particles with velocity
$v\in\R^3$ at time $t\geq 0$ in position $x \in \mathbb{T}^{3}$. Here $\mathbb{T}^{3} \colonequals  [0,1]^{3}$ is the periodic box with unit volume $|\mathbb{T}^{3}|=1$.
The  quantum Boltzmann operator $Q_{\Phi, \hbar}$ acting only on velocity variable
$v$ is defined by
\ben \label{U-U-operator}
Q_{\Phi, \hbar}(g,h)(v) \colonequals   \int_{{\mathbb S}^2 \times {\mathbb R}^3} B_{\Phi,\hbar}(v- v_{*},\sigma)
 \mathrm{D}\big(g_{*}^{\prime} h^{\prime}(1 +  \hbar^{3}g_{*})(1 +  \hbar^{3}h)\big)
\mathrm{d}\sigma \mathrm{d}v_{*},
\een
where according to \cite{erdHos2004quantum} and \cite{benedetto2005weak} the Boltzmann kernel $B_{\Phi,\hbar}(v- v_{*},\sigma)$ has the following form
\ben \label{scaling-Boltzmann-kernel}
B_{\Phi,\hbar}(v- v_{*},\sigma) \colonequals   \hbar^{-4} |v-v_{*}| \big(
\hat{\Phi} (\hbar^{-1} |v-v^{\prime}|)
+ \hat{\Phi} (\hbar^{-1} |v-v^{\prime}_{*}|)
\big)^{2}.
\een
Here $\hbar$ is the Plank constant. In \eqref{scaling-Boltzmann-kernel}, the radial function $\hat{\Phi}(|\xi|) \colonequals   \hat{\Phi}(\xi) = \int_{\mathbb{R}^{3}} e^{-\mathrm{i}x\cdot\xi}\Phi(x)\mathrm{d}x$ is the Fourier transform of a radial potential function  $\Phi$.

In \eqref{U-U-operator} and the rest of the article,
we use the convenient  shorthand $h=h(v)$, $g_*=g(v_*)$,
$h'=h(v')$, $g'_*=g(v'_*)$ where $v'$, $v_*'$ are given by
\ben\label{v-prime-v-prime-star}
v'=\frac{v+v_{*}}{2}+\frac{|v-v_{*}|}{2}\sigma, \quad  v'_{*}=\frac{v+v_{*}}{2}-\frac{|v-v_{*}|}{2}\sigma, \quad \sigma\in\SS^{2}.
\een
Now it remains to see the notation $\mathrm{D}(\cdot)$ in \eqref{U-U-operator}. For $n=1$ or $n=2$, we denote
\ben \label{shorthand-D}
\mathrm{D}^{n}(f(v,v_{*},v^{\prime},v^{\prime}_{*})) \colonequals \left(f(v,v_{*},v^{\prime},v^{\prime}_{*}) - f(v^{\prime},v^{\prime}_{*},v,v_{*})\right)^{n}.
\een
If $n=1$, we write
$\mathrm{D}(\cdot) = \mathrm{D}^{1} (\cdot)$.
The term $\mathrm{D}$ is interpreted as ``difference'' before and after collision.
We will consider a singular kernel and always need some function difference to remove the singularity.
The notation $\mathrm{D}(\cdot)$ or $\mathrm{D}^{2}(\cdot)$ indicates which function is
offering help to cancel the singularity.

By the following scaling
\ben \label{scaling-tranform}
\tilde{F}(t,x,v) = \hbar^{3}F(\hbar^{3}t,x,\hbar^{-3}v), \quad \phi(|x|) = \hbar^{4}\Phi(\hbar^{4}|x|),
\een
we can normalize the Plank constant $\hbar$. Indeed, it is easy to check $F$ is a solution to \eqref{quantum-Boltzmann-UU} if and only if $\tilde{F}$ is a solution of the following normalized equation
\ben \label{quantum-Boltzmann-UU-scaling}
\partial _t F +  v \cdot \nabla_{x} F=Q_{\phi}(F,F), ~~t > 0, x \in \mathbb{T}^{3}, v \in \R^3 ; \quad
F|_{t=0}(x,v) = \hbar^{3}F_{0}(x,\hbar^{-3}v),
\een
where the operator $Q_{\phi}$ is defined by
\ben \label{U-U-operator-scaling}
Q_{\phi}(g,h)(v) \colonequals   \int_{{\mathbb S}^2 \times {\mathbb R}^3} B_{\phi}(v- v_{*},\sigma)
\mathrm{D}\big(g_{*}^{\prime} h^{\prime}(1 +  g_{*})(1 +  h)\big)
\mathrm{d}\sigma \mathrm{d}v_{*}
\een
with the kernel $B_{\phi}(v- v_{*},\sigma)$ given by
\ben \label{scaling-Boltzmann-kernel-scaling}
B_{\phi}(v- v_{*},\sigma) \colonequals   |v-v_{*}| \big(
\hat{\phi} (|v-v^{\prime}|)
+ \hat{\phi} (|v-v^{\prime}_{*}|)
\big)^{2}.
\een

In this article, we will take the inverse power law potential $\phi(x) = |x|^{-p}$. Before that, we give a short review on some relevant mathematical research of the quantum Boltzmann equation.

\subsection{Existing results} We recall some existing mathematical results on the quantum Boltzmann equation in this subsection.
The quantum Boltzmann equations include two models: one for bosons or Bose-Einstein(B-E) particles, the other for fermions or Fermi-Dirac(F-D) particles. The quantum Boltzmann equation for B-E particles is written by \eqref{quantum-Boltzmann-UU-scaling}, \eqref{U-U-operator-scaling} and \eqref{scaling-Boltzmann-kernel-scaling} and usually referred as ``BBE equation''. If the four ``$+$'' in \eqref{U-U-operator-scaling} and the one ``$+$'' in \eqref{scaling-Boltzmann-kernel-scaling} are replaced by minus sign ``$-$'', we will get the quantum Boltzmann equation for F-D particles which is usually referred as ``BFD equation''.

We begin with the mathematical results on BFD equations. In spatially homogeneous case $F=F(t,v)$, we refer the monograph Escobedo-Mischler-Valle \cite{escobedo2003homogeneous} for global existence and weak  convergence to equilibrium.
We also refer Lu \cite{lu2001spatially} and Lu-Wennberg \cite{lu2003stability} for weak  and strong  convergence to equilibrium of mild solutions. On the whole space $x\in\mathbb{R}^{3}$,
Dolbeault \cite{dolbeault1994kinetic} established
global existence and uniqueness of mild solutions for globally integrable kernels. Lions \cite{lions1994compactness} proved global existence of distributional solutions for locally integrable kernels. Alexandre \cite{alexandre2000some} proved global existence of H-solutions for non-cutoff kernels that allow angular singularity given by inverse power law potentials.
Based on the averaging compactness result in \cite{lu2006boltzmann}, on the torus $x\in\mathbb{T}^{3}$,
Lu \cite{lu2008boltzmann} proved global existence of weak solutions for kernels with very soft potentials.
Recently in a perturbation framework, on the whole space $x\in\mathbb{R}^{3}$, Jiang-Xiong-Zhou \cite{jiang2021incompressible} went further to consider the incompressible Navier-Stokes-Fourier limit from the BFD equation with hard sphere collisions.

Unlike Fermi-Dirac particles whose density has a natural bound $0 \leq F \leq 1$ due to Pauli's exclusion principle, density of Bose-Einstein particles may blow up in finite time, which corresponds to the intriguing phenomenon: Bose-Einstein condensation(BEC) in low temperature. As a consequence, the mathematical study of BBE equations is much more difficult than BFD equations. As a result, most of existing results on BBE equations are concerned with isotropic solutions $F(t,v)=F(t,|v|)$ in the homogeneous case with non-singular kernel, for instance hard sphere model or hard potentials with some cutoff.
Lu \cite{lu2000modified} proved global existence of $L^{1}$ solutions for some kernels with strong cutoff assumption (not satisfied by the hard sphere model). Based on the results in \cite{lu2000modified},
Lu \cite{lu2004isotropic} further established global existence of conservative distributional (measure-valued) isotropic solutions for some kernels including the hard sphere model. Also see
Escobedo-Mischler-Vel{\'a}zquez \cite{escobedo2007fundamental} and Escobedo-Mischler-Vel{\'a}zquez \cite{escobedo2008singular} for some singular (near $|v|=0$) solutions.

In terms of long time behavior of the conservative measure-valued isotropic solutions, we refer
Lu \cite{lu2005boltzmann} for strong convergence to equilibrium and single point concentration. With some local condition on the initial datum, long time strong convergence to equilibrium was proved in \cite{lu2016long} and \cite{lu2018long}. Remarkably for any temperature and without any local condition on the initial datum,
Cai-Lu \cite{cai2019spatially} provided algebraic rate of strong convergence to equilibrium.

BEC is an interesting physical phenomenon that deserves deep mathematical understanding. It is fortunate to see many excellent mathematical results (for instance, Spohn \cite{spohn2010kinetics}, Lu \cite{lu2013boltzmann, lu2014boltzmann}, Escobedo-Vel{\'a}zquez \cite{escobedo2014blow, escobedo2015finite}) on this topic.

Note that all the above results on BBE equations are concerned with isotropic solutions $F(t,v)=F(t,|v|)$.
In the homogeneous case,  local existence of anisotropic solutions was first proved in Briant-Einav \cite{briant2016cauchy}, while global existence was first proved in Li-Lu \cite{li2019global} for very high temperature.

Note that all the above results on BBE equations are obtained in the homogeneous case. In the inhomogeneous case,
there are at least two recent results:  Bae-Jang-Yun \cite{bae2021relativistic} and  Ouyang-Wu \cite{ouyang2021quantum}. These two works have different focuses and will be revisited later.

\subsection{Potential function and Boltzmann kernel} In weak-coupling regime for bosons, the Boltzmann kernel $B_{\phi}$ depends on the potential function $\phi$ via \eqref{scaling-Boltzmann-kernel-scaling}. As mentioned before, existing results on BBE equations are mostly concerned with the hard sphere model. This amounts to taking the Dirac delta function as the potential function $\phi(x) = \delta(x)$ and thus $B_{\phi} = (v- v_{*},\sigma) = C |v-v_{*}|$ for some universal constant $C>0$. Observe that the kernel in the quantum case is the same as that in the classical case.
Another physically relevant potential is the inverse power law $\phi(x) = |x|^{-p}$ which has been extensively studied in the research of classical Boltzmann equations but has rarely been considered in the quantum case. To our best knowledge,
this article seems to be the first to study BBE equation with inverse power law.

Note that $p=1$ corresponds to the famous Coulomb potential which is the critical case for Boltzmann equation to be meaningful.
In this article, we work in dimension 3 where $p=3$ is the critical value such that $|x|^{-p}$ is locally integrable near $|x|=0$. Fix $1<p<3$, then the Fourier transform of $\phi(x) = |x|^{-p}$ is
$
\hat{\phi}(\xi) = \hat{\phi}(|\xi|) = C |\xi|^{p-3}.
$ Let $\theta$ be the angle between $v-v_{*}$ and $\sigma$, then
$|v-v_{*}| \sin \frac{\theta}{2} = |v-v^{\prime}|, |v-v_{*}| \cos \frac{\theta}{2} = |v-v^{\prime}_{*}|$.
As a result, the Boltzmann kernel given by \eqref{scaling-Boltzmann-kernel-scaling} is
\beno
B_{\phi}(v- v_{*},\sigma) = C |v-v_{*}|^{2p-5} (\sin^{p-3} \frac{\theta}{2} + \cos^{p-3} \frac{\theta}{2})^{2}.
\eeno
Here and in the rest of the article $C$ denotes a constant that depends only on fixed parameters and could
change from line to line.
Due the symmetry structure of \eqref{scaling-Boltzmann-kernel-scaling}, we can always assume $0 \leq \theta \leq \pi/2$. Then $\sin \frac{\theta}{2} \leq \cos \frac{\theta}{2}$, since $p<3$, we get
\beno
B_{\phi}(v- v_{*},\sigma) \leq C |v-v_{*}|^{2p-5} \sin^{2p-6} \frac{\theta}{2}.
\eeno
Since $p>1$, the following integral over $\mathbb{S}^{2}$ is bounded,
\ben \label{ub-ipl-mean-momentum-transfer}
\int B_{\phi}(v- v_{*},\sigma) \sin^{2} \frac{\theta}{2} \mathrm{d}\sigma \leq C |v-v_{*}|^{2p-5}.
\een
As in \eqref{ub-ipl-mean-momentum-transfer}, in the rest of the article, when taking integration, the range of some frequently
used variables will be omitted if there is no ambiguity. For instance,
the usual ranges $\sigma \in \mathbb{S}^{2}, x \in \mathbb{T}^{3}, v, v_{*} \in \mathbb{R}^{3}$ will be consistently used unless otherwise specified. Whenever a new variable appears, we will specify its range once and then omit it thereafter.

By \eqref{ub-ipl-mean-momentum-transfer}, for $1<p<3$ the kernel $B_{\phi}$ has finite momentum transfer which is a basic condition for the classical Boltzmann equation to be well-posed.
The constant $C$ in \eqref{ub-ipl-mean-momentum-transfer} blows up as $p \rightarrow 1^{+}$ since $\int_{0}^{\sqrt{2}/{2}} t^{2p-3} \mathrm{d}t \sim \frac{1}{p-1}$ for $1<p<3$.
If $p > 2$, the angular function $\sin^{2p-6} \frac{\theta}{2}$ satisfies Grad's angular cutoff assumption
\beno
\int \sin^{2p-6} \frac{\theta}{2} \mathrm{d}\sigma \lesssim  \frac{1}{p-2}.
\eeno
Note that $p=2$ is the critical value such that Grad's cutoff assumption fails. To summarize,
$2<p<3$ can be seen as angular cutoff while $1<p<2$ can be seen as angular non-cutoff.
We consider the much harder case $1<p<2$ in this article. By taking $s=2-p, \gamma=2p-5$, it suffices to consider the more general kernel,
\ben \label{Boltzmann-kernel-general-and-old}
B(v- v_{*},\sigma) \colonequals   C |v-v_{*}|^{\gamma} (\sin^{-1-s} \frac{\theta}{2} + \cos^{-1-s} \frac{\theta}{2})^{2} \mathrm{1}_{0 \leq \theta \leq \pi/2}, \quad -3<\gamma < 0<s<1, \gamma+2s \leq 0.
\een
The parameter pair $(\gamma, s)$ is commonly used in the study of classical Boltzmann equations with inverse power law. We find the resulting kernel \eqref{Boltzmann-kernel-general-and-old} has close relation (see \eqref{Boltzmann-kernel-general-ub} for details) to the Boltzmann kernel
$B^{ipl}(v- v_{*},\sigma)$ defined in \eqref{inverse-power-law-kernel}. The condition $\gamma+2s \leq 0$ is referred as (very) soft potentials. From now on, the notation $B(v- v_{*},\sigma)$ or $B$ stands for the kernel in \eqref{Boltzmann-kernel-general-and-old} unless otherwise specified.

To reiterate,
in this article we will study the following Cauchy problem
\ben \label{quantum-Boltzmann-CP}
\partial _t F +  v \cdot \nabla_{x} F=Q(F,F), ~~t > 0, x \in \mathbb{T}^{3}, v \in \R^3 ; \quad
F|_{t=0}(x,v) = F_{0}(x,v),
\een
where the operator $Q$ is defined through \eqref{U-U-operator-scaling} with $B_{\phi}$ replaced by $B$ in \eqref{Boltzmann-kernel-general-and-old}. Here $F_{0}$ is a given initial datum satisfying some high temperature condition which will be given in the next subsection.

\subsection{Temperature and initial datum} Temperature plays an important role in the study of quantum Boltzmann equation. For example, for B-E particles, BEC will happen in low temperature. We now introduce some basic knowledge about temperature in the quantum context.
Let us consider a homogeneous density $f=f(v) \geq 0$ with zero mean $\int v f(v) \mathrm{d}v=0$. For $k \geq 0$, we recall the moment function
\beno
M_{k}(f) \colonequals   \int |v|^{k} f(v)  \mathrm{d}v.
\eeno
Let $M_{0} = M_{0}(f), M_{2}=M_{2}(f)$ for simplicity. Let $m$ be the mass of a particle, then
 $m M_{0}$ and $\f{1}{2}m M_{2}$ are the total mass and kinetic energy per unit space volume.
Referring \cite{lu2004isotropic},
the kinetic temperature $\bar{T}$ and the critical temperature $\bar{T}_{c}$ of the particle system are defined by
\ben \label{kinetic-temperature-and-critical}
\bar{T} = \frac{1}{3 k_{B}} \frac{m M_{2}}{M_{0}}, \quad \bar{T}_{c} = \frac{m\zeta(5/2)}{2 \pi k_{B} \zeta(3/2)}
 \big( \frac{M_{0} }{\zeta(3/2)} \big)^{\f23},
\een
where $k_{B}$ is the Boltzmann constant and $\zeta(s) = \sum_{n=1}^{\infty} \frac{1}{n^{s}}$ is the Riemann zeta function. 

We now recall some known facts about equilibrium distribution. The equilibrium of the classical Boltzmann equation is the Maxwellian distribution with density $\mu_{\rho, v_{0}, T}$ function defined by
\beno
\mu_{\rho, v_{0}, T}(v) \colonequals   \rho(2\pi T)^{-\frac{3}{2}}e^{-\frac{1}{2T}|v-v_{0}|^{2}}, \quad \rho, T>0, v_{0} \in \mathbb{R}^{3}.
\eeno
Here  $\rho$ is density, $v_{0}$ is mean velocity  and $T$ is temperature. The famous Bose-Einstein distribution has density
function
\ben \label{equilibrium}
\mathcal{M}_{\rho, v_{0}, T} \colonequals   \frac{\mu_{\rho, v_{0}, T}}{1 - \mu_{\rho, v_{0}, T}}, \quad \rho(2\pi T)^{-\frac{3}{2}} \leq 1. \een
Now $\rho, v_{0}$ and $T$ do not represent density, mean velocity and temperature anymore, but only three parameters. The ratio $\bar{T}/\bar{T}_{c}$ quantifies high and low temperature.
In high temperature  $\bar{T}/\bar{T}_{c} > 1$,  the equilibrium of BBE equation is the Bose-Einstein distribution \eqref{equilibrium} with
$\rho(2\pi T)^{-\frac{3}{2}} < 1$. In low temperature  $\bar{T}/\bar{T}_{c} < 1$,  the equilibrium of BBE equation is the Bose-Einstein distribution \eqref{equilibrium} with
$\rho(2\pi T)^{-\frac{3}{2}} = 1$ plus some Dirac delta function. That is, the equilibrium contains a Dirac measure. In the critical case  $\bar{T}/\bar{T}_{c} = 1$, the equilibrium is \eqref{equilibrium} with
$\rho(2\pi T)^{-\frac{3}{2}} = 1$. One can refer to \cite{lu2005boltzmann} for the classification of equilibria.

In this article, we work with high temperature and thus the equilibrium of \eqref{quantum-Boltzmann-UU}
is $\mathcal{M}_{\rho, v_{0}, T}$ defined in \eqref{equilibrium}.
For perturbation around equilibrium, we define
\ben \label{multiplier-function} \mathcal{N}_{\rho, v_{0}, T} \colonequals   \mathcal{M}_{\rho, v_{0}, T}^{\f{1}{2}} (1 +  \mathcal{M}_{\rho, v_{0}, T})^{\f{1}{2}} = \frac{\mu_{\rho, v_{0}, T}^{\f{1}{2}}}{1 - \mu_{\rho, v_{0}, T}}. \een
We remark that the function $\mathcal{N}_{\rho, v_{0}, T}$ serves as the multiplier in the expansion $F = \mathcal{M}_{\rho, v_{0}, T} + \mathcal{N}_{\rho, v_{0}, T} f$.

Recall that the solution of \eqref{quantum-Boltzmann-UU} conserves mass, momentum and energy. That is, for any $t \geq 0$,
\ben \label{conversation-mass-momentum-energy}
\int (1, v, |v|^{2}) F(t,x,v) \mathrm{d}x \mathrm{d}v = \int (1, v, |v|^{2}) F_{0}(x,v) \mathrm{d}x \mathrm{d}v.
\een
Once $F_{0}$ is appropriately given, the constants $\rho,T>0, v_{0} \in \mathbb{R}^{3}$ in \eqref{equilibrium} are uniquely determined through 
\ben \label{F0-determines-equilibrium}
\int (1, v, |v|^{2}) F_{0}(x,v) \mathrm{d}x \mathrm{d}v = \int (1, v, |v|^{2}) \mathcal{M}_{\rho, v_{0}, T}(v) \mathrm{d}x \mathrm{d}v.
\een
Without any loss of generality, we assume that $F_{0}$ has zero mean and thus gives $v_{0}=0$ from now on. Also without any loss of generality, we assume that $F_{0}$ gives $T=1$ in this article. Indeed, we can make the transform $f(v) \to f(T^{1/2}v)$ to reduce the general $T \neq 1$ case to the special case $T=1$.
As a result, we only keep $\rho$ as a parameter. That is, we only consider those initial data with $T=1, v_{0}=0$  according to \eqref{F0-determines-equilibrium}. Taking $\mu (v) \colonequals   (2\pi)^{-\frac{3}{2}} e^{-\f{1}{2}|v|^{2}},$ the equilibrium $\mathcal{M}_{\rho, v_{0}, T}$ and the multiplier function $\mathcal{N}_{\rho, v_{0}, T}$ reduce to
\ben \label{equilibrium-rho}
\mathcal{M}_{\rho} \colonequals   \frac{\rho \mu}{1 - \rho \mu}, \quad \mathcal{N}_{\rho} \colonequals   \frac{\rho^{\f{1}{2}} \mu^{\f{1}{2}}}{1 - \rho \mu}.  \een

When $\rho$ is small, it is easy to see
$
M_{2}(\mathcal{M}_{\rho}) \sim \rho, M_{0}(\mathcal{M}_{\rho}) \sim \rho$
and so by recalling \eqref{kinetic-temperature-and-critical},
\beno
\frac{\bar{T}}{\bar{T}_{c}} \sim \rho^{-\frac{2}{3}}.
\eeno
In our main result Theorem  \ref{global-well-posedness}, we will assume  $0< \rho \leq \rho_{*}$ for some small constant $0<\rho_{*} \ll 1$ which means
\ben \label{high-temperature-condition}
\frac{\bar{T}}{\bar{T}_{c}} \gtrsim \rho_{*}^{-\frac{2}{3}} \gg 1.
\een
That is, we need high temperature assumption. Note that high temperature assumption is also imposed in \cite{li2019global} to prove global well-posedness of homogeneous BBE equation with (slightly general than) hard sphere collisions.

\subsection{Perturbation around equilibrium and main result}
For simplicity, let $\mathcal{M}\colonequals  \mathcal{M}_{\rho}, \mathcal{N}\colonequals  \mathcal{N}_{\rho}$. With the expansion $F = \mathcal{M} + \mathcal{N} f$,
the linearized quantum Boltzmann equation corresponding to \eqref{quantum-Boltzmann-CP} reads
\ben\label{linearized-quantum-Boltzmann-eq} \left\{ \begin{aligned}
&\partial _t f +  v \cdot \nabla_{x} f + \mathcal{L}^{\rho}f = \Gamma_{2}^{\rho}(f,f) + \Gamma_{3}^{\rho}(f,f,f), ~~t > 0, x \in \mathbb{T}^{3}, v \in \R^3 ;\\
&f|_{t=0} = f_{0} = \frac{(1-  \rho\mu)F_{0} - \rho\mu}{\rho^{\f{1}{2}} \mu^{\f{1}{2}}}.
\end{aligned} \right.
\een
Here the linearized quantum Boltzmann operator $\mathcal{L}^{\rho}$ is define by
\ben \label{linearized-quantum-Boltzmann-operator-UU}
(\mathcal{L}^{\rho}f)(v)  \colonequals   \int  B \mathcal{N}_{*} \mathcal{N}^{\prime} \mathcal{N}^{\prime}_{*} \mathrm{S}(\mathcal{N}^{-1}f)
\mathrm{d}\sigma \mathrm{d}v_{*},
\een
where $\mathrm{S}(\cdot)$ is defined by
\ben  \label{symmetry-operator}
\mathrm{S}(g) \colonequals     g + g_{*} - g^{\prime} - g^{\prime}_{*}.
\een
The bilinear term $\Gamma_{2}^{\rho}(\cdot,\cdot)$ and the trilinear term $\Gamma_{3}^{\rho}(\cdot,\cdot,\cdot)$ are defined by
\ben \label{definition-Gamma-2-epsilon}
\Gamma_{2}^{\rho}(g,h) &\colonequals&   \mathcal{N}^{-1}\int
B \Pi_{2}(g,h)
\mathrm{d}\sigma \mathrm{d}v_{*}.
\\ \label{definition-Gamma-3-epsilon}
\Gamma_{3}^{\rho}(g,h,\varrho)  &\colonequals&     \mathcal{N}^{-1}\int
B  \mathrm{D} \big( (\mathcal{N}g)_{*}^{\prime} (\mathcal{N}h)^{\prime} ((\mathcal{N}\varrho)_{*} + \mathcal{N}\varrho) \big) \mathrm{d}\sigma \mathrm{d}v_{*}. \quad
\quad
\een
The notation $\Pi_{2}$ in \eqref{definition-Gamma-2-epsilon} is defined by
 \ben   \label{definition-A-2}
\Pi_{2}(g,h) &\colonequals  & \mathrm{D} \big((\mathcal{N}g)^{\prime}_{*}(\mathcal{N}h)^{\prime}\big)
\\ \label{line-1} &&+
\mathrm{D}\big((\mathcal{N}g)^{\prime}_{*}(\mathcal{N}h)^{\prime}(\mathcal{M} +  \mathcal{M}_{*})\big)
\\ \label{line-2} &&+ \mathrm{D}\big( (\mathcal{N}g)_{*}(\mathcal{N}h)^{\prime}( \mathcal{M}^{\prime}_{*} -  \mathcal{M}) \big)
\\ \label{line-3} &&+ \big( (\mathcal{N}g)^{\prime}(\mathcal{N}h)\mathrm{D}(\mathcal{M}^{\prime}_{*}) + (\mathcal{N}g)^{\prime}_{*}(\mathcal{N}h)_{*}\mathrm{D}(\mathcal{M}^{\prime})\big).
\een
Remark that the three operators $\mathcal{L}^{\rho}$, $\Gamma_{2}^{\rho}(\cdot,\cdot)$ and $\Gamma_{3}^{\rho}(\cdot,\cdot,\cdot)$ depends on $\rho$ through $\mathcal{M} =  \mathcal{M}_{\rho},  \mathcal{N} =  \mathcal{N}_{\rho}$.

Our goal is to prove global well-posedness of \eqref{linearized-quantum-Boltzmann-eq} in some weighted Sobolev space.
More precisely, we use the following  energy and dissipation functional
\ben \label{definition-energy-and-dissipation}
\mathcal{E}_{N}(f) \colonequals   \sum_{|\alpha|+|\beta| \leq N} \|W_{l_{|\alpha|,|\beta|}} \partial^{\alpha}_{\beta} f \|_{L^{2}_{x}L^{2}}^{2}, \quad \mathcal{D}_{N}(f) \colonequals   \sum_{|\alpha|+|\beta| \leq N}
\|W_{l_{|\alpha|,|\beta|}} \partial^{\alpha}_{\beta} f\|_{L^{2}_{x}\mathcal{L}^{s}_{\gamma/2}}^{2},
\een
where $\partial^{\alpha}_{\beta}\colonequals   \partial^{\alpha}_{x}\partial^{\beta}_{v}$.
See subsection \ref{notation} for the definition of $\|\cdot\|_{L^{2}_{x}L^{2}}$ and $\|\cdot\|_{L^{2}_{x}\mathcal{L}^{s}_{\gamma/2}}$. See \eqref{definition-of-norm-L-epsilon-gamma} for the definition of $|\cdot|_{\mathcal{L}^{s}_{\gamma/2}}$. For the moment, just keep in mind $\|\cdot\|_{L^{2}_{x}\mathcal{L}^{s}_{\gamma/2}}$ is the dissipation corresponding to the energy $\|\cdot\|_{L^{2}_{x}L^{2}}$.
Here $W_{l}(v)\colonequals   (1+|v|^{2})^{\frac{l}{2}}$ is a polynomial weight on the velocity variable $v$. The weight order $l_{|\alpha|,|\beta|} \geq 0$ depends on the derivative order $|\alpha|,|\beta|$
and the sequence  $\{l_{|\alpha|, |\beta|}\}_{|\alpha|+|\beta| \leq N}$ verifies
\ben \label{weight-order-condition-1}
 l_{|\alpha|, |\beta|} + |\gamma+2s| \leq l_{|\alpha|+1, |\beta|-1}, \quad \text{ for } |\alpha| \leq N-1,  1 \leq  |\beta| \leq N, |\alpha|+|\beta| \leq N.
 \\
 \label{weight-order-condition-2}
 l_{|\alpha|,0} \leq l_{0,|\alpha|-1}, \quad \text{ for } 1 \leq |\alpha| \leq N.
\een

The condition \eqref{weight-order-condition-1} is used to deal with linear streaming term $v \cdot \nabla_{x} f$ as $\gamma+2s \leq 0$ (see \eqref{streaming-term-analysis} below). The two conditions \eqref{weight-order-condition-1} and \eqref{weight-order-condition-2} together ensure that $l_{|\alpha|,|\beta|}$ increases as $|\alpha|+|\beta|$ decreases (see \eqref{weight-order-condition-on-total-regularity}).

Now we are ready to give global well-posedness of \eqref{linearized-quantum-Boltzmann-eq} in the following theorem.
\begin{thm}\label{global-well-posedness} Let $N \geq 9$. There exist two universal constants  $\rho_{*}, \delta_{*}>0$ such that the following global well-posedness is valid. Let $0< \rho \leq \rho_{*}$.
If
\ben \label{condition-on-initial-data}
\frac{\rho \mu}{1 - \rho \mu} + \frac{\rho^{\f{1}{2}}\mu^{\f{1}{2}}}{1 - \rho \mu}f_{0} \geq 0, \quad \mathcal{E}_{N}(f_{0}) \leq \delta_{*} \rho^{2N+1},
\een
then the Cauchy problem \eqref{linearized-quantum-Boltzmann-eq} has a unique global solution $f^{\rho} \in L^{\infty}([0,\infty); \mathcal{E}_{N})$ satisfying for some universal constant $C$,
\ben \label{uniform-estimate-global} \sup_{t \geq 0}\mathcal{E}_{N}(f^{\rho}(t)) +  \rho \int_{0}^{\infty}\mathcal{D}_{N}(f^{\rho}(t)) \mathrm{d}t \leq C \rho^{-2N} \mathcal{E}_{N}(f_{0}),\een
and	for all $t \geq 0$,
\beno
\frac{\rho \mu}{1 - \rho \mu} + \frac{\rho^{\f{1}{2}}\mu^{\f{1}{2}}}{1 - \rho \mu}f^{\rho}(t) \geq 0.
\eeno
\end{thm}

We give some explanations and comments on Theorem \ref{global-well-posedness} in the rest of this subsection and the next two subsections.

We emphasize that the constants $\rho_{*}, \delta_{*}, C$ in Theorem \ref{global-well-posedness} could depend on $N, l_{0,0}, s, \gamma$. In Theorem \ref{global-well-posedness} and the rest of the article,  if a constant depends only on the given parameters $N, l_{0,0}, s, \gamma$, we say it is  a universal constant. In particular, the two constants $\delta_{*}, C$ are independent of $0< \rho \leq \rho_{*}$.

The condition $N \geq 9$ owes to Theorem \ref{Gamma-3-energy-estimate}. Understandably, a certain high order is required to deal with the singular kernel \eqref{Boltzmann-kernel-general-and-old} and the trilinear operator $\Gamma_{3}^{\rho}(\cdot,\cdot,\cdot)$. Recall that the hard sphere model in \cite{ouyang2021quantum} needs $N \geq 8$.

Note that Theorem \ref{global-well-posedness} gives global well-posedness under the condition
\beno
0 <\rho \ll 1, \quad \mathcal{M}_{\rho} = \frac{\rho \mu}{1 - \rho \mu} \sim \rho \mu, \quad \mathcal{N}_{\rho} = \frac{\rho^{\f{1}{2}} \mu^{\f{1}{2}}}{1 - \rho \mu} \sim \rho^{\f{1}{2}} \mu^{\f{1}{2}}, \quad
\mathcal{E}_{N} \big(\frac{(1-  \rho\mu)F_{0} - \rho\mu}{\rho^{\f{1}{2}} \mu^{\f{1}{2}}} \big) \ll \rho^{2N+1}.
\eeno
Roughly speaking, these conditions mean
\beno
F_{0} \approx \rho \mu + \mathrm{o}(\rho^{N+1}) \mu^{\f{1}{2}}.
\eeno
Note that this $\rho$-related smallness did not appear in previous works where the parameter $\rho$ was regarded as a fixed constant.
Other effects of viewing $\rho$ as a variable parameter will be pointed out later.

By \eqref{condition-on-initial-data} and \eqref{uniform-estimate-global}, we have $\mathcal{E}_{N}(f^{\rho}(t)) \lesssim \rho$ for any $t \geq 0$. This estimate is consistent with the trivial case $\rho=0$ where $F(t) \equiv 0$ is the solution to problem \eqref{quantum-Boltzmann-CP} starting with a zero initial datum $F_{0} \equiv 0$.

In terms of physical relevance of Theorem \ref{Gamma-3-energy-estimate}, we give the following two remarks.
\begin{rmk} Note that $\|f^{\rho}(t)\|_{L^{\infty}_{x,v}}^{2} \lesssim \mathcal{E}_{N}(f^{\rho}(t)) \lesssim \rho$ and thus
$\|\frac{\rho \mu}{1 - \rho \mu} + \frac{\rho^{\f{1}{2}}\mu^{\f{1}{2}}}{1 - \rho \mu}f^{\rho}(t)\|_{L^{\infty}_{x,v}}\lesssim \rho$ for any $t \geq 0$.
Theorem \ref{global-well-posedness} shows that BEC will not happen if the initial datum is close enough to
an equilibrium with high temperature, which is reasonable and consistent with physical observation.
\end{rmk}
\begin{rmk}
By the scaling \eqref{scaling-tranform}, Theorem \ref{global-well-posedness} ensures the global well-posedness of \eqref{quantum-Boltzmann-UU} with $\Phi(x) = \hbar^{4p-4}|x|^{-p}$ and initial datum $F_{0}$ close enough to the equilibrium $\mathcal{M}_{\rho, \hbar}$ where $\mathcal{M}_{\rho, \hbar}(v) = \hbar^{-3}\mathcal{M}_{\rho}(\hbar^{3} v)$. Simply looking at the equilibrium $\mathcal{M}_{\rho, \hbar}$, we have $\|\mathcal{M}_{\rho, \hbar}\|_{L^{\infty}_{x,v}} \sim \hbar^{-3} \rho, \|\mathcal{M}_{\rho, \hbar}\|_{L^{2}_{x}L^{2}} \sim \hbar^{-\f{15}{2}} \rho.$ Since $\hbar$ is a relatively small constant, the magnitude of the solution to the problem \eqref{quantum-Boltzmann-UU} can be relatively large.
\end{rmk}

\subsection{Feature of our result} In terms of condition and conclusion, Theorem \ref{global-well-posedness} is the closest to the main result (Theorem 1.4) of \cite{li2019global}. More precisely, both of these two results validate global existence of anisotropic solutions under very high temperature condition, see \eqref{high-temperature-condition} and Remark 1.6 of \cite{li2019global} for more details. The main differences of these two results are also obvious. To reiterate, \cite{li2019global} considers spatially homogeneous case and (slightly general than) hard sphere model, while we work with the inhomogeneous case and non-cutoff kernels.

In terms of mathematical methods, this article is closer to \cite{bae2021relativistic}, \cite{ouyang2021quantum} and \cite{jiang2021incompressible} since all of these works fall into the close-to-equilibrium framework well established for the classical Boltzmann equation. There are many great works that contribute to this mathematically satisfactory theory for global well-posedness of the classical Boltzmann equation. For reference, we mention
\cite{ukai1974existence,guo2003classical} for angular cutoff kernels and \cite{gressman2011global,alexandre2012boltzmann} for non-cutoff kernels.


Each of  \cite{bae2021relativistic}, \cite{ouyang2021quantum} and \cite{jiang2021incompressible} has its own features and focuses. The work \cite{bae2021relativistic} is the first to investigate both relativistic and quantum effects.
The article \cite{jiang2021incompressible} studies the hydro dynamic limit from BFD (but not BBE) to incompressible Navier-Stokes-Fourier equation.  The work \cite{ouyang2021quantum} includes three cases: torus near equilibrium, whole space near equilibrium and whole space near vacuum. These novel works contribute to the literature of quantum Boltzmann equation from different aspects.

Like the above works, our main result Theorem \ref{global-well-posedness} also has some unique features that may better our understanding of quantum Boltzmann equation. In particular, this article may be the first in the spatially inhomogeneous case to
\begin{itemize}
\item study the quantum Boltzmann equation with inverse power law potentials;
\item incorporate high temperature condition into global well-posedness;
\item rule out BEC globally in time.
\end{itemize}
Different from the three works (\cite{bae2021relativistic}, \cite{ouyang2021quantum} and \cite{jiang2021incompressible}),
we keep the parameter $\rho$ along our derivation throughout the article. This intentional choice enables us to relate the high temperature condition to the smallness of $\rho$ quantitatively in \eqref{high-temperature-condition}.
As a result, BEC is rigorously ruled out globally.

\subsection{Possible improvements} \label{possible}
Compared to the well-established global well-posedness theory for
classical Boltzmann equations with inverse power law potentials, our Theorem \ref{global-well-posedness} has much room to improve. In this subsection, some possible improvements of Theorem \ref{global-well-posedness} are given based on the author's limited knowledge. To keep the present article in a reasonable length, we leave these improvements in future works.

Recall $\mathcal{E}_{N}$ is defined in \eqref{definition-energy-and-dissipation} with
the weight order $\{l_{|\alpha|, |\beta|}\}_{|\alpha|+|\beta| \leq N}$ satisfying \eqref{weight-order-condition-1} and \eqref{weight-order-condition-2}. Let $\mathcal{E}_{N}^{*}$ be the space with the minimal weight order $\{l^{*}_{|\alpha|, |\beta|}\}_{|\alpha|+|\beta| \leq N}$ where $l^{*}_{0, N}=0$ and $\{l^{*}_{|\alpha|, |\beta|}\}_{|\alpha|+|\beta| \leq N}$ verifies the two conditions \eqref{weight-order-condition-1} and \eqref{weight-order-condition-2} with identities. To illustrate, we give an example of the weight order $\{l^{*}_{|\alpha|, |\beta|}\}_{|\alpha|+|\beta| \leq N}$.
Recall for the inverse power law $\phi(x)=|x|^{-p}$, $\gamma=2-p, s= 2p-5$ and thus $\gamma+2s=-1$.
Let us take $N=9$, the corresponding $\{l^{*}_{|\alpha|, |\beta|}\}_{|\alpha|+|\beta| \leq 9}$ with $\gamma+2s=-1$ is
 shown in Table \ref{weight-order-an-example}. In Table \ref{weight-order-an-example}, the column index represents derivative order $|\alpha|$ of space variable $x$ and the row index represents derivative order $|\beta|$ of velocity variable $v$. For example $l_{0,9} = 0, l_{1,8}=1, l_{0,8}=9, l_{0,1}=44, l_{1,0}=l_{0,0}=45$.
\begin{table}[!htbp]
\centering
\caption{Minimal weight order when $N=9, \gamma+2s=-1$}\label{weight-order-an-example}
\begin{tabular}{c|c|c|c|c|c|c|c|c|c|c|c}
\hline
9 & 0 & $\times$ & $\times$ & $\times$ & $\times$ & $\times$ & $\times$ & $\times$ & $\times$ & $\times$ & $\times$\\
\hline
8 & 9 & 1 & $\times$ & $\times$ & $\times$ & $\times$ & $\times$ & $\times$ & $\times$ & $\times$ & $\times$\\
\hline
7 & 17 & 10 & 2 & $\times$ & $\times$ & $\times$ & $\times$ & $\times$ & $\times$ & $\times$ & $\times$\\
\hline
6 & 24 & 18 & 11 & 3 & $\times$ & $\times$ & $\times$ & $\times$ & $\times$ & $\times$ & $\times$\\
\hline
5 & 30 & 25 & 19 & 12 & 4 & $\times$ & $\times$ & $\times$ & $\times$ & $\times$ & $\times$\\
\hline
4 & 35 & 31 & 26 & 20 & 13 & 5 & $\times$ & $\times$ & $\times$ & $\times$ & $\times$\\
\hline
3 & 39 & 36 & 32 & 27 & 21 & 14 & 6 & $\times$ & $\times$ & $\times$ & $\times$\\
\hline
2 & 42 & 40 & 37 & 33 & 28 & 22 & 15 & 7 & $\times$ & $\times$ & $\times$\\
\hline
1 & 44 & 43 & 41 & 38 & 34 & 29 & 23 & 16 & 8 & $\times$ & $\times$\\
\hline
0 & 45 & 45 & 44 & 42 & 39 & 35 & 30 & 24 & 17 & 9 & $\times$\\
\hline
$(|\beta|/|\alpha|)$ & 0 & 1 & 2 & 3 & 4 & 5 & 6 & 7 & 8 & 9 & 10\\
\hline
\end{tabular}
\end{table}

Note that in Theorem \ref{global-well-posedness}, smallness of $\mathcal{E}_{N}(f_{0})$ is imposed to prove well-posedness in $L^{\infty}([0,\infty); \mathcal{E}_{N})$.
We may try to prove well-posedness in $L^{\infty}([0,\infty); \mathcal{E}_{N})$ under the condition
$\mathcal{E}_{9}^{*}(f_{0})\ll \rho^{19}, \mathcal{E}_{N}(f_{0}) < \infty$. That is,
smallness assumption is only imposed on $\mathcal{E}_{9}^{*}$. Such kind better result was derived by Guo in \cite{guo2012vlasov} for the Vlasov-Possion-Landau system.

Theorem \ref{global-well-posedness} does not consider relaxation to equilibrium. Using the methods in \cite{strain2006almost} and \cite{duan2021global} for classical Boltzmann equations, it is very promising to derive similar long time behaviors for the solution $f^{\rho}$ obtained in Theorem \ref{global-well-posedness}.
For example, we can derive almost exponential decay like in \cite{strain2006almost} by recalling  \eqref{definition-of-norm-L-epsilon-gamma} $|\cdot|_{\mathcal{L}^{s}_{\gamma/2}} \geq |\cdot|_{L^{2}_{s+\gamma/2}}$ and using some interpolation inequality to deal with soft potential $\gamma+2s<0$. In this regard, a possible result may be as follows.
Let $N \geq 9, l_{2} > l_{1} \geq 0$ and assume  $\mathcal{E}_{9}^{*}(f_{0})\ll \rho^{19}, \mathcal{E}_{N}(W_{l_{2}}f_{0}) < \infty$. Try to prove for any $t \geq 0$,
\beno
\mathcal{E}_{N}(W_{l_{2}}f^{\rho}(t)) \lesssim \rho^{-2N}C(\mathcal{E}_{N}(W_{l_{2}}f_{0})), \quad \mathcal{E}_{N}(W_{l_{1}}f^{\rho}(t)) \lesssim (1+t)^{-q}\rho^{-2N}\mathcal{E}_{N}(W_{l_{1}}f_{0}).
\eeno
The first upper bound estimate on $\mathcal{E}_{N}(W_{l_{2}}f^{\rho}(t))$ is used together with interpolation method to derive the polynomial decay of  $\mathcal{E}_{N}(W_{l_{1}}f^{\rho}(t))$. By interpolation, we should have $q = -\frac{l_{2} - l_{1}}{ s+\gamma/2}$.

We may try to derive sub-exponential decay rate under more assumption on initial data like in \cite{duan2021global}. Roughly speaking, we could get something as follows.  Let $\lambda>0$ be small enough and
assume  $\mathcal{E}_{N}(e^{2 \lambda \langle v \rangle} f_{0}) \ll \rho^{2N+1}$. Try to prove for any $t \geq 0$,
\beno
 \mathcal{E}_{N}(e^{ 2 \lambda \langle v \rangle} f^{\rho}(t))\lesssim \rho^{-2N}\mathcal{E}_{N}(e^{2 \lambda \langle v \rangle} f_{0}), \quad \mathcal{E}_{N}(f^{\rho}(t)) \lesssim e^{ -\lambda t^{\kappa}}\rho^{-2N}\mathcal{E}_{N}(e^{ 2 \lambda \langle v \rangle} f_{0}),
\eeno
where $\kappa = \frac{1}{1+|\gamma+2s|}$.

Another topic is to
prove global well-posedness in a larger space than $\mathcal{E}_{9}^{*}$. Note that for classical Boltzmann equation with inverse power law potential, see \cite{duan2021global} for global well-posedness in
the up-to-date largest space  $L^{1}_{k}L^{2}$ (containing $\supset H^{\f{3}{2}+\delta}_{x}L^{2}$ for any $\delta>0$). One may try to establish global well-posedness in such kind low regularity space.

\subsection{Strategy of proof}
The proof of  Theorem \ref{global-well-posedness} will be given in subsection \ref{global-proof} by a rigorous continuity argument based on the local well-posedness result in Theorem \ref{local-well-posedness-LBE} and the {\it a priori} estimate in Theorem \ref{a-priori-estimate-LBE}. We spend this subsection to outline the procedure and give the key points in deriving Theorem \ref{local-well-posedness-LBE} and Theorem \ref{a-priori-estimate-LBE}. {\it The best way is to glance over this subsection  first and come back to read it carefully when appropriate.}

\subsubsection{An auxiliary Cauchy problem}
In order to prove local well-posedness of \eqref{quantum-Boltzmann-CP} or \eqref{linearized-quantum-Boltzmann-eq}, we study the following linear problem
\ben\label{quantum-Boltzmann-UU-linear}
\partial _t F +  v \cdot \nabla_{x} F=\tilde{Q}(G,F), ~~t > 0, x \in \mathbb{T}^{3}, v \in \R^3; \quad
F|_{t=0} = F_{0}.
\een
Here $G$ is a given function and $F$ is unknown. The operator $\tilde{Q}(\cdot,\cdot)$ is defined by
\ben \label{U-U-operator-linear}
\tilde{Q}(g,h)(v) \colonequals    \int  B(v- v_{*},\sigma)
 \mathrm{D}\big(g_{*}^{\prime} h^{\prime}(1 +  g_{*} +  g)\big)
\mathrm{d}\sigma \mathrm{d}v_{*}.
\een
Since $G$ is a given function, the operator $\tilde{Q}(G,F)$  is linear in $F$. Using the expansion $F = \mathcal{M} + \mathcal{N} f, G = \mathcal{M} + \mathcal{N} g$, the equation \eqref{quantum-Boltzmann-UU-linear} is equivalent to
\ben \label{linearized-quantum-Boltzmann-eq-linear} \left\{ \begin{aligned}
&\partial _t f +  v \cdot \nabla_{x} f + \mathcal{L}^{\rho}f =  \mathcal{L}^{\rho}_{r}f + \mathcal{C}^{\rho}f+
\rho^{\f{1}{2}}\Gamma_{2,m}^{\rho}(g, f + \rho^{\f{1}{2}}\mu^{\f{1}{2}}) + \Gamma_{3}^{\rho}(g, f + \rho^{\f{1}{2}}\mu^{\f{1}{2}},g);\\
&f|_{t=0} = f_{0} = \frac{(1-  \rho\mu)F_{0} - \rho\mu}{\rho^{\f{1}{2}} \mu^{\f{1}{2}}}.
\end{aligned} \right.
\een
where the linear operators $\mathcal{L}^{\rho}_{r}, \mathcal{C}^{\rho}$ are defined by
\ben \label{correction-operator-linear}
(\mathcal{C}^{\rho}f)(v)  \colonequals    \rho^{\f{1}{2}}  \int  B \mathcal{N}_{*} \mathcal{N}^{\prime} \mathcal{N}^{\prime}_{*} \mathrm{D}\big((\mu^{\f{1}{2}}f)^{\prime} \big)
\mathrm{d}\sigma \mathrm{d}v_{*},
\\ \label{L-epsilon-pm-tau-rest-part-original}
(\mathcal{L}^{\rho}_{r}f)(v)  \colonequals    \int  B \mathcal{N}_{*} \mathcal{N}^{\prime} \mathcal{N}^{\prime}_{*}
\mathrm{D}\big( (\mathcal{N}^{-1}f)_{*} \big)
\mathrm{d}\sigma \mathrm{d}v_{*},
\een
and the operator $\Gamma_{2,m}^{\rho}$ is defined in \eqref{definition-Gamma-main}. In this article, the subscripts ``$m$'' and
``$r$'' are referred to ``main'' and ``remaining'' respectively. For example,  $\Gamma_{2,m}^{\rho}$ is the main part of $\Gamma_{2}^{\rho}$. Corresponding to the remaining part $\mathcal{L}^{\rho}_{r}$ of $\mathcal{L}^{\rho}$, we also have the main part $\mathcal{L}^{\rho}_{m}$ defined in \eqref{definition-of-lm}. Note that $\mathcal{L}^{\rho}=\mathcal{L}^{\rho}_{m}+\mathcal{L}^{\rho}_{r}$.

Local well-posedness of \eqref{linearized-quantum-Boltzmann-eq} is proved by iterating on equation \eqref{linearized-quantum-Boltzmann-eq-linear}. In order to implement energy method on equations \eqref{linearized-quantum-Boltzmann-eq} and \eqref{linearized-quantum-Boltzmann-eq-linear}, we first give necessary estimates on the linear operators
$\mathcal{L}^{\rho}, \mathcal{L}^{\rho}_{r}, \mathcal{C}^{\rho}, \tilde{Q}(g,\cdot)$, the bilinear operators $\Gamma_{2,m}^{\rho}, \Gamma_{2}^{\rho}$ and the trilinear operator $\Gamma_{3}^{\rho}$ in Section \ref{linear}, \ref{bilinear} and \ref{trilinear} respectively.
We now give some keys ideas of these estimates.

\subsubsection{Linear operators} In this sequel, we consider the four linear operators: $\mathcal{L}^{\rho}, \mathcal{L}^{\rho}_{r}, \mathcal{C}^{\rho}, \tilde{Q}(g,\cdot)$.

{\it Coercivity estimate of $\mathcal{L}^{\rho}$.} Coercivity estimate of $\mathcal{L}^{\rho}$ plays a central role in the close-to-equilibrium framework. Note that $\mathcal{L}^{\rho}$ is a self-joint operator. Indeed,
\ben \label{self-joint-operator}
\langle \mathcal{L}^{\rho}g, h \rangle = \frac{1}{4} \int B
\mathcal{N} \mathcal{N}_{*}  \mathcal{N}^{\prime} \mathcal{N}^{\prime}_{*} \mathrm{S}( \mathcal{N}^{-1} g ) \mathrm{S}(\mathcal{N}^{-1} h) \mathrm{d}\sigma \mathrm{d}v_{*} \mathrm{d}v = \langle g,  \mathcal{L}^{\rho}h \rangle.
\een
The null space of $\mathcal{L}^{\rho}$ is
\ben \label{null-space-of-L}
\mathrm{Null}^{\rho} =  \mathrm{span} \{ \mathcal{N}_{\rho}, \mathcal{N}_{\rho}v_{1}, \mathcal{N}_{\rho}v_{2}, \mathcal{N}_{\rho}v_{3}, \mathcal{N}_{\rho}|v|^{2} \}.
\een
Let $\mathbb{P}_{\rho}$(see \eqref{definition-projection-operator} for its precise definition) be the projection operator on the null space $\mathrm{Null}^{\rho}$.
Then $\mathcal{L}^{\rho}$ enjoys the following estimate (see Theorem \ref{main-theorem} for the precise statement)
\ben \label{coercivity-result}
 \langle \mathcal{L}^{\rho} f , f \rangle \sim \rho|(\mathbb{I}-\mathbb{P}_{\rho})f|_{\mathcal{L}^{s}_{\gamma/2}}^{2},
\een
where the norm $|\cdot|_{\mathcal{L}^{s}_{\gamma/2}}$ is defined in \eqref{definition-of-norm-L-epsilon-gamma}. We use $\mathbb{I}$ to denote the identity operator.
Note that $\mathcal{L}^{\rho}$ vanishes with order-1 as $\rho \rightarrow 0$.  There are two key points in deriving \eqref{coercivity-result}.
\begin{itemize}
\item The first one is to reduce quantum case to classical case.
Let us recall the classical (non-quantum) linearized Boltzmann operator with kernel $B$ in \eqref{Boltzmann-kernel-general-and-old},
\ben \label{linearized-Boltzmann-operator-UU}
(\mathcal{L}_{c}f)(v)  \colonequals    \int  B \mu_{*} \mu^{\f{1}{2}}
\mathrm{S}(\mu^{-\f{1}{2}}f)
\mathrm{d}\sigma \mathrm{d}v_{*},
\een
where we use $\mathrm{S}(\cdot)$  defined in \eqref{symmetry-operator}. In this article, the subscript
``$c$'' is referred to ``classical''. When $\rho$ is small, we reduce $\langle \mathcal{L}^{\rho} f , f \rangle$ to
$\langle \rho\mathcal{L}_{c}f , f \rangle$ with some small correction term.
See Lemma  \ref{tau-reduce-to-0} for details.
\item The second one is to give an estimate of $\mathbb{P}_{\rho}-\mathbb{P}_{0}$ where $\mathbb{P}_{0}$ is
the projection operator on the null space ($\mathrm{Null}^{0}$ defined in \eqref{null-based-on-N-limit-0}) of classical operator $\mathcal{L}_{c}$. Roughly speaking, $|\mathbb{P}_{\rho}-\mathbb{P}_{0}| = O(\rho)$,
see Lemma \ref{projection-close} for details.
\end{itemize}
With these two key observations, using the coercivity estimate of $\mathcal{L}_{c}$ in Theorem \ref{classical-linearized-operator-dissipation},
we get Theorem \ref{main-theorem} for the coercivity estimate of $\mathcal{L}^{\rho}$.

{\it Some preliminary formulas.} We collect some preliminary formulas for Boltzmann type integrals in subsection \ref{preliminary}.  We first give two commonly used changes of variable in \eqref{change-v-and-v-star} and \eqref{change-v-and-v-prime}. When estimating integrals, we have to deal with singularity in the kernel $B$. To cancel angular singularity, we give Lemma \ref{cancel-singularity-in-two-type-integral} (based on Taylor expansion) and Lemma  \ref{cancelation-lemma-result}(based on the cancelation lemma in \cite{alexandre2000entropy}) for  integrals that have particular structures. To cancel velocity  singularity, we prepare Lemma \ref{cancel-velocity-regularity}. In order to retain negative exponential weight, we give Lemma \ref{v-square-sum-bounded} and Remark \ref{optimal-constant}.

{\it Upper bound of $\mathcal{L}^{\rho}_{r}$.}  We make several rearrangements to $\langle \mathcal{L}^{\rho}_{r}g, h\rangle$ and then use the formulas in subsection \ref{preliminary} to get $|\langle \mathcal{L}^{\rho}_{r}g, h\rangle| \lesssim \rho|\mu^{\frac{1}{64}}g|_{L^{2}}|\mu^{\frac{1}{64}}h|_{L^{2}}$  in Proposition \ref{ub-linearized-L2}. Note that we succeed in obtaining an upper bound only involving $L^{2}$-norm
with negative exponential weight.
This result is comparable to the classical case, see Lemma 2.15 in \cite{alexandre2012boltzmann}.

{\it Upper bound of $\mathcal{C}^{\rho}$.} Using the upper bound for functionals $\mathcal{N}^{ipl}(\mu^{a}, f)$ and $\mathcal{N}^{ipl}(f, \mu^{a})$ involving classical Boltzmann kernel $B^{ipl}$ defined in \eqref{inverse-power-law-kernel},
we prove in Theorem \ref{two-functional-ub-by-norm} that the two functionals $\mathcal{N}(\mu^{a}, f)$ and $\mathcal{N}(f, \mu^{a})$ are bounded from up by $|f|_{\mathcal{L}^{s}_{\gamma/2}}^{2}$. See \eqref{funcation-N-pm} for the definition of these functionals. We make several rearrangements to $\langle \mathcal{C}^{\rho}g, h\rangle$ such that it can be controlled by the two functionals, which enables us to apply Theorem \ref{two-functional-ub-by-norm} to get $|\langle \mathcal{C}^{\rho} h, f\rangle| \lesssim \rho^{2} |\mu^{\frac{1}{4}}h|_{\mathcal{L}^{s}_{\gamma/2}}|\mu^{\frac{1}{4}}f|_{\mathcal{L}^{s}_{\gamma/2}}$. See Proposition \ref{ub-for-correction-term} for details.

{\it Upper bound of $\langle \tilde{Q}(g, f), f\rangle$.} We prove $\langle \tilde{Q}(g, f), f\rangle \lesssim |\mu^{-\f14} g|_{H^{4}}(1+|\mu^{-\f14} g|_{H^{4}}) |f|_{L^{2}}^{2}$ in Proposition \ref{for-positivity} which is used in Proposition \ref{positivity} to derive non-negativity of solutions to the linear problem \eqref{quantum-Boltzmann-UU-linear}. The key idea in the proof of Proposition \ref{for-positivity} is to put all regularity and weight on $g$ through those preliminary formulas mentioned before.

\subsubsection{Nonlinear operators}
To analyze the nonlinear operators $\Gamma_{2,m}^{\rho}, \Gamma_{2}^{\rho}, \Gamma_{3}^{\rho}$,
we will make use of the classical (non-quantum) Boltzmann operator with kernel $B$ in \eqref{Boltzmann-kernel-general-and-old},
\ben \label{classical-Boltzmann-operator-UU}
Q_{c}(g, h)(v)  \colonequals    \int  B(v- v_{*},\sigma)
 \mathrm{D}(g_{*}^{\prime} h^{\prime}) \mathrm{d}\sigma \mathrm{d}v_{*}.
\een
By some rearrangement, we find
\beno
\Gamma_{2}^{\rho}(g,h) = \rho^{\f{1}{2}} \big(Q_{c} (N_{\rho} g, h) + I^{\rho}(g,h) \big) + \rho^{\frac{3}{2}}\Gamma_{2,r}^{\rho}(g,h), \quad \Gamma_{2,m}^{\rho}(g,h) = Q_{c} (N_{\rho}g, h) + I^{\rho}(g,h),
\eeno
where $I^{\rho}$ is defined in \eqref{definition-I-epsilon} and $\Gamma_{2,r}^{\rho}(g,h)$ is defined in \eqref{Gamma-remaining-into-three-terms}. Here $N_{\rho}= (1-\rho \mu)^{-1} \mu^{\f12}$.
In order to estimate $\Gamma_{2}^{\rho}(g, h)$, we need to consider $Q_{c} (N_{\rho} g, h), I^{\rho}(g,h), \Gamma_{2,r}^{\rho}(g,h)$. We now give some keys points.

{\it Upper bound of $\langle Q_{c} (N_{\rho} g, h), f\rangle$.} Observing \eqref{Boltzmann-kernel-general-ub}, we first get an upper bound of $\langle Q_{c} (g, h), f\rangle$ in Corollary  \ref{Q-up-bound-full-on-first-or-second} by using some known estimates of classical Boltzmann operator with $B^{ipl}$ defined in \eqref{inverse-power-law-kernel}. Then noting that $N_{\rho} \sim \mu^{\f12}$ and using the elementary Lemma \ref{product-take-out}, we get the upper bound of $\langle Q_{c} (N_{\rho} g, h), f\rangle$  in Proposition \ref{Q-pm-Ng-h-f}.

{\it Upper bound of $\langle I^{\rho}(g,h), f\rangle$.} See from \eqref{I-into-main-and-rest}, \eqref{definition-Im-epsilon} and \eqref{definition-Ir-epsilon}  that $\langle I^{\rho}(g,h), f\rangle$ contains differences
$\mathrm{D}(\mu^{\f{1}{2}}_{*})$ and $\mathrm{D}(\mu^{\f{1}{2}})$. To deal with functionals containing such differences, we introduce Lemma \ref{functional-X-g-h} and Remark \ref{mu-to-N-or-M-X-Y}. We then make suitable rearrangements and use some basic tools (such as Cauchy-Schwartz inequality, changes of variable in \eqref{change-v-and-v-star} and \eqref{change-v-and-v-prime}) to bound $\langle I^{\rho}(g,h), f\rangle$ from up by the functionals in Lemma \ref{functional-X-g-h} and Theorem \ref{two-functional-ub-by-norm}. See Proposition \ref{I-pm-rho-g-h-f} for more details.

{\it Upper bounds of $\Gamma_{2,r}^{\rho}, \Gamma_{3}^{\rho}$.} Lemma \ref{optimal-weight-estimate} and Remark \ref{star-prime-or-prime-version} are introduced to deal with various functionals containing $\mathrm{D}^{2}(g_{*})$. For the upper bound of $\Gamma_{2,r}^{\rho}$ and $\Gamma_{3}^{\rho}$,
we make a variety of rearrangements and use some basic tools (such as Cauchy-Schwartz inequality, changes of variable in \eqref{change-v-and-v-star} and \eqref{change-v-and-v-prime}, weight retainment \eqref{mu-weight-result}, the imbedding $H^{2} \hookrightarrow L^{\infty}$, the usual change of variable $v \to v^{\prime}$ or $v_{*} \to v^{\prime}_{*}$ ) to bound $\langle \Gamma_{2,r}^{\rho}(g,h), f\rangle$ and $\langle \Gamma_{3}^{\rho}(g,h,\varrho), f\rangle$ from up by $\langle Q_{c} (g, h), f\rangle, \langle I^{\rho}(g,h), f\rangle$ and the functionals in Lemma \ref{functional-X-g-h}(Remark \ref{mu-to-N-or-M-X-Y}), Lemma \ref{optimal-weight-estimate}(Remark \ref{star-prime-or-prime-version}) and Theorem \ref{two-functional-ub-by-norm}. keep in mind that the ending upper bounds will be used in later energy estimate of the equations  \eqref{linearized-quantum-Boltzmann-eq} and \eqref{linearized-quantum-Boltzmann-eq-linear} and so they must be sharp enough. 
We illustrate this point by looking at $\Gamma_{3}^{\rho}$.
We use $\partial_{\beta} \colonequals  \partial^{\beta}_{v}$ to denote $v$-derivative.
In later energy estimate of $\partial_{\beta}f$,
we will encounter
\beno
\langle \Gamma_{3}^{\rho}(\partial_{\beta_{1}}f,\partial_{\beta_{2}}f,\partial_{\beta_{3}}f),  \partial_{\beta}f\rangle,
\eeno
where $\beta = \beta_{1} + \beta_{2} + \beta_{3}$. The most dangerous term appears when all the derivatives fall on a single  function. There are three cases
\beno
\langle \Gamma_{3}^{\rho}(\partial_{\beta}f, f, f),  \partial_{\beta}f\rangle, \quad
\langle \Gamma_{3}^{\rho}(f,\partial_{\beta}f, f),  \partial_{\beta}f\rangle, \quad
\langle \Gamma_{3}^{\rho}(f, f,\partial_{\beta}f),  \partial_{\beta}f\rangle. \quad
\eeno
To deal with these terms, the upper bound of $\langle \Gamma_{3}^{\rho}(g,h,\varrho), f\rangle$ must only involve  $L^{2}$ or $\mathcal{L}^{s}_{\gamma/2}$ for at least one of $g,h,\varrho$. To achieve such flexibility, we try hard to figure out suitable rearrangements. See Theorem \ref{Gamma-3-ub} for the three final estimates of $\langle \Gamma_{3}^{\rho}(g,h,\varrho), f\rangle$. Applicable estimates of $\langle \Gamma_{2,r}^{\rho}(g,h), f\rangle$ are given in
Propositions \ref{upper-bound-of-Gamma-2-r-1-g-h-f},
\ref{upper-bound-of-Gamma-2-r-2-g-h-f} and \ref{upper-bound-of-Gamma-2-r-3} corresponding to the three terms in \eqref{Gamma-remaining-into-three-terms}.

\subsubsection{Commutator estimates} To implement energy method in weighted Sobolev space, we need to consider commutators between the weight function $W_{l}$ and the above operators.  Such commutators always contain the difference $\mathrm{D}(W_{l})$. To deal with functionals containing such difference, we introduce Lemma \ref{full-integral-Wl-difference-g-h-f} and Remark \ref{still-true-h-prime}. Then estimates of commutators like
$[W_{l},\mathcal{L}^{\rho}], [W_{l},\Gamma_{2}^{\rho}(g, \cdot)], [W_{l},\Gamma_{3}^{\rho}(g, \cdot,\varrho)]$ are derived by using Lemma \ref{full-integral-Wl-difference-g-h-f}(Remark \ref{still-true-h-prime}) and some other known results. With these commutator estimates and the operator estimates in Section \ref{linear}, \ref{bilinear} and \ref{trilinear}, we give weighted inner product estimates first in $L^{2}(\R^{3})$ space and then in $L^{2}(\mathbb{T}^{3} \times \R^{3})$ space. Finally, some useful energy estimates are derived in Theorem \ref{L-energy-estimate}, \ref{Gamma-2-energy-estimate-label} and \ref{Gamma-3-energy-estimate}. See Section \ref{commutator} for details.

\subsubsection{Local well-posedness}
We first prove well-posedness of the equation \eqref{linearized-quantum-Boltzmann-eq-linear} in Proposition \ref{a-priori-estimate-of-linear-equation} under suitable smallness assumption on the given function $g$. In addition, some continuity (see \eqref{solution-also-verifies-estimate}) to initial datum is provided. Based on these results on equation \eqref{linearized-quantum-Boltzmann-eq-linear},  we construct a function sequence through iteration. More concretely, we start with $f^{0} \equiv 0$ and take $(g, f) = (f^{n-1}, f^{n})$ in \eqref{linearized-quantum-Boltzmann-eq-linear} to construct a function sequence $\{f^{n}\}_{n \geq 0}$. Then using the continuity result \eqref{solution-also-verifies-estimate},
the sequence $\{f^{n}\}_{n \geq 0}$ is proved to be a Cauchy sequence generating a local solution to the nonlinear equation \eqref{linearized-quantum-Boltzmann-eq}. See Theorem \ref{local-well-posedness-LBE} and its proof for more details.

\subsubsection{A priori estimate}
The standard macro-micro decomposition method is used in this step. The decomposition reads
$f = \mathbb{P}_{\rho}f + (f - \mathbb{P}_{\rho}f)$. Recall that $\mathbb{P}_{\rho}f$ and $(f - \mathbb{P}_{\rho}f)$ are referred as ``macroscopic'' and ``microscopic'' parts respectively. For the ``macroscopic'' part, we first derive a system of macroscopic equations \eqref{macroscopic-system} and some local conservation laws \eqref{local-conservation-laws}. With these equations and some other elementary estimates, the dissipation on $\mathbb{P}_{\rho}f$ can be derived as in \cite{duan2008cauchy}. See Lemma \ref{estimate-for-highorder-abc} and \ref{estimate-for-highorder-abc-full} for the precise results.
Full dissipation functional $\mathcal{D}_{N}(f)$ in \eqref{definition-energy-and-dissipation} is derived
in Proposition \ref{essential-estimate-of-micro-macro} for the equation $\partial_{t}f + v\cdot \nabla_{x} f + \mathcal{L}^{\rho}f= g$ where  $g$ is a general source term.

Let us see a key point in the proof of Proposition \ref{essential-estimate-of-micro-macro}. The most difficult term is the free streaming term $v \cdot \nabla_{x} f$ when taking $v$-derivative. More precisely, we need to deal with the commutator $[v \cdot \nabla_{x}, \partial^{\alpha}_{\beta}]$. By the condition \eqref{weight-order-condition-1}, in Lemma \ref{transport-commutator} for any $\eta>0$ we get
\ben \label{streaming-term-analysis}
|([v \cdot \nabla_{x}, \partial^{\alpha}_{\beta}]f, W_{2l_{|\alpha|,|\beta|}} \partial^{\alpha}_{\beta}f)|
\leq \eta \|W_{l_{|\alpha|,|\beta|}} \partial^{\alpha}_{\beta}f\|_{L^{2}_{x}\mathcal{L}^{s}_{\gamma/2}}^{2} + \frac{1}{\eta } \sum_{j=1}^{3} |\beta^{j}|^{2} \|W_{l_{|\alpha|+1,|\beta|-1}} \partial^{\alpha+e^{j}}_{\beta-e^{j}}f\|_{L^{2}_{x}\mathcal{L}^{s}_{\gamma/2}}^{2}.
\een
Here $\beta = (\beta^{1},\beta^{2},\beta^{3}), e^{1}=(1,0,0), e^{2}=(0,1,0), e^{3}=(0,0,1)$.
Observe that there is a factor $\rho$ before the dissipation (see \eqref{coercivity-result} or Theorem \ref{main-theorem}) for the corresponding energy. Here we have
\beno
\|W_{l_{|\alpha|,|\beta|}} \partial^{\alpha}_{\beta} f \|_{L^{2}_{x}L^{2}}^{2}  \to \rho \|W_{l_{|\alpha|,|\beta|}} \partial^{\alpha}_{\beta}f\|_{L^{2}_{x}\mathcal{L}^{s}_{\gamma/2}}^{2},
\\ \|W_{l_{|\alpha|+1,|\beta|-1}} \partial^{\alpha+e^{j}}_{\beta-e^{j}}f \|_{L^{2}_{x}L^{2}}^{2} \to \rho \|W_{l_{|\alpha|+1,|\beta|-1}} \partial^{\alpha+e^{j}}_{\beta-e^{j}}f\|_{L^{2}_{x}\mathcal{L}^{s}_{\gamma/2}}^{2}.
\eeno
For this reason,  in \eqref{streaming-term-analysis} we should take $\eta=\delta \rho$ for some sufficiently small $\delta>0$. The resulting latter term $ \delta^{-1} \rho^{-1}\sum_{j=1}^{3} |\beta^{j}|^{2} \|W_{l_{|\alpha|+1,|\beta|-1}} \partial^{\alpha+e^{j}}_{\beta-e^{j}}f\|_{L^{2}_{x}\mathcal{L}^{s}_{\gamma/2}}^{2}$ needs the dissipation in the energy estimate of $M \delta^{-1} \rho^{-2}\sum_{j=1}^{3} \|W_{l_{|\alpha|+1,|\beta|-1}} \partial^{\alpha+e^{j}}_{\beta-e^{j}}f \|_{L^{2}_{x}L^{2}}^{2}$ for some large constant $M$. Such kind of treatment results in the following combination
\ben \label{combination-of-energy}
\sum_{j=0}^{N}K_{j}\rho^{2j-2N}\sum_{|\alpha|\leq N-j, |\beta|=j}\|W_{l_{|\alpha|,|\beta|}} \partial^{\alpha}_{\beta} f\|^{2}_{L^{2}_{x}L^{2}}
\een
for some constants $\{K_{j}\}_{0 \leq j \leq N}$. Note that the power $2j-2N$ of $\rho$ depends on the $v$-derivative order $|\beta| = j$. See the proof of Proposition \ref{essential-estimate-of-micro-macro} for a detailed and rigourous derivation.
A comparison of \eqref{combination-of-energy} and the energy functional $\mathcal{E}_{N}$ in  \eqref{definition-energy-and-dissipation} explains  the factor $\rho^{-2N}$ in Theorem \ref{global-well-posedness}.

{\it {A priori}} estimate of the problem \eqref{linearized-quantum-Boltzmann-eq} is proved in
Theorem \ref{a-priori-estimate-LBE} by applying Proposition \ref{essential-estimate-of-micro-macro} with $g = \Gamma_{2}^{\rho}(f,f) + \Gamma_{3}^{\rho}(f,f,f)$
and using the energy estimates (\eqref{Gamma-2-energy-estimate} for $\Gamma_{2}^{\rho}$ and Theorem \ref{Gamma-3-energy-estimate} for $\Gamma_{3}^{\rho}$) prepared in Section \ref{commutator}.

\subsubsection{Global well-posedness and smallness of parameter $\rho$}
Global well-posedness (Theorem \ref{global-well-posedness}) is proved in subsection \ref{global-proof} by a continuity argument based on Theorem \ref{local-well-posedness-LBE} and Theorem \ref{a-priori-estimate-LBE} for $0<\rho \leq \rho_{*}$. Let us see the places where smallness condition appears.
\begin{enumerate}
\item In Lemma \ref{M-N-mu}, we need
$\rho \leq \f{1}{2}(2 \pi)^{\f32}$ in order to bound from below the denominator $1 - \rho  \mu \geq \f{1}{2}$.
\item In Lemma \ref{projection-close}, we need $\rho \leq \frac{1}{160}$ to estimate the operator difference $\mathbb{P}_{\rho} - \mathbb{P}_{0}$.
\item In Theorem \ref{main-theorem}, we need $\rho \leq \rho_{0}$ to get the coercivity estimate of $\mathcal{L}^{\rho}$.
\item In Proposition \ref{a-priori-estimate-of-linear-equation}, we need $\rho \leq \rho_{1}$ to prove well-posedness of the linear equation \eqref{linearized-quantum-Boltzmann-eq-linear}. More precisely, smallness is used to cancel out the linear term $\mathcal{C}^{\rho}$.
\item In Theorem \ref{local-well-posedness-LBE}, we need $\rho \leq \rho_{2}$ to prove the local well-posedness of the nonlinear equation \eqref{linearized-quantum-Boltzmann-eq}.
\end{enumerate}
As $0<\rho_{2} \leq \rho_{1} \leq \rho_{0} \leq  \frac{1}{160} \leq \f{1}{2}(2 \pi)^{\f32}$, we finally take $\rho_{*} = \rho_{2}$ in Theorem \ref{global-well-posedness}.

\subsection{Notations} \label{notation} In this subsection, we give a list of notations.

\noindent $\bullet$ Given a set $A$,  $\mathrm{1}_A$ is the characteristic function of $A$.

\noindent $\bullet$ Given two operator $T_{1},T_{2}$, their commutator is denoted by $[T_{1},T_{2}]\colonequals T_{1}T_{2}-T_{2}T_{1}$.

\noindent $\bullet$ The notation $a\lesssim b$  means that  there is a universal constant $C$
such that $a\leq Cb$.  The constant $C$ could depend on the kernel parameters $\gamma, s$ and the energy space index $N, l_{0,0}$.

\noindent $\bullet$ If both $a\lesssim b$ and $b \lesssim a$, we write $a\sim b$.

\noindent $\bullet$ We denote $C(\lambda_1,\lambda_2,\cdots, \lambda_n)$ or $C_{\lambda_1,\lambda_2,\cdots, \lambda_n}$  by a constant depending on $\lambda_1,\lambda_2,\cdots, \lambda_n$.

\noindent $\bullet$ The bracket $\langle \cdot\rangle$ is defined by $\langle v \rangle \colonequals   (1+|v|^2)^{\f{1}{2}}$. The weight function  $W_l(v)\colonequals    \langle v\rangle^l $.

\noindent $\bullet$ For $f,g \in L^{2}({\R^3})$,  $\langle f,g\rangle\colonequals    \int_{\R^3}f(v)g(v) \mathrm{d}v$ and $|f|_{L^{2}}^{2}\colonequals   \langle f,f\rangle$.

\noindent $\bullet$ For $f,g \in L^{2}({\TT^3})$,   $\langle f,g\rangle_{x}\colonequals    \int_{\TT^3}f(x)g(x) \mathrm{d}x$ and $|f|_{L^{2}_{x}}^{2}\colonequals   \langle f,f\rangle_{x}$.

\noindent $\bullet$ For $f,g \in L^{2}({\TT^3 \times \R^3})$,   $(f,g)\colonequals    \int_{\TT^3 \times \R^3} f(x,v)g(x,v) \mathrm{d}x\mathrm{d}v$ and $\|f\|_{L^{2}_{x}L^{2}}^{2}\colonequals   (f, f)$.

\noindent $\bullet$ For a multi-index
$\alpha =(\alpha_1,\alpha_2,\alpha_3) \in \mathbb{N}^{3}$, define
$|\alpha|\colonequals   \alpha_1+\alpha_2+\alpha_3$.

\noindent $\bullet$ For  $\alpha, \beta \in \mathbb{N}^{3}$,
denote $\partial^{\alpha}\colonequals   \partial^{\alpha}_{x}, \partial_{\beta}\colonequals   \partial^{\beta}_{v}, \partial^{\alpha}_{\beta}\colonequals   \partial^{\alpha}_{x}\partial^{\beta}_{v}$.

We now introduce some norm.

\noindent $\bullet$ For $n \in \mathbb{N}, l \in \mathbb{R}$ and a function $f(v)$ on $\mathbb{R}^{3}$, define
\ben \label{Sobolev-norm}
|f|_{H^{n}_{l}}^{2} \colonequals    \sum_{|\beta| \leq n} |W_{l} \partial_{\beta}f|_{L^{2}}^{2}, \quad |f|_{L^{2}_{l}}\colonequals   |f|_{H^{0}_{l}}, \quad |f|_{L^{\infty}_{l}} \colonequals     \esssup_{v \in \mathbb{R}^{3}}  |W_{l}(v)f(v)|.
\een
Note that $|f|_{L^{2}}=|f|_{L^{2}_{0}}$.

\noindent $\bullet$ For $n \in \mathbb{N}$ and a function $f(x)$ on $\mathbb{T}^{3}$, define
\beno
|f|_{H^{n}_{x}}^{2} \colonequals    \sum_{|\alpha| \leq n} |\partial^{\alpha}f|_{L^{2}_{x}}^{2}, \quad |f|_{L^{\infty}_{x}} :=  \esssup_{x \in \mathbb{T}^{3}}  |f(x)|.
\eeno
Note that $|f|_{L^{2}_{x}}=|f|_{H^{0}_{x}}$.

\noindent $\bullet$
For $m, n \in \mathbb{N}, l \in \mathbb{R}$ and a function $f(x,v)$ on $\mathbb{T}^{3}\times \mathbb{R}^{3}$, define
\ben \label{not-mix-x-v-norm-energy}
\|f\|_{H^{m}_{x}H^{n}_{l}}^{2} \colonequals     \sum_{|\alpha| \leq m, |\beta| \leq n} \|W_{l}\partial^{\alpha}_{\beta} f\|_{L^{2}_{x}L^{2}}^{2}, \quad \|f\|_{L^{2}_{x}L^{2}_{l}} \colonequals   \|f\|_{H^{0}_{x}H^{0}_{l}}, \quad \|f\|_{H^{m}_{x}H^{n}} \colonequals   \|f\|_{H^{m}_{x}H^{n}_{0}},
\\ \label{mix-x-v-norm-energy}
 \|f\|_{H^{n}_{x,v}}^{2}  \colonequals     \sum_{|\alpha| + |\beta| \leq n} \|\partial^{\alpha}_{\beta} f\|_{L^{2}_{x}L^{2}}^{2}, \quad \|f\|_{L^{\infty}_{x,v}} \colonequals    \esssup_{x \in \mathbb{T}^{3}, v \in \mathbb{R}^{3}}  |f(x,v)|.
\een

\noindent $\bullet$ For $n, l \in \mathbb{R}$ and a function $f(v)$ on $\mathbb{R}^{3}$, define
\ben \label{Sobolev-regularity-norm-real-index}
|f|_{H^{n}}^{2} \colonequals    \int_{\mathbb{R}^{3}} (1+|\xi|^{2})^{n} |\hat{f}(\xi)|^{2} \mathrm{d}\xi, \quad |f|_{H^{n}_{l}}^{2}\colonequals   |W_{l}f|_{H^{n}}^{2}, \een
where $\hat{f}$ is the Fourier transform.

\noindent $\bullet$ For $n \in \mathbb{R}$ and a function $f(v)$ on $\mathbb{R}^{3}$, define
\ben \label{anisotropic-regularity-norm-real-index}
|f|_{A^{n}}^{2}\colonequals    \sum_{l=0}^\infty\sum_{m=-l}^{l} \int_{0}^{\infty}  \big( 1+l(l+1) \big)^{n} (f^{m}_{l}(r))^{2} r^{2} \mathrm{d}r. \een
where 
$l \in \mathbb{N}, m \in \mathbb{Z}, -l\leq m \leq l,  f^{m}_{l}(r) = \int Y^{m}_{l}(\sigma) f(r \sigma) \mathrm{d}\sigma$. Here  $Y_{l}^{m}$ are the real spherical harmonics verifying that
$(-\triangle_{\mathbb{S}^2})Y_{l}^{m}=l(l+1)Y_{l}^{m}$ where $-\triangle_{\mathbb{S}^2}$ is the Laplacian operator on the unit sphere $\mathbb{S}^2$. Note that $\{Y_{l}^{m}\}_{l \geq 0, -l \leq m \leq l}$ is an orthonormal basis of $L^{2}(\mathbb{S}^2)$. Here the notation $A$ refers to ``anisotropic regularity''.

\noindent $\bullet$
For $l \in \mathbb{R}, 0<s<1$, we define
\ben \label{definition-of-norm-L-epsilon-gamma}
|f|_{\mathcal{L}^{s}_{l}}^{2} \colonequals    |W_{l}f|^{2}_{L^{2}_{s}} + |W_{l}f|_{H^{s}}^{2} + |W_{l}f|_{A^{s}}^{2}.
\een
Recalling \eqref{Sobolev-norm}, \eqref{Sobolev-regularity-norm-real-index} and \eqref{anisotropic-regularity-norm-real-index}, the three norms on the right-hand of \eqref{definition-of-norm-L-epsilon-gamma} share a common weight function $W_{s}(r) = (1+|r|^{2})^{\frac{s}{2}}$ for $|v|, |\xi|, (l(l+1))^{\f12}$.

\noindent $\bullet$
For $n \in \mathbb{N}, l \in \mathbb{R}$ and a function $f(v)$ on $\mathbb{R}^{3}$, define
\ben \label{definition-of-L-n}
|f|_{\mathcal{L}^{n,s}_{l}}^{2}\colonequals     \sum_{|\beta| \leq n} |\partial_{\beta} f|_{\mathcal{L}^{s}_{l}}^{2}.
\een
Note that $|f|_{\mathcal{L}^{s}_{l}} = |f|_{\mathcal{L}^{0,s}_{l}}$.

\noindent $\bullet$ For $m, n \in \mathbb{N}, l \in \mathbb{R}$ and a function $f(x,v)$ on $\mathbb{T}^{3}\times \mathbb{R}^{3}$, define
\ben \label{not-mix-x-v-norm-dissipation}
\|f\|_{L^{2}_{x}\mathcal{L}^{s}_{l}}^{2} \colonequals   \int |f(x,\cdot)|_{\mathcal{L}^{s}_{l}}^{2} \mathrm{d}x, \quad \|g\|_{H^{n}_{x}\mathcal{L}^{m,s}_{l}}^{2} \colonequals     \sum_{|\alpha| \leq n, |\beta| \leq m} \|\partial^{\alpha}_{\beta} g\|_{L^{2}_{x}\mathcal{L}^{s}_{l}}^{2}.
\een

\subsection{Plan of the article} \label{plan}
Section \ref{linear} contains estimates of linear operators, including coercivity estimate of $\mathcal{L}^{\rho}$ and upper bounds of $\mathcal{L}^{\rho}_{r}, \mathcal{C}^{\rho}, \tilde{Q}(g,\cdot)$. Section \ref{bilinear} and \ref{trilinear} are devoted to upper bounds of the bilinear operator $\Gamma_{2}^{\rho}$ and the trilinear operator $\Gamma_{3}^{\rho}$ respectively. In Section \ref{commutator},
various functionals that will appear in
later energy method are estimated after necessary commutator estimates. In Section \ref{local}, we derive local well-posedness. In Section \ref{global}, we first prove a priori estimate and then establish global well-posedness. Section \ref{appendix} is an appendix in which we put some elementary proof for the sake of completeness.

\section{Linear operator estimate} \label{linear}
In the rest of the article, in the various functional estimates, the involved functions $g, h, \varrho, f$ are assumed to be functions on $\R^{3}$ or  $\mathbb{T}^{3} \times \R^{3}$ such that the corresponding norms of them
are well-defined. For simplicity, we use the notation $\mathrm{d}V \colonequals \mathrm{d}\sigma \mathrm{d}v_{*} \mathrm{d}v$.

\subsection{Coercivity estimate of $\mathcal{L}^{\rho}$} Recall \eqref{equilibrium-rho}, \eqref{linearized-quantum-Boltzmann-operator-UU} and \eqref{linearized-Boltzmann-operator-UU}.
The operator $\mathcal{L}^{\rho}$ will vanish as $\rho \rightarrow 0$. However, formally it is easy to see $\rho^{-1}\mathcal{L}^{\rho} \rightarrow \mathcal{L}_{c}$ as $\rho \rightarrow 0$.
We define
\ben \label{equilibrium-rho-not-vanish}
M_{\rho} \colonequals    \frac{\mu}{1 - \rho \mu}, \quad N_{\rho} \colonequals    \frac{\mu^{\f{1}{2}}}{1 - \rho \mu}.
\een
Recalling \eqref{equilibrium-rho}, we can see that
\ben \label{M--with-mathcal-M}
\mathcal{M}_{\rho} = \rho M_{\rho}, \quad \mathcal{N}_{\rho} = \rho^{\f{1}{2}} N_{\rho}.
\een

When $\rho$ is small and close to 0, it is obvious that $M_{\rho} \sim \mu, N_{\rho} \sim \mu^{\f{1}{2}}$. More precisely, we have
\begin{lem} \label{M-N-mu}
If $0 \leq \rho \leq \f{1}{2}(2 \pi)^{\frac{3}{2}}$, then
\ben \label{M-rho-mu-N}
\mu \leq M_{\rho} \leq 2 \mu, \quad  \mu^{\f{1}{2}} \leq N_{\rho} \leq 2 \mu^{\f{1}{2}}.
\een
As a direct result, there holds
\ben \label{K-2-mu}
 \mu \mu_{*} \leq  N_{\rho} (N_{\rho})_{*} (N_{\rho})^{\prime} (N_{\rho})^{\prime}_{*} \leq 2^{4} \mu \mu_{*}.
\een
\end{lem}
\begin{proof} If $0 \leq \rho \leq \f{1}{2}(2 \pi)^{\frac{3}{2}}$, then $0 \leq \rho \mu \leq \f{1}{2}$ and so $\f{1}{2} \leq 1 - \rho \mu \leq 1$, which gives \eqref{M-rho-mu-N} by recalling the definition of  $M_{\rho}$ and  $N_{\rho}$ in \eqref{equilibrium-rho-not-vanish}. As a direct result, we get \eqref{K-2-mu} since $\mu\mu_{*}=\mu^{\prime}\mu^{\prime}_{*}$.
\end{proof}

{\it{In the rest of the article, we always assume $0 \leq \rho \leq \f{1}{2}(2 \pi)^{\frac{3}{2}}$ which enables us to use the results in Lemma \ref{M-N-mu}.}}
 Other smallness condition on $\rho$ will be specified as we go further.

We relate $\mathcal{L}^{\rho}$ and $\mathcal{L}_{c}$ in the following lemma.
\begin{lem}\label{tau-reduce-to-0} It holds that
\beno
 \rho (\f{1}{2}\langle \mathcal{L}_{c}f, f \rangle - \rho^{2}\langle \mathcal{L}_{c}\mu f, \mu f \rangle) \leq \langle \mathcal{L}^{\rho} f, f \rangle \leq 2^{5} \rho (\langle \mathcal{L}_{c}f, f \rangle +  \rho^{2}\langle \mathcal{L}_{c}\mu f, \mu f \rangle).
\eeno
\end{lem}
\begin{proof} Recalling \eqref{self-joint-operator} and \eqref{M--with-mathcal-M}, we have
\beno
\langle \mathcal{L}^{\rho}f, f \rangle = \frac{\rho}{4} \int B N_{\rho} (N_{\rho})_{*} (N_{\rho})^{\prime} (N_{\rho})^{\prime}_{*}  \mathrm{S}^{2}( N_{\rho}^{-1} f )  \mathrm{d}V.
\eeno
Thanks to \eqref{K-2-mu}, we have
\ben \label{reduce-to-mu-mu-star}
\frac{\rho}{4}   \int B  \mu \mu_{*}  \mathrm{S}^{2}( N_{\rho}^{-1} f )  \mathrm{d}V \leq \langle \mathcal{L}^{\rho}f, f \rangle \leq 2^{4} \times \frac{\rho}{4}   \int B \mu \mu_{*}  \mathrm{S}^{2}( N_{\rho}^{-1} f )  \mathrm{d}V.
\een
Recalling \eqref{equilibrium-rho-not-vanish}, we have $N_{\rho}^{-1} =  \mu^{-\f{1}{2}}  - \rho \mu^{\f{1}{2}}$ and thus
\beno
\f{1}{2}\mathrm{S}^{2}( \mu^{-\f{1}{2}} f ) - \rho^{2} \mathrm{S}^{2}( \mu^{\f{1}{2}} f )
\leq
 \mathrm{S}^{2}( N_{\rho}^{-1} f )   =\big(\mathrm{S}( \mu^{-\f{1}{2}} f ) - \rho\mathrm{S}( \mu^{\f{1}{2}} f ) \big)^{2}
\leq 2 \mathrm{S}^{2}( \mu^{-\f{1}{2}} f ) + 2 \rho^{2} \mathrm{S}^{2}( \mu^{\f{1}{2}} f ) .
\eeno
Plugging which into \eqref{reduce-to-mu-mu-star}, we have
\beno
\rho (\f{1}{2}\mathcal{I}_{1} - \rho^{2}\mathcal{I}_{2}) \leq \langle \mathcal{L}^{\rho}f, f \rangle \leq 2^{5}\rho (\mathcal{I}_{1} +  \rho^{2}\mathcal{I}_{2}),
\eeno
where
\beno
\mathcal{I}_{1} \colonequals    \frac{1}{4} \int B \mu \mu_{*} \mathrm{S}^{2}( \mu^{-\f{1}{2}} f )  \mathrm{d}V, \quad
\mathcal{I}_{2} \colonequals    \frac{1}{4} \int B \mu \mu_{*} \mathrm{S}^{2}( \mu^{\f{1}{2}} f )   \mathrm{d}V.
\eeno
Similar to \eqref{self-joint-operator}, there holds
\ben \label{self-joint-operator-classical-case}
\langle \mathcal{L}_{c}g, h \rangle = \frac{1}{4} \int B \mu\mu_{*} \mathrm{S}( \mu^{-\f{1}{2}} g ) \mathrm{S}(\mu^{-\f{1}{2}} h) \mathrm{d}V = \langle g,  \mathcal{L}_{c}h \rangle,
\een
By \eqref{self-joint-operator-classical-case}, we observe $\mathcal{I}_{1}=\langle \mathcal{L}_{c}f, f \rangle, \mathcal{I}_{2}=\langle \mathcal{L}_{c}\mu f, \mu f \rangle$ and get the desired result.
\end{proof}

Recalling \eqref{null-space-of-L}, since $\mathcal{N}_{\rho} = \rho^{\f{1}{2}} N_{\rho}$, for $\rho>0$, we can also write
\ben \label{null-based-on-N-not-vanish}
\mathrm{Null}^{\rho} =  \mathrm{span} \{ N_{\rho}, N_{\rho}v_{1}, N_{\rho}v_{2}, N_{\rho}v_{3}, N_{\rho}|v|^{2} \}.
\een
Note that based on \eqref{null-based-on-N-not-vanish}, $\mathrm{Null}^{0}$ is well-defined by
\ben \label{null-based-on-N-limit-0}
\mathrm{Null}^{0} \colonequals     \mathrm{span} \{ \mu^{\f{1}{2}}, \mu^{\f{1}{2}}v_{1}, \mu^{\f{1}{2}}v_{2}, \mu^{\f{1}{2}}v_{3}, \mu^{\f{1}{2}}|v|^{2} \}.
\een
Observe that $\mathrm{Null}^{0}$ is the null space of $\mathcal{L}_{c}$.

We construct an orthogonal basis for $\mathrm{Null}^{\rho}$ for $\rho \geq 0$ as follows
\ben \label{definition-of-di}
\{ d^{\rho}_{i} \}_{1 \leq i \leq 5} \colonequals    \{ N_{\rho}, N_{\rho} v_{1}, N_{\rho} v_{2}, N_{\rho} v_{3}, N_{\rho}|v|^{2} -
\langle N_{\rho}|v|^{2} , N_{\rho}\rangle|N_{\rho}|^{-2}_{L^{2}}N_{\rho} \}.
\een
Note that $\langle N_{\rho}|v|^{2} , N_{\rho}\rangle|N_{\rho}|^{-2}_{L^{2}}N_{\rho}$ is the projection of $N_{\rho}|v|^{2}$ on $N_{\rho}$.
We denote the coefficient $C^{\rho}_{5,1} \colonequals   \langle N_{\rho}|v|^{2} , N_{\rho}\rangle|N_{\rho}|^{-2}_{L^{2}}$ which depends only on $\rho$.
Then $ d^{\rho}_{5} = d^{\rho}_{1}|v|^{2} - C^{\rho}_{5,1} d^{\rho}_{1}$.
By normalizing $\{ d^{\rho}_{i} \}_{1 \leq i \leq 5}$, an orthonormal basis of $\mathrm{Null}^{\rho}$ can be obtained as
\ben \label{definition-of-e-pm-rho}
\{ e^{\rho}_{i} \}_{1 \leq i \leq 5} \colonequals    \{ \frac{d^{\rho}_{i}}{|d^{\rho}_{i}|_{L^{2}}} \}_{1 \leq i \leq 5}.
\een
With this orthonormal basis, the projection operator $\mathbb{P}_{\rho}$ on the null space  $\mathrm{Null}^{\rho}$ is defined by
\ben \label{definition-projection-operator}
\mathbb{P}_{\rho}f \colonequals    \sum_{i=1}^{5} \langle f, e^{\rho}_{i}\rangle e^{\rho}_{i}.
\een

Let us see $\mathbb{P}_{\rho}$ more clearly.
Let us define
\ben \label{definition-of-m-i}
m_{\rho,0} \colonequals    |d^{\rho}_{1}|_{L^{2}}, \quad  m_{\rho,1} \colonequals    |d^{\rho}_{2}|_{L^{2}} = |d^{\rho}_{3}|_{L^{2}} =|d^{\rho}_{4}|_{L^{2}} , \quad m_{\rho,2} \colonequals    |d^{\rho}_{5}|_{L^{2}}.
\een
For simplicity, we set $m_{i}=m_{\rho,i}$ for $i=0,1,2$ and $N=N_{\rho}$. Then by direct derivation and rearrangement, we have
\ben \nonumber
\mathbb{P}_{\rho}f &=& \langle f, \frac{N}{m_{0}} \rangle \frac{N}{m_{0}} + \sum_{i=1}^{3} \langle f, \frac{N v_{i}}{m_{1}} \rangle \frac{N v_{i}}{m_{1}} + \langle f, \frac{N(|v|^{2}-C^{\rho}_{5,1})}{m_{2}} \rangle \frac{N(|v|^{2}-C^{\rho}_{5,1})}{m_{2}}
\\ \label{linear-combination-of-basis} &=& (a^{f}_{\rho} + b^{f}_{\rho} \cdot v + c^{f}_{\rho}|v|^{2})N,
\een
where
\beno
a^{f}_{\rho} \colonequals    \langle f, (\frac{1}{m_{0}^{2}} + \frac{(C^{\rho}_{5,1})^{2}}{m_{2}^{2}})N -  \frac{C^{\rho}_{5,1}}{m_{2}^{2}}N|v|^{2} \rangle, \quad
b^{f}_{\rho} \colonequals    \langle f, \frac{N v}{m_{1}^{2}} \rangle,
\quad
c^{f}_{\rho} \colonequals    \langle f, \frac{1}{m_{2}^{2}} N|v|^{2} - \frac{C^{\rho}_{5,1}}{m_{2}^{2}} N \rangle.
\eeno
Note that $b^{f}_{\rho}$ is a vector of length 3.  Let us define
\ben \label{defintion-of-l-i}
l_{\rho,1}\colonequals   \frac{1}{m_{0}^{2}} + \frac{(C^{\rho}_{5,1})^{2}}{m_{2}^{2}}, \quad
l_{\rho,2}\colonequals   \frac{C^{\rho}_{5,1}}{m_{2}^{2}}, \quad
l_{\rho,3}\colonequals   \frac{1}{m_{1}^{2}}, \quad
l_{\rho,4}\colonequals   \frac{1}{m_{2}^{2}}.
\een
For simplicity, let $l_{i} = l_{\rho, i}$.
Then there holds
\ben \label{explicit-defintion-of-abc}
a^{f}_{\rho} =  \langle f, l_{1}N -  l_{2}N|v|^{2} \rangle, \quad
b^{f}_{\rho} = \langle f, l_{3}N v \rangle, \quad
c^{f}_{\rho} = \langle f, l_{4} N|v|^{2} - l_{2} N \rangle.
\een


The next lemma shows that $\mathbb{P}_{\rho} - \mathbb{P}_{0}$ is of order $O(\rho)$ when $\rho$ is small.
\begin{lem} \label{projection-close} Let $ m \geq 0, l \in \mathbb{R}, 0 \leq \rho \leq \frac{1}{160}$.
There holds
\beno
|\mathbb{P}_{\rho}f - \mathbb{P}_{0}f|_{H^{m}_{l}} \leq C_{m,l}  \rho |\mu^{\frac{1}{4}}f|_{L^{2}},
\eeno
where $C_{m,l}$ is defined in \eqref{constant-C-s-l}.
\end{lem}
The proof of Lemma \ref{projection-close} will be given in the appendix.

We now give the coercivity estimate of $\mathcal{L}_{c}$ based on some
known result on the classical linearized Boltzmann operator with inverse power law potential. For inverse power law potential, it suffices to consider the following Boltzmann kernel
\ben \label{inverse-power-law-kernel}
B^{ipl}(v-v_{*},\sigma) \colonequals    |v-v_{*}|^{\gamma} \sin^{-2-2s}\frac{\theta}{2} \mathrm{1}_{0 \leq \theta \leq \pi/2}.
\een
Note that the superscript ``$ipl$'' is short for ``inverse power law potential''.
Let $\mathcal{L}^{ipl}_{c}$ be the associated classical linearized Boltzmann operator, i.e., $\mathcal{L}^{ipl}_{c}$ is defined through \eqref{linearized-Boltzmann-operator-UU} by replacing
$B$ with $B^{ipl}$. By \cite{he2022asymptotic}, it turns out that
\ben \label{classical-coercivity-result-rough-1}
\langle \mathcal{L}^{ipl}_{c} f , f \rangle \sim |(\mathbb{I}-\mathbb{P}_{0})f|_{\mathcal{L}^{s}_{\gamma/2}}^{2},
\een
where the norm $|\cdot|_{\mathcal{L}^{s}_{\gamma/2}}$ is defined in \eqref{definition-of-norm-L-epsilon-gamma}. We remind that the norm $|\cdot|_{\mathcal{L}^{s}_{\gamma/2}}$ is equivalent to $\mathcal{N}^{s,\gamma}$ in \cite{gressman2011global} and $|||\cdot|||$ in \cite{alexandre2012boltzmann}.

Since $\sin \frac{\theta}{2} \leq \cos \frac{\theta}{2}$ for $0 \leq \theta \leq \frac{\pi}{2}$, recalling \eqref{inverse-power-law-kernel} and \eqref{Boltzmann-kernel-general-and-old},
there holds
\ben \label{Boltzmann-kernel-general-ub}
C B^{ipl}(v-v_{*},\sigma) \leq B(v- v_{*},\sigma) \leq 4 C B^{ipl}(v-v_{*},\sigma) .
\een
Recalling \eqref{self-joint-operator-classical-case} and using \eqref{Boltzmann-kernel-general-ub}, we get
\beno
\langle \mathcal{L}_{c} f , f \rangle = \langle \mathcal{L}_{c} (\mathbb{I}-\mathbb{P}_{0})f , (\mathbb{I}-\mathbb{P}_{0})f \rangle \sim \langle \mathcal{L}^{ipl}_{c} (\mathbb{I}-\mathbb{P}_{0})f , (\mathbb{I}-\mathbb{P}_{0})f \rangle = \langle \mathcal{L}^{ipl}_{c} f , f \rangle.
\eeno
By using \eqref{classical-coercivity-result-rough-1},  we get the following theorem.
\begin{thm}\label{classical-linearized-operator-dissipation} It holds that
\beno
\langle \mathcal{L}_{c} f , f \rangle \sim |(\mathbb{I}-\mathbb{P}_{0})f|_{\mathcal{L}^{s}_{\gamma/2}}^{2}.
\eeno
For later reference, let  $\lambda_{*}, C_{*}$ be the two optimal universal constants such that
\ben \label{classical-coercivity-result}
\lambda_{*} |(\mathbb{I}-\mathbb{P}_{0})f|_{\mathcal{L}^{s}_{\gamma/2}}^{2} \leq \langle \mathcal{L}_{c} f , f \rangle \leq C_{*}|(\mathbb{I}-\mathbb{P}_{0})f|_{\mathcal{L}^{s}_{\gamma/2}}^{2}.
\een
\end{thm}

We now prove the following the following coercivity estimate of $\mathcal{L}^{\rho}$  by using Theorem \ref{classical-linearized-operator-dissipation}, Lemma \ref{tau-reduce-to-0} and Lemma \ref{projection-close}.
\begin{thm}\label{main-theorem}
There are three universal constants $\rho_{0}, \lambda_{0}, C_{0}>0$ such that for any $\rho$ verifying $0 \leq \rho \leq \rho_{0}$, it holds that
\beno
\lambda_{0}\rho |(\mathbb{I}-\mathbb{P}_{\rho})f|_{\mathcal{L}^{s}_{\gamma/2}}^{2} \leq \langle \mathcal{L}^{\rho} f , f \rangle \leq C_{0}\rho|(\mathbb{I}-\mathbb{P}_{\rho})f|_{\mathcal{L}^{s}_{\gamma/2}}^{2}.
\eeno
The constants $\rho_{0}, \lambda_{0}, C_{0}>0$ are explicitly given in \eqref{constant-rho0-lambda0-C0}.
\end{thm}
\begin{proof} Since $\langle \mathcal{L}^{\rho} f, f \rangle =
\langle \mathcal{L}^{\rho} (\mathbb{I}-\mathbb{P}_{\rho})f, (\mathbb{I}-\mathbb{P}_{\rho})f \rangle$, it suffices to prove
$\langle \mathcal{L}^{\rho} f, f \rangle \sim \rho |f|_{\mathcal{L}^{s}_{\gamma/2}}^{2}$ for $f$ verifying $\mathbb{P}_{\rho}f=0$. From now on, we assume $\mathbb{P}_{\rho}f=0$.
By Lemma \ref{tau-reduce-to-0}, for $0 \leq \rho \leq \f{1}{2}(2 \pi)^{\frac{3}{2}}$, we have
\ben \label{epsilon-reduce-to-0}
\rho(\f{1}{2}\langle \mathcal{L}_{c}f, f \rangle - \rho^{2}\langle \mathcal{L}_{c}\mu f, \mu f \rangle) \leq \langle \mathcal{L}^{\rho}f, f \rangle \leq 2^{5}\rho (\langle \mathcal{L}_{c}f, f \rangle +  \rho^{2}\langle \mathcal{L}_{c}\mu f, \mu f \rangle).
\een
By \eqref{classical-coercivity-result} in Theorem \ref{classical-linearized-operator-dissipation}, we have
\ben \label{key-estimate-tau=0}
\lambda_{*} |(\mathbb{I}-\mathbb{P}_{0})f|_{\mathcal{L}^{s}_{\gamma/2}}^{2} \leq \langle \mathcal{L}_{c} f, f \rangle \leq C_{*} |(\mathbb{I}-\mathbb{P}_{0})f|_{\mathcal{L}^{s}_{\gamma/2}}^{2} \leq C_{*} C_{1}|f|_{\mathcal{L}^{s}_{\gamma/2}}^{2}.
\\
\label{L-tau-0-muf-ub}
 \langle \mathcal{L}_{c} \mu f , \mu f \rangle \leq C_{*} |(\mathbb{I}-\mathbb{P}_{0}) \mu f|_{\mathcal{L}^{s}_{\gamma/2}}^{2}   \leq C_{*} C_{2}|f|_{\mathcal{L}^{s}_{\gamma/2}}^{2}.
\een
Here $C_{1}$ is the optimal constant such that $|(\mathbb{I}-\mathbb{P}_{0})f|_{\mathcal{L}^{s}_{\gamma/2}}^{2} \leq C_{1}|f|_{\mathcal{L}^{s}_{\gamma/2}}^{2}$ for any $f$
and $C_{2}$ is the optimal constant such $|(\mathbb{I}-\mathbb{P}_{0}) \mu f|_{\mathcal{L}^{s}_{\gamma/2}}^{2} \leq C_{2}|f|_{\mathcal{L}^{s}_{\gamma/2}}^{2}$ for any $f$. The existence of $C_{2}$ is ensured by \eqref{taking-out-with-a-weight}.
We remark that $C_{1},C_{2}$ are universal constants independent of
$\rho, s, \gamma$.

Then plugging \eqref{key-estimate-tau=0} and \eqref{L-tau-0-muf-ub} into \eqref{epsilon-reduce-to-0}, we first have
\ben \label{ub-direction}
 \langle \mathcal{L}^{\rho} f, f \rangle \leq 2^{5} C_{*}  (C_{1}+C_{2}\rho^{2}) \rho |f|_{\mathcal{L}^{s}_{\gamma/2}}^{2}.
\een

Since $\mathbb{P}_{\rho}f=0$, then by Lemma \ref{projection-close} if $ \rho \leq \frac{1}{160}$, we have
\beno
|\mathbb{P}_{0}f|_{\mathcal{L}^{s}_{\gamma/2}} \leq C_{3} |\mathbb{P}_{0}f|_{H^{s}_{s+\gamma/2}} \leq C_{3} |\mathbb{P}_{0}f|_{H^{1}_{1}} \leq C_{3}C_{1,1} \rho |\mu^{\f14}f|_{L^{2}} \leq C_{3}C_{1,1} |\mu^{\f14}|_{L^{\infty}_{3/2}} \rho |f|_{\mathcal{L}^{s}_{\gamma/2}}.
\eeno
Here $C_{3}$ is the optimal constant such that $|f|_{\mathcal{L}^{s}_{\gamma/2}} \leq C_{3}|f|_{H^{s}_{s+\gamma/2}}$ for any $f$. Note that $|f|_{\mathcal{L}^{s}_{\gamma/2}} \geq |f|_{L^{2}_{-3/2}}$ is used.
Let $C_{4}\colonequals    C_{3}C_{1,1} |\mu^{\f14}|_{L^{\infty}_{3/2}}$. If  $C^{2}_{4}\rho^{2} \leq \frac{1}{4}$, then
\ben \label{I-P0-is-dominated-by-I}
|(\mathbb{I}-\mathbb{P}_{0})f|_{\mathcal{L}^{s}_{\gamma/2}}^{2} \geq \f{1}{2} |f|_{\mathcal{L}^{s}_{\gamma/2}}^{2} - |\mathbb{P}_{0}f|_{\mathcal{L}^{s}_{\gamma/2}}^{2} \geq \f{1}{2} |f|_{\mathcal{L}^{s}_{\gamma/2}}^{2} - C^{2}_{4}\rho^{2} |f|_{\mathcal{L}^{s}_{\gamma/2}}^{2} \geq \frac{1}{4} |f|_{\mathcal{L}^{s}_{\gamma/2}}^{2}.
\een
Plugging \eqref{key-estimate-tau=0}, \eqref{L-tau-0-muf-ub} and \eqref{I-P0-is-dominated-by-I} into \eqref{epsilon-reduce-to-0}, if $C_{*}C_{2}\rho^{2} \leq \frac{1}{16}\lambda_{*}$, we have
\ben \nonumber
 \langle \mathcal{L}^{\rho} f, f \rangle &\geq& \rho(\f{1}{2}\lambda_{*} |(\mathbb{I}-\mathbb{P}_{0})f|_{\mathcal{L}^{s}_{\gamma/2}}^{2} - C_{*}C_{2}\rho^{2} |f|_{\mathcal{L}^{s}_{\gamma/2}}^{2})
 \\ \label{lower-bound-L-tau} &\geq& \rho(\frac{1}{8}\lambda_{*} - C_{*}C_{2}\rho^{2}) |f|_{\mathcal{L}^{s}_{\gamma/2}}^{2} \geq \frac{\lambda_{*}}{16}\rho|f|_{\mathcal{L}^{s}_{\gamma/2}}^{2}.
\een
Define
\ben \label{constant-rho0-lambda0-C0}
\rho_{0} \colonequals    \min\{\frac{1}{160}, \frac{1}{2C_{4}}, \sqrt{\frac{\lambda_{*}}{16C_{*}C_{2}}}\}, \quad \lambda_{0} \colonequals    \frac{1}{16}\lambda_{*}, \quad C_{0}\colonequals   2^{5} C_{*}  (C_{1}+\frac{1}{4}(2 \pi)^{3}C_{2}).
\een
Patching together \eqref{ub-direction} and \eqref{lower-bound-L-tau}, we finish the proof.
\end{proof}

For the upper bound estimate in Theorem \ref{main-theorem} we only need $\rho \leq \f{1}{2}(2 \pi)^{\frac{3}{2}}$. Indeed,
based on the proof, we only need Lemma \ref{tau-reduce-to-0} and Theorem \ref{classical-linearized-operator-dissipation} to get \eqref{ub-direction}.  With the constant  $C_{0}$ defined in \eqref{constant-rho0-lambda0-C0}, we have
\begin{rmk}\label{up-no-small-on-rho} For $0 \leq \rho \leq \f{1}{2}(2 \pi)^{\frac{3}{2}}$, we have
$\langle \mathcal{L}^{\rho} f , f \rangle \leq C_{0}\rho|(\mathbb{I}-\mathbb{P}_{\rho})f|_{\mathcal{L}^{s}_{\gamma/2}}^{2}$.
\end{rmk}

\subsection{Some preliminary formulas} \label{preliminary}
In this subsection,
we recall some useful formulas for the computation of integrals involving $B=B(v-v_{*}, \sigma)$ defined in \eqref{Boltzmann-kernel-general-and-old}.
It is obvious that the change of variable $(v, v_{*}, \sigma) \rightarrow (v_{*}, v, -\sigma)$ has unit Jacobian and thus
\ben \label{change-v-and-v-star}
\int B(v-v_{*}, \sigma) f(v,v_{*},v^{\prime},v^{\prime}_{*}) \mathrm{d}V   =
\int B(v-v_{*}, \sigma) f(v_{*},v,v^{\prime}_{*},v^{\prime}) \mathrm{d}V  ,
\een
where $f$ is a general function such that the integral exists. Thanks to the symmetry of elastic collision formula \eqref{v-prime-v-prime-star}, the change of variable $(v, v_{*}, \sigma) \rightarrow (v^{\prime}, v^{\prime}_{*}, \frac{v-v_{*}}{|v-v_{*}|})$ has unit Jacobian and thus
\ben \label{change-v-and-v-prime}
\int B(v-v_{*}, \sigma) f(v,v_{*},v^{\prime},v^{\prime}_{*}) \mathrm{d}V   =
\int B(v-v_{*}, \sigma) f(v^{\prime},v^{\prime}_{*},v,v_{*}) \mathrm{d}V  .
\een

From now on, we will frequently use the notation \eqref{shorthand-D}.
By \eqref{shorthand-D} and the shorthand $f=f(v), f^{\prime} = f(v^{\prime}), f_{*}=f(v_{*}), f^{\prime}_{*} = f(v^{\prime}_{*})$, it is easy to see
\beno
\mathrm{D}(f) = f(v) - f(v^{\prime}) = - \mathrm{D}(f^{\prime}),
\quad
\mathrm{D}(f_{*}) = f(v_{*}) - f(v^{\prime}_{*}) = - \mathrm{D}(f^{\prime}_{*}),
\\
 \mathrm{D}^{2}(f) = \mathrm{D}^{2}(f^{\prime}), \quad \mathrm{D}^{2}(f_{*}) = \mathrm{D}^{2}(f^{\prime}_{*}).
\eeno
Similar to \eqref{shorthand-D}, we introduce
\ben \label{shorthand-A}
\mathrm{A}(f(v,v_{*},v^{\prime},v^{\prime}_{*})) \colonequals f(v,v_{*},v^{\prime},v^{\prime}_{*}) + f(v^{\prime},v^{\prime}_{*},v,v_{*}).
\een
The term $\mathrm{A}$ is interpreted as ``addition'' before and after collision.

In upper bound estimate, we frequently encounter quantities like $\int B(v-v_{*},\sigma) g_{*} h \mathrm{D}(f)  \mathrm{d}V  $ and $\int B(v-v_{*},\sigma) g_{*} h^{\prime} \mathrm{D}(f)  \mathrm{d}V  $. Thanks to \eqref{Boltzmann-kernel-general-ub},
it suffices to consider the kernel $B^{ipl}$ defined
\eqref{inverse-power-law-kernel} for upper bound estimates involving $B$.
To cancel the angular singularity of $\sin^{-2-2s}\frac{\theta}{2}$ near $\theta=0$, one usually relies on the order-2 factor $\sin^{2}\frac{\theta}{2}$ and the following fact for $0<s<1$,
\ben \label{cancel-angular-singularity}
\int \sin^{-2s}\frac{\theta}{2} \mathrm{d}\sigma \lesssim \frac{1}{1-s}, \quad \int B(v-v_{*},\sigma) \sin^{2}\frac{\theta}{2} \mathrm{d}\sigma \lesssim \frac{|v-v_{*}|^{\gamma}}{1-s}.
\een
Thanks to the symmetry of $\sigma$-integral,
the factor $\sin^{2}\frac{\theta}{2}$ will appear if we appropriately apply Taylor expansion to $\mathrm{D}(f)$. By Taylor expansion, we have the following two candidates.
\ben \label{Taylor-at-v}
-\mathrm{D}(f)=\nabla f (v) \cdot (v^{\prime}-v) + \int_{0}^{1}(1-\kappa) \nabla^{2} f (v(\kappa)): (v^{\prime}-v)  \otimes (v^{\prime}-v)  \mathrm{d}\kappa.
\\ \label{Taylor-at-v-prime}
-\mathrm{D}(f)=   \nabla f  (v^{\prime}) \cdot  (v^{\prime}-v)  -  \int_{0}^{1} \kappa  \nabla^{2} f (v(\kappa)): (v^{\prime}-v)  \otimes (v^{\prime}-v)  \mathrm{d}\kappa.
\een
Here $v(\kappa)\colonequals    v + \kappa(v^{\prime}-v)$. Since $|v^{\prime}-v|=|v-v_{*}|\sin\frac{\theta}{2}$, the second order contains
$\sin^{2}\frac{\theta}{2}$. So it remains to deal with the first order term.
Fortunately we have the following two identities,
\ben \label{theta-squre-out}
\int B(v-v_{*}, \sigma) (v^{\prime} - v)  \mathrm{d}\sigma = \int B(v-v_{*}, \sigma)\sin^{2} \frac{\theta}{2}  (v_{*} - v) \mathrm{d}\sigma.
\\ \label{for-fix-v-star-vanish}
\int B(v-v_{*}, \sigma) h^{\prime} (v^{\prime} - v) \mathrm{d}v  \mathrm{d}\sigma = 0.
\een
We remark that \eqref{for-fix-v-star-vanish} holds for any fixed $v_{*}$.

Note that the right-hand side of \eqref{theta-squre-out} contains $\sin^{2}\frac{\theta}{2}$ and so we can
 use \eqref{Taylor-at-v} and \eqref{theta-squre-out} to cancel the angular singularity  in $\int B g_{*} h \mathrm{D}(f)  \mathrm{d}V  $.
We can use \eqref{Taylor-at-v-prime} and \eqref{for-fix-v-star-vanish} to cancel the angular singularity in $\int B g_{*} h^{\prime} \mathrm{D}(f)  \mathrm{d}V  $ since the first order term vanishes by \eqref{for-fix-v-star-vanish}. As a result, we have
\begin{lem}\label{cancel-singularity-in-two-type-integral} The following two estimates are valid.
\ben \label{type-1-cancel}
 |\int B g_{*} h \mathrm{D}(f)  \mathrm{d}V  | &\lesssim& \int |v-v_{*}|^{\gamma+1}|g_{*}h \nabla f|
 \mathrm{d}v_{*} \mathrm{d}v
 \\ \nonumber &&+ \int |v-v_{*}|^{\gamma+2} \mathrm{1}_{0 \leq \theta \leq \pi/2}\sin^{-2s}\frac{\theta}{2}|g_{*}h \nabla^{2} f(v(\kappa))| \mathrm{d}\kappa \mathrm{d}V.
 \\ \label{type-2-cancel}
 |\int B g_{*} h^{\prime} \mathrm{D}(f)   \mathrm{d}V  | &\lesssim& \int |v-v_{*}|^{\gamma+2} \mathrm{1}_{0 \leq \theta \leq \pi/2} \sin^{-2s}\frac{\theta}{2}|g_{*}h^{\prime} \nabla^{2} f(v(\kappa))| \mathrm{d}\kappa \mathrm{d}V.
\een
\end{lem}

The following lemma is a direct consequence of Cancelation Lemma 1 in \cite{alexandre2000entropy} and \eqref{cancel-angular-singularity}.
\begin{lem}\label{cancelation-lemma-result} There holds
$
 |\int B g_{*}\mathrm{D}(h) \mathrm{d}V| \lesssim \int |v-v_{*}|^{\gamma}|g_{*}h| \mathrm{d}v_{*} \mathrm{d}v.
$
\end{lem}

The following Lemma \ref{cancel-velocity-regularity} is used to cancel singularity near $|v-v_{*}|=0$ when $\gamma<0$.
\begin{lem} \label{cancel-velocity-regularity}
For $\gamma>-3, a >0$, there holds
$ 
\int |v-v_{*}|^{\gamma}  \mu^{a}(v) \mathrm{d}v \leq C_{\gamma,a} \langle v_{*} \rangle^{\gamma}.
$ 
\end{lem}
One can refer to Lemma 2.5 in \cite{alexandre2012boltzmann} for a general version and a short proof to Lemma \ref{cancel-velocity-regularity}.

Note that the points $v(\kappa)$ for $0 \leq \kappa \leq 1$ connecting $v$ and $v^{\prime}$
appear
in the formulas \eqref{Taylor-at-v} and \eqref{Taylor-at-v-prime}
. If we apply
Taylor expansion to $\mathrm{D}(f_{*})$, the points $v_{*}(\iota) \colonequals    v_{*} + \iota(v^{\prime}_{*}-v_{*})$ for $0 \leq \iota \leq 1$ will appear. Thanks to the symmetry of elastic collision formula \eqref{v-prime-v-prime-star}, we have
\begin{lem}\label{v-square-sum-bounded} For any $0 \leq \kappa, \iota \leq 1$, there holds
\beno
(1-\frac{\sqrt{2}}{2})(|v|^{2}+|v_{*}|^{2}) \leq |v(\kappa)|^{2}+|v_{*}(\iota)|^{2} \leq (1+\frac{\sqrt{2}}{2})(|v|^{2}+|v_{*}|^{2}).
\eeno
\end{lem}
We omit the proof of Lemma \ref{v-square-sum-bounded} as it is elementary.
We next give a remark on Lemma \ref{v-square-sum-bounded}.
\begin{rmk}\label{optimal-constant}
Since
\beno
\frac{1}{4} \leq 1-\frac{\sqrt{2}}{2} \approx 0.293 \leq \frac{1}{3}, \quad \frac{3}{2} \leq 1+\frac{\sqrt{2}}{2} \approx 1.707 \leq 2,
\eeno
for simplicity in the rest of the article we will use
\beno
\frac{1}{4}(|v|^{2}+|v_{*}|^{2}) \leq |v(\kappa)|^{2}+|v_{*}(\iota)|^{2} \leq 2(|v|^{2}+|v_{*}|^{2}).
\eeno
As a result, recalling $\mu (v) =   (2\pi)^{-\frac{3}{2}} e^{-\f{1}{2}|v|^{2}}$, for any $0 \leq \kappa, \iota \leq 1$, there holds
\ben \label{mu-weight-result}
\mu^{2}(v) \mu^{2}(v_{*}) \leq \mu(v(\kappa)) \mu(v_{*}(\iota)) \leq \mu^{\frac{1}{4}}(v) \mu^{\frac{1}{4}}(v_{*}).
\een
\end{rmk}
We will frequently use \eqref{mu-weight-result} to retain the good negative exponential ($\mu$-type) weight.

\subsection{Estimate of $\mathcal{L}^{\rho}_{r}$}
After the preliminary preparation in subsection \ref{preliminary},
we are ready to give upper bound estimate of the operator $\mathcal{L}^{\rho}_{r}$.
\begin{prop}\label{ub-linearized-L2} There holds
$
|\langle \mathcal{L}^{\rho}_{r}g, h\rangle| \lesssim \rho|\mu^{\frac{1}{64}}g|_{L^{2}}|\mu^{\frac{1}{64}}h|_{L^{2}}.
$
\end{prop}
\begin{proof} Recalling \eqref{L-epsilon-pm-tau-rest-part-original} and \eqref{M--with-mathcal-M}, we have
\ben \label{L-epsilon-pm-tau-rest-part}
(\mathcal{L}^{\rho}_{r}f)(v)  = \rho \int  B N_{*} N^{\prime} N^{\prime}_{*}
\mathrm{D}((N^{-1}f)_{*})
\mathrm{d}\sigma \mathrm{d}v_{*}.
\een
By \eqref{L-epsilon-pm-tau-rest-part} and \eqref{change-v-and-v-prime}, using the identity $ \mathrm{D}(N^{-1}h) = N^{-1}\mathrm{D}(h) + h^{\prime}\mathrm{D}(N^{-1})$, we have
\beno
\langle \mathcal{L}^{\rho}_{r}g, h\rangle = \rho \int B N_{*} N^{\prime} N^{\prime}_{*}
\mathrm{D}((N^{-1}g)_{*})  h
\mathrm{d}V
= \rho \int B N N^{\prime} N^{\prime}_{*} g_{*} \mathrm{D}(N^{-1}h) \mathrm{d}V = \mathcal{I}_{1}+\mathcal{I}_{2},
\\
\mathcal{I}_{1} \colonequals   \rho \int B N^{\prime} N^{\prime}_{*} g_{*} \mathrm{D}(h) \mathrm{d}V,
\quad
\mathcal{I}_{2} \colonequals   \rho \int B N^{\prime}_{*} g_{*} h^{\prime}\mathrm{D}(N^{\prime}) \mathrm{d}V.
\eeno

We now go to see $\mathcal{I}_{1}$. Note that
\beno
N^{\prime} N^{\prime}_{*} \mathrm{D}(h) &=&  \mathrm{D}(N^{\prime}_{*}) N^{\prime} \mathrm{D}(h) + N_{*} N^{\prime} \mathrm{D}(h)
\\&=& - \mathrm{D}(N^{\prime}_{*})(Nh)^{\prime} + \mathrm{D}(N^{\prime}_{*}) \mathrm{D}(N^{\prime}) h + \mathrm{D}(N^{\prime}_{*}) N h
+ N_{*} \mathrm{D}(N^{\prime})h  + N_{*}\mathrm{D}(Nh),
\eeno
which gives
$
\mathcal{I}_{1} =  \mathcal{I}_{1,1} + \mathcal{I}_{1,2} + \mathcal{I}_{1,3} + \mathcal{I}_{1,4} + \mathcal{I}_{1,5},
$
where
\ben \label{I-1-1-new-expression}
\mathcal{I}_{1,1} \colonequals
- \rho \int B \mathrm{D}(N^{\prime}_{*})g_{*}(Nh)^{\prime} \mathrm{d}V =
 \rho \int B \mathrm{D}(N^{\prime}) g^{\prime}(Nh)_{*} \mathrm{d}V,
\\ \nonumber
\mathcal{I}_{1,2} \colonequals    \rho \int B \mathrm{D}(N^{\prime}_{*}) \mathrm{D}(N^{\prime}) g_{*}h \mathrm{d}V, \quad
\mathcal{I}_{1,3} \colonequals    \rho \int B \mathrm{D}(N^{\prime}_{*}) g_{*} N h \mathrm{d}V,
\\ \nonumber
\mathcal{I}_{1,4} \colonequals    \rho \int B N_{*}g_{*} \mathrm{D}(N^{\prime})h \mathrm{d}V, \quad
\mathcal{I}_{1,5} \colonequals    \rho \int B N_{*}g_{*}\mathrm{D}(Nh) \mathrm{d}V.
\een
Here we use \eqref{change-v-and-v-star} and \eqref{change-v-and-v-prime} in \eqref{I-1-1-new-expression}.

In what follows, we derive the estimates of $\mathcal{I}_{1,1}$ and $\mathcal{I}_{1,2}$ with full details because the tricks used here will be used later.
Recalling \eqref{equilibrium-rho-not-vanish},
it is easy to see
\ben \label{derivative-N-M-bounded-by-mu}
N \lesssim \mu^{\f{1}{2}}, \quad
|\nabla N| \lesssim \mu^{\frac{1}{4}}, \quad |\nabla^{2} N| \lesssim \mu^{\frac{1}{4}}. \een
Using \eqref{type-2-cancel} and \eqref{derivative-N-M-bounded-by-mu}, we have
\beno
|\mathcal{I}_{1,1}| \lesssim
 \rho \int |v-v_{*}|^{\gamma+2} \mathrm{1}_{0 \leq \theta \leq \pi/2} \sin^{-2s}\frac{\theta}{2} \mu^{\frac{1}{4}}(v(\kappa)) \mu^{\f{1}{2}}_{*} |g^{\prime}h_{*}| \mathrm{d}\kappa \mathrm{d}V.
\eeno
By \eqref{mu-weight-result}, we have
$
|v-v_{*}|^{2} \mu^{\frac{1}{4}}(v(\kappa)) \mu^{\f{1}{2}}_{*} \lesssim |v-v_{*}|^{2} \mu^{\frac{1}{16}} \mu^{\frac{1}{16}}_{*} \mu^{\frac{1}{4}}_{*} \lesssim \mu^{\frac{1}{32}} \mu^{\frac{1}{32}}_{*} \mu^{\frac{1}{4}}_{*} \lesssim (\mu^{\frac{1}{32}})^{\prime} \mu^{\frac{1}{4}}_{*}
$ and thus
\ben \nonumber
|\mathcal{I}_{1,1}| &\lesssim&
 \rho \int |v-v_{*}|^{\gamma} \mathrm{1}_{0 \leq \theta \leq \pi/2} \sin^{-2s}\frac{\theta}{2} (\mu^{\frac{1}{32}})^{\prime} \mu^{\frac{1}{4}}_{*} |g^{\prime}h_{*}| \mathrm{d}V
\\ \nonumber &\lesssim&
 \rho \int |v^{\prime}-v_{*}|^{\gamma} \mathrm{1}_{0 \leq \theta^{\prime} \leq \pi/4} \sin^{-2s}\theta^{\prime} (\mu^{\frac{1}{32}})^{\prime} \mu^{\frac{1}{4}}_{*} |g^{\prime}h_{*}| \mathrm{d}V^{\prime}
\\ \label{last-line-I-11} &\lesssim&
 \rho \int |v^{\prime}-v_{*}|^{\gamma}  (\mu^{\frac{1}{32}})^{\prime} \mu^{\frac{1}{4}}_{*} |g^{\prime}h_{*}|
 \mathrm{d}v_{*} \mathrm{d}v^{\prime} \lesssim \rho |\mu^{\frac{1}{64}}g|_{L^{2}}|\mu^{\frac{1}{64}}h|_{L^{2}},
\een
where from the first to the second line we use the following change of variable (see Lemma 1 in \cite{alexandre2000entropy} for details)
\ben \label{change-v-to-v-prime}
v \to v^{\prime}, \quad |\frac{\partial v}{\partial v^{\prime}}| \lesssim 1, \quad \theta^{\prime} = \frac{\theta}{2}, \quad
\cos \theta^{\prime} = \frac{v^{\prime}-v_{*}}{|v^{\prime}-v_{*}|} \cdot \sigma.
\een
From the second to the third line we use the fact $\int \mathrm{1}_{0 \leq \theta^{\prime} \leq \pi/4} \sin^{-2s}\theta^{\prime} \mathrm{d}\sigma = 2 \pi \int_{0}^{\pi/4} \sin^{1-2s}\theta^{\prime} d \theta^{\prime} \lesssim \frac{1}{1-s}$. The last inequality in \eqref{last-line-I-11} follows Cauchy-Schwartz inequality and Lemma \ref{cancel-velocity-regularity}. More precisely, we deduce that
\ben \label{cancel-velocity-singularity}
&&\int |v^{\prime}-v_{*}|^{\gamma}  (\mu^{\frac{1}{32}})^{\prime} \mu^{\frac{1}{4}}_{*} |g^{\prime}h_{*}| \mathrm{d}v_{*} \mathrm{d}v^{\prime}
\\ \nonumber &\lesssim& \big( \int |v^{\prime}-v_{*}|^{\gamma}  (\mu^{\frac{1}{32}})^{\prime} \mu^{\frac{1}{4}}_{*} |g^{\prime}|^{2} \mathrm{d}v_{*} \mathrm{d}v^{\prime} \big)^{\f12}
\big(\int |v^{\prime}-v_{*}|^{\gamma}  (\mu^{\frac{1}{32}})^{\prime} \mu^{\frac{1}{4}}_{*} |h_{*}|^{2} \mathrm{d}v_{*} \mathrm{d}v^{\prime}\big)^{\f12}
\\ \nonumber &\lesssim& \big(\int \langle v^{\prime} \rangle^{\gamma}  (\mu^{\frac{1}{32}})^{\prime}  |g^{\prime}|^{2}  \mathrm{d}v^{\prime}\big)^{\f12}
\big(\int \langle v_{*} \rangle^{\gamma}  \mu^{\frac{1}{4}}_{*} |h_{*}|^{2} \mathrm{d}v_{*} \big)^{\f12}
\lesssim  |\mu^{\frac{1}{64}}g|_{L^{2}}|\mu^{\frac{1}{64}}h|_{L^{2}},
\een
where $\gamma \leq 0$ is used to in the last inequality.

We now go to see $\mathcal{I}_{1,2}$.
Recalling $N \sim \mu^{\f{1}{2}}$ by \eqref{M-rho-mu-N}, recalling \eqref{shorthand-D} and \eqref{shorthand-A},
using $|\nabla N^{\f{1}{2}}| \lesssim 1$,
we get
\ben \label{N-prime-N-difference}
|\mathrm{D}(N^{\prime})| = |\mathrm{D}(N)| = |\mathrm{D}(N^{\f{1}{2}}) \mathrm{A}(N^{\f{1}{2}})| \lesssim
|v-v_{*}|\sin\frac{\theta}{2} \mathrm{A}(\mu^{\frac{1}{4}}).
\een
Here the role of $\mathrm{D}(N^{\f{1}{2}})$ is to produce $\sin\frac{\theta}{2}$ to cancel angular singularity later. The factor $\mathrm{A}(\mu^{\frac{1}{4}})$ is kept to retain $\mu$-type weight. Such a treatment will be used frequently in the rest of the article. Similar to \eqref{N-prime-N-difference}, we have
\ben \label{N-prime-star-N-star-difference}
|\mathrm{D}(N^{\prime}_{*})| \lesssim
|v-v_{*}|\sin\frac{\theta}{2} \mathrm{A}(\mu^{\frac{1}{4}}_{*}).
\een
Patching together \eqref{N-prime-N-difference} and \eqref{N-prime-star-N-star-difference}, using
\eqref{mu-weight-result}, we have
\beno
|\mathrm{D}(N^{\prime}_{*}) \mathrm{D}(N^{\prime})| \lesssim \mu^{\frac{1}{16}}(v) \mu^{\frac{1}{16}}(v_{*}) |v-v_{*}|^{2}\sin^{2}\frac{\theta}{2}  \lesssim \mu^{\frac{1}{32}}(v) \mu^{\frac{1}{32}}(v_{*}) \sin^{2}\frac{\theta}{2}.
\eeno
Note that we get the factor  $\sin^{2}\frac{\theta}{2}$ to cancel singularity and also the good $\mu$-type weight factor
$\mu^{\frac{1}{32}}(v) \mu^{\frac{1}{32}}(v_{*})$. As a result, by \eqref{cancel-angular-singularity}, using \eqref{cancel-velocity-singularity}, we get
\beno
|\mathcal{I}_{1,2}|  \lesssim \rho \int |v-v_{*}|^{\gamma} \mu^{\frac{1}{32}}\mu^{\frac{1}{32}}_{*} |g_{*} h| \mathrm{d}v_{*} \mathrm{d}v \lesssim \rho |\mu^{\frac{1}{64}}g|_{L^{2}}|\mu^{\frac{1}{64}}h|_{L^{2}}.
\eeno

By Lemma \ref{cancelation-lemma-result}, using \eqref{cancel-velocity-singularity}, we get
\beno
|\mathcal{I}_{1,5}|  \lesssim \rho \int |v-v_{*}|^{\gamma} \mu^{\f{1}{2}}\mu^{\f{1}{2}}_{*} |g_{*} h| \mathrm{d}v_{*} \mathrm{d}v \lesssim \rho |\mu^{\frac{1}{64}}g|_{L^{2}}|\mu^{\frac{1}{64}}h|_{L^{2}}.
\eeno

Applying \eqref{type-1-cancel} to $\mathcal{I}_{1,4}$, using \eqref{derivative-N-M-bounded-by-mu}, \eqref{mu-weight-result} and \eqref{cancel-angular-singularity}, we get
\beno
|\mathcal{I}_{1,4}|  \lesssim \rho \int |v-v_{*}|^{\gamma} \mu^{\frac{1}{32}}\mu^{\frac{1}{32}}_{*} |g_{*} h| \mathrm{d}v_{*} \mathrm{d}v \lesssim \rho |\mu^{\frac{1}{64}}g|_{L^{2}}|\mu^{\frac{1}{64}}h|_{L^{2}}.
\eeno
Using \eqref{change-v-and-v-star}, we have
$
\mathcal{I}_{1,3} = \rho \int B \mathrm{D}(N^{\prime}) g (N h)_{*} \mathrm{d}V
$
which is the same as $\mathcal{I}_{1,4}$ if we exchange $g$ and $h$. Therefore we also have $|\mathcal{I}_{1,3}| \lesssim \rho |\mu^{\frac{1}{64}}g|_{L^{2}}|\mu^{\frac{1}{64}}h|_{L^{2}}.$

Patching together the above estimates of $\mathcal{I}_{1,i}$ for $1 \leq i \leq 5$, we get
\ben \label{ub-I-1-by-mu-weight}
|\mathcal{I}_{1}|  \lesssim \rho |\mu^{\frac{1}{64}}g|_{L^{2}}|\mu^{\frac{1}{64}}h|_{L^{2}}.
\een

We next go to see $\mathcal{I}_{2}$. It is further decomposed into two terms
$
\mathcal{I}_{2} = \mathcal{I}_{2,1} + \mathcal{I}_{2,2},
$
where
\beno
\mathcal{I}_{2,1} \colonequals
\rho \int B \mathrm{D}(N^{\prime}_{*})g_{*} h^{\prime}\mathrm{D}(N^{\prime}) \mathrm{d}V,
\quad
\mathcal{I}_{2,2} \colonequals   \rho \int B N_{*} g_{*} h^{\prime}\mathrm{D}(N^{\prime}) \mathrm{d}V.
\eeno
Comparing $\mathcal{I}_{2,1}$ and  $\mathcal{I}_{1,2}$, we can use the same arguments for $\mathcal{I}_{1,2}$ and \eqref{change-v-to-v-prime} to get
$
|\mathcal{I}_{2,1}|  \lesssim \rho|\mu^{\frac{1}{64}}g|_{L^{2}}|\mu^{\frac{1}{64}}h|_{L^{2}}.
$
Note that $\mathcal{I}_{2,2}$ is the same as $\mathcal{I}_{1,1}$ if we exchange $g$ and $h$. Therefore we also have $|\mathcal{I}_{2,2}| \lesssim \rho |\mu^{\frac{1}{64}}g|_{L^{2}}|\mu^{\frac{1}{64}}h|_{L^{2}}.$
Patching together the above estimates of $\mathcal{I}_{2,1}$ and $\mathcal{I}_{2,2}$, we get
\ben \label{ub-I-2-by-mu-weight}
|\mathcal{I}_{2}|  \lesssim \rho|\mu^{\frac{1}{64}}g|_{L^{2}}|\mu^{\frac{1}{64}}h|_{L^{2}}.
\een

Patching together \eqref{ub-I-1-by-mu-weight} and \eqref{ub-I-2-by-mu-weight}, we finish the proof.
\end{proof}

\subsection{Estimate of $\mathcal{C}^{\rho}$}
Based on  $B$(see \eqref{Boltzmann-kernel-general-and-old}) and $B^{ipl}$(see \eqref{inverse-power-law-kernel}),
we define
\ben \label{funcation-N-pm}
\mathcal{N}(g, h)\colonequals     \int B g_{*}^{2}\mathrm{D}^{2}(h) \mathrm{d}V , \quad
\mathcal{N}^{ipl}(g, h)\colonequals     \int B^{ipl} g_{*}^{2}\mathrm{D}^{2}(h) \mathrm{d}V .
\een

\begin{thm} \label{two-functional-ub-by-norm} For $\frac{1}{16} \leq a\leq 1$, there holds
\beno
\mathcal{N}(\mu^{a}, f) + \mathcal{N}(f, \mu^{a}) \lesssim |f|_{\mathcal{L}^{s}_{\gamma/2}}^{2}.
\eeno
\end{thm}
\begin{proof} By \eqref{Boltzmann-kernel-general-ub}, we have $\mathcal{N}(\mu^{a}, f) + \mathcal{N}(f, \mu^{a}) \lesssim \mathcal{N}^{ipl}(\mu^{a}, f) + \mathcal{N}^{ipl}(f, \mu^{a}).$ Then the desired estimate is a direct result of $\mathcal{N}^{ipl}(\mu^{a}, f) + \mathcal{N}^{ipl}(f, \mu^{a}) \lesssim |f|_{\mathcal{L}^{s}_{\gamma/2}}^{2}$ by \cite{he2022asymptotic}(Theorem 1.1 and the proof of Theorem 2.1). We remark that notations in \cite{he2022asymptotic} and this article are different.
\end{proof}

Recalling \eqref{M-rho-mu-N} and following the proof of Theorem 1.1 and Theorem 2.1 in \cite{he2022asymptotic},  we give the following remark.
\begin{rmk} \label{mu-to-N-or-M-functional-N} If  $\mu$ is replaced by $N^{2}$ or $M$ in Theorem \ref{two-functional-ub-by-norm},
the result is still valid.
\end{rmk}

We now derive an upper bound estimate for the linear operator $\mathcal{C}^{\rho}$ defined in \eqref{correction-operator-linear}.
\begin{prop}\label{ub-for-correction-term} There holds
$
|\langle \mathcal{C}^{\rho} h, f\rangle| \lesssim \rho^{2} |\mu^{\frac{1}{4}}h|_{\mathcal{L}^{s}_{\gamma/2}}|\mu^{\frac{1}{4}}f|_{\mathcal{L}^{s}_{\gamma/2}}.
$
\end{prop}
\begin{proof} Recalling \eqref{equilibrium-rho-not-vanish} and \eqref{M--with-mathcal-M}, we rewrite the operator $\mathcal{C}^{\rho}$ in \eqref{correction-operator-linear} as
\beno
\mathcal{C}^{\rho} h  =  -\rho^{2} \int  B N_{*} N^{\prime} N^{\prime}_{*} \mathrm{D}(\mu^{\f{1}{2}}h)
\mathrm{d}\sigma \mathrm{d}v_{*}.
\eeno
Taking inner product with $f$, we have
\beno
\langle \mathcal{C}^{\rho} h, f\rangle  =  -\rho^{2} \int  B N N_{*} N^{\prime} N^{\prime}_{*} \mathrm{D}(\mu^{\f{1}{2}}h) N^{-1}f
 \mathrm{d}V
=  -\f{1}{2} \rho^{2} \int  B N N_{*} N^{\prime} N^{\prime}_{*} \mathrm{D}(\mu^{\f{1}{2}}h)
\mathrm{D}(N^{-1}f) \mathrm{d}V  ,
\eeno
where we use \eqref{change-v-and-v-prime} in the second line.
Then by Cauchy-Schwartz inequality and \eqref{K-2-mu}, we get
\beno
|\langle \mathcal{C}^{\rho} h, f\rangle|  &\lesssim& \rho^{2} \big(\int  B \mu^{\f{1}{2}} \mu^{\f{1}{2}}_{*} \mathrm{D}^{2}(\mu^{\f{1}{2}}h)  \mathrm{d}V   \big)^{\f12} \big(\int  B \mu^{\frac{3}{2}} \mu^{\frac{3}{2}}_{*} \mathrm{D}^{2}(N^{-1}f) \mathrm{d}V  \big)^{\f12}
\\ &\lesssim& \rho^{2}|\mu^{\f{1}{2}}h|_{\mathcal{L}^{s}_{\gamma/2}} \big(\int  B \mu^{\frac{3}{2}} \mu^{\frac{3}{2}}_{*} \mathrm{D}^{2}(N^{-1}f) \mathrm{d}V  \big)^{\f12},
\eeno
where we use Theorem \ref{two-functional-ub-by-norm} in the second line. Noting that $N^{-1}=\mu^{-\f{1}{2}} - \rho \mu^{\f{1}{2}}$ by \eqref{equilibrium-rho-not-vanish}, we get
\beno
\mathrm{D}(N^{-1}f) = \mathrm{D}(\mu^{-\f{1}{2}}f) - \rho \mathrm{D}(\mu^{\f{1}{2}}f),
\eeno
which gives
$
\int  B \mu^{\frac{3}{2}} \mu^{\frac{3}{2}}_{*} \mathrm{D}^{2}(N^{-1}f) \mathrm{d}V   \lesssim \mathcal{I}_{1} + \mathcal{I}_{2},$ where
\beno
\mathcal{I}_{1} \colonequals  \int  B \mu^{\frac{3}{2}} \mu^{\frac{3}{2}}_{*} \mathrm{D}^{2}(\mu^{-\f{1}{2}}f) \mathrm{d}V = \int  B  \mathrm{D}^{2}(\mu^{\frac{3}{4}}_{*}\mu^{\frac{1}{4}}f) \mathrm{d}V, \quad
\mathcal{I}_{2} \colonequals \int  B \mu^{\frac{3}{2}} \mu^{\frac{3}{2}}_{*} \mathrm{D}^{2}(\mu^{\f{1}{2}}f) \mathrm{d}V  .
\eeno
By Theorem \ref{two-functional-ub-by-norm}, we have
$
\mathcal{I}_{2} \lesssim |\mu^{\f{1}{2}}f|_{\mathcal{L}^{s}_{\gamma/2}}^{2}.
$
Since $\mathrm{D}(\mu^{\frac{3}{4}}_{*}\mu^{\frac{1}{4}}f)
=\mu^{\frac{3}{4}}_{*}\mathrm{D}(\mu^{\frac{1}{4}}f)+
\mathrm{D}(\mu^{\frac{3}{4}}_{*})(\mu^{\frac{1}{4}}f)^{\prime}$, we get
\beno
\mathcal{I}_{1}
\leq 2 \int  B  \mu^{\frac{3}{2}}_{*}\mathrm{D}^{2}(\mu^{\frac{1}{4}}f)\mathrm{d}V   +
2\int  B  \mathrm{D}^{2}(\mu^{\frac{3}{4}}_{*}) (\mu^{\f{1}{2}}f^{2})^{\prime} \mathrm{d}V   \lesssim |\mu^{\frac{1}{4}}f|_{\mathcal{L}^{s}_{\gamma/2}}^{2},
\eeno
where we use \eqref{change-v-and-v-star}, \eqref{change-v-and-v-prime} and Theorem \ref{two-functional-ub-by-norm}. Patching together the above estimates, we finish the proof.
\end{proof}

\subsection{Estimate of $\tilde{Q}$} \label{estimate-of-tilde-Q}
We now derive a result about $\tilde{Q}$ which will be used in
Proposition \ref{positivity} to derive non-negativity of solutions to the linear equation \eqref{quantum-Boltzmann-UU-linear}.
\begin{prop}\label{for-positivity} If $g \geq  0$, then
$
\langle \tilde{Q}(g, f), f\rangle \lesssim |\mu^{-\f14} g|_{H^{4}}(1+|\mu^{-\f14} g|_{H^{4}}) |f|_{L^{2}}^{2}.
$
\end{prop}
\begin{proof}
Recalling \eqref{U-U-operator-linear}, using $f\mathrm{D}(f^{\prime}) = -\f12 \mathrm{D}^{2}(f)
- \f12 \mathrm{D}(f^{2})$ and $g_{*}(1 + g_{*}^{\prime} +  g^{\prime}) \geq 0$,
we have
\beno
\langle \tilde{Q}(g, f), f\rangle &=& \int B
 g_{*} f (1 + g_{*}^{\prime} +  g^{\prime})\mathrm{D}(f^{\prime}) \mathrm{d}V
\\&=& - \f12 \int B
 g_{*}  (1 + g_{*}^{\prime} +  g^{\prime}) \mathrm{D}^{2}(f) \mathrm{d}V
 -  \f12 \int B
 g_{*}  (1 + g_{*}^{\prime} +  g^{\prime})\mathrm{D}(f^{2}) \mathrm{d}V
\\&\leq& - \f12 \int B
 g_{*}  (1 + g_{*}^{\prime} +  g^{\prime})\mathrm{D}(f^{2}) \mathrm{d}V = - \f12(\mathcal{I}_{1} + \mathcal{I}_{2}+ \mathcal{I}_{3}),
\eeno
where
\beno
\mathcal{I}_{1} \colonequals   \int B
 g_{*}  \mathrm{D}(f^{2}) \mathrm{d}V, \quad
\mathcal{I}_{2} \colonequals   \int B
 g_{*}  g_{*}^{\prime} \mathrm{D}(f^{2}) \mathrm{d}V, \quad
\mathcal{I}_{3} \colonequals   \int B
 g_{*}   g^{\prime} \mathrm{D}(f^{2}) \mathrm{d}V.
\eeno
By Lemma \ref{cancelation-lemma-result}, the imbedding $H^{2} \hookrightarrow L^{\infty}$ and
Lemma \ref{cancel-velocity-regularity}, we get
\beno
|\mathcal{I}_{1}| \lesssim  \int |v-v_{*}|^{\gamma} g_{*} f^{2}
 \mathrm{d}v_{*} \mathrm{d}v \lesssim |\mu^{-\f14} g|_{L^{\infty}} \int |v-v_{*}|^{\gamma} \mu^{\f14}_{*} f^{2}
 \mathrm{d}v_{*} \mathrm{d}v \lesssim |\mu^{-\f14} g|_{H^{2}} |f|_{L^{2}_{\gamma/2}}^{2} \lesssim |\mu^{-\f14} g|_{H^{2}} |f|_{L^{2}}^{2}.
\eeno
By \eqref{change-v-and-v-prime}, we have $\mathcal{I}_{2}=0$ and $\mathcal{I}_{3} = \int B
 \mathrm{D}(g_{*}g^{\prime})  f^{2} \mathrm{d}V
  = \int B  \mathrm{D}(g_{*}) g f^{2} \mathrm{d}V + \int B
g_{*} \mathrm{D}(g^{\prime}) f^{2} \mathrm{d}V$.
By \eqref{change-v-and-v-star}, Lemma \ref{cancelation-lemma-result}, the imbedding $H^{2} \hookrightarrow L^{\infty}$ and
Lemma \ref{cancel-velocity-regularity}, we get
\beno
|\int B  \mathrm{D}(g_{*}) g f^{2} \mathrm{d}V| = |\int B  \mathrm{D}(g) (g f^{2})_{*} \mathrm{d}V|
\lesssim \int |v-v_{*}|^{\gamma} |g (g f^{2})_{*}|
 \mathrm{d}v_{*} \mathrm{d}v \lesssim |\mu^{-\f14} g|_{H^{2}}^{2} |f|_{L^{2}}^{2}.
\eeno
By \eqref{type-1-cancel}, we get
\beno
|\int B
g_{*} \mathrm{D}(g^{\prime}) f^{2} \mathrm{d}V| \lesssim \int |v-v_{*}|^{\gamma+1}|g_{*}f^{2} \nabla g| \mathrm{d}v_{*} \mathrm{d}v
+ \int |v-v_{*}|^{\gamma+2} \sin^{-2s}\frac{\theta}{2}|g_{*}f^{2} \nabla^{2} g(v(\kappa))| \mathrm{d}\kappa \mathrm{d}V.
\eeno
Noting the following fact
\beno
|\nabla g| = |\nabla (\mu^{\f14}\mu^{-\f14}g)| \lesssim | \mu^{-\f14}g \nabla \mu^{\f14}| + |\mu^{\f14} \nabla (\mu^{-\f14}g)|
\lesssim \mu^{\f18} (| \mu^{-\f14}g|_{L^{\infty}} + |\nabla (\mu^{-\f14}g)|_{L^{\infty}}) \lesssim \mu^{\f18} | \mu^{-\f14}g|_{H^{3}},
\eeno
and similarly $|\nabla^{2} g|  \lesssim \mu^{\f18} |\mu^{-\f14}g|_{H^{4}}$, using \eqref{mu-weight-result}, \eqref{cancel-angular-singularity} and Lemma \ref{cancel-velocity-regularity},
we get
\beno
|\int B
g_{*} \mathrm{D}(g^{\prime}) f^{2} \mathrm{d}V| \lesssim  |\mu^{-\f14}g|_{H^{4}}^{2} \int (|v-v_{*}|^{\gamma+1}+|v-v_{*}|^{\gamma+2}) \mu^{\f18}_{*} \mu^{\f{1}{32}} \mu^{\f{1}{32}}_{*}  f^{2} \mathrm{d}v_{*} \mathrm{d}v \lesssim |\mu^{-\f14} g|_{H^{2}}^{4} |f|_{L^{2}}^{2}.
\eeno
Patching together the above estimates, we  finish the proof.
\end{proof}

\section{Bilinear operator estimate} \label{bilinear}

Recall \eqref{equilibrium-rho-not-vanish}. For simplicity we continue to write $N=N_{\rho}, M=M_{\rho}$. Recalling the relation \eqref{M--with-mathcal-M} between
$N, M$ and $\mathcal{N}, \mathcal{M}$, the definition of $\Gamma_{2}^{\rho}$ in \eqref{definition-Gamma-2-epsilon} and $\Pi_{2}$ in \eqref{definition-A-2}, we have
\ben \label{Gamma-main-remaining}
\Gamma_{2}^{\rho}(g,h) = \rho^{\f{1}{2}}\Gamma_{2,m}^{\rho}(g,h) + \rho^{\frac{3}{2}}\Gamma_{2,r}^{\rho}(g,h).
\een
Here $\Gamma_{2,m}^{\rho}$ stands for the main term (``$m$'' is referred to ``main'') defined by
\ben \label{definition-Gamma-main}
\Gamma_{2,m}^{\rho}(g,h) \colonequals   N^{-1} Q_{c} (Ng, Nh) = Q_{c} (Ng, h) + I^{\rho}(g,h),
\een
where $Q_{c}$ is defined in \eqref{classical-Boltzmann-operator-UU} and $I^{\rho}$ is defined by
\ben \label{definition-I-epsilon}
I^{\rho}(g,h) \colonequals     \int
B \mathrm{D}(N^{-1})(Ng)^{\prime}_{*}(Nh)^{\prime}
\mathrm{d}\sigma \mathrm{d}v_{*}.
\een
Here $\Gamma_{2,r}^{\rho}$ represents the remaining term (``$r$'' is referred to ``remaining''). This term consists of three parts $\Gamma_{2,r,1}^{\rho},\Gamma_{2,r,2}^{\rho},\Gamma_{2,r,3}^{\rho}$ corresponding to \eqref{line-1}, \eqref{line-2}, \eqref{line-3} respectively,
\ben
 \label{Gamma-remaining-into-three-terms} \Gamma_{2,r}^{\rho}(g,h) &\colonequals  &  \Gamma_{2,r,1}^{\rho}(g,h)
+\Gamma_{2,r,2}^{\rho}(g,h)+\Gamma_{2,r,3}^{\rho}(g,h).
\\
\label{definition-Gamma-remaining-1}
\Gamma_{2,r,1}^{\rho}(g,h) &\colonequals    & N^{-1} \int
B \mathrm{D} \big((Ng)^{\prime}_{*}(Nh)^{\prime}(M + M_{*}) \big)
\mathrm{d}\sigma \mathrm{d}v_{*}.
\\ \label{definition-Gamma-remaining-2}
\Gamma_{2,r,2}^{\rho}(g,h) &\colonequals    & N^{-1} \int
B \mathrm{D}\big((Ng)_{*}(Nh)^{\prime}(M^{\prime}_{*} - M) \big)
\mathrm{d}\sigma \mathrm{d}v_{*}.
\\ \label{definition-Gamma-remaining-3}
\Gamma_{2,r,3}^{\rho}(g,h) &\colonequals    & N^{-1} \int
B \big((Ng)^{\prime}(Nh) \mathrm{D}(M^{\prime}_{*}) + (Ng)^{\prime}_{*}(Nh)_{*} \mathrm{D}(M^{\prime}) \big)
\mathrm{d}\sigma \mathrm{d}v_{*}.
\een

We call $\Gamma_{2,m}^{\rho}$ the main term for two reasons. First, the factor before $\Gamma_{2,m}^{\rho}$ is $\rho^{\f{1}{2}}$, while the factor before $\Gamma_{2,r}^{\rho}(g,h)$ is $\rho^{\frac{3}{2}}$. Second, when $\rho=0$, the term $\Gamma_{2,m}^{\rho}$ corresponds to the nonlinear term $\mu^{-\f{1}{2}} Q_{c} (\mu^{\f{1}{2}}g, \mu^{\f{1}{2}} h)$ in the classical linearized Boltzmann equation.

\subsection{The main operator $\Gamma_{2, m}^{\rho}$}



As a result of Theorem 2.2 in \cite{he2022asymptotic} and \eqref{Boltzmann-kernel-general-ub},
 we have
\begin{prop} \label{Q-up-bound-full-on-first-or-second-prop} It holds that
\ben \label{Q-g-h-f-up-full-on-g}
 |  \langle  Q_{c}(g, h), f \rangle   |
&\lesssim&  |g|_{H^{2}_{7}}  |h|_{\mathcal{L}^{s}_{\gamma/2}} |f|_{\mathcal{L}^{s}_{\gamma/2}}, \\ \label{Q-g-h-f-up-full-on-h}
 |  \langle  Q_{c}(g, h), f \rangle   |  &\lesssim&  |g|_{L^{2}_{7}}
(|h|_{\mathcal{L}^{s}_{\gamma/2}}+|h|_{H^{s+2}_{\gamma/2}}) |f|_{\mathcal{L}^{s}_{\gamma/2}}.
\een
\end{prop}

For the convenience of later reference, we write Proposition \ref{Q-up-bound-full-on-first-or-second-prop} as
\begin{col}\label{Q-up-bound-full-on-first-or-second} There holds
\beno
 |  \langle  Q_{c}(g, h), f \rangle   |  \lesssim  \min\{|g|_{H^{2}_{7}}  |h|_{\mathcal{L}^{s}_{\gamma/2}} , |g|_{L^{2}_{7}}(|h|_{\mathcal{L}^{s}_{\gamma/2}}+|h|_{H^{s+2}_{\gamma/2}})\} |f|_{\mathcal{L}^{s}_{\gamma/2}}.
\eeno
\end{col}

The following lemma is used to deal with the norm of the product of two functions.

\begin{lem} \label{product-take-out} Let $0 \leq a \leq 2 \leq b$.
Then
\ben \label{deal-with-product}
|g h|_{H^{a}} \lesssim |g|_{H^{2}} |h|_{H^{a}}, \quad |g h|_{H^{b}} \lesssim |g|_{H^{b}} |h|_{H^{b}}.
\een
Let $ n \geq 0, l \in \mathbb{R}, c>0$. There exists a constant $C(n,l,c)$ such that
\ben \label{deal-with-polynomial-weight}
|\mu^{c} W_{l}|_{H^{n}} \leq C(n,l,c), \quad |\mu^{2c} f|_{H^{n}_{l}} \leq C(n,l,c) |\mu^{c} f|_{H^{n}}.
\een
Recall $-3<\gamma<0<s<1$. There holds
\ben \label{taking-out-with-a-weight}
|g h|_{\mathcal{L}^{s}_{\gamma/2}} \lesssim |g|_{H^{2}_{1}} |h|_{\mathcal{L}^{s}_{\gamma/2}}.
\een

\end{lem}
\begin{proof} The two results \eqref{deal-with-product} and  \eqref{deal-with-polynomial-weight} are standard. We now prove \eqref{taking-out-with-a-weight}. Recalling \eqref{definition-of-norm-L-epsilon-gamma} and using \eqref{deal-with-product}, we have
\beno
|g h|_{\mathcal{L}^{s}_{\gamma/2}} \lesssim |g h|_{H^{s}_{s+\gamma/2}} \lesssim |g|_{H^{2}_{s}}|h|_{H^{s}_{\gamma/2}}
\lesssim |g|_{H^{2}_{1}} |h|_{\mathcal{L}^{s}_{\gamma/2}}.
\eeno
\end{proof}

As a direct application of Corollary \ref{Q-up-bound-full-on-first-or-second}, using \eqref{deal-with-product} and the fact
$|\mu^{-\frac{1}{8}} N|_{H^{2}_{7}} \lesssim 1$, we have
\begin{prop} \label{Q-pm-Ng-h-f} There holds
\beno
 |  \langle  Q_{c}(Ng, h), f \rangle   |  \lesssim  \min\{|\mu^{\frac{1}{8}}g|_{H^{2}}  |h|_{\mathcal{L}^{s}_{\gamma/2}} , |\mu^{\frac{1}{8}}g|_{L^{2}}(|h|_{\mathcal{L}^{s}_{\gamma/2}}+|h|_{H^{s+2}_{\gamma/2}})\} |f|_{\mathcal{L}^{s}_{\gamma/2}}.
\eeno
\end{prop}
Note that Proposition \ref{Q-pm-Ng-h-f} has the flexibility to balance the regularity between $g$ and $h$. Such flexibility allows us to close energy estimate in high order Sobolev spaces.

Now we set to consider $I^{\rho}(g,h).$ We first prepare some intermediate estimates.

\begin{lem} \label{functional-X-g-h} Let $\frac{1}{16} \leq a,b \leq 1$. Let $s_{1},s_{2},s_{3}\geq 0,s_{1}+s_{2}+s_{3}=\f{1}{2}$,
then the following four estimates are valid.
\ben \label{estimate-of-X-g-h-f-total-one-half}
|\int
 B \mathrm{D}(\mu^{a}_{*}) g_{*} h f \mathrm{d}V| &\lesssim& |g|_{L^{1}} |h|_{L^{2}_{\gamma/2+s}} |f|_{L^{2}_{\gamma/2+s}} + |\mu^{\frac{1}{256}}g|_{H^{s_{1}}} |\mu^{\frac{1}{256}}h|_{H^{s_{2}}} |\mu^{\frac{1}{256}}f|_{H^{s_{3}}}.
\\ \label{estimate-of-X-g-h-f-total-one-half-with-mu}
|\int
 B \mathrm{D}(\mu^{a}_{*})\mu^{b} g_{*} h f \mathrm{d}V| &\lesssim& |\mu^{\frac{1}{256}}g|_{H^{s_{1}}} |\mu^{\frac{1}{256}}h|_{H^{s_{2}}}|\mu^{\frac{1}{256}}f|_{H^{s_{3}}}.
\\ \label{estimate-of-Y-g-h-f-total-one-half-with-mu}
|\int
 B \mathrm{D}(\mu^{a}_{*})\mu^{b} g (h f)_{*} \mathrm{d}V| &\lesssim& |\mu^{\frac{1}{256}}g|_{H^{s_{1}}} |\mu^{\frac{1}{256}}h|_{H^{s_{2}}}|\mu^{\frac{1}{256}}f|_{H^{s_{3}}}.
\\ \label{estimate-of-D-g-h-f-total-one-half-with-mu}
|\int
 B \mathrm{D}(\mu^{a}_{*})\mu^{b} g h f \mathrm{d}V| &\lesssim& |\mu^{\frac{1}{256}}g|_{H^{2}} |\mu^{\frac{1}{256}}h|_{L^{2}}|\mu^{\frac{1}{256}}f|_{L^{2}}.
\een
Let $s_{1},s_{2}\geq 0,s_{1}+s_{2}=\f{1}{2}$,
then the following two estimates are valid.
\ben \label{estimate-of-Z-g2-h2-total-one-half}
\int
B  \mathrm{D}^{2}(\mu^{a}_{*}) g_{*}^{2} h^{2} \mathrm{d}V
&\lesssim& |g|_{H^{s_{1}}}^{2}|h|_{H^{s_{2}}_{\gamma/2+s}}^{2}.
\\ \label{estimate-of-Z-g2-h2-total-one-half-with-mu}
\int
B  \mathrm{D}^{2}(\mu^{a}_{*}) \mu^{b} g_{*}^{2} h^{2} \mathrm{d}V &\lesssim& |\mu^{\frac{1}{256}}g|_{H^{s_{1}}}^{2}|\mu^{\frac{1}{256}}h|_{H^{s_{2}}}^{2}.
\een
Let $\{a_{1},a_{2},a_{3}\} = \{0, \f12, 2\}$(which means the two sets are equal) and $s_{1},s_{2} \geq 0,s_{1}+s_{2}=\f{1}{2}$, then the following four estimates are valid.
\ben \label{estimate-of-A-g-h-varrho-f-total-one-half-plus-2}
|\int
B  \mathrm{D}(\mu^{a}_{*}) g_{*}  h \varrho_{*} f \mathrm{d}V| &\lesssim& |g|_{H^{a_{1}}} |h|_{H^{a_{2}}_{\gamma/2+s}} |\varrho|_{H^{a_{3}}} |f|_{L^{2}_{\gamma/2+s}}.
\\ \label{estimate-of-A-g-h-varrho-f-with-mu}
|\int
B  \mathrm{D}(\mu^{a}_{*})\mu^{b} g_{*}  h \varrho_{*} f \mathrm{d}V| &\lesssim&
|\mu^{\frac{1}{256}}g|_{H^{a_{1}}} |\mu^{\frac{1}{256}}h|_{H^{a_{2}}} |\mu^{\frac{1}{256}}\varrho|_{H^{a_{3}}} |\mu^{\frac{1}{256}}f|_{L^{2}}.
\\ \label{estimate-of-B-g-h-varrho-f-with-mu}
|\int
B  \mathrm{D}(\mu^{a}_{*}) \mu^{b} g_{*}  h \varrho f \mathrm{d}V| &\lesssim&
|\mu^{\frac{1}{256}}g|_{H^{s_{1}}} |\mu^{\frac{1}{256}}h|_{H^{s_{2}}} |\mu^{\frac{1}{256}}\varrho|_{H^{2}} |\mu^{\frac{1}{256}}f|_{L^{2}}.
\\ \label{estimate-of-C-g-h-varrho-f-with-mu}
|\int
B  \mathrm{D}(\mu^{a}_{*}) \mu^{b} g  (h \varrho f)_{*} \mathrm{d}V| &\lesssim&
|\mu^{\frac{1}{256}}g|_{H^{s_{1}}} |\mu^{\frac{1}{256}}h|_{H^{s_{2}}} |\mu^{\frac{1}{256}}\varrho|_{H^{2}} |\mu^{\frac{1}{256}}f|_{L^{2}}.
\een
\end{lem}
The proof of Lemma \ref{functional-X-g-h} is given in the Appendix \ref{appendix}.
Thanks to Remark \ref{mu-to-N-or-M}, we have
\begin{rmk} \label{mu-to-N-or-M-X-Y} If we replace $\mu$ with $N^{2}$ or $M$ in Lemma \ref{functional-X-g-h},
all the results are still valid. Let $P_{1}, P_{2}$ be two polynomials on $\mathbb{R}^{3}$. If we replace $\mu^{a}$ and $\mu^{b}$ with  $P_{1}\mu^{a}$ and $P_{2}\mu^{b}$ respectively in Lemma \ref{functional-X-g-h},
all the results are still valid. Moreover, if we replace $\mu^{a}$ and $\mu^{b}$ with  $P_{1}N^{2a}$(or $P_{1}M^{a}$) and $P_{2}N^{2b}$(or $P_{2}M^{b}$) respectively in Lemma \ref{functional-X-g-h},
all the results are still valid. Since $|v-v_{*}| \sim |v^{\prime}-v_{*}|$, if we replace $h^{2}$ with $(h^{2})^{\prime}$ in
\eqref{estimate-of-Z-g2-h2-total-one-half}, the result is still valid.
\end{rmk}
Remark \ref{mu-to-N-or-M-X-Y} says that Lemma \ref{functional-X-g-h} have a more general version. This flexibility allows us to deal with various similar integrals. For example, we will use Remark \ref{mu-to-N-or-M-X-Y} to deal with integrals involving derivatives of $\mu$ in subsection \ref{derivative}.

Now we are ready to prove the following upper bound estimate for the operator $I^{\rho}$.
\begin{prop} \label{I-pm-rho-g-h-f} Let $s_{1},s_{2}\geq 0,s_{1}+s_{2}=\f{1}{2}$.
There holds
$
 |  \langle  I^{\rho}(g,h), f \rangle   |  \lesssim |g|_{H^{s_{1}}}|h|_{H^{s_{2}}_{\gamma/2+s}} |f|_{\mathcal{L}^{s}_{\gamma/2}}.
$
\end{prop}
\begin{proof} Note that
$\mathrm{D}(N^{-1}) = \mathrm{D}(\mu^{-\f{1}{2}})- \rho \mathrm{D}(\mu^{\f{1}{2}})
$ by recalling \eqref{equilibrium-rho-not-vanish}.
From which together with $\mu^{\f{1}{2}}\mu^{\f{1}{2}}_{*}=(\mu^{\f{1}{2}})^{\prime}(\mu^{\f{1}{2}})^{\prime}_{*}$,
we have
\ben \label{N-power-minus-1-difference}
\mathrm{D}(N^{-1})N^{\prime}_{*}N^{\prime} = \big(\mathrm{D}(\mu^{\f{1}{2}}_{*}) - \rho \mathrm{D}(\mu^{\f{1}{2}}) \mu^{\f{1}{2}}\mu^{\f{1}{2}}_{*}\big) (\frac{1}{1 - \rho \mu})^{\prime}_{*}(\frac{1}{1 - \rho \mu})^{\prime}.
\een
Plugging \eqref{N-power-minus-1-difference} into \eqref{definition-I-epsilon}, we have
\ben \label{I-into-main-and-rest}
I^{\rho}(g,h) &=& I^{\rho}_{m}(g,h) + I^{\rho}_{r}(g,h),
\\ \label{definition-Im-epsilon}
I^{\rho}_{m}(g,h) &\colonequals    &
\int
B \mathrm{D}(\mu^{\f{1}{2}}_{*})
(\frac{g}{1 - \rho \mu})^{\prime}_{*}(\frac{h}{1 - \rho \mu})^{\prime}
\mathrm{d}\sigma \mathrm{d}v_{*} ,
\\ \label{definition-Ir-epsilon}
I^{\rho}_{r}(g,h) &\colonequals    &
-\rho \int
B \mathrm{D}(\mu^{\f{1}{2}}) \mu^{\f{1}{2}}\mu^{\f{1}{2}}_{*} (\frac{g}{1 - \rho \mu})^{\prime}_{*}(\frac{h}{1 - \rho \mu})^{\prime}
\mathrm{d}\sigma \mathrm{d}v_{*}.
\een
Now it suffices to estimate $ \langle  I^{\rho}_{m}(g,h), f \rangle  $ and $ \langle  I^{\rho}_{r}(g,h), f \rangle  $.

{\it Estimate of $ \langle  I^{\rho}_{m}(g,h), f \rangle  $.} Recalling \eqref{definition-Im-epsilon} and using \eqref{change-v-and-v-prime}, we have
\beno
 \langle  I^{\rho}_{m}(g,h), f \rangle    = \int
B \mathrm{D}((\mu^{\f{1}{2}})^{\prime}_{*})
(\frac{g}{1 - \rho \mu})_{*}(\frac{h}{1 - \rho \mu}) f^{\prime} \mathrm{d}V.
\eeno
Using $f^{\prime} = \mathrm{D}(f^{\prime}) + f$ and
$
\mathrm{D}((\mu^{\f{1}{2}})^{\prime}_{*})
= \mathrm{D}^{2}(\mu^{\f{1}{4}}_{*}) + 2\mu_{*}^{\frac{1}{4}}\mathrm{D}((\mu^{\f{1}{4}})^{\prime}_{*}),
$
we get $\langle  I^{\rho}_{m}(g,h), f \rangle   = \mathcal{I}^{\rho}_{m,1} + \mathcal{I}^{\rho}_{m,2} + \mathcal{I}^{\rho}_{m,3}
$ where
\ben \nonumber
\mathcal{I}^{\rho}_{m,1} &\colonequals    & \int
B \mathrm{D}^{2}(\mu^{\f{1}{4}}_{*})
(\frac{g}{1 - \rho \mu})_{*}(\frac{h}{1 - \rho \mu}) f^{\prime} \mathrm{d}V,
\\ \label{definition-I-m2}
\mathcal{I}^{\rho}_{m,2} &\colonequals    & 2 \int
B \mu_{*}^{\frac{1}{4}}\mathrm{D}((\mu^{\f{1}{4}})^{\prime}_{*})
(\frac{g}{1 - \rho \mu})_{*}(\frac{h}{1 - \rho \mu}) \mathrm{D}(f^{\prime}) \mathrm{d}V,
\\ \label{definition-I-m3}
\mathcal{I}^{\rho}_{m,3} &\colonequals    & 2 \int
B \mu_{*}^{\frac{1}{4}}\mathrm{D}((\mu^{\f{1}{4}})^{\prime}_{*})
(\frac{g}{1 - \rho \mu})_{*}(\frac{h}{1 - \rho \mu}) f \mathrm{d}V.
\een
Since $0 \leq \rho \leq \f{1}{2}(2 \pi)^{\frac{3}{2}}$,
there holds $1 \leq \frac{1}{1 - \rho \mu} \leq 2$. By Cauchy-Schwartz inequality, \eqref{change-v-and-v-star} and \eqref{change-v-and-v-prime},  recalling \eqref{funcation-N-pm}, using \eqref{estimate-of-Z-g2-h2-total-one-half} and Theorem \ref{two-functional-ub-by-norm},
we get
\ben \label{to-X-g-h-Y-f}
|\mathcal{I}^{\rho}_{m,1}| \lesssim
 \big( \int
B \mathrm{D}^{2}(\mu^{\f{1}{4}}_{*})
g_{*}^{2} h^{2} \mathrm{d}V \big)^{\f{1}{2}} \big( \int
B f_{*}^{2}\mathrm{D}^{2}(\mu^{\f{1}{4}})
\mathrm{d}V \big)^{\f{1}{2}}  \lesssim |g|_{H^{s_{1}}}|h|_{H^{s_{2}}_{\gamma/2+s}} |f|_{\mathcal{L}^{s}_{\gamma/2}}.
\een
Recalling \eqref{definition-I-m2}, by Cauchy-Schwartz inequality, recalling \eqref{funcation-N-pm}, using \eqref{estimate-of-Z-g2-h2-total-one-half} and Theorem \ref{two-functional-ub-by-norm},
we have
\ben \label{I-pm-rho-m-2}
|\mathcal{I}^{\rho}_{m,2}| \lesssim  \big(\int
B \mathrm{D}^{2}(\mu^{\f{1}{4}}_{*})
g_{*}^{2} h^{2}  \mathrm{d}V \big)^{\f{1}{2}}
\big(\int B \mu_{*}^{\f{1}{2}} \mathrm{D}^{2}(f) \mathrm{d}V \big)^{\f{1}{2}}
\lesssim |g|_{H^{s_{1}}}|h|_{H^{s_{2}}_{\gamma/2+s}} |f|_{\mathcal{L}^{s}_{\gamma/2}}.
\een
Recalling \eqref{definition-I-m3}, applying \eqref{estimate-of-X-g-h-f-total-one-half} and Lemma \ref{product-take-out}, we get
\beno
|\mathcal{I}^{\rho}_{m,3}|
\lesssim (|\frac{\mu^{\frac{1}{4}}}{1 - \rho \mu}g|_{L^{1}} |\frac{h}{1 - \rho \mu}|_{L^{2}_{\gamma/2+s}} + |\frac{\mu^{\frac{1}{256}}\mu^{\frac{1}{4}}}{1 - \rho \mu}g|_{H^{s_{1}}} |\frac{\mu^{\frac{1}{256}}}{1 - \rho \mu}h|_{H^{s_{2}}}) |f|_{L^{2}_{\gamma/2+s}} \lesssim |g|_{H^{s_{1}}}|h|_{H^{s_{2}}_{\gamma/2+s}} |f|_{L^{2}_{\gamma/2+s}}.
\eeno
Patching together the estimates of $\mathcal{I}^{\rho}_{m,i}$ for $1 \leq i \leq 3$, we get
\ben \label{up-Im-g-h-f}
 |  \langle  I^{\rho}_{m}(g,h), f \rangle   |  \lesssim |g|_{H^{s_{1}}}|h|_{H^{s_{2}}_{\gamma/2+s}} |f|_{\mathcal{L}^{s}_{\gamma/2}}.
\een

Recalling \eqref{definition-Ir-epsilon} and using \eqref{change-v-and-v-prime}, we get
\beno
 \langle  I^{\rho}_{r}(g,h), f \rangle   &=& \rho \int
B \mathrm{D}(\mu^{\f{1}{2}}) \mu^{\f{1}{2}}\mu^{\f{1}{2}}_{*}
(\frac{g}{1 - \rho \mu})_{*}(\frac{h}{1 - \rho \mu}) f^{\prime} \mathrm{d}V = \mathcal{I}^{\rho}_{r,1} + \mathcal{I}^{\rho}_{r,2},
\\
\mathcal{I}^{\rho}_{r,1} &\colonequals    & \rho \int
B \mathrm{D}(\mu^{\f{1}{2}}) \mu^{\f{1}{2}}\mu^{\f{1}{2}}_{*}
(\frac{g}{1 - \rho \mu})_{*}(\frac{h}{1 - \rho \mu}) \mathrm{D}(f^{\prime}) \mathrm{d}V,
\\
\mathcal{I}^{\rho}_{r,2} &\colonequals    & \rho \int
B \mathrm{D}(\mu^{\f{1}{2}}) \mu^{\f{1}{2}}\mu^{\f{1}{2}}_{*}
(\frac{g}{1 - \rho \mu})_{*}(\frac{h}{1 - \rho \mu}) f \mathrm{d}V.
\eeno
By Cauchy-Schwartz inequality and \eqref{change-v-and-v-star}, using \eqref{estimate-of-Z-g2-h2-total-one-half-with-mu} and Theorem \ref{two-functional-ub-by-norm},
we have
\beno
|\mathcal{I}^{\rho}_{r,1}| \lesssim
\cs{\int B \mathrm{D}^{2}(\mu^{\f{1}{2}}_{*}) \mu^{\f{1}{2}}
g^{2}h_{*}^{2}  \mathrm{d}V}
{\int B \mu_{*}^{\f{1}{2}} \mathrm{D}^{2}(f) \mathrm{d}V}
\lesssim   |\mu^{\frac{1}{256}}g|_{H^{s_{1}}}|\mu^{\frac{1}{256}}h|_{H^{s_{2}}} |f|_{\mathcal{L}^{s}_{\gamma/2}}
\eeno
By \eqref{change-v-and-v-star}, using the estimate \eqref{estimate-of-Y-g-h-f-total-one-half-with-mu} and Lemma \ref{product-take-out},
 we have
\beno
|\mathcal{I}^{\rho}_{r,2}| &=&| \rho
\int
B \mathrm{D}(\mu^{\f{1}{2}}_{*}) \mu^{\frac{1}{4}}
\frac{\mu^{\frac{1}{4}} g}{1 - \rho \mu} (\frac{\mu^{\f{1}{2}}h}{1 - \rho \mu})_{*} f_{*} \mathrm{d}V|
\\&\lesssim& |\frac{\mu^{\frac{1}{256}}\mu^{\frac{1}{4}} }{1 - \rho \mu}g|_{H^{s_{1}}} | \frac{\mu^{\frac{1}{256}}\mu^{\f{1}{2}}}{1 - \rho \mu}h|_{H^{s_{2}}}|\mu^{\frac{1}{256}}f|_{L^{2}} \lesssim |\mu^{\frac{1}{256}}g|_{H^{s_{1}}}|\mu^{\frac{1}{256}}h|_{H^{s_{2}}} |f|_{\mathcal{L}^{s}_{\gamma/2}}.
\eeno
Patching together the estimates of $\mathcal{I}^{\rho}_{r,1}$ and $\mathcal{I}^{\rho}_{r,2}$, using Lemma \ref{product-take-out}, we have
\ben \label{up-Ir-g-h-f}
 |  \langle  I^{\rho}_{r}(g,h), f \rangle   |  \lesssim  |\mu^{\frac{1}{256}}g|_{H^{s_{1}}}|\mu^{\frac{1}{256}}h|_{H^{s_{2}}} |f|_{\mathcal{L}^{s}_{\gamma/2}} \lesssim |g|_{H^{s_{1}}}|h|_{H^{s_{2}}_{\gamma/2+s}} |f|_{\mathcal{L}^{s}_{\gamma/2}}.
\een

Patching together \eqref{up-Im-g-h-f} and \eqref{up-Ir-g-h-f}, we finish the proof.
\end{proof}

Recalling \eqref{definition-Gamma-main}, patching together Proposition \ref{Q-pm-Ng-h-f} and  Proposition \ref{I-pm-rho-g-h-f}, we get
\begin{thm}\label{upper-bound-of-Gamma-2-m}  It holds that
\beno
 |  \langle  \Gamma_{2,m}^{\rho}(g,h), f \rangle   |  \lesssim \min \{|g|_{H^{2}}  |h|_{\mathcal{L}^{s}_{\gamma/2}}, |g|_{L^{2}} ( |h|_{\mathcal{L}^{s}_{\gamma/2}} +|h|_{H^{2+s}_{\gamma/2}} + |h|_{H^{\f{1}{2}}_{\gamma/2+s}} ) \} |f|_{\mathcal{L}^{s}_{\gamma/2}}.
\eeno
\end{thm}

\subsection{Some remark on derivative before operator} \label{derivative}
In this subsection, we address the issue of taking derivative $\partial_{\beta}$ w.r.t. variable $v$ of various Boltzmann type operators appeared in this article. The conclusion is that we can focus on derivatives of the involved functions $g, h, \varrho, f$ and safely ignore the derivatives on those good functions like $\mu, N, M$ appearing in the operators.
We begin with the following
 fact about Boltzmann type operator
\ben \label{distribute-derivatives}
\partial_{\beta}  \int B g_{*} h \varrho^{\prime}_{*} f^{\prime}  \mathrm{d}\sigma \mathrm{d}v_{*}
= \sum_{\beta_{1}+\beta_{2}+\beta_{3}+\beta_{4}=\beta} C^{\beta}_{\beta_{1},\beta_{2},\beta_{3},\beta_{4}} \int B (\partial_{\beta_{1}}g)_{*} (\partial_{\beta_{2}}h) (\partial_{\beta_{3}}\varrho)^{\prime}_{*} (\partial_{\beta_{4}}f)^{\prime} \mathrm{d}\sigma \mathrm{d}v_{*}.
\een
That is, when taking derivative w.r.t. $v$, we can use the binomial formula for the product $g h \varrho f$ first and then take integral.

Recalling \eqref{definition-Gamma-main} and \eqref{I-into-main-and-rest}, we have
$
\Gamma_{2,m}^{\rho}(g,h) = Q_{c} (Ng, h) + I^{\rho}_{m}(g,h) + I^{\rho}_{r}(g,h).
$
Therefore
\ben \label{derivatives-on-Gamma-2m}
\partial_{\beta}\Gamma_{2,m}^{\rho}(g,h) = \partial_{\beta} Q_{c} (Ng, h) + \partial_{\beta} I^{\rho}_{m}(g,h) +\partial_{\beta} I^{\rho}_{r}(g,h).
\een
Thanks to \eqref{distribute-derivatives}, we have
\beno
\partial_{\beta} Q_{c} (Ng, h) = \sum_{\beta_{1}+\beta_{2}+\beta_{3}=\beta} C^{\beta}_{\beta_{1},\beta_{2},\beta_{3}} Q_{c} (\partial_{\beta_{3}}N \partial_{\beta_{1}}g, \partial_{\beta_{2}}h).
\eeno
Note that $|\mu^{-\frac{1}{8}}\partial_{\beta_{3}}N |_{H^{2}_{7}} \lesssim 1$. Then by Corollary \ref{Q-up-bound-full-on-first-or-second} and Lemma \ref{product-take-out}, we have
\ben \label{partial-beta-3-on-N}
 |  \langle  Q_{c}(\partial_{\beta_{3}}N \partial_{\beta_{1}}g, \partial_{\beta_{2}}h), f \rangle   |  \lesssim  \min\{|\mu^{\frac{1}{8}}\partial_{\beta_{1}}g|_{H^{2}}  |\partial_{\beta_{2}}h|_{\mathcal{L}^{s}_{\gamma/2}} , |\mu^{\frac{1}{8}}\partial_{\beta_{1}}g|_{L^{2}}(|\partial_{\beta_{2}}h|_{\mathcal{L}^{s}_{\gamma/2}}
+|\partial_{\beta_{2}}h|_{H^{s+2}_{\gamma/2}})\} |f|_{\mathcal{L}^{s}_{\gamma/2}}.
\een
Note that if we use Proposition \ref{Q-pm-Ng-h-f} to estimate $ \langle  Q_{c} (N \partial_{\beta_{1}}g, \partial_{\beta_{2}}h), f \rangle  $, we will get the same result as \eqref{partial-beta-3-on-N}. Therefore we conclude that to estimate $ \langle  \partial_{\beta} Q_{c} (Ng, h), f \rangle  $, it suffices to estimate
$\langle Q_{c}(N \partial_{\beta_{1}}g, \partial_{\beta_{2}}h), f \rangle$ using Proposition \ref{Q-pm-Ng-h-f} for all $\beta_{1},\beta_{2}$ such that $\beta_{1}+\beta_{2} \leq \beta$. That is, we can safely regard $\beta_{3} \neq 0$ as $\beta_{3} = 0$.

Now we visit $\partial_{\beta} I^{\rho}_{m}(g,h)$. Recalling \eqref{definition-Im-epsilon} and using \eqref{distribute-derivatives}, we have
\beno
\partial_{\beta} I^{\rho}_{m}(g,h) &=& \sum_{\beta_{1}+\beta_{2}+\beta_{3}+\beta_{4}+\beta_{5}=\beta} C^{\beta}_{\beta_{1},\beta_{2},\beta_{3},\beta_{4},\beta_{5}} I^{\rho}_{m}(\partial_{\beta_{1}}g, \partial_{\beta_{2}}h; \beta_{3},\beta_{4},\beta_{5}),
\\
I^{\rho}_{m}(G, H; \beta_{3},\beta_{4},\beta_{5}) &\colonequals&     \int
 B  \mathrm{D}(\partial_{\beta_{3}}\mu^{\f{1}{2}})_{*})
 (G \partial_{\beta_{4}}(1 - \rho \mu)^{-1})^{\prime}_{*}(H \partial_{\beta_{5}}(1 - \rho \mu)^{-1})^{\prime}
\mathrm{d}\sigma \mathrm{d}v_{*}.
\eeno
Note that $I^{\rho}_{m}(G, H; 0,0,0) = I^{\rho}_{m}(G, H)$ when $\beta_{3}=\beta_{4}=\beta_{5}=0$.
For index $\beta \in \mathbb{N}^{3}$, there is a polynomial $P_{\beta}$ such that $\partial_{\beta}\mu^{\f{1}{2}} = \mu^{\f{1}{2}} P_{\beta}$. Observe
\ben \label{nice-decomposition}
 \mathrm{D}((\mu^{\f{1}{2}}P_{\beta})_{*}^{\prime})=
 \mathrm{D}((\mu^{\frac{1}{4}})^{\prime}_{*})\mathrm{D}((\mu^{\frac{1}{4}}P_{\beta})_{*}^{\prime})
 + \mu^{\frac{1}{4}}_{*}\mathrm{D}((\mu^{\frac{1}{4}}P_{\beta})_{*}^{\prime}) +  (\mu^{\frac{1}{4}}P_{\beta})_{*}\mathrm{D}((\mu^{\frac{1}{4}})^{\prime}_{*}).
\een
Using \eqref{change-v-and-v-prime} and \eqref{nice-decomposition}, setting $\mathcal{K}(G,H,f)\colonequals  (G \partial_{\beta_{4}}(1 - \rho \mu)^{-1})_{*}(H \partial_{\beta_{5}}(1 - \rho \mu)^{-1}) f^{\prime}$ for simplicity,
we have
\beno
 \langle  I^{\rho}_{m}(G, H; \beta_{3},\beta_{4},\beta_{5}), f \rangle
&=& \int
B \mathrm{D}((\mu^{\f{1}{2}} P_{\beta_{3}})^{\prime}_{*})
\mathcal{K}(G,H,f) \mathrm{d}V
\\&=& \mathcal{I}_{1}(G, H, f) + \mathcal{I}_{2}(G, H, f) + \mathcal{I}_{3}(G, H, f),
\\
\mathcal{I}_{1}(G, H, f) &\colonequals&     \int B
\mathrm{D}((\mu^{\frac{1}{4}})^{\prime}_{*})
\mathrm{D}((\mu^{\frac{1}{4}}P_{\beta_{3}})_{*}^{\prime})
\mathcal{K}(G,H,f) \mathrm{d}V,
\\
\mathcal{I}_{2}(G, H, f) &\colonequals&     \int B
\mu^{\frac{1}{4}}_{*} \mathrm{D}((\mu^{\frac{1}{4}}P_{\beta_{3}})_{*}^{\prime})
\mathcal{K}(G,H,f) \mathrm{d}V,
\\
\mathcal{I}_{3}(G, H, f) &\colonequals&     \int B
 (\mu^{\frac{1}{4}}P_{\beta_{3}})_{*}\mathrm{D}((\mu^{\frac{1}{4}})^{\prime}_{*})
\mathcal{K}(G,H,f) \mathrm{d}V.
\eeno
Recall $0 \leq \rho \leq \f{1}{2}(2 \pi)^{\frac{3}{2}}$. Therefore
 for any index $\beta \in \mathbb{N}^{3}$, there exists some constant $C_{\beta}$ such that
\ben \label{derivative-bounded}
|\partial_{\beta}(1 - \rho \mu)^{-1}| \leq C_{\beta}.
\een
Then by Cauchy-Schwartz inequality, using \eqref{change-v-and-v-star} and \eqref{change-v-and-v-prime}, we have
\beno
|\mathcal{I}_{1}(G, H, f)| \lesssim
\cs
{\int B \mathrm{D}^{2}((\mu^{\frac{1}{4}}P_{\beta_{3}})_{*})
 G^{2}_{*} H^{2}
\mathrm{d}V} {\int B
 \mathrm{D}^{2}(\mu^{\frac{1}{4}})
f^{2}_{*} \mathrm{d}V},
\eeno
which can be handled like in \eqref{to-X-g-h-Y-f} thanks to Remark \ref{mu-to-N-or-M-X-Y}. With the identity
$f^{\prime} = \mathrm{D}(f^{\prime}) + f$, both $\mathcal{I}_{2}(G, H, f)$ and $\mathcal{I}_{3}(G, H, f)$ can be split into two terms like \eqref{definition-I-m2}  and \eqref{definition-I-m3}. Then by following the estimate  of $ \mathcal{I}^{\rho}_{m,2}$ and $\mathcal{I}^{\rho}_{m,3}$ in Proposition  \ref{I-pm-rho-g-h-f}, thanks to
Remark \ref{mu-to-N-or-M-X-Y}, using \eqref{derivative-bounded}, we will get the same upper bound. In a word, $ \langle  I^{\rho}_{m}(G, H; \beta_{3},\beta_{4},\beta_{5}), f \rangle  $ shares the same upper bound  as \eqref{up-Im-g-h-f} for $ \langle  I^{\rho}_{m}(G, H), f \rangle  $.  Therefore we conclude that
  to estimate $\partial_{\beta} I^{\rho}_{m}(g,h)$, it suffices to consider $ I^{\rho}_{m}(\partial_{\beta_{1}}g, \partial_{\beta_{2}}h)$ for all $\beta_{1},\beta_{2}$ such that $\beta_{1}+\beta_{2}\leq\beta$. That is, we can safely regard general $\beta_{3},\beta_{4},\beta_{5}$ as $\beta_{3}=\beta_{4}=\beta_{5}=0$.  By nearly the same analysis, the same conclusion holds for
 $\partial_{\beta} I^{\rho}_{r}(g,h)$. Recalling \eqref{derivatives-on-Gamma-2m}, we arrive at the following remark.

\begin{rmk} \label{sufficient-to-consider}
In order to estimate $\langle  \partial_{\beta} \Gamma_{2,m}^{\rho}(g,h), f \rangle$, it suffices to consider $ \langle  \Gamma_{2,m}^{\rho}(\partial_{\beta_{1}}g, \partial_{\beta_{2}}h), f \rangle  $ for all $\beta_{1},\beta_{2}$ such that $\beta_{1}+\beta_{2} \leq \beta$.
\end{rmk}

Following the proof of upper bound estimates in Section \ref{linear}, \ref{bilinear} and \ref{trilinear},
we can go further to conclude that
\begin{rmk} \label{sufficient-to-consider-123}
To estimate upper bound concerning $\partial_{\beta} \Gamma_{2}^{\rho}(g,h)$, it suffices to consider $\Gamma_{2}^{\rho}(\partial_{\beta_{1}}g, \partial_{\beta_{2}}h)$ for all $\beta_{1},\beta_{2}$ such that $\beta_{1}+\beta_{2} \leq \beta$. To estimate upper bound concerning $\partial_{\beta} \Gamma_{3}^{\rho}(g,h,\varrho)$, it suffices to consider $\Gamma_{2}^{\rho}(\partial_{\beta_{1}}g, \partial_{\beta_{2}}h, \partial_{\beta_{3}}\varrho)$ for all $\beta_{1},\beta_{2},\beta_{3}$ such that $\beta_{1}+\beta_{2}+\beta_{3} \leq \beta$.
To estimate upper bound concerning $\partial_{\beta} \mathcal{L}^{\rho}g, \partial_{\beta} \mathcal{L}^{\rho}_{r} g$ and  $\partial_{\beta}\mathcal{C}^{\rho} g$, it suffices to consider $\mathcal{L}^{\rho} \partial_{\beta_{1}} g, \mathcal{L}^{\rho}_{r} \partial_{\beta_{1}} g$ and $\mathcal{C}^{\rho} \partial_{\beta_{1}} g$ for all $\beta_{1}$ such that $\beta_{1} \leq \beta$.
\end{rmk}

\subsection{The remaining operator $\Gamma_{2,r}^{\rho}$}
We first prepare some intermediate estimates.
\begin{lem} \label{optimal-weight-estimate} Let $s_{1},s_{2}\geq 0,s_{1}+s_{2}=\f{1}{2}$. The following four estimates are valid.
\ben \label{mu-star-g-star-difference-h2}
\int B \mu^{\f{1}{16}}_{*}\mathrm{D}^{2}(g_{*}) h^{2} \mathrm{d}V &\lesssim& |g|_{H^{s_{1}+1}}^{2} |h|_{H^{s_{2}}_{\gamma/2+s}}^{2}.
\\ \label{mu-star-g-star-difference-h2-f2-star}
\int B \mu^{\f{1}{16}}_{*}\mathrm{D}^{2}(g_{*}) f_{*}^{2} h^{2}  \mathrm{d}V &\lesssim& |g|_{H^{3}}^{2} |\mu^{\f{1}{128}}f|_{H^{s_{1}}}^{2}|h|_{H^{s_{2}}_{\gamma/2+s}}^{2}.
\\ \label{mu-star-g-star-difference-h2-f2-star-another}
\int B \mu^{\f{1}{16}}_{*} \mu^{\f{1}{16}} \mathrm{D}^{2}(g_{*}) f_{*}^{2} h^{2}  \mathrm{d}V &\lesssim& |\mu^{\f{1}{256}}g|_{H^{3}}^{2} |\mu^{\f{1}{256}}f|_{H^{s_{1}}}^{2}|\mu^{\f{1}{256}}h|_{H^{s_{2}}}^{2}.
\\ \label{key-estimate-only-star}
\int B  \mathrm{D}^{2}(g_{*}) \mu^{\f{1}{16}} h^{2} \mathrm{d}V &\lesssim&
|g|_{H^{s_{1}+1}_{\gamma/2+s}}^{2}|\mu^{\frac{1}{64}}h|_{H^{s_{2}}}^{2}.
\een
\end{lem}

We will only use \eqref{key-estimate-only-star} in this subsection. The proof of Lemma \ref{optimal-weight-estimate} is given in the Appendix \ref{appendix}. Appropriately revising the proof of Lemma \ref{optimal-weight-estimate}, we have
\begin{rmk}\label{star-prime-or-prime-version} Since $|v-v_{*}| \sim |v-v_{*}^{\prime}| \sim |v-v_{*}(\iota)|$, if we replace $\mu^{\f{1}{16}}_{*}$ with $(\mu^{\f{1}{16}})_{*}^{\prime}$ in
\eqref{mu-star-g-star-difference-h2}, the result is still valid. Since $|v-v_{*}| \sim |v^{\prime}-v_{*}|$, if we replace $h^{2}$ with $(h^{2})^{\prime}$ in
\eqref{mu-star-g-star-difference-h2-f2-star}, the result is still valid.
\end{rmk}

In the rest of the article, $s_{1},s_{2}$ are two constants verifying $s_{1},s_{2}\geq 0, s_{1}+s_{2}=\f{1}{2}$ unless otherwise specified.

Recalling \eqref{Gamma-remaining-into-three-terms}, there holds $\Gamma_{2,r}^{\rho} = \Gamma_{2,r,1}^{\rho} + \Gamma_{2,r,2}^{\rho} + \Gamma_{2,r,3}^{\rho}$. We will give estimates of $\Gamma_{2,r,1}^{\rho}, \Gamma_{2,r,2}^{\rho}, \Gamma_{2,r,3}^{\rho}$ in
Propositions \ref{upper-bound-of-Gamma-2-r-1-g-h-f},
\ref{upper-bound-of-Gamma-2-r-2-g-h-f}, \ref{upper-bound-of-Gamma-2-r-3} respectively.

\begin{prop}\label{upper-bound-of-Gamma-2-r-1-g-h-f}  It holds that
\beno
 |  \langle  \Gamma_{2,r,1}^{\rho}(g, h), f \rangle   |  \lesssim  \min \{|g|_{H^{2}}  |h|_{\mathcal{L}^{s}_{\gamma/2}}, |g|_{L^{2}} ( |h|_{\mathcal{L}^{s}_{\gamma/2}} +|h|_{H^{2+s}_{\gamma/2}} + |h|_{H^{\f{1}{2}}_{\gamma/2+s}} ) \} |f|_{\mathcal{L}^{s}_{\gamma/2}}.
\eeno
\end{prop}
\begin{proof} Recalling \eqref{definition-Gamma-remaining-1}, we have
\ben \label{inner-product-remaining-1}
 \langle  \Gamma_{2,r,1}^{\rho}(g, h), f \rangle   = \int
B  (Nh)(Ng)_{*}\mathrm{D}((N^{-1}f)^{\prime}) (M^{\prime} + M^{\prime}_{*}) \mathrm{d}V.
\een
We need to make rearrangement in order to use previous results. By the following identity
\ben \label{usual-decomposition-into-N-and-f}
\mathrm{D}((N^{-1}f)^{\prime}) = \mathrm{D}((N^{-1})^{\prime})f^{\prime} + N^{-1} \mathrm{D}(f^{\prime}),
\een
we have
\ben \label{G-2-r-1-into-1-2-3-4}
 \langle  \Gamma_{2,r,1}^{\rho}(g, h), f \rangle   = \mathcal{T}_{1} + \mathcal{T}_{2} + \mathcal{T}_{3} + \mathcal{T}_{4},
\\ \nonumber
\mathcal{T}_{1}\colonequals      \int B  (Ng)_{*}(Nh)(Mf)^{\prime}\mathrm{D}((N^{-1})^{\prime}) \mathrm{d}V,
\quad 
\mathcal{T}_{2}\colonequals      \int B  (Ng)_{*}(Nh)f^{\prime}M^{\prime}_{*}\mathrm{D}((N^{-1})^{\prime}) \mathrm{d}V,
\\ \nonumber
\mathcal{T}_{3}\colonequals      \int B  (Ng)_{*} h \mathrm{D}(f^{\prime})M^{\prime} \mathrm{d}V, \quad
\mathcal{T}_{4}\colonequals      \int B  (Ng)_{*} h \mathrm{D}(f^{\prime})M^{\prime}_{*} \mathrm{d}V.
\een

{\it {Estimate of $\mathcal{T}_{1}$.}}
By \eqref{change-v-and-v-prime} and recalling \eqref{definition-I-epsilon}, using Proposition \ref{I-pm-rho-g-h-f} and Lemma \ref{product-take-out},
we have
\ben \nonumber
|\mathcal{T}_{1}| &=& | \int
B  (Ng)^{\prime}_{*}(Nh)^{\prime}Mf \mathrm{D}(N^{-1})
\mathrm{d}V |
= |  \langle  I^{\rho}(g, h), Mf \rangle   |
\\ \label{U-1-result} &\lesssim& |g|_{H^{s_{1}}}|h|_{H^{s_{2}}_{\gamma/2+s}} |Mf|_{\mathcal{L}^{s}_{\gamma/2}} \lesssim |g|_{H^{s_{1}}}|h|_{H^{s_{2}}_{\gamma/2+s}} |f|_{\mathcal{L}^{s}_{\gamma/2}}.
\een

{\it {Estimate of $\mathcal{T}_{2}$.}} Using the identity $M^{\prime}_{*} = \mathrm{D}(M^{\prime}_{*}) + M_{*}$, we have
$\mathcal{T}_{2} = \mathcal{T}_{2,1} + \mathcal{T}_{2,2}$ where
\beno
\mathcal{T}_{2,1} \colonequals     \int B  (MNg)_{*}(Nh)f^{\prime}\mathrm{D}((N^{-1})^{\prime}) \mathrm{d}V,
\quad
\mathcal{T}_{2,2} \colonequals     \int B  (Ng)_{*}(Nh)f^{\prime}\mathrm{D}(M^{\prime}_{*})\mathrm{D}((N^{-1})^{\prime}) \mathrm{d}V.
\eeno
Recalling \eqref{definition-I-epsilon}, using Proposition \ref{I-pm-rho-g-h-f} and Lemma \ref{product-take-out},
we have
\ben \label{U-2-1-result}
|\mathcal{T}_{2,1}| = | \langle  I^{\rho}(M g, h), f \rangle  | \lesssim |g|_{H^{s_{1}}}|h|_{H^{s_{2}}_{\gamma/2+s}} |f|_{\mathcal{L}^{s}_{\gamma/2}}.
\een
By Cauchy-Schwartz inequality, have
\ben \label{U-2-2-Cauchy-Schwartz}
|\mathcal{T}_{2,2}| \lesssim
\cs{\int B  |(Ng)_{*}(Nh)|^{2}\mathrm{D}^{2}(N^{-1})
\mathrm{d}V}{\int B  (f^{2})^{\prime}\mathrm{D}^{2}(M_{*})
\mathrm{d}V}.
\een
By \eqref{change-v-and-v-star} and \eqref{change-v-and-v-prime}, we observe
\ben  \label{f-square-with-M-difference}
\int B  (f^{2})^{\prime}\mathrm{D}^{2}(M_{*})
\mathrm{d}V = \int B  f^{2}_{*}\mathrm{D}^{2}(M) \mathrm{d}V
= \mathcal{N}(f,M^{\f{1}{2}}) \lesssim |f|_{\mathcal{L}^{s}_{\gamma/2}}^{2},
\een
where we use Theorem \ref{two-functional-ub-by-norm} and Remark \ref{mu-to-N-or-M-functional-N}.
Similar to \eqref{N-power-minus-1-difference}, we have
\beno
\mathrm{D}(N^{-1})N_{*}N = \big(\mathrm{D}(\mu^{\f{1}{2}}_{*}) - \rho \mathrm{D}(\mu^{\f{1}{2}}) \mu^{\f{1}{2}}\mu^{\f{1}{2}}_{*}\big) (\frac{1}{1 - \rho \mu})_{*}(\frac{1}{1 - \rho \mu}),
\eeno
which yields
$
|\mathrm{D}(N^{-1}) N_{*}N|^{2} \lesssim
\mathrm{D}^{2}(\mu^{\f{1}{2}}_{*}) + \mu \mu_{*} \mathrm{D}^{2}(\mu^{\f{1}{2}})
$
and thus
\ben \label{U-2-2-part-g2-h2}
&& \int B  |(Ng)_{*}(Nh)|^{2}\mathrm{D}^{2}(N^{-1}) \mathrm{d}V
\\ \nonumber  &\lesssim&
\int B   \mathrm{D}^{2}(\mu^{\f{1}{2}}_{*}) g_{*}^{2}h^{2} \mathrm{d}V
+ \int B   \mathrm{D}^{2}(\mu^{\f{1}{2}}) g_{*}^{2}h^{2}\mu \mu_{*} \mathrm{d}V
\lesssim|g|_{H^{s_{1}}}^{2}|h|_{H^{s_{2}}_{\gamma/2+s}}^{2},
\een
where we use \eqref{change-v-and-v-star} for the latter integral and then get the final estimate by using \eqref{estimate-of-Z-g2-h2-total-one-half}.
Plugging \eqref{f-square-with-M-difference} and \eqref{U-2-2-part-g2-h2} into \eqref{U-2-2-Cauchy-Schwartz}, we get
\ben \label{U-2-2-result}
|\mathcal{T}_{2,2}| \lesssim |g|_{H^{s_{1}}}|h|_{H^{s_{2}}_{\gamma/2+s}} |f|_{\mathcal{L}^{s}_{\gamma/2}}.
\een
Patching together \eqref{U-2-1-result} and \eqref{U-2-2-result}, recalling $\mathcal{T}_{2} = \mathcal{T}_{2,1} + \mathcal{T}_{2,2}$, we have
\ben \label{U-2-result}
|\mathcal{T}_{2}| \lesssim |g|_{H^{s_{1}}}|h|_{H^{s_{2}}_{\gamma/2+s}} |f|_{\mathcal{L}^{s}_{\gamma/2}}.
\een

{\it {Estimate of $\mathcal{T}_{3}$.}} Using the identity $M^{\prime} = \mathrm{D}(M^{\prime}) + M$, we have $\mathcal{T}_{3} = \mathcal{T}_{3,1} + \mathcal{T}_{3,2}$ where
\beno
\mathcal{T}_{3,1}\colonequals      \int B  (Ng)_{*} M h \mathrm{D}(f^{\prime}) \mathrm{d}V, \quad
\mathcal{T}_{3,2}\colonequals      \int B  (Ng)_{*} h \mathrm{D}(f^{\prime}) \mathrm{D}(M^{\prime}) \mathrm{d}V.
\eeno
Observe $\mathcal{T}_{3,1} =   \langle  Q_{c}(Ng, Mh), f \rangle  $. Then  Corollary \ref{Q-up-bound-full-on-first-or-second} and Lemma \ref{product-take-out} yield
\ben \nonumber
|\mathcal{T}_{3,1}| &\lesssim& \min\{|N g|_{H^{2}_{7}}  |M h|_{\mathcal{L}^{s}_{\gamma/2}} , |N g|_{L^{2}_{7}}(|M h|_{\mathcal{L}^{s}_{\gamma/2}}+|M h|_{H^{s+2}_{\gamma/2}})\} |f|_{\mathcal{L}^{s}_{\gamma/2}}
\\ \label{U-3-1-result} &\lesssim& \min\{|g|_{H^{2}}  |h|_{\mathcal{L}^{s}_{\gamma/2}} , |g|_{L^{2}}(|h|_{\mathcal{L}^{s}_{\gamma/2}}+| h|_{H^{s+2}_{\gamma/2}})\} |f|_{\mathcal{L}^{s}_{\gamma/2}}.
\een
By Cauchy-Schwartz inequality, \eqref{M-rho-mu-N} and \eqref{change-v-and-v-star}, have
\ben \label{U-3-2-result}
|\mathcal{T}_{3,2}| \lesssim
\cs{ \int B   \mu^{\f12} g^{2} h_{*}^{2}\mathrm{D}^{2}(M_{*})
\mathrm{d}V}{\int B  \mu^{\f12}_{*}\mathrm{D}^{2}(f)
\mathrm{d}V}
 \lesssim |\mu^{\frac{1}{256}}g|_{H^{s_{1}}}|\mu^{\frac{1}{256}}h|_{H^{s_{2}}} |f|_{\mathcal{L}^{s}_{\gamma/2}}.
\een
where we use \eqref{estimate-of-Z-g2-h2-total-one-half-with-mu}, Remark \ref{mu-to-N-or-M-X-Y} and Theorem \ref{two-functional-ub-by-norm}.
Patching together \eqref{U-3-1-result} and \eqref{U-3-2-result}, we have
\ben \label{U-3-result}
|\mathcal{T}_{3}| \lesssim \min\{|g|_{H^{2}}  |h|_{\mathcal{L}^{s}_{\gamma/2}} , |g|_{L^{2}}(|h|_{\mathcal{L}^{s}_{\gamma/2}}+| h|_{H^{s+2}_{\gamma/2}})\} |f|_{\mathcal{L}^{s}_{\gamma/2}}.
\een

{\it {Estimate of $\mathcal{T}_{4}$.}} Using the identity $M^{\prime}_{*} = \mathrm{D}(M^{\prime}_{*}) + M_{*}$, we have $\mathcal{T}_{4} = \mathcal{T}_{4,1} + \mathcal{T}_{4,2}$ where
\ben \label{U42-definition}
\mathcal{T}_{4,1} \colonequals      \int B  (MNg)_{*} h \mathrm{D}(f^{\prime}) \mathrm{d}V, \quad
\mathcal{T}_{4,2} \colonequals      \int B  (Ng)_{*} h \mathrm{D}(f^{\prime})\mathrm{D}(M^{\prime}_{*}) \mathrm{d}V.
\een
Observe $\mathcal{T}_{4,1} =   \langle  Q_{c}(MNg, h), f \rangle  $. Then  Corollary \ref{Q-up-bound-full-on-first-or-second} and Lemma \ref{product-take-out} yield
\beno
|\mathcal{T}_{4,1}|  \lesssim \min\{|g|_{H^{2}}  |h|_{\mathcal{L}^{s}_{\gamma/2}} , |g|_{L^{2}}(|h|_{\mathcal{L}^{s}_{\gamma/2}}+| h|_{H^{s+2}_{\gamma/2}})\} |f|_{\mathcal{L}^{s}_{\gamma/2}}.
\eeno
Note that $\mathcal{T}_{4,2}$ has the same structure as $\mathcal{I}^{\rho}_{m,2}$ defined in
\eqref{definition-I-m2}. Then similar to \eqref{I-pm-rho-m-2}, we get
\ben \label{result-U42}
|\mathcal{T}_{4,2}| \lesssim |g|_{H^{s_{1}}}|h|_{H^{s_{2}}_{\gamma/2+s}} |f|_{\mathcal{L}^{s}_{\gamma/2}}.
\een
Patching together the previous two estimates,  we get
\ben \label{U-4-result}
|\mathcal{T}_{4}| \lesssim \min\{|g|_{H^{2}}  |h|_{\mathcal{L}^{s}_{\gamma/2}} , |g|_{L^{2}}(|h|_{\mathcal{L}^{s}_{\gamma/2}}+| h|_{H^{s+2}_{\gamma/2}}+|h|_{H^{\f{1}{2}}_{\gamma/2+s}})\} |f|_{\mathcal{L}^{s}_{\gamma/2}}.
\een

Patching together \eqref{U-1-result}, \eqref{U-2-result}, \eqref{U-3-result} and \eqref{U-4-result},
recalling \eqref{G-2-r-1-into-1-2-3-4}, we finish the proof.
\end{proof}

\begin{prop}\label{upper-bound-of-Gamma-2-r-2-g-h-f}  It holds that
\beno
 |  \langle  \Gamma_{2,r,2}^{\rho}(g, h), f \rangle   |  \lesssim   \min\{|g|_{H^{2}}  |h|_{\mathcal{L}^{s}_{\gamma/2}} , |g|_{L^{2}}(|h|_{\mathcal{L}^{s}_{\gamma/2}}+| h|_{H^{s+2}_{\gamma/2}}+| h|_{H^{\frac{3}{2}}_{\gamma/2+s}})\} |f|_{\mathcal{L}^{s}_{\gamma/2}}.
\eeno
\end{prop}

\begin{proof}
Recalling the definition of $\Gamma_{2,r,2}^{\rho}$ in \eqref{definition-Gamma-remaining-2}, we have
\ben \label{inner-product-remaining-3}
 \langle  \Gamma_{2,r,2}^{\rho}(g, h), f \rangle   = \int
B  (Nh)(Ng)^{\prime}_{*}\mathrm{D}((N^{-1}f)^{\prime})(M_{*}-M^{\prime})
\mathrm{d}V.
\een
Using the identity \eqref{usual-decomposition-into-N-and-f}, we get
\ben  \nonumber
 \langle  \Gamma_{2,r,2}^{\rho}(g, h), f \rangle   &=&  \int B (Ng)^{\prime}_{*}(Nh)f^{\prime}(M_{*}-M^{\prime})\mathrm{D}((N^{-1})^{\prime}) \mathrm{d}V
\\ \label{G-2-r-2-into-1-2-3-4} &&+ \int B (Ng)^{\prime}_{*} h \mathrm{D}(f^{\prime})(M_{*}-M^{\prime}) \mathrm{d}V
 = \mathcal{U}_{1} + \mathcal{U}_{2} + \mathcal{U}_{3} + \mathcal{U}_{4},
\\ \label{G-2-r-2-U-1}
\mathcal{U}_{1} &\colonequals&      \int B  (Ng)^{\prime}_{*}(Nh)f^{\prime}M_{*}\mathrm{D}((N^{-1})^{\prime}) \mathrm{d}V,
\\ \label{G-2-r-2-U-2}
\mathcal{U}_{2} &\colonequals&      \int B  (Ng)^{\prime}_{*}(Nh) (Mf)^{\prime}\mathrm{D}(N^{-1}) \mathrm{d}V,
\\ \label{G-2-r-2-U-3}
\mathcal{U}_{3} &\colonequals&      \int B  (Ng)^{\prime}_{*} (h-h^{\prime}) \mathrm{D}(f^{\prime})(M_{*}-M^{\prime}) \mathrm{d}V,
\\ \label{G-2-r-2-U-4}
\mathcal{U}_{4} &\colonequals&      \int B  (Ng)^{\prime}_{*} h^{\prime} \mathrm{D}(f^{\prime})(M_{*}-M^{\prime}) \mathrm{d}V.
\een

{\it {Estimate of $\mathcal{U}_{1}$.}} Using \eqref{change-v-and-v-prime} to \eqref{G-2-r-2-U-1}, we have
\ben \label{a-very-hard-term}
\mathcal{U}_{1} = \int B  (Ng)_{*}(Nh)^{\prime}f M^{\prime}_{*}\mathrm{D}(N^{-1}) \mathrm{d}V.
\een
Recalling \eqref{equilibrium-rho-not-vanish}, there holds $M = \mu^{\f{1}{2}} N.$
Similar to \eqref{N-power-minus-1-difference}, we have
\ben \label{N-prime-M-prime-star}
N^{\prime}M^{\prime}_{*}\mathrm{D}(N^{-1}) = \big(\mathrm{D}(\mu^{\f{1}{2}}_{*}) - \rho \mathrm{D}(\mu^{\f{1}{2}}) \mu^{\f{1}{2}}\mu^{\f{1}{2}}_{*}\big) N^{\prime}_{*}(\frac{1}{1 - \rho \mu})^{\prime}.
\een
Plugging \eqref{N-prime-M-prime-star} into \eqref{a-very-hard-term},
using $N^{\prime}_{*} = \mathrm{D}(N^{\prime}_{*}) + N_{*}$ and $(\frac{h}{1 - \rho \mu})^{\prime} =
\mathrm{D}((\frac{h}{1 - \rho \mu})^{\prime}) +\frac{h}{1 - \rho \mu}$, we have
\ben \label{U1-into-1-2-3}
\mathcal{U}_{1} &=& \mathcal{U}_{1,1} + \mathcal{U}_{1,2} + \mathcal{U}_{1,3},
\\ \nonumber
\mathcal{U}_{1,1} &\colonequals&      \int B  (Ng)_{*} (\frac{h}{1 - \rho \mu})^{\prime} f \mathrm{D}(N^{\prime}_{*}) \big(\mathrm{D}(\mu^{\f{1}{2}}_{*}) - \rho \mathrm{D}(\mu^{\f{1}{2}}) \mu^{\f{1}{2}}\mu^{\f{1}{2}}_{*}\big) \mathrm{d}V,
\\ \nonumber
\mathcal{U}_{1,2} &\colonequals&      \int B  (N^{2}g)_{*} \mathrm{D}((\frac{h}{1 - \rho \mu})^{\prime}) f \big(\mathrm{D}(\mu^{\f{1}{2}}_{*}) - \rho \mathrm{D}(\mu^{\f{1}{2}}) \mu^{\f{1}{2}}\mu^{\f{1}{2}}_{*}\big)
\mathrm{d}V,
\\ \nonumber
\mathcal{U}_{1,3} &\colonequals&      \int B  (N^{2}g)_{*} \frac{h}{1 - \rho \mu} f \big(\mathrm{D}(\mu^{\f{1}{2}}_{*}) - \rho \mathrm{D}(\mu^{\f{1}{2}}) \mu^{\f{1}{2}}\mu^{\f{1}{2}}_{*}\big)
\mathrm{d}V.
\een

{\it {Estimate of $\mathcal{U}_{1,1}$.}}
By Cauchy-Schwartz inequality and \eqref{M-rho-mu-N}, we have
\beno
|\mathcal{U}_{1,1}| &\lesssim&
\cs{
\int B  (\mu^{\f{1}{2}}g)^{2}_{*} (h^{2})^{\prime} \mathrm{D}^{2}(N_{*})
\mathrm{d}V}{\int
B   f^{2} \mathrm{D}^{2}(\mu^{\f{1}{2}}_{*})
\mathrm{d}V + \int
B  \mu \mu_{*} f^{2} \mathrm{D}^{2}(\mu^{\f{1}{2}})
\mathrm{d}V}.
\eeno
By \eqref{estimate-of-Z-g2-h2-total-one-half} and Remark \ref{mu-to-N-or-M-X-Y}, we have
\beno
\int B  (\mu^{\f{1}{2}}g)^{2}_{*} (h^{2})^{\prime} \mathrm{D}^{2}(N_{*})
\mathrm{d}V \lesssim |g|_{H^{s_{1}}}^{2}|h|_{H^{s_{2}}_{\gamma/2+s}}^{2}.
\eeno
By \eqref{change-v-and-v-star} and Theorem \ref{two-functional-ub-by-norm}, we get
\beno
\int
B   f^{2} \mathrm{D}^{2}(\mu^{\f{1}{2}}_{*})
\mathrm{d}V = \int
B   f_{*}^{2} \mathrm{D}^{2}(\mu^{\f{1}{2}})
\mathrm{d}V =  \mathcal{N}(f,\mu^{\f{1}{2}}) \lesssim |f|_{\mathcal{L}^{s}_{\gamma/2}}^{2}.
\eeno
By \eqref{change-v-and-v-star} and the estimate \eqref{estimate-of-Z-g2-h2-total-one-half-with-mu}, we get
\beno
\int
B  \mu \mu_{*} f^{2} \mathrm{D}^{2}(\mu^{\f{1}{2}})
\mathrm{d}V =\int
B  \mu \mu_{*} f_{*}^{2} \mathrm{D}^{2}(\mu^{\f{1}{2}}_{*})
\mathrm{d}V \lesssim |\mu^{\f{1}{2}}f|_{L^{2}}^{2}.
\eeno
Patching together the previous three estimates, we have
\ben \label{result-U-1-1}
|\mathcal{U}_{1,1}| \lesssim |g|_{H^{s_{1}}}|h|_{H^{s_{2}}_{\gamma/2+s}}|f|_{\mathcal{L}^{s}_{\gamma/2}}.
\een

{\it {Estimate of $\mathcal{U}_{1,2}$.}}
By Cauchy-Schwartz inequality and \eqref{M-rho-mu-N}, we have
\ben \label{U12-by-Cauchy-Schwartz-inequality}
|\mathcal{U}_{1,2}| \lesssim \big(\int
B  \mu_{*}g^{2}_{*} f^{2} \big(\mathrm{D}^{2}(\mu^{\f{1}{2}}_{*}) + \mathrm{D}^{2}(\mu^{\f{1}{2}}) \mu \mu_{*}\big)
\mathrm{d}V \big)^{\f{1}{2}}
 \big(\int
B  \mu_{*} \mathrm{D}^{2}(\frac{h}{1 - \rho \mu})
\mathrm{d}V \big)^{\f{1}{2}}.
\een
By \eqref{estimate-of-Z-g2-h2-total-one-half}, we have
\ben \label{term1-on-g2-f2-by-norm}
\int B  \mu_{*}g^{2}_{*} f^{2} \mathrm{D}^{2}(\mu^{\f{1}{2}}_{*}) \mathrm{d}V
\lesssim |g|_{H^{\f{1}{2}}}^{2}|f|_{L^{2}_{\gamma/2+s}}^{2}.
\een
By \eqref{change-v-and-v-star} and the estimate \eqref{estimate-of-Z-g2-h2-total-one-half-with-mu}, we get
\ben \label{term2-on-g2-f2-by-norm}
\int B  \mu_{*}g^{2}_{*} f^{2}  \mathrm{D}^{2}(\mu^{\f{1}{2}}) \mu \mu_{*}
\mathrm{d}V = \int B  \mu g^{2} f_{*}^{2}  \mathrm{D}^{2}(\mu^{\f{1}{2}}_{*}) \mu \mu_{*}
\mathrm{d}V
\lesssim |\mu^{\frac{1}{256}} g|_{H^{\f{1}{2}}}^{2}|\mu^{\frac{1}{256}}f|_{L^{2}}^{2}.
\een
Observe
$
\mathrm{D}((\frac{h}{1 - \rho \mu})^{\prime}) = \frac{\mathrm{D}(h^{\prime})}{(1 - \rho \mu)^{\prime}} - \frac{\rho\mathrm{D}(\mu)h}{(1 - \rho \mu)^{\prime}(1 - \rho \mu)}
$
and thus
\ben \label{difference-h-over}
\mathrm{D}^{2}(\frac{h}{1 - \rho \mu}) \lesssim \mathrm{D}^{2}(h) + \mathrm{D}^{2}(\mu)h^{2},
\een
which gives
\ben \label{term3-on-h2-by-norm}
\int B  \mu_{*} \mathrm{D}^{2}(\frac{h}{1 - \rho \mu})
\mathrm{d}V
 \lesssim \int B  \mu_{*} \mathrm{D}^{2}(h) \mathrm{d}V +
\int B  \mu \mathrm{D}^{2}(\mu_{*})h_{*}^{2} \mathrm{d}V \lesssim |h|_{\mathcal{L}^{s}_{\gamma/2}}^{2},
\een
where we use \eqref{change-v-and-v-star}, Theorem \ref{two-functional-ub-by-norm} and the estimate \eqref{estimate-of-Z-g2-h2-total-one-half-with-mu}.
Plugging \eqref{term1-on-g2-f2-by-norm}, \eqref{term2-on-g2-f2-by-norm} and
\eqref{term3-on-h2-by-norm} into \eqref{U12-by-Cauchy-Schwartz-inequality}, we get
\ben \label{result-U-1-2}
|\mathcal{U}_{1,2}| \lesssim |g|_{H^{\f{1}{2}}} |h|_{\mathcal{L}^{s}_{\gamma/2}}|f|_{\mathcal{L}^{s}_{\gamma/2}}.
\een

{\it {Another estimate of $\mathcal{U}_{1,2}$.}}
We now give another estimate of $\mathcal{U}_{1,2}$ to put more regularity on $h$. By Cauchy-Schwartz inequality and \eqref{M-rho-mu-N}, we have
\ben \label{U12-by-Cauchy-Schwartz-inequality-2}
|\mathcal{U}_{1,2}| \lesssim \big(\int
B  \mu_{*} f^{2} \big(\mathrm{D}^{2}(\mu^{\f{1}{2}}_{*}) +
\mathrm{D}^{2}(\mu^{\f{1}{2}}) \mu \mu_{*}\big)
\mathrm{d}V \big)^{\f{1}{2}}
\big(\int
B   \mu_{*} g^{2}_{*}\mathrm{D}^{2}(\frac{h}{1 - \rho \mu})
\mathrm{d}V \big)^{\f{1}{2}}.
\een
By \eqref{estimate-of-Z-g2-h2-total-one-half}, we have
$
\int B  \mu_{*} f^{2} \mathrm{D}^{2}(\mu^{\f{1}{2}}_{*}) \mathrm{d}V
\lesssim |f|_{L^{2}_{\gamma/2+s}}^{2}.
$
By \eqref{change-v-and-v-star} and the estimate \eqref{estimate-of-Z-g2-h2-total-one-half-with-mu}, we get
\beno
\int B   \mu_{*} f^{2}  \mathrm{D}^{2}(\mu^{\f{1}{2}}) \mu \mu_{*}
\mathrm{d}V = \int B  \mu f_{*}^{2}  \mathrm{D}^{2}(\mu^{\f{1}{2}}_{*}) \mu \mu_{*}
\mathrm{d}V
\lesssim |\mu^{\frac{1}{256}}f|_{L^{2}}^{2}.
\eeno
By \eqref{difference-h-over} and \eqref{change-v-and-v-star}, using \eqref{key-estimate-only-star} and \eqref{estimate-of-Z-g2-h2-total-one-half-with-mu},
we have
\beno
\int B   \mu_{*} g^{2}_{*}\mathrm{D}^{2}(\frac{h}{1 - \rho \mu})
\mathrm{d}V
&\lesssim& \int B  \mu g^{2} \mathrm{D}^{2}(h_{*}) \mathrm{d}V +
\int B  \mu g^{2} \mathrm{D}^{2}(\mu_{*}) h^{2}_{*} \mathrm{d}V
\\&\lesssim& |\mu^{\frac{1}{64}} g|_{H^{s_{1}}}^{2}|h|_{H^{s_{2}+1}_{\gamma/2+s}}^{2} + |\mu^{\frac{1}{256}} g|_{H^{s_{1}}}^{2}|\mu^{\frac{1}{256}}h|_{H^{s_{2}}}^{2}.
\eeno
Plugging the previous three estimates into \eqref{U12-by-Cauchy-Schwartz-inequality-2}, we get
\ben \label{result-U-1-2-another}
|\mathcal{U}_{1,2}| \lesssim |\mu^{\frac{1}{256}} g|_{H^{s_{1}}} |h|_{H^{s_{2}+1}_{\gamma/2+s}}|f|_{L^{2}_{\gamma/2+s}}.
\een

{\it {Estimate of $\mathcal{U}_{1,3}$.}}
Now we set to estimate $\mathcal{U}_{1,3}$. By \eqref{change-v-and-v-star},
using \eqref{estimate-of-X-g-h-f-total-one-half} and \eqref{estimate-of-Y-g-h-f-total-one-half-with-mu}, we have
\ben \nonumber
|\mathcal{U}_{1,3}| &\leq& |\int B  (N^{2}g)_{*} \frac{h}{1 - \rho \mu} f \mathrm{D}(\mu^{\f{1}{2}}_{*})
\mathrm{d}V|
\\ \label{result-U-1-3}
&&+|\rho \int B  N^{2}g (\frac{\mu^{\f{1}{2}} h f}{1 - \rho \mu})_{*} \mathrm{D}(\mu^{\f{1}{2}}_{*}) \mu^{\f{1}{2}}
\mathrm{d}V|
 \lesssim |g|_{H^{s_{1}}}|h|_{H^{s_{2}}_{\gamma/2+s}} |f|_{L^{2}_{\gamma/2+s}}.
\een
Patching together \eqref{result-U-1-1}, \eqref{result-U-1-2} and \eqref{result-U-1-3}, recalling \eqref{U1-into-1-2-3}, we get
\ben \label{result-U-1}
|\mathcal{U}_{1}| \lesssim |g|_{H^{\f{1}{2}}} |h|_{\mathcal{L}^{s}_{\gamma/2}} |f|_{\mathcal{L}^{s}_{\gamma/2}}.
\een
Patching together \eqref{result-U-1-1}, \eqref{result-U-1-2-another} and \eqref{result-U-1-3}, recalling \eqref{U1-into-1-2-3}, we get
\ben \label{result-U-1-another}
|\mathcal{U}_{1}| \lesssim |g|_{H^{s_{1}}} |h|_{H^{s_{2}+1}_{\gamma/2+s}}|f|_{\mathcal{L}^{s}_{\gamma/2}}.
\een

{\it {Estimate of $\mathcal{U}_{2}$.}} Recalling \eqref{G-2-r-2-U-2} and $M = \mu^{\f{1}{2}} N$, using \eqref{change-v-and-v-prime}, we have
\ben \label{U-2-into-U21-and-U22}
\mathcal{U}_{2}=  \int
B  (Ng)^{\prime}_{*} h (\mu^{\f{1}{2}}f)^{\prime} \mathrm{D}(N^{\prime})
\mathrm{d}V
=  \int
B  (Ng)_{*} h^{\prime} \mu^{\f{1}{2}}f \mathrm{D}(N)
\mathrm{d}V = \mathcal{U}_{2,1} + \mathcal{U}_{2,2},
 \\  \nonumber
\mathcal{U}_{2,1} \colonequals      \int B  (Ng)_{*} \mathrm{D}(h^{\prime}) \mu^{\f{1}{2}}f \mathrm{D}(N) \mathrm{d}V,
\quad
\mathcal{U}_{2,2} \colonequals      \int B  (Ng)_{*} h \mu^{\f{1}{2}}f \mathrm{D}(N) \mathrm{d}V.
\een
By Cauchy-Schwartz inequality, \eqref{M-rho-mu-N} and \eqref{change-v-and-v-star},
we have
\ben \nonumber
|\mathcal{U}_{2,1}| &\lesssim& \cs{\int
B  \mu^{\f12}_{*} \mathrm{D}^{2}(h)
\mathrm{d}V} {\int
B  \mu^{\f12} g^{2} (\mu f^{2})_{*} \mathrm{D}^{2}(N_{*})
\mathrm{d}V}
\\ \label{result-U-2-1} &\lesssim&  |\mu^{\frac{1}{256}}g|_{H^{\f{1}{2}}} |h|_{\mathcal{L}^{s}_{\gamma/2}} |\mu^{\frac{1}{256}}f|_{L^{2}}  \lesssim |g|_{H^{\f{1}{2}}} |h|_{\mathcal{L}^{s}_{\gamma/2}} |f|_{\mathcal{L}^{s}_{\gamma/2}},
\een
where we use Theorem \ref{two-functional-ub-by-norm}, the estimate \eqref{estimate-of-Z-g2-h2-total-one-half-with-mu} and Remark \ref{mu-to-N-or-M-X-Y}.

We now give another estimate of $\mathcal{U}_{2,1}$.
By Cauchy-Schwartz inequality, \eqref{M-rho-mu-N} and \eqref{change-v-and-v-star},
we have
\ben \nonumber
|\mathcal{U}_{2,1}| &\lesssim& \cs{\int
B  \mu^{\f12} g^{2} \mathrm{D}^{2}(h_{*})
\mathrm{d}V} {\int
B  \mu^{\f12}  (\mu f^{2})_{*} \mathrm{D}^{2}(N_{*})
\mathrm{d}V}
\\ \label{result-U-2-1-another} &\lesssim& |\mu^{\frac{1}{64}} g|_{H^{s_{1}}}|h|_{H^{s_{2}+1}_{\gamma/2+s}}|\mu^{\frac{1}{256}}f|_{L^{2}} \lesssim
|g|_{H^{s_{1}}}|h|_{H^{s_{2}+1}_{\gamma/2+s}}|f|_{\mathcal{L}^{s}_{\gamma/2}},
\een
where we use \eqref{key-estimate-only-star}, the estimate \eqref{estimate-of-Z-g2-h2-total-one-half-with-mu} and Remark \ref{mu-to-N-or-M-X-Y}.

By \eqref{change-v-and-v-star}, the estimate \eqref{estimate-of-Y-g-h-f-total-one-half-with-mu} and Remark \ref{mu-to-N-or-M-X-Y},
 we have
\ben \label{result-U-2-2}
|\mathcal{U}_{2,2}| = |\int B  N g (h \mu^{\f{1}{2}}f)_{*} \mathrm{D}(N_{*}) \mathrm{d}V| \lesssim |\mu^{\frac{1}{256}}g|_{H^{s_{1}}}|\mu^{\frac{1}{256}}h|_{H^{s_{2}}}|\mu^{\frac{1}{256}}f|_{L^{2}}.
\een
Patching together \eqref{result-U-2-1} and \eqref{result-U-2-2}, recalling \eqref{U-2-into-U21-and-U22}, we get
\ben \label{result-U-2}
|\mathcal{U}_{2}| \lesssim |g|_{H^{\f{1}{2}}} |h|_{\mathcal{L}^{s}_{\gamma/2}} |f|_{\mathcal{L}^{s}_{\gamma/2}}.
\een
Patching together \eqref{result-U-2-1-another} and \eqref{result-U-2-2}, recalling \eqref{U-2-into-U21-and-U22}, we get
\ben \label{result-U-2-another}
|\mathcal{U}_{2}| \lesssim |g|_{H^{s_{1}}}|h|_{H^{s_{2}+1}_{\gamma/2+s}}|f|_{\mathcal{L}^{s}_{\gamma/2}}.
\een

{\it {Estimate of $\mathcal{U}_{3}$.}} Recalling \eqref{G-2-r-2-U-3},
by Cauchy-Schwartz inequality, the imbedding $H^{2} \hookrightarrow L^{\infty}$, \eqref{M-rho-mu-N}, \eqref{change-v-and-v-prime} and Theorem \ref{two-functional-ub-by-norm}, we have
\ben \nonumber
|\mathcal{U}_{3}|  &\lesssim& |g|_{L^{\infty}} \cs{\int
B  N^{\prime}_{*} \mathrm{D}^{2}(h)
\mathrm{d}V} {\int
B  N^{\prime}_{*} \mathrm{D}^{2}(f)
\mathrm{d}V}
\\ \label{result-U-3} &\lesssim& |g|_{H^{2}} \cs{\int B  \mu^{\f12}_{*} \mathrm{D}^{2}(h) \mathrm{d}V}
{\int B  \mu^{\f12}_{*} \mathrm{D}^{2}(f) \mathrm{d}V}
\lesssim |g|_{H^{2}}|h|_{\mathcal{L}^{s}_{\gamma/2}}|f|_{\mathcal{L}^{s}_{\gamma/2}}. \quad
\een
By Cauchy-Schwartz inequality, \eqref{M-rho-mu-N}, \eqref{change-v-and-v-star}, \eqref{change-v-and-v-prime}, the estimate \eqref{key-estimate-only-star} and Theorem \ref{two-functional-ub-by-norm}, we have
\ben \nonumber
|\mathcal{U}_{3}|  &\lesssim& \cs{ \int
B  N^{\prime}_{*} (g^{2})^{\prime}_{*} \mathrm{D}^{2}(h)
\mathrm{d}V} {\int
B  N^{\prime}_{*} \mathrm{D}^{2}(f)
\mathrm{d}V}
\\ \label{result-U-3-another} &\lesssim&  \cs{ \int B  \mu^{\f12} g^{2} \mathrm{D}^{2}(h_{*}) \mathrm{d}V}
{\int B  \mu^{\f12}_{*} \mathrm{D}^{2}(f) \mathrm{d}V}
\lesssim |g|_{H^{s_{1}}}|h|_{H^{s_{2}+1}_{\gamma/2+s}}|f|_{\mathcal{L}^{s}_{\gamma/2}}. \quad \quad
\een

{\it {Estimate of $\mathcal{U}_{4}$.}} Recalling \eqref{G-2-r-2-U-4} and using \eqref{change-v-and-v-prime}, we have
\beno
\mathcal{U}_{4} &=&  \int B  (Ng)_{*} h \mathrm{D}(f)(M_{*}^{\prime}-M) \mathrm{d}V
\\&=&   \int B  (Ng)_{*} Mh \mathrm{D}(f^{\prime}) \mathrm{d}V + \int B  (M Ng)_{*} h \mathrm{D}(f) \mathrm{d}V
+ \int B  (Ng)_{*} h \mathrm{D}(f)\mathrm{D}(M_{*}^{\prime}) \mathrm{d}V
\\&=&    \langle  Q_{c}(Ng, Mh), f \rangle   -  \langle  Q_{c}(MNg, h), f \rangle
- \mathcal{T}_{4,2},
\eeno
where we recall \eqref{U42-definition} for the notation $\mathcal{T}_{4,2}$.
By Corollary \ref{Q-up-bound-full-on-first-or-second} and Lemma \ref{product-take-out}, we have
\beno
| \langle  Q_{c}(Ng, Mh), f \rangle  | + | \langle  Q_{c}(MNg, h), f \rangle  | \lesssim
\min\{|g|_{H^{2}}  |h|_{\mathcal{L}^{s}_{\gamma/2}} , |g|_{L^{2}}(|h|_{\mathcal{L}^{s}_{\gamma/2}}+| h|_{H^{s+2}_{\gamma/2}})\} |f|_{\mathcal{L}^{s}_{\gamma/2}}.
\eeno
From which together with \eqref{result-U42}, we get
\ben \label{result-U-4}
|\mathcal{U}_{4}| \lesssim \min\{|g|_{H^{2}}  |h|_{\mathcal{L}^{s}_{\gamma/2}} , |g|_{L^{2}}(|h|_{\mathcal{L}^{s}_{\gamma/2}}+| h|_{H^{s+2}_{\gamma/2}}+| h|_{H^{\f{1}{2}}_{\gamma/2+s}})\} |f|_{\mathcal{L}^{s}_{\gamma/2}}.
\een

Patching together \eqref{result-U-1}, \eqref{result-U-2}, \eqref{result-U-3} and \eqref{result-U-4}, recalling \eqref{G-2-r-2-into-1-2-3-4},
we get
\ben \label{result-U-final-1}
 |  \langle  \Gamma_{2,r,2}^{\rho}(g, h), f \rangle   |  \lesssim |g|_{H^{2}}|h|_{\mathcal{L}^{s}_{\gamma/2}}|f|_{\mathcal{L}^{s}_{\gamma/2}}.
\een
Patching together \eqref{result-U-1-another}, \eqref{result-U-2-another}, \eqref{result-U-3-another} and \eqref{result-U-4}, by taking $s_{1}=0,s_{2}=\f{1}{2}$, recalling \eqref{G-2-r-2-into-1-2-3-4},
we get
\ben \label{result-U-final-2}
 |  \langle  \Gamma_{2,r,2}^{\rho}(g, h), f \rangle   |  \lesssim |g|_{L^{2}}(|h|_{\mathcal{L}^{s}_{\gamma/2}}+| h|_{H^{s+2}_{\gamma/2}}+| h|_{H^{\frac{3}{2}}_{\gamma/2+s}})|f|_{\mathcal{L}^{s}_{\gamma/2}}.
\een
Patching together \eqref{result-U-final-1} and \eqref{result-U-final-2}, we finish the proof.
\end{proof}

\begin{prop}\label{upper-bound-of-Gamma-2-r-3}  It holds that
\beno
| \langle  \Gamma_{2,r,3}^{\rho}(g, h), f \rangle  | \lesssim  \min\{|g|_{H^{2}}  |h|_{\mathcal{L}^{s}_{\gamma/2}}, |g|_{L^{2}}  (|\mu^{\frac{1}{256}}h|_{H^{2}}+|h|_{H^{\frac{3}{2}}_{\gamma/2+s}})\} |f|_{\mathcal{L}^{s}_{\gamma/2}}.
\eeno
\end{prop}
\begin{proof} Recalling the definition of $\Gamma_{2,r,3}^{\rho}$ in \eqref{definition-Gamma-remaining-3}, we have
\ben \label{Gamma-2-r-3-into-1-2}
\Gamma_{2,r,3}^{\rho}(g, h) &=& \Gamma_{2,r,3,1}^{\rho}(g,h) + \Gamma_{2,r,3,2}^{\rho}(g,h),
\\ \label{definition-Gamma-remaining-3-1}
\Gamma_{2,r,3,1}^{\rho}(g,h) &\colonequals    & N^{-1} \int
B (Ng)^{\prime}(Nh)\mathrm{D}(M^{\prime}_{*}) \mathrm{d}\sigma \mathrm{d}v_{*},
\\ \label{definition-Gamma-remaining-3-2}
\Gamma_{2,r,3,2}^{\rho}(g,h) &\colonequals    & N^{-1} \int
B  (Ng)^{\prime}_{*}(Nh)_{*}\mathrm{D}(M^{\prime})
\mathrm{d}\sigma \mathrm{d}v_{*}.
\een
Then it suffices to consider $ \langle  \Gamma_{2,r,3,1}^{\rho}(g, h), f \rangle  $ and $ \langle  \Gamma_{2,r,3,2}^{\rho}(g, h), f \rangle  $.

{\it Estimate of $ \langle  \Gamma_{2,r,3,1}^{\rho}(g, h), f \rangle  $.} Recalling \eqref{definition-Gamma-remaining-3-1}, using \eqref{change-v-and-v-prime}, we have
\ben
 \langle  \Gamma_{2,r,3,1}^{\rho}(g, h), f \rangle   =  \int
B (Ng)^{\prime} h f \mathrm{D}(M^{\prime}_{*}) \mathrm{d}V
 \label{Gamma-2-r-3-1-ghf} =  \int
B (Ng) (hf)^{\prime} \mathrm{D}(M_{*})
\mathrm{d}V = \mathcal{V}_{1} + \mathcal{V}_{2} + \mathcal{V}_{3},
\\  
\mathcal{V}_{1} \colonequals      \int
B (Ng) h^{\prime}\mathrm{D}(f^{\prime}) \mathrm{D}(M_{*})
\mathrm{d}V,
\quad 
\mathcal{V}_{2} \colonequals      \int
B (Ng) \mathrm{D}(h^{\prime})f \mathrm{D}(M_{*})
\mathrm{d}V,
\quad \label{g-no-difference-second} \mathcal{V}_{3} \colonequals      \int
B (Ng) h f \mathrm{D}(M_{*})
\mathrm{d}V.
\een
By Cauchy-Schwartz inequality, the imbedding $H^{2} \hookrightarrow L^{\infty}$, \eqref{change-v-and-v-prime}, \eqref{M-rho-mu-N} and \eqref{mu-weight-result}, we get
\ben \nonumber
|\mathcal{V}_{1}| &\lesssim& |g|_{L^{\infty}} \cs{\int B N^{2} (h^{2})^{\prime}
\mathrm{D}^{2}(M^{\f{1}{2}}_{*})
\mathrm{d}V}  {\int B  \mathrm{D}^{2}(f) \mathrm{A}(M_{*})
\mathrm{d}V}
\\ &\lesssim&|g|_{H^{2}} \cs{\int B  \mu^{\frac{1}{8}} \mu^{\frac{1}{8}}_{*} h^{2}
\mathrm{D}^{2}(M^{\f{1}{4}}_{*})
\mathrm{d}V}  {\int B  \mu_{*} \mathrm{D}^{2}(f)
\mathrm{d}V}
 \label{V1-result-g-highest} \lesssim |g|_{H^{2}} |\mu^{\frac{1}{256}}h|_{L^{2}} |f|_{\mathcal{L}^{s}_{\gamma/2}},
\een
where we use Theorem \ref{two-functional-ub-by-norm}, the estimate \eqref{estimate-of-Z-g2-h2-total-one-half-with-mu} and Remark \ref{mu-to-N-or-M-X-Y} in the last line.
By the same derivation, we have
\ben
|\mathcal{V}_{2}| \lesssim |g|_{L^{\infty}} \cs{\int B \mathrm{D}^{2}(h) \mathrm{A}(M_{*})
\mathrm{d}V}  {\int B N^{2} f^{2} \mathrm{D}^{2}(M^{\f{1}{2}}_{*})
\mathrm{d}V}
\label{V2-result-g-highest} \lesssim |g|_{H^{2}} |h|_{\mathcal{L}^{s}_{\gamma/2}} |\mu^{\frac{1}{256}}f|_{L^{2}} .
\een
By \eqref{estimate-of-D-g-h-f-total-one-half-with-mu}(in which $g, h, f$ play the same role) and Remark \ref{mu-to-N-or-M-X-Y},
we have
\ben \label{V3-result-g-h-flexiable}
|\mathcal{V}_{3}| \lesssim  |g|_{H^{s_{3}}}|\mu^{\frac{1}{256}}h|_{H^{s_{4}}}|\mu^{\frac{1}{256}}f|_{L^{2}}.
\een
Here and in the rest of the article, $(s_{3}, s_{4})=(2, 0)$ or $(s_{3}, s_{4})=(0, 2)$ unless otherwise specified.

Patching together \eqref{V1-result-g-highest}, \eqref{V2-result-g-highest} and \eqref{V3-result-g-h-flexiable}, recalling \eqref{Gamma-2-r-3-1-ghf}, by taking $(s_{3}, s_{4})=(2, 0)$,
we get
\ben \label{Gamma-2-r-3-1-result-g-highest}
| \langle  \Gamma_{2,r,3,1}^{\rho}(g, h), f \rangle  | \lesssim |g|_{H^{2}} |h|_{\mathcal{L}^{s}_{\gamma/2}} |f|_{\mathcal{L}^{s}_{\gamma/2}}.
\een

{\it Another estimate of $ \langle  \Gamma_{2,r,3,1}^{\rho}(g, h), f \rangle  $.}
Recalling  \eqref{g-no-difference-second} for $\mathcal{V}_{3}$, we have
\ben \label{Gamma-2-r-3-1-ghf-into-WV}
 \langle  \Gamma_{2,r,3,1}^{\rho}(g, h), f \rangle   =  \int
B (Ng) (hf)^{\prime} \mathrm{D}(M_{*})
\mathrm{d}V = \mathcal{W}_{1} + \mathcal{W}_{2} + \mathcal{V}_{3},
\\ \nonumber 
\mathcal{W}_{1} \colonequals      \int
B \mu^{-\frac{1}{4}} N g \mathrm{D}(\mu^{\frac{1}{4}})(h f)^{\prime} \mathrm{D}(M_{*})
\mathrm{d}V,
\quad 
\mathcal{W}_{2} \colonequals      \int
B \mu^{-\frac{1}{4}} N g \mathrm{D}((\mu^{\frac{1}{4}} h f)^{\prime}) \mathrm{D}(M_{*})
\mathrm{d}V.
\een
By Cauchy-Schwartz inequality, \eqref{M-rho-mu-N},  \eqref{change-v-and-v-prime}, \eqref{change-v-and-v-star}, the estimate \eqref{mu-weight-result}, we get
\ben \nonumber
|\mathcal{W}_{1}| &\lesssim& \cs{ \int B \mu^{\frac{1}{4}} g^{2} \mathrm{D}^{2}(M^{\f{1}{2}}_{*})
\mathrm{d}V}  {\int B  \mu^{\frac{1}{4}} \mathrm{D}^{2}(\mu^{\frac{1}{4}}) (h^{2} f^{2})^{\prime} \mathrm{A}(M_{*})
\mathrm{d}V}
\\ \nonumber &\lesssim& \cs{\int B  \mu^{\frac{1}{16}} \mu^{\frac{1}{16}}_{*} g^{2}\mathrm{D}^{2}(M^{\f{1}{4}}_{*})
\mathrm{d}V} {\int B \mu^{\frac{1}{16}} \mu^{\frac{1}{16}}_{*} h^{2}_{*} f^{2}_{*}
\mathrm{D}^{2}(\mu^{\frac{1}{4}}_{*})
\mathrm{d}V}
\\ \label{W1-result-h-highest} &\lesssim& |\mu^{\frac{1}{256}}g|_{L^{2}} |\mu^{\frac{1}{256}}h|_{H^{2}} |\mu^{\frac{1}{256}}f|_{L^{2}},
\een
where we use the imbedding $H^{2} \hookrightarrow L^{\infty}$, the estimate \eqref{estimate-of-Z-g2-h2-total-one-half-with-mu} and Remark \ref{mu-to-N-or-M-X-Y} in the  last line.
By Cauchy-Schwartz inequality, \eqref{change-v-and-v-prime}, the estimate \eqref{mu-weight-result},  Theorem \ref{two-functional-ub-by-norm},
the estimate \eqref{estimate-of-Z-g2-h2-total-one-half-with-mu} and Remark \ref{mu-to-N-or-M-X-Y}, we get
\ben \nonumber
|\mathcal{W}_{2}| &\lesssim& \cs{\int B \mu^{\f{1}{2}} g^{2} \mathrm{D}^{2}(M^{\f{1}{2}}_{*})
\mathrm{d}V} {\int B   \mathrm{D}^{2}(\mu^{\frac{1}{4}} h f) \mathrm{A}(M_{*})
\mathrm{d}V}
\\ \nonumber &\lesssim& \cs{ \int B  \mu^{\frac{1}{8}} \mu^{\frac{1}{8}}_{*} g^{2}\mathrm{D}^{2}(M^{\f{1}{4}}_{*})
\mathrm{d}V} {\int B \mu_{*}\mathrm{D}^{2}(\mu^{\frac{1}{4}} h f)
\mathrm{d}V}
\\ \label{W2-result} &\lesssim& |\mu^{\frac{1}{256}}g|_{L^{2}} |\mu^{\frac{1}{4}} h f|_{\mathcal{L}^{s}_{\gamma/2}} \lesssim  |\mu^{\frac{1}{256}}g|_{L^{2}} |\mu^{\frac{1}{16}} h |_{H^{2}} |f|_{\mathcal{L}^{s}_{\gamma/2}},
\een
where we use \eqref{taking-out-with-a-weight} in the last inequality.
Patching together \eqref{W1-result-h-highest}, \eqref{W2-result} and \eqref{V3-result-g-h-flexiable}, recalling \eqref{Gamma-2-r-3-1-ghf-into-WV}, by taking $(s_{3}, s_{4})=(0, 2)$, we get
\ben \label{Gamma-2-r-3-1-result-h-highest}
| \langle  \Gamma_{2,r,3,1}^{\rho}(g, h), f \rangle  | \lesssim |g|_{L^{2}} |\mu^{\frac{1}{256}}h|_{H^{2}} |f|_{\mathcal{L}^{s}_{\gamma/2}}.
\een

Patching together \eqref{Gamma-2-r-3-1-result-g-highest} and \eqref{Gamma-2-r-3-1-result-h-highest},  we get
\ben \label{G-2-r-3-1-result}
| \langle  \Gamma_{2,r,3,1}^{\rho}(g, h), f \rangle  | \lesssim \min\{ |g|_{H^{2}}  |h|_{\mathcal{L}^{s}_{\gamma/2}}, |g|_{L^{2}} |\mu^{\frac{1}{256}}h|_{H^{2}} \}|f|_{\mathcal{L}^{s}_{\gamma/2}}.
\een

{\it Estimate of $ \langle  \Gamma_{2,r,3,2}^{\rho}(g, h), f \rangle  $.}
Recalling \eqref{definition-Gamma-remaining-3-2}, taking inner product, we have
\ben \label{inner-product-remaining-2-g-g-star}
 \langle  \Gamma_{2,r,3,2}^{\rho}(g, h), f \rangle   = \int
B  (Ng)^{\prime}_{*}(Nh)_{*}\mathrm{D}(M^{\prime}) N^{-1} f
\mathrm{d}V.
\een
Recalling $M = \mu^{\f{1}{2}} N$, we have
$
\mathrm{D}(M^{\prime}) = \mathrm{D}((\mu^{\f{1}{2}}N)^{\prime}) = (\mu^{\f{1}{2}})^{\prime}\mathrm{D}(N^{\prime})+ \mathrm{D}((\mu^{\f{1}{2}})^{\prime})N
$
which gives
\beno
N^{\prime}_{*}N_{*}\mathrm{D}(M^{\prime}) N^{-1} &=& N^{\prime}_{*}N_{*}(\mu^{\f{1}{2}})^{\prime}\mathrm{D}(N^{\prime}) N^{-1}
+N^{\prime}_{*}N_{*}\mathrm{D}((\mu^{\f{1}{2}})^{\prime})
\\&=&  (\frac{1}{1 - \rho \mu})^{\prime}_{*} M_{*} \mathrm{D}(N^{\prime}) (1 - \rho \mu)
+N^{\prime}_{*}N_{*}\mathrm{D}((\mu^{\f{1}{2}})^{\prime}).
\eeno
Plugging which into \eqref{inner-product-remaining-2-g-g-star}, we get
\ben \label{Gamma-2-r-3-2-into-2-terms}
 \langle  \Gamma_{2,r,3,2}^{\rho}(g, h), f \rangle   &=& \mathcal{X}_{1}+\mathcal{X}_{2},
\\ \label{defintion-X1}
\mathcal{X}_{1} &\colonequals    &  \int
B (\frac{g}{1 - \rho \mu})^{\prime}_{*} (Mh)_{*} (1 - \rho \mu)f \mathrm{D}(N^{\prime})
\mathrm{d}V,
\\ \label{defintion-X2}
\mathcal{X}_{2} &\colonequals    &  \int
B (Ng)^{\prime}_{*}(Nh)_{*}f\mathrm{D}((\mu^{\f{1}{2}})^{\prime})
\mathrm{d}V.
\een
Since the two quantities have a similar structure, we only consider $\mathcal{X}_{2}$. Indeed, we will not use the factor $N$ before $g$ in $\mathcal{X}_{2}$ and so the estimate of $\mathcal{X}_{1}$ is similar. Using \eqref{change-v-and-v-prime},
note that
\ben \label{X2-into-4-terms}
\mathcal{X}_{2} =  \int
B (Ng)_{*}(Nh)^{\prime}_{*}f^{\prime}\mathrm{D}(\mu^{\f{1}{2}})
\mathrm{d}V = \mathcal{X}_{2,1} + \mathcal{X}_{2,2} + \mathcal{X}_{2,3} + \mathcal{X}_{2,4},
\\ \nonumber
\mathcal{X}_{2,1} \colonequals      \int
B (Ng)_{*}(Nh)^{\prime}_{*}\mathrm{D}(f^{\prime})\mathrm{D}(\mu^{\f{1}{2}})
\mathrm{d}V,
\quad
\mathcal{X}_{2,2} \colonequals      \int
B (Ng)_{*}N^{\prime}_{*}\mathrm{D}(h^{\prime}_{*})f\mathrm{D}(\mu^{\f{1}{2}})
\mathrm{d}V,
\\ \nonumber
\mathcal{X}_{2,3} \colonequals      \int
B (Ng)_{*}\mathrm{D}(N^{\prime}_{*})h_{*}f\mathrm{D}(\mu^{\f{1}{2}})
\mathrm{d}V,
\quad
\mathcal{X}_{2,4} \colonequals      \int
B (Ng)_{*}(Nh)_{*}f\mathrm{D}(\mu^{\f{1}{2}})
\mathrm{d}V.
\een
By Cauchy-Schwartz inequality, \eqref{M-rho-mu-N}, \eqref{change-v-and-v-prime}, \eqref{change-v-and-v-star}, the estimate \eqref{mu-weight-result}, the imbedding $H^{2} \hookrightarrow L^{\infty}$, the estimate \eqref{estimate-of-Z-g2-h2-total-one-half-with-mu}, Remark \ref{mu-to-N-or-M-X-Y} and Theorem \ref{two-functional-ub-by-norm},
we have
\beno
|\mathcal{X}_{2,1}| &\lesssim& \cs{ \int
B g^{2}_{*}(h^{2})^{\prime}_{*} (\mu^{\f12})^{\prime}_{*}\mathrm{D}^{2}(\mu^{\f{1}{2}})
\mathrm{d}V} { \int
B (\mu^{\f12})^{\prime}_{*}\mathrm{D}^{2}(f)
\mathrm{d}V}
\\&\lesssim& \cs{ \int
B g^{2}(h^{2})^{\prime} \mu^{\frac{1}{8}}\mu^{\frac{1}{8}}_{*}\mathrm{D}^{2}(\mu^{\frac{1}{4}}_{*})
\mathrm{d}V} { \int
B \mu^{\f12}_{*}\mathrm{D}^{2}(f)
\mathrm{d}V}
\lesssim |\mu^{\frac{1}{256}}g|_{H^{s_{3}}} |\mu^{\frac{1}{256}}h|_{H^{s_{4}}} |f|_{\mathcal{L}^{s}_{\gamma/2}}.
\eeno

We will deal with $\mathcal{X}_{2,2}$ in two ways. On one hand, by Cauchy-Schwartz inequality, \eqref{M-rho-mu-N}, \eqref{change-v-and-v-prime}, \eqref{change-v-and-v-star}, the estimate \eqref{mu-weight-result}, Theorem \ref{two-functional-ub-by-norm} and  the estimate \eqref{estimate-of-Z-g2-h2-total-one-half-with-mu},
we have
\beno
|\mathcal{X}_{2,2}| &\lesssim& \cs{ \int
B \mathrm{D}^{2}(h_{*})\mathrm{A}(\mu^{\f{1}{2}})
\mathrm{d}V} {\int
B \mu^{\prime}_{*} g^{2}_{*} f^{2}\mathrm{D}^{2}(\mu^{\frac{1}{4}})
\mathrm{d}V}
\\&\lesssim& \cs{ \int
B \mu^{\f{1}{2}}_{*}\mathrm{D}^{2}(h)
\mathrm{d}V} {\int
B \mu^{\frac{1}{16}}\mu^{\frac{1}{16}}_{*} g^{2} f^{2}_{*}\mathrm{D}^{2}(\mu^{\frac{1}{8}}_{*})
\mathrm{d}V}
\lesssim |\mu^{\frac{1}{256}}g|_{H^{\f{1}{2}}} |h|_{\mathcal{L}^{s}_{\gamma/2}} |\mu^{\frac{1}{256}}f|_{L^{2}} .
\eeno
On the other hand, by Cauchy-Schwartz inequality, \eqref{M-rho-mu-N}, \eqref{change-v-and-v-star}, the estimate \eqref{mu-weight-result}, the estimate \eqref{estimate-of-Z-g2-h2-total-one-half-with-mu} and the estimate \eqref{key-estimate-only-star},
we have
\beno
|\mathcal{X}_{2,2}| &\lesssim& \cs{ \int
B g^{2}_{*} (\mu^{\f12})^{\prime}_{*} \mathrm{D}^{2}(\mu^{\frac{1}{4}})
\mathrm{d}V}{ \int
B (\mu^{\f12})^{\prime}_{*} \mathrm{D}^{2}(h_{*})f^{2}\mathrm{A}(\mu^{\f{1}{2}})
\mathrm{d}V}
\\&\lesssim& \cs{ \int
B g^{2}\mu^{\frac{1}{16}}\mu^{\frac{1}{16}}_{*}\mathrm{D}^{2}(\mu^{\frac{1}{8}}_{*})
\mathrm{d}V} {\int
B \mu^{\frac{1}{8}}\mu^{\frac{1}{8}}_{*}\mathrm{D}^{2}(h_{*})f^{2}
\mathrm{d}V}
\lesssim |\mu^{\frac{1}{256}}g|_{L^{2}} |h|_{H^{\frac{3}{2}}_{\gamma/2+s}} |\mu^{\frac{1}{64}}f|_{L^{2}} .
\eeno

Now we set to estimate $\mathcal{X}_{2,3}$. By Cauchy-Schwartz inequality, \eqref{M-rho-mu-N}, \eqref{change-v-and-v-star}, the estimate \eqref{mu-weight-result}, the imbedding $H^{2} \hookrightarrow L^{\infty}$, the estimate \eqref{estimate-of-Z-g2-h2-total-one-half-with-mu}, Theorem \ref{two-functional-ub-by-norm} and Remark \ref{mu-to-N-or-M-functional-N},
we have
\beno
|\mathcal{X}_{2,3}| &\lesssim& \cs{ \int
B g^{2}_{*}h^{2}_{*}
\mathrm{A}(\mu^{\f{1}{2}}_{*})\mathrm{D}^{2}(\mu^{\f{1}{2}})
\mathrm{d}V} { \int
B \mathrm{D}^{2}(N^{\frac{1}{2}}_{*})f^{2}
\mathrm{d}V}
\\&\lesssim& \cs{ \int
B g^{2}h^{2}\mu^{\frac{1}{8}}\mu^{\frac{1}{8}}_{*}\mathrm{D}^{2}(\mu^{\frac{1}{4}}_{*})
\mathrm{d}V} { \int
B f^{2}_{*}\mathrm{D}^{2}(N^{\frac{1}{2}})
\mathrm{d}V}
\lesssim |\mu^{\frac{1}{256}}g|_{H^{s_{3}}} |\mu^{\frac{1}{256}}h|_{H^{s_{4}}} |f|_{\mathcal{L}^{s}_{\gamma/2}} .
\eeno

By \eqref{change-v-and-v-star}, using  \eqref{estimate-of-X-g-h-f-total-one-half-with-mu} and Remark \ref{mu-to-N-or-M-X-Y}, we have
\beno
|\mathcal{X}_{2,4}| = |\int
B g h f_{*} N^{2} (\mu^{\f{1}{2}}_{*}-(\mu^{\f{1}{2}})^{\prime}_{*})
\mathrm{d}V| \lesssim |\mu^{\frac{1}{256}}g|_{H^{s_{1}}} |\mu^{\frac{1}{256}}h|_{H^{s_{2}}}|\mu^{\frac{1}{256}}f|_{L^{2}}.
\eeno

Patching together the above estimates of $\mathcal{X}_{2,1}, \mathcal{X}_{2,2}, \mathcal{X}_{2,3}, \mathcal{X}_{2,4}$, recalling \eqref{X2-into-4-terms}, we have
\ben \label{G-2-r-3-2-X-result}
|\mathcal{X}_{2}| \lesssim \min\{|g|_{H^{2}}  |h|_{\mathcal{L}^{s}_{\gamma/2}}, |g|_{L^{2}} ( |\mu^{\frac{1}{256}}h|_{H^{2}} + |h|_{H^{\frac{3}{2}}_{\gamma/2+s}})\} |f|_{\mathcal{L}^{s}_{\gamma/2}}.
\een
We emphasize that $\mathcal{X}_{1}$ enjoys the same upper bound as that of $\mathcal{X}_{2}$. By recalling \eqref{Gamma-2-r-3-2-into-2-terms}, we have
\ben \label{G-2-r-3-2-result}
| \langle  \Gamma_{2,r,3,2}^{\rho}(g, h), f \rangle  | \lesssim  \min\{|g|_{H^{2}}  |h|_{\mathcal{L}^{s}_{\gamma/2}}, |g|_{L^{2}} ( |\mu^{\frac{1}{256}}h|_{H^{2}} + |h|_{H^{\frac{3}{2}}_{\gamma/2+s}})\} |f|_{\mathcal{L}^{s}_{\gamma/2}}.
\een

Patching together \eqref{G-2-r-3-1-result} and \eqref{G-2-r-3-2-result}, we finish the proof.
\end{proof}

Recalling \eqref{Gamma-main-remaining} and \eqref{Gamma-remaining-into-three-terms},
patching together Theorem \ref{upper-bound-of-Gamma-2-m}, Propositions \ref{upper-bound-of-Gamma-2-r-1-g-h-f},
\ref{upper-bound-of-Gamma-2-r-2-g-h-f} and \ref{upper-bound-of-Gamma-2-r-3}, since $|\mu^{\frac{1}{256}}h|_{H^{2}} \lesssim |h|_{H^{2}_{\gamma/2}}$,
we get
\begin{thm} \label{Gamma-2-ub}
It holds that
\beno
 |\langle \Gamma_{2}^{\rho}(g,h), f \rangle|
\lesssim \rho^{\f12}\min\{|g|_{H^{2}}  |h|_{\mathcal{L}^{s}_{\gamma/2}} , |g|_{L^{2}}(|h|_{\mathcal{L}^{s}_{\gamma/2}}+| h|_{H^{s+2}_{\gamma/2}}+| h|_{H^{\frac{3}{2}}_{\gamma/2+s}})\} |f|_{\mathcal{L}^{s}_{\gamma/2}}.
\eeno
\end{thm}

\section{Trilinear operator estimate} \label{trilinear}
In this section, we estimate the trilinear operator $\Gamma_{3}^{\rho}$ defined in \eqref{definition-Gamma-3-epsilon}. Recalling the relation \eqref{M--with-mathcal-M} between
$N$ and $\mathcal{N}$, we have
\ben \label{Gamma3-into-Gamma31-and-Gamma32}
\Gamma_{3}^{\rho}(g,h,\varrho) &=& \rho \Gamma_{3,1}^{\rho}(g,h,\varrho) + \rho \Gamma_{3,2}^{\rho}(g,h,\varrho),
\\ \label{defintion-Gamma31}
\Gamma_{3,1}^{\rho}(g,h,\varrho) &\colonequals&      N^{-1}\int
B  \mathrm{D}\big( (Ng)_{*}^{\prime} (Nh)^{\prime} (N\varrho)_{*} \big) \mathrm{d}\sigma \mathrm{d}v_{*},
\\ \label{defintion-Gamma32}
\Gamma_{3,2}^{\rho}(g,h,\varrho) &\colonequals&      N^{-1}\int
B  \mathrm{D}\big( (Ng)_{*}^{\prime} (Nh)^{\prime} N\varrho \big) \mathrm{d}\sigma \mathrm{d}v_{*}.
\een

The whole section is devoted to prove
\begin{thm} \label{Gamma-3-ub}  The following functional estimates are valid.
\ben \label{rho-highest-derivative}
|\langle \Gamma_{3}^{\rho}(g,h,\varrho), f\rangle| &\lesssim& \rho|g|_{H^{2}}|\mu^{\frac{1}{256}}h|_{H^{2}}|\varrho|_{\mathcal{L}^{s}_{\gamma/2}}|f|_{\mathcal{L}^{s}_{\gamma/2}}
\\ \nonumber &&+ \rho|g|_{H^{3}}(|h|_{\mathcal{L}^{s}_{\gamma/2}}+|h|_{H^{s+2}_{\gamma/2}}
+|h|_{H^{\frac{3}{2}}_{\gamma/2+s}})|\mu^{\frac{1}{256}}\varrho|_{L^{2}}|f|_{\mathcal{L}^{s}_{\gamma/2}}.
\\ \label{h-highest-derivative}
|\langle \Gamma_{3}^{\rho}(g,h,\varrho), f\rangle| &\lesssim& \rho|g|_{H^{2}}|h|_{\mathcal{L}^{s}_{\gamma/2}}|\mu^{\frac{1}{256}}\varrho|_{H^{3}}|f|_{\mathcal{L}^{s}_{\gamma/2}}.
\\ \label{g-highest-derivative}
|\langle \Gamma_{3}^{\rho}(g,h,\varrho), f\rangle| &\lesssim& \rho|g|_{L^{2}}(|h|_{\mathcal{L}^{s}_{\gamma/2}}+|h|_{H^{s+2}_{\gamma/2}}
+|h|_{H^{\f{1}{2}}_{\gamma/2+s}})|\mu^{\frac{1}{256}}\varrho|_{H^{3}}|f|_{\mathcal{L}^{s}_{\gamma/2}}.
\een
\end{thm}

\begin{proof} [Proof of \eqref{rho-highest-derivative}.]
Recalling \eqref{defintion-Gamma31}, using \eqref{usual-decomposition-into-N-and-f}, we have
\ben \label{cubic1-inner-product}
\langle \Gamma_{3,1}^{\rho}(g,h,\varrho), f\rangle &=&  \int
B (Ng)_{*} (Nh) (N\varrho)^{\prime}_{*}  \mathrm{D}((N^{-1}f)^{\prime}) \mathrm{d}V = \mathcal{A}_{1} + \mathcal{A}_{2},
\\  \label{Gamma3-1-defintion-A1}
\mathcal{A}_{1} &\colonequals    &  \int
B (Ng)_{*} h (N\varrho)^{\prime}_{*}  \mathrm{D}(f^{\prime})  \mathrm{d}V,
\\  \label{Gamma3-1-defintion-A2}
\mathcal{A}_{2} &\colonequals    &  \int
B (Ng)_{*} (Nh) (N\varrho)^{\prime}_{*} f^{\prime} \mathrm{D}((N^{-1})^{\prime})  \mathrm{d}V.
\een

{\it Estimate of $\mathcal{A}_{1}$.} Using \eqref{change-v-and-v-prime}, $(Ng)^{\prime}_{*} = \mathrm{D}((Ng)^{\prime}_{*}) + (Ng)_{*}$ and $h^{\prime} = \mathrm{D}(h^{\prime}) + h$,
we get
\ben \label{Z1-into-Z11-12-13}
\mathcal{A}_{1} &=&  \int
B (Ng)^{\prime}_{*} h^{\prime} (N\varrho)_{*}  \mathrm{D}(f)  \mathrm{d}V = \mathcal{A}_{1,1} + \mathcal{A}_{1,2} + \mathcal{A}_{1,3},
\\ \nonumber 
\mathcal{A}_{1,1} &\colonequals&     \int
B \mathrm{D}((Ng)^{\prime}_{*}) h^{\prime} (N\varrho)_{*}  \mathrm{D}(f)  \mathrm{d}V,
\\ \nonumber \mathcal{A}_{1,2} &\colonequals&     \int
B (N^{2}g\varrho)_{*} \mathrm{D}(h^{\prime})  \mathrm{D}(f)  \mathrm{d}V,
\\ \nonumber 
\mathcal{A}_{1,3} &\colonequals&     \int
B (N^{2}g\varrho)_{*} h  \mathrm{D}(f)  \mathrm{d}V.
\een
By Cauchy-Schwartz inequality and \eqref{M-rho-mu-N}, using \eqref{mu-star-g-star-difference-h2-f2-star}, Remark \ref{star-prime-or-prime-version} and Theorem \ref{two-functional-ub-by-norm},
we have
\ben \nonumber
|\mathcal{A}_{1,1}| &\lesssim& \cs{\int
B \mathrm{D}^{2}((Ng)_{*}) (h^{2})^{\prime} (\mu^{\f12} \varrho^{2})_{*}
\mathrm{d}V } {\int
B  \mu^{\f12}_{*} \mathrm{D}^{2}(f)
\mathrm{d}V }
\\ \label{result-Z11}  &\lesssim& |N g|_{H^{3}} |h|_{H^{\f{1}{2}}_{\gamma/2+s}} |\mu^{\frac{1}{128}}\varrho|_{L^{2}} |f|_{\mathcal{L}^{s}_{\gamma/2}} \lesssim |g|_{H^{3}} |h|_{H^{\f{1}{2}}_{\gamma/2+s}} |\mu^{\frac{1}{128}}\varrho|_{L^{2}} |f|_{\mathcal{L}^{s}_{\gamma/2}}.
\een
By Cauchy-Schwartz inequality,  the imbedding $H^{2} \hookrightarrow L^{\infty}$, \eqref{M-rho-mu-N}, \eqref{change-v-and-v-star}, the estimate \eqref{key-estimate-only-star}
and Theorem \ref{two-functional-ub-by-norm}, we have
\ben
|\mathcal{A}_{1,2}|\lesssim |g|_{L^{\infty}} \cs{\int
B \mu \varrho^{2}\mathrm{D}^{2}(h_{*})
\mathrm{d}V}{\int
B  \mu_{*} \mathrm{D}^{2}(f)
\mathrm{d}V}
\label{result-Z12} \lesssim |g|_{H^{2}} |h|_{H^{\frac{3}{2}}_{\gamma/2+s}} |\mu^{\frac{1}{256}}\varrho|_{L^{2}} |f|_{\mathcal{L}^{s}_{\gamma/2}}.
\een
Note that $|\mathcal{A}_{1,3}| = | \langle  Q_{c}(N^{2}g\varrho, h), f \rangle  |$. Then by Corollary \ref{Q-up-bound-full-on-first-or-second} and Lemma \ref{product-take-out}, we have
\ben \label{result-Z13}
|\mathcal{A}_{1,3}|  \lesssim |N^{2}g\varrho|_{L^{2}_{7}}(|h|_{\mathcal{L}^{s}_{\gamma/2}}+|h|_{H^{s+2}_{\gamma/2}})|f|_{\mathcal{L}^{s}_{\gamma/2}} \lesssim
|g|_{H^{2}}(|h|_{\mathcal{L}^{s}_{\gamma/2}}+|h|_{H^{s+2}_{\gamma/2}})|\mu^{\f{1}{2}}\varrho|_{L^{2}}
|f|_{\mathcal{L}^{s}_{\gamma/2}}.
\een
Recalling \eqref{Z1-into-Z11-12-13}, patching together \eqref{result-Z11}, \eqref{result-Z12} and \eqref{result-Z13}, we get
\ben \label{result-Z1}
|\mathcal{A}_{1}| \lesssim
|g|_{H^{3}}(|h|_{\mathcal{L}^{s}_{\gamma/2}}+|h|_{H^{s+2}_{\gamma/2}}
+|h|_{H^{\frac{3}{2}}_{\gamma/2+s}})|\mu^{\frac{1}{256}}\varrho|_{L^{2}}|f|_{\mathcal{L}^{s}_{\gamma/2}}.
\een

{\it Estimate of $\mathcal{A}_{2}$.} Using \eqref{change-v-and-v-prime}, we get
\ben
 \label{definition-of-Z2}
\mathcal{A}_{2} =   \int
B (Ng)^{\prime}_{*} (Nh)^{\prime} (N\varrho)_{*}f  \mathrm{D}(N^{-1})  \mathrm{d}V.
\een
Plugging the following two identities
\ben \nonumber
(\frac{g}{1 - \rho \mu})^{\prime}_{*} (\frac{h}{1 - \rho \mu})^{\prime} &=& (\frac{g}{1 - \rho \mu})^{\prime}_{*}
\mathrm{D}((\frac{h}{1 - \rho \mu})^{\prime})
+ \mathrm{D}((\frac{g}{1 - \rho \mu})^{\prime}_{*})\frac{h}{1 - \rho \mu} + (\frac{g}{1 - \rho \mu})_{*} \frac{h}{1 - \rho \mu},
\\ \label{identity-to-cancel-N-negative-diff}
\mu^{\f{1}{2}}_{*}\mu^{\f{1}{2}}\mathrm{D}(N^{-1}) &=&
 \mathrm{D}(\mu^{\f{1}{2}}_{*}) - \rho \mu^{\f{1}{2}}_{*}\mu^{\f{1}{2}}\mathrm{D}(\mu^{\f{1}{2}}),
\een
into \eqref{definition-of-Z2}, we have
\ben \label{Z2-into-6-terms}
\mathcal{A}_{2} &=&   \int
B (\frac{g}{1 - \rho \mu})^{\prime}_{*} (\frac{h}{1 - \rho \mu})^{\prime} (N\varrho)_{*}f  \mu^{\f{1}{2}}\mu^{\f{1}{2}}_{*}\mathrm{D}(N^{-1})  \mathrm{d}V
= \sum_{i=1}^{6} \mathcal{A}_{2,i},
\een
where
\ben
\nonumber
\mathcal{A}_{2,1} &\colonequals    &  \int
B (\frac{g}{1 - \rho \mu})^{\prime}_{*} \mathrm{D}((\frac{h}{1 - \rho \mu})^{\prime}) (N\varrho)_{*}f  \mathrm{D}(\mu^{\f{1}{2}}_{*})  \mathrm{d}V,
\\  \nonumber
\mathcal{A}_{2,2} &\colonequals    &  -\rho \int
B (\frac{g}{1 - \rho \mu})^{\prime}_{*} \mathrm{D}((\frac{h}{1 - \rho \mu})^{\prime}) (N\varrho)_{*}f  \mu^{\f{1}{2}}_{*}\mu^{\f{1}{2}}\mathrm{D}(\mu^{\f{1}{2}})  \mathrm{d}V,
\\  \nonumber
\mathcal{A}_{2,3} &\colonequals    &  \int
B \mathrm{D}((\frac{g}{1 - \rho \mu})^{\prime}_{*})\frac{h}{1 - \rho \mu} (N\varrho)_{*}f
\mathrm{D}(\mu^{\f{1}{2}}_{*})  \mathrm{d}V,
\\  \nonumber
\mathcal{A}_{2,4} &\colonequals    &  -\rho \int
B \mathrm{D}((\frac{g}{1 - \rho \mu})^{\prime}_{*})\frac{h}{1 - \rho \mu} (N\varrho)_{*}f  \mu^{\f{1}{2}}_{*}\mu^{\f{1}{2}}\mathrm{D}(\mu^{\f{1}{2}})   \mathrm{d}V,
\\ \label{definition-Z25}
\mathcal{A}_{2,5} &\colonequals    &  \int
B (\frac{gN\varrho}{1 - \rho \mu})_{*} \frac{hf}{1 - \rho \mu}
\mathrm{D}(\mu^{\f{1}{2}}_{*})  \mathrm{d}V,
\\ \label{definition-Z26}
\mathcal{A}_{2,6} &\colonequals    &  -\rho \int
B (\frac{gN\varrho}{1 - \rho \mu})_{*} \frac{hf}{1 - \rho \mu}  \mu^{\f{1}{2}}_{*}\mu^{\f{1}{2}}\mathrm{D}(\mu^{\f{1}{2}})   \mathrm{d}V.
\een
By Cauchy-Schwartz inequality and \eqref{M-rho-mu-N}, we have
\beno
|\mathcal{A}_{2,1}| \lesssim
\cs{ \int
B (g^{2})^{\prime}_{*} f^{2}  \mathrm{D}^{2}(\mu^{\f{1}{2}}_{*})
\mathrm{d}V}{ \int
B  \mathrm{D}^{2}(\frac{h}{1 - \rho \mu}) \mu_{*}\varrho^{2}_{*}
\mathrm{d}V}.
\eeno
By the change of variable $v_{*} \rightarrow v^{\prime}_{*}$ and the estimate \eqref{estimate-of-Z-g2-h2-total-one-half}, we get
\beno
\int
B (g^{2})^{\prime}_{*} f^{2}  \mathrm{D}^{2}(\mu^{\f{1}{2}}_{*})
\mathrm{d}V \lesssim |g|_{H^{\f{1}{2}}}^{2} |f|_{L^{2}_{\gamma/2+s}}^{2}.
\eeno
Recalling \eqref{difference-h-over} and \eqref{change-v-and-v-star}, using  \eqref{key-estimate-only-star} and \eqref{estimate-of-Z-g2-h2-total-one-half-with-mu}
we get
\beno
\int
B  \mathrm{D}^{2}(\frac{h}{1 - \rho \mu}) \mu_{*}\varrho^{2}_{*}
\mathrm{d}V
\lesssim \int
B  \mathrm{D}^{2}(h_{*}) \mu \varrho^{2}
\mathrm{d}V  + \int
B  \mathrm{D}^{2}(\mu_{*})h^{2}_{*} \mu \varrho^{2}
\mathrm{d}V
\lesssim |h|_{H^{\frac{3}{2}}_{\gamma/2+s}}^{2} |\mu^{\frac{1}{256}}\varrho|_{L^{2}}^{2}.
\eeno
Patching together the previous two estimates, we get
\beno
|\mathcal{A}_{2,1}| \lesssim
 |g|_{H^{\f{1}{2}}} |h|_{H^{\frac{3}{2}}_{\gamma/2+s}} |\mu^{\frac{1}{256}}\varrho|_{L^{2}} |f|_{L^{2}_{\gamma/2+s}}.
\eeno
By nearly the same argument as that for $\mathcal{A}_{2,1}$, thanks to the factor $\mu^{\f{1}{2}}_{*}\mu^{\f{1}{2}}$,
we can also get
\beno
|\mathcal{A}_{2,2}| \lesssim
 |\mu^{\frac{1}{256}}g|_{H^{\f{1}{2}}} |\mu^{\frac{1}{256}}h|_{H^{\frac{3}{2}}} |\mu^{\frac{1}{256}}\varrho|_{L^{2}} |\mu^{\frac{1}{256}}f|_{L^{2}}.
\eeno
By Cauchy-Schwartz inequality and \eqref{M-rho-mu-N}, we have
\beno
|\mathcal{A}_{2,3}| \lesssim
\cs{ \int
B \mathrm{D}^{2}((\frac{g}{1 - \rho \mu})_{*}) \mu^{\f12}_{*}f^{2}
\mathrm{d}V}{\int
B  h^{2} (\mu^{\f12} \varrho^{2})_{*}\mathrm{D}^{2}(\mu^{\f{1}{2}}_{*})
\mathrm{d}V}.
\eeno
By the estimate \eqref{estimate-of-Z-g2-h2-total-one-half}, we get
\ben \label{Z23-t1}
\int
B h^{2} (\mu^{\f12} \varrho^{2})_{*}  \mathrm{D}^{2}(\mu^{\f{1}{2}}_{*})
\mathrm{d}V \lesssim  |h|_{H^{\f{1}{2}}_{\gamma/2+s}}^{2} |\mu^{\frac{1}{4}}\varrho|_{L^{2}}^{2}.
\een
Recalling \eqref{difference-h-over}, we have
\ben \label{difference-g-over}
\mathrm{D}^{2}((\frac{g}{1 - \rho \mu})_{*}) \lesssim (g^{\prime}_{*}-g_{*})^{2} + \mathrm{D}^{2}(\mu_{*})g^{2}_{*}.
\een
By  \eqref{difference-g-over}, using \eqref{mu-star-g-star-difference-h2} and \eqref{estimate-of-Z-g2-h2-total-one-half}, we get
\ben
 \int
B \mathrm{D}^{2}((\frac{g}{1 - \rho \mu})_{*}) \mu^{\f12}_{*}f^{2}
\mathrm{d}V
 \label{Z23-t2} \lesssim \int
B (g^{\prime}_{*}-g_{*})^{2} \mu^{\f12}_{*} f^{2}
\mathrm{d}V  + \int
B  \mathrm{D}^{2}(\mu_{*})g^{2}_{*} \mu^{\f12}_{*} f^{2}
\mathrm{d}V
\lesssim |g|_{H^{\frac{3}{2}}}^{2} |f|_{L^{2}_{\gamma/2+s}}^{2}.
\een
Patching together \eqref{Z23-t1} and \eqref{Z23-t2}, we get
$
|\mathcal{A}_{2,3}| \lesssim
 |g|_{H^{\frac{3}{2}}} |h|_{H^{\f{1}{2}}_{\gamma/2+s}} |\mu^{\frac{1}{4}}\varrho|_{L^{2}} |f|_{L^{2}_{\gamma/2+s}}.
$
By nearly the same argument as that for $\mathcal{A}_{2,3}$, we can also get
$
|\mathcal{A}_{2,4}| \lesssim
|g|_{H^{\frac{3}{2}}} |h|_{H^{\f{1}{2}}_{\gamma/2+s}} |\mu^{\frac{1}{4}}\varrho|_{L^{2}} |f|_{L^{2}_{\gamma/2+s}}.
$
Using \eqref{estimate-of-A-g-h-varrho-f-total-one-half-plus-2}, for $\{ a_{1}, a_{2}, a_{3}\} = \{ 2, \f{1}{2}, 0\}$, we get
\ben \label{general-result-without-mu}
|\mathcal{A}_{2,5}| \lesssim  |\mu^{\frac{1}{8}}g|_{H^{a_{1}}} |h|_{H^{a_{2}}_{\gamma/2+s}} |\mu^{\frac{1}{8}}\varrho|_{H^{a_{3}}} |f|_{L^{2}_{\gamma/2+s}}.
\een
By \eqref{change-v-and-v-star} and \eqref{estimate-of-A-g-h-varrho-f-with-mu},
we have
\ben \label{general-result-with-mu}
|\mathcal{A}_{2,6}|
\lesssim  |\mu^{\frac{1}{256}}g|_{H^{a_{1}}} |\mu^{\frac{1}{256}}h|_{H^{a_{2}}} |\mu^{\frac{1}{256}}\varrho|_{H^{a_{3}}}  |\mu^{\frac{1}{256}}f|_{L^{2}}.
\een

Recalling \eqref{Z2-into-6-terms}, patching together the above estimates of $\mathcal{A}_{2,i}(1 \leq i\leq 6)$, taking $( a_{1}, a_{2}, a_{3}) = ( 2, \f{1}{2}, 0)$, we get
\ben \label{result-Z2}
|\mathcal{A}_{2}| \lesssim  |g|_{H^{2}} |h|_{H^{\frac{3}{2}}_{\gamma/2+s}} |\mu^{\frac{1}{256}}\varrho|_{L^{2}} |f|_{L^{2}_{\gamma/2+s}}.
\een
Recalling \eqref{cubic1-inner-product}, patching together \eqref{result-Z1} and \eqref{result-Z2}, we have
\ben \label{Gamma-31-rho-highest-derivative}
|\langle \Gamma_{3,1}^{\rho}(g,h,\varrho), f\rangle| \lesssim |g|_{H^{3}}(|h|_{\mathcal{L}^{s}_{\gamma/2}}+|h|_{H^{s+2}_{\gamma/2}}
+|h|_{H^{\frac{3}{2}}_{\gamma/2+s}})|\mu^{\frac{1}{256}}\varrho|_{L^{2}}|f|_{\mathcal{L}^{s}_{\gamma/2}}.
\een

We set to estimate $\langle \Gamma_{3,2}^{\rho}(g,h,\varrho), f\rangle$. Recalling \eqref{defintion-Gamma32}, using the following identity
\ben \label{type2-decomposition}
\mathrm{D}((N^{-1}f)^{\prime}) =  (N^{-1})^{\prime} \mathrm{D}(f^{\prime}) + \mathrm{D}((N^{-1})^{\prime})f,
\een
we have
\ben \label{cubic2-inner-product}
\langle \Gamma_{3,2}^{\rho}(g,h,\varrho), f\rangle =  \int
B (Ng)_{*} (Nh) (N\varrho)^{\prime}  \mathrm{D}((N^{-1}f)^{\prime}) \mathrm{d}V = \mathcal{B}_{1} + \mathcal{B}_{2},
\\  \label{Gamma3-2-defintion-B1}
\mathcal{B}_{1}\colonequals      \int
B (Ng)_{*} (Nh) \varrho^{\prime} \mathrm{D}(f^{\prime})  \mathrm{d}V,
\quad
\mathcal{B}_{2}\colonequals      \int
B (Ng)_{*} (Nh) (N\varrho)^{\prime} \mathrm{D}((N^{-1})^{\prime})f  \mathrm{d}V.
\een

We now divide $\mathcal{B}_{1}$ into two terms.
Using $\varrho^{\prime}= \mathrm{D}(\varrho^{\prime}) + \varrho$, we get $\mathcal{B}_{1} = \mathcal{B}_{1,1} + \mathcal{B}_{1,2}$ where
\ben \label{definition-of-B11}
\mathcal{B}_{1,1} \colonequals     \int
B (Ng)_{*} (Nh) \mathrm{D}(\varrho^{\prime}) \mathrm{D}(f^{\prime})  \mathrm{d}V,
\quad
\mathcal{B}_{1,2} \colonequals     \int
B (Ng)_{*} (Nh\varrho) \mathrm{D}(f^{\prime})  \mathrm{d}V.
\een
By Cauchy-Schwartz inequality, \eqref{M-rho-mu-N}, the imbedding $H^{2} \hookrightarrow L^{\infty}$ and Theorem \ref{two-functional-ub-by-norm},
we have
\beno
|\mathcal{B}_{1,1}| \lesssim |g|_{L^{\infty}}|\mu^{\f{1}{2}}h|_{L^{\infty}} \cs{\int
B \mu^{\f{1}{2}}_{*}  \mathrm{D}^{2}(\varrho)
\mathrm{d}V}{\int
B  \mu^{\f{1}{2}}_{*}  \mathrm{D}^{2}(f)
\mathrm{d}V}
\lesssim |g|_{H^{2}}|\mu^{\f{1}{2}} h|_{H^{2}} |\varrho|_{\mathcal{L}^{s}_{\gamma/2}}|f|_{\mathcal{L}^{s}_{\gamma/2}}.
\eeno
Using Corollary \ref{Q-up-bound-full-on-first-or-second} and Lemma \ref{product-take-out}, we have
\beno
|\mathcal{B}_{1,2}| = |  \langle  Q_{c}(Ng, Nh\varrho), f \rangle  | \lesssim
|N g|_{H^{2}_{7}}|N h \varrho|_{\mathcal{L}^{s}_{\gamma/2}}|f|_{\mathcal{L}^{s}_{\gamma/2}} \lesssim |g|_{H^{2}}|\mu^{\frac{1}{8}}h|_{H^{2}} |\varrho|_{\mathcal{L}^{s}_{\gamma/2}}|f|_{\mathcal{L}^{s}_{\gamma/2}}.
\eeno

We now divide $\mathcal{B}_{2}$ into two terms. Using $N N^{\prime} \mathrm{D}((N^{-1})^{\prime}) =  \mathrm{D}(N)$ and
$\varrho^{\prime}= \mathrm{D}(\varrho^{\prime}) + \varrho$, we get $\mathcal{B}_{2} = \mathcal{B}_{2,1}+ \mathcal{B}_{2,2}$ where
\beno   
\mathcal{B}_{2,1} \colonequals   \int
B (Ng)_{*} h \mathrm{D}(\varrho^{\prime}) f\mathrm{D}(N)  \mathrm{d}V,
\quad 
\mathcal{B}_{2,2} \colonequals   \int
B (Ng)_{*} h \varrho f\mathrm{D}(N)  \mathrm{d}V.
\eeno
By Cauchy-Schwartz inequality, \eqref{M-rho-mu-N}, \eqref{change-v-and-v-star}, the estimate \eqref{mu-weight-result}, the imbedding $H^{2} \hookrightarrow L^{\infty}$, Theorem \ref{two-functional-ub-by-norm}, the estimate \eqref{estimate-of-Z-g2-h2-total-one-half-with-mu} and Remark \ref{mu-to-N-or-M-X-Y},
we have
\beno
|\mathcal{B}_{2,1}| &\lesssim&
\cs{ \int
B \mu^{\f{1}{2}}_{*} g_{*}^{2}h^{2}\mathrm{A}(N) \mathrm{D}^{2}(\varrho)
\mathrm{d}V}{\int
B  \mu^{\f{1}{2}}_{*}  f^{2} \mathrm{D}^{2}(N^{\f{1}{2}})
\mathrm{d}V}
\\&\lesssim& \cs{ \int
B g_{*}^{2}h^{2}\mu^{\frac{1}{8}}\mu^{\frac{1}{8}}_{*} \mathrm{D}^{2}(\varrho)
\mathrm{d}V}{\int
B  \mu^{\f{1}{2}}  f^{2}_{*} \mathrm{D}^{2}(N^{\f{1}{2}}_{*})
\mathrm{d}V}
\lesssim |g|_{H^{2}}|\mu^{\frac{1}{16}}h|_{H^{2}} |\varrho|_{\mathcal{L}^{s}_{\gamma/2}}|\mu^{\frac{1}{256}}f|_{L^{2}}.
\eeno
By \eqref{change-v-and-v-star}, using \eqref{estimate-of-C-g-h-varrho-f-with-mu}(in which $h, \varrho, f$ play the same role),
we have
\beno
|\mathcal{B}_{2,2}| =
|\int
B Ng (h \varrho f)_{*} \mathrm{D}(N_{*})  \mathrm{d}V|
\lesssim |\mu^{\frac{1}{256}} g|_{H^{\f{1}{2}}}|\mu^{\frac{1}{256}}h|_{H^{2}} |\mu^{\frac{1}{256}}\varrho|_{L^{2}}|\mu^{\frac{1}{256}}f|_{L^{2}}.
\eeno
Recalling \eqref{cubic2-inner-product}, patching together the above estimates of $\mathcal{B}_{1,1}, \mathcal{B}_{1,2}, \mathcal{B}_{2,1}, \mathcal{B}_{2,2}$, we have
\ben  \label{Gamma-32-rho-highest-derivative}
|\langle \Gamma_{3,2}^{\rho}(g,h,\varrho), f\rangle| = |\mathcal{B}_{1,1}+\mathcal{B}_{1,2}+ \mathcal{B}_{2,1} + \mathcal{B}_{2,2}| \lesssim |g|_{H^{2}}|\mu^{\frac{1}{256}}h|_{H^{2}} |\varrho|_{\mathcal{L}^{s}_{\gamma/2}} |f|_{\mathcal{L}^{s}_{\gamma/2}}.
\een

Recalling \eqref{Gamma3-into-Gamma31-and-Gamma32},
patching together \eqref{Gamma-31-rho-highest-derivative} and \eqref{Gamma-32-rho-highest-derivative},
we arrive at \eqref{rho-highest-derivative}.
\end{proof}

\begin{proof} [Proof of \eqref{h-highest-derivative} and \eqref{g-highest-derivative}.]
{\it{Estimate of $\langle \Gamma_{3,1}^{\rho}(g,h,\varrho), f\rangle.$}} Recalling \eqref{cubic1-inner-product}, \eqref{Gamma3-1-defintion-A1} and \eqref{Gamma3-1-defintion-A2}. We need give new estimates of $\mathcal{A}_{1}$ and $\mathcal{A}_{2}$.

We first estimate $\mathcal{A}_{1}$. Recalling \eqref{Gamma3-1-defintion-A1},
using $(N\varrho)^{\prime}_{*} = \mathrm{D}((N\varrho)^{\prime}_{*}) + (N\varrho)_{*}$, we get
\ben \label{A1-into-A11-12}
\mathcal{A}_{1} = \tilde{\mathcal{A}}_{1,1} + \tilde{\mathcal{A}}_{1,2},
\quad 
\tilde{\mathcal{A}}_{1,1} \colonequals   \int
B (Ng)_{*} h \mathrm{D}((N\varrho)^{\prime}_{*})  \mathrm{D}(f^{\prime})  \mathrm{d}V,
\quad  
\tilde{\mathcal{A}}_{1,2} \colonequals   \int
B (N^{2}g\varrho)_{*} h  \mathrm{D}(f^{\prime})  \mathrm{d}V.
\een
By Cauchy-Schwartz inequality, \eqref{M-rho-mu-N}, \eqref{mu-star-g-star-difference-h2-f2-star} and Theorem \ref{two-functional-ub-by-norm}, we have
\ben \nonumber
|\tilde{\mathcal{A}}_{1,1}| &\lesssim& \cs{ \int
B (\mu^{\f12}g^{2})_{*} h^{2} \mathrm{D}^{2}((N\varrho)_{*})
\mathrm{d}V}{\int
B  \mu^{\f12}_{*} \mathrm{D}^{2}(f)
\mathrm{d}V}
\\ \label{result-A11} &\lesssim& |\mu^{\frac{1}{128}}g|_{H^{s_{1}}} |h|_{H^{s_{2}}_{\gamma/2+s}} |\mu^{\frac{1}{4}}\varrho|_{H^{3}} |f|_{\mathcal{L}^{s}_{\gamma/2}}.
\een
Observe $|\tilde{\mathcal{A}}_{1,2}| = | \langle  Q_{c}(N^{2}g\varrho, h), f \rangle  |$. Then
by Corollary \ref{Q-up-bound-full-on-first-or-second} and Lemma \ref{product-take-out}, we have
\ben \nonumber
|\tilde{\mathcal{A}}_{1,2}| &\lesssim&
\min\{|N^{2}g\varrho|_{H^{2}_{7}}  |h|_{\mathcal{L}^{s}_{\gamma/2}} , |N^{2}g\varrho|_{L^{2}_{7}}(|h|_{\mathcal{L}^{s}_{\gamma/2}}+|h|_{H^{s+2}_{\gamma/2}})\} |f|_{\mathcal{L}^{s}_{\gamma/2}}.
\\ \label{result-A12-for-g-or-h} && \min\{|\mu^{\frac{1}{4}}g|_{H^{2}}  |h|_{\mathcal{L}^{s}_{\gamma/2}} , |\mu^{\frac{1}{4}}g|_{L^{2}}(|h|_{\mathcal{L}^{s}_{\gamma/2}}+|h|_{H^{s+2}_{\gamma/2}})\} |\mu^{\frac{1}{4}}\varrho|_{H^{2}} |f|_{\mathcal{L}^{s}_{\gamma/2}}.
\een
Recalling \eqref{A1-into-A11-12}, patching together \eqref{result-A11} and \eqref{result-A12-for-g-or-h}, we get
\ben \label{result-A1-for-g-or-h}
|\mathcal{A}_{1}| \lesssim
\min\{|g|_{H^{2}}  |h|_{\mathcal{L}^{s}_{\gamma/2}} , |g|_{L^{2}}(|h|_{\mathcal{L}^{s}_{\gamma/2}}+|h|_{H^{s+2}_{\gamma/2}}+|h|_{H^{\f{1}{2}}_{\gamma/2+s}})\} |\mu^{\frac{1}{4}}\varrho|_{H^{3}} |f|_{\mathcal{L}^{s}_{\gamma/2}}.\een

We then estimate $\mathcal{A}_{2}$. Recalling \eqref{Gamma3-1-defintion-A2}, using \eqref{identity-to-cancel-N-negative-diff} and
$
(N\varrho)^{\prime}_{*}f^{\prime} = \mathrm{D}((N\varrho)^{\prime}_{*})f^{\prime} + (N\varrho)_{*}\mathrm{D}(f^{\prime}) + (N\varrho)_{*}f,
$
we have
\ben \label{A2-into-6-terms}
\mathcal{A}_{2} &=&   -\int
B (\frac{g}{1 - \rho \mu})_{*} \frac{h}{1 - \rho \mu} (N\varrho)^{\prime}_{*} f^{\prime}  \mu^{\f{1}{2}}\mu^{\f{1}{2}}_{*}\mathrm{D}(N^{-1})  \mathrm{d}V
= \sum_{i=1}^{6} \tilde{\mathcal{A}}_{2,i},
\\ \nonumber \tilde{\mathcal{A}}_{2,1} &\colonequals    &  -\int
B (\frac{g}{1 - \rho \mu})_{*} \frac{h}{1 - \rho \mu} \mathrm{D}((N\varrho)^{\prime}_{*})f^{\prime}
\mathrm{D}(\mu^{\f{1}{2}}_{*})  \mathrm{d}V,
\\ \nonumber \tilde{\mathcal{A}}_{2,2} &\colonequals    &  \rho \int
B (\frac{g}{1 - \rho \mu})_{*} \frac{h}{1 - \rho \mu} \mathrm{D}((N\varrho)^{\prime}_{*})f^{\prime}  \mu^{\f{1}{2}}_{*}\mu^{\f{1}{2}}\mathrm{D}(\mu^{\f{1}{2}})  \mathrm{d}V,
\\ \nonumber \tilde{\mathcal{A}}_{2,3} &\colonequals    &  -\int
B (\frac{gN\varrho}{1 - \rho \mu})_{*} \frac{h}{1 - \rho \mu} \mathrm{D}(f^{\prime})
\mathrm{D}(\mu^{\f{1}{2}}_{*})  \mathrm{d}V,
\\ \nonumber \tilde{\mathcal{A}}_{2,4} &\colonequals    &  \rho \int
B (\frac{gN\varrho}{1 - \rho \mu})_{*} \frac{h}{1 - \rho \mu} \mathrm{D}(f^{\prime}) \mu^{\f{1}{2}}_{*}\mu^{\f{1}{2}}\mathrm{D}(\mu^{\f{1}{2}})   \mathrm{d}V,
\\ \nonumber \tilde{\mathcal{A}}_{2,5} &\colonequals    &  -\int
B (\frac{gN\varrho}{1 - \rho \mu})_{*} \frac{hf}{1 - \rho \mu}
\mathrm{D}(\mu^{\f{1}{2}}_{*})  \mathrm{d}V,
\\ \nonumber \tilde{\mathcal{A}}_{2,6} &\colonequals    &  \rho \int
B (\frac{gN\varrho}{1 - \rho \mu})_{*} \frac{hf}{1 - \rho \mu}  \mu^{\f{1}{2}}_{*}\mu^{\f{1}{2}}\mathrm{D}(\mu^{\f{1}{2}})
\mathrm{d}V.
\een
By Cauchy-Schwartz inequality, \eqref{estimate-of-Z-g2-h2-total-one-half}, \eqref{mu-star-g-star-difference-h2} and Remark \ref{star-prime-or-prime-version},
we have
\beno
|\tilde{\mathcal{A}}_{2,1}| \lesssim
\cs{ \int
B g^{2}_{*} h^{2}  \mathrm{D}^{2}(\mu^{\frac{1}{4}}_{*})
\mathrm{d}V}{\int
B \mathrm{D}^{2}((N\varrho)_{*})f^{2} \mathrm{A}(\mu^{\f{1}{2}}_{*})
\mathrm{d}V}
\lesssim |g|_{H^{s_{1}}} |h|_{H^{s_{2}}_{\gamma/2+s}} |\mu^{\frac{1}{4}}\varrho|_{H^{\frac{3}{2}}} |f|_{L^{2}_{\gamma/2+s}}.
\eeno
By nearly the same argument as that for $\tilde{\mathcal{A}}_{2,1}$, thanks to the factor $\mu^{\f{1}{2}}_{*}\mu^{\f{1}{2}}$, we can get
\beno
|\tilde{\mathcal{A}}_{2,2}|  \lesssim |\mu^{\frac{1}{256}}g|_{H^{s_{1}}} |\mu^{\frac{1}{256}}h|_{H^{s_{2}}} |\mu^{\frac{1}{4}}\varrho|_{H^{\frac{3}{2}}} |\mu^{\frac{1}{256}}f|_{L^{2}}.
\eeno
By Cauchy-Schwartz inequality, \eqref{M-rho-mu-N}, the imbedding $H^{2} \hookrightarrow L^{\infty}$,  the estimate \eqref{estimate-of-Z-g2-h2-total-one-half}
and Theorem \ref{two-functional-ub-by-norm},
we have
\beno
|\tilde{\mathcal{A}}_{2,3}| \lesssim
\cs{\int
B (g^{2}\mu^{\f12}\varrho^{2})_{*} h^{2}  \mathrm{D}^{2}(\mu^{\f{1}{2}}_{*})
\mathrm{d}V}{\int
B \mu^{\f12}_{*}\mathrm{D}^{2}(f)
\mathrm{d}V}
\lesssim |\mu^{\frac{1}{8}}g|_{H^{s_{1}}} |h|_{H^{s_{2}}_{\gamma/2+s}} |\mu^{\frac{1}{8}}\varrho|_{H^{2}} |f|_{\mathcal{L}^{s}_{\gamma/2}}.
\eeno
By nearly the same argument as that for $\tilde{\mathcal{A}}_{2,3}$, thanks to the factor $\mu^{\f{1}{2}}_{*}\mu^{\f{1}{2}}$, we get
\beno
|\tilde{\mathcal{A}}_{2,4}|  \lesssim |\mu^{\frac{1}{256}}g|_{H^{s_{1}}} |\mu^{\frac{1}{256}}h|_{H^{s_{2}}} |\mu^{\frac{1}{8}}\varrho|_{H^{2}} |f|_{\mathcal{L}^{s}_{\gamma/2}}.
\eeno
Note that $|\tilde{\mathcal{A}}_{2,5}| = |\mathcal{A}_{2,5}|$ where $\mathcal{A}_{2,5}$ is defined in \eqref{definition-Z25}. Then by taking $( a_{1}, a_{2}, a_{3}) = (\f{1}{2}, 0, 2)$ or $( a_{1}, a_{2}, a_{3}) = (0, \f{1}{2}, 2)$ in \eqref{general-result-without-mu}, we get
\beno
|\tilde{\mathcal{A}}_{2,5}| \lesssim |\mu^{\frac{1}{4}}g|_{H^{s_{1}}} |h|_{H^{s_{2}}_{\gamma/2+s}} |\mu^{\frac{1}{4}}\varrho|_{H^{2}} |f|_{L^{2}_{\gamma/2+s}}.
\eeno
Note that $|\tilde{\mathcal{A}}_{2,6}| = |\mathcal{A}_{2,6}|$ where $\mathcal{A}_{2,6}$ is defined in \eqref{definition-Z26}. Then by taking $( a_{1}, a_{2}, a_{3}) = (\f{1}{2}, 0, 2)$ or $( a_{1}, a_{2}, a_{3}) = (0, \f{1}{2}, 2)$ in \eqref{general-result-with-mu}, we get
\beno
|\tilde{\mathcal{A}}_{2,6}| \lesssim |\mu^{\frac{1}{256}}g|_{H^{s_{1}}} |\mu^{\frac{1}{256}}h|_{H^{s_{2}}} |\mu^{\frac{1}{256}}\varrho|_{H^{2}}  |\mu^{\frac{1}{256}}f|_{L^{2}}.
\eeno
Recalling \eqref{A2-into-6-terms}, patching together the above estimates of $\tilde{\mathcal{A}}_{2,i}(1\leq i\leq 6)$,  we get
\ben \label{result-A2-for-h}
|\mathcal{A}_{2}| \lesssim  |g|_{H^{s_{1}}} |h|_{H^{s_{2}}_{\gamma/2+s}} |\mu^{\frac{1}{256}}\varrho|_{H^{2}} |f|_{\mathcal{L}^{s}_{\gamma/2}}.
\een

Recalling \eqref{cubic1-inner-product},
patching together \eqref{result-A1-for-g-or-h} and \eqref{result-A2-for-h}, we have
\ben \label{Gamma3-1-for-g-or-h}
|\langle \Gamma_{3,1}^{\rho}(g,h,\varrho), f\rangle| \lesssim \min\{|g|_{H^{2}}  |h|_{\mathcal{L}^{s}_{\gamma/2}} , |g|_{L^{2}}(|h|_{\mathcal{L}^{s}_{\gamma/2}}+|h|_{H^{s+2}_{\gamma/2}}+|h|_{H^{\f{1}{2}}_{\gamma/2+s}})\} |\mu^{\frac{1}{256}}\varrho|_{H^{3}} |f|_{\mathcal{L}^{s}_{\gamma/2}}.
\een

{\it{Estimate of $\langle \Gamma_{3,2}^{\rho}(g,h,\varrho), f\rangle.$}} Recall from \eqref{cubic2-inner-product} and
\eqref{Gamma3-2-defintion-B1} that
$\langle \Gamma_{3,2}^{\rho}(g,h,\varrho), f\rangle = \mathcal{B}_{1} + \mathcal{B}_{2}$.

We first deal with $\mathcal{B}_{1}$.
Recall $\mathcal{B}_{1} = \mathcal{B}_{1,1} + \mathcal{B}_{1,2}$ where $\mathcal{B}_{1,1}$ and $\mathcal{B}_{1,2}$ are defined in \eqref{definition-of-B11}. We now give new estimates to $\mathcal{B}_{1,1}$ and $\mathcal{B}_{1,2}$.
By Cauchy-Schwartz inequality, \eqref{M-rho-mu-N}, \eqref{change-v-and-v-star},
\eqref{mu-star-g-star-difference-h2-f2-star-another} and
Theorem \ref{two-functional-ub-by-norm},  we have
\beno
|\mathcal{B}_{1,1}| \lesssim \cs{\int
B \mu^{\f12}g^{2} (\mu h^{2})_{*} \mathrm{D}^{2}(\varrho_{*})
\mathrm{d}V }{\int
B  \mu^{\f12}_{*} \mathrm{D}^{2}(f)
\mathrm{d}V}
\lesssim |\mu^{\frac{1}{256}}g|_{H^{s_{1}}} |\mu^{\frac{1}{256}}h|_{H^{s_{2}}} |\mu^{\frac{1}{256}}\varrho|_{H^{3}}|f|_{\mathcal{L}^{s}_{\gamma/2}}.
\eeno
Observe $|\mathcal{B}_{1,2}| = | \langle  Q_{c}(Ng, Nh\varrho), f \rangle  |$. Then
by Corollary \ref{Q-up-bound-full-on-first-or-second} and Lemma \ref{product-take-out}, we have
\beno
|\mathcal{B}_{1,2}| &\lesssim&
\min\{|Ng|_{H^{2}_{7}}  |Nh\varrho|_{\mathcal{L}^{s}_{\gamma/2}} , |Ng|_{L^{2}_{7}}(|Nh\varrho|_{\mathcal{L}^{s}_{\gamma/2}}+|Nh\varrho|_{H^{s+2}_{\gamma/2}})\} |f|_{\mathcal{L}^{s}_{\gamma/2}}.
 \\&\lesssim& \min\{|g|_{H^{2}}  |h|_{\mathcal{L}^{s}_{\gamma/2}} , |g|_{L^{2}}(|h|_{\mathcal{L}^{s}_{\gamma/2}}+|h|_{H^{s+2}_{\gamma/2}})\} |\mu^{\frac{1}{64}}\varrho|_{H^{3}} |f|_{\mathcal{L}^{s}_{\gamma/2}}.
\eeno

We then consider $\mathcal{B}_{2}$. Recalling \eqref{Gamma3-2-defintion-B1}, using
$
(N\varrho)^{\prime} = N^{\prime}\mathrm{D}(\varrho^{\prime})+\mathrm{D}(N^{\prime})\varrho + N \varrho
$
and the identity \eqref{identity-to-cancel-N-negative-diff},
we have
\ben \label{B2-into-6-terms}
\mathcal{B}_{2} &=&   -\int
B (\frac{g}{1 - \rho \mu})_{*} \frac{hf}{1 - \rho \mu} (N\varrho)^{\prime}
\mu^{\f{1}{2}}\mu^{\f{1}{2}}_{*}\mathrm{D}(N^{-1})  \mathrm{d}V
= \sum_{i=1}^{6}  \tilde{\mathcal{B}}_{2,i},
\\ \nonumber
\tilde{\mathcal{B}}_{2,1} &\colonequals    &  -\int
B (\frac{g}{1 - \rho \mu})_{*} \frac{hf}{1 - \rho \mu} N^{\prime}\mathrm{D}(\varrho^{\prime})
\mathrm{D}(\mu^{\f{1}{2}}_{*})  \mathrm{d}V,
\\ \nonumber \tilde{\mathcal{B}}_{2,2} &\colonequals    &  \rho \int
B (\frac{g}{1 - \rho \mu})_{*} \frac{hf}{1 - \rho \mu} N^{\prime}\mathrm{D}(\varrho^{\prime})   \mu^{\f{1}{2}}_{*}\mu^{\f{1}{2}}\mathrm{D}(\mu^{\f{1}{2}})  \mathrm{d}V,
\\ \nonumber \tilde{\mathcal{B}}_{2,3} &\colonequals    &  -\int
B (\frac{g}{1 - \rho \mu})_{*} \frac{h \varrho f}{1 - \rho \mu} \mathrm{D}(N^{\prime})
\mathrm{D}(\mu^{\f{1}{2}}_{*})  \mathrm{d}V,
\\ \nonumber \tilde{\mathcal{B}}_{2,4} &\colonequals    &  \rho \int
B (\frac{g}{1 - \rho \mu})_{*} \frac{h \varrho f}{1 - \rho \mu} \mathrm{D}(N^{\prime})  \mu^{\f{1}{2}}_{*}\mu^{\f{1}{2}}\mathrm{D}(\mu^{\f{1}{2}})   \mathrm{d}V,
\\ \nonumber \tilde{\mathcal{B}}_{2,5} &\colonequals    &  -\int
B (\frac{g}{1 - \rho \mu})_{*} \frac{N h\varrho f}{1 - \rho \mu}
\mathrm{D}(\mu^{\f{1}{2}}_{*})  \mathrm{d}V,
\\ \nonumber \tilde{\mathcal{B}}_{2,6} &\colonequals    &  \rho \int
B (\frac{g}{1 - \rho \mu})_{*} \frac{N h\varrho f}{1 - \rho \mu}  \mu^{\f{1}{2}}_{*}\mu^{\f{1}{2}}\mathrm{D}(\mu^{\f{1}{2}})
\mathrm{d}V.
\een
By Cauchy-Schwartz inequality, \eqref{M-rho-mu-N}, \eqref{change-v-and-v-star}, \eqref{mu-weight-result}, using  \eqref{estimate-of-Z-g2-h2-total-one-half-with-mu} and \eqref{mu-star-g-star-difference-h2-f2-star-another}(by regarding $\mu^{\frac{1}{8}} = \mu^{\frac{1}{16}} \mu^{\frac{1}{16}}$),
we have
\beno
|\tilde{\mathcal{B}}_{2,1}| &\lesssim&
\cs{ \int
B g^{2}_{*} h^{2} (\mu^{\f{1}{2}})^{\prime} \mathrm{D}^{2}(\mu^{\frac{1}{4}}_{*})
\mathrm{d}V}{\int
B \mathrm{D}^{2}(\varrho)f^{2} (\mu^{\f{1}{2}})^{\prime} \mathrm{A}(\mu^{\f{1}{2}}_{*})
\mathrm{d}V}
\\ &\lesssim&
\cs{ \int
B g^{2}_{*} h^{2} \mu^{\frac{1}{16}}_{*} \mu^{\frac{1}{16}} \mathrm{D}^{2}(\mu^{\f{1}{8}}_{*})
\mathrm{d}V}{\int
B \mathrm{D}^{2}(\varrho_{*}) f^{2}_{*} \mu^{\frac{1}{8}}_{*} \mu^{\frac{1}{8}}
\mathrm{d}V}
\\ &\lesssim& |\mu^{\frac{1}{256}}g|_{H^{s_{1}}} |\mu^{\frac{1}{256}}h|_{H^{s_{2}}} |\mu^{\frac{1}{256}}\varrho|_{H^{3}} |\mu^{\frac{1}{256}}f|_{L^{2}}.
\eeno
By nearly the same argument as that for $\tilde{\mathcal{B}}_{2,1}$, using the factor $\mu^{\f{1}{2}}_{*}\mu^{\f{1}{2}}$, we get
\beno
|\tilde{\mathcal{B}}_{2,2}|  \lesssim |\mu^{\frac{1}{256}}g|_{H^{s_{1}}} |\mu^{\frac{1}{256}}h|_{H^{s_{2}}} |\mu^{\frac{1}{256}}\varrho|_{H^{3}} |\mu^{\frac{1}{256}}f|_{L^{2}}.
\eeno
By Cauchy-Schwartz inequality, \eqref{M-rho-mu-N}, \eqref{change-v-and-v-star}, \eqref{mu-weight-result}, the imbedding $H^{2} \hookrightarrow L^{\infty}$ for $\mu^{\frac{1}{16}}\rho$, Theorem \ref{two-functional-ub-by-norm},
the estimate \eqref{estimate-of-Z-g2-h2-total-one-half-with-mu}
and Remark \ref{mu-to-N-or-M-X-Y},
we have
\beno
|\tilde{\mathcal{B}}_{2,3}| &\lesssim&
\big(\int
B g^{2}_{*} h^{2} \varrho^{2} \mathrm{D}^{2}(N^{\f{1}{2}}) \mathrm{A}(N)\mathrm{A}(\mu^{\f{1}{2}}_{*})
\mathrm{d}V\big)^{\f{1}{2}}
\big(\int
B f^{2}  \mathrm{D}^{2}(\mu^{\frac{1}{4}}_{*})
\mathrm{d}V\big)^{\f{1}{2}}
\\ &\lesssim&
\cs{\int
B g^{2} h^{2}_{*} \varrho^{2}_{*} \mu^{\frac{1}{8}}_{*} \mu^{\frac{1}{8}} \mathrm{D}^{2}(N^{\frac{1}{2}}_{*})
\mathrm{d}V}{\int
B f^{2}_{*}  \mathrm{D}^{2}(\mu^{\frac{1}{4}})
\mathrm{d}V}
\lesssim |\mu^{\frac{1}{256}}g|_{H^{s_{1}}} |\mu^{\frac{1}{256}}h|_{H^{s_{2}}} |\mu^{\frac{1}{16}}\varrho|_{H^{2}} |f|_{\mathcal{L}^{s}_{\gamma/2}}.
\eeno
By nearly the same argument as that for $\tilde{\mathcal{B}}_{2,3}$, using the factor $\mu^{\f{1}{2}}_{*}\mu^{\f{1}{2}}$, we get
\beno
|\tilde{\mathcal{B}}_{2,4}|  \lesssim |\mu^{\frac{1}{256}}g|_{H^{s_{1}}} |\mu^{\frac{1}{256}}h|_{H^{s_{2}}} |\mu^{\frac{1}{256}}\varrho|_{H^{2}} |f|_{\mathcal{L}^{s}_{\gamma/2}}.
\eeno
By \eqref{estimate-of-B-g-h-varrho-f-with-mu},
we have
$
|\tilde{\mathcal{B}}_{2,5}| \lesssim
|g|_{H^{s_{1}}} |\mu^{\frac{1}{256}}h|_{H^{s_{2}}} |\mu^{\frac{1}{256}}\varrho|_{H^{2}}  |\mu^{\frac{1}{256}}f|_{L^{2}}.
$
Using \eqref{change-v-and-v-star} and \eqref{estimate-of-C-g-h-varrho-f-with-mu}, we get
$
|\tilde{\mathcal{B}}_{2,6}| \lesssim |g|_{H^{s_{1}}} |\mu^{\frac{1}{256}}h|_{H^{s_{2}}} |\mu^{\frac{1}{256}}\varrho|_{H^{2}}  |\mu^{\frac{1}{256}}f|_{L^{2}}.
$
Recalling $\langle \Gamma_{3,2}^{\rho}(g,h,\varrho), f\rangle = \mathcal{B}_{1} + \mathcal{B}_{2} = \mathcal{B}_{1,1} + \mathcal{B}_{1,2} + \sum_{i=1}^{6}  \tilde{\mathcal{B}}_{2,i}$. Patching together the above estimates of $\mathcal{B}_{1,1}, \mathcal{B}_{1,2}, \tilde{\mathcal{B}}_{2,i}(1 \leq i \leq 6)$, we have
\ben \label{Gamma3-2-for-g-or-h}
|\langle \Gamma_{3,2}^{\rho}(g,h,\varrho), f\rangle| \lesssim \min\{|g|_{H^{2}}  |h|_{\mathcal{L}^{s}_{\gamma/2}} , |g|_{L^{2}}(|h|_{\mathcal{L}^{s}_{\gamma/2}}+|h|_{H^{s+2}_{\gamma/2}})\} |\mu^{\frac{1}{256}}\varrho|_{H^{3}} |f|_{\mathcal{L}^{s}_{\gamma/2}}.
\een

Recalling \eqref{Gamma3-into-Gamma31-and-Gamma32},
patching together \eqref{Gamma3-1-for-g-or-h} and  \eqref{Gamma3-2-for-g-or-h}, we get \eqref{h-highest-derivative} and \eqref{g-highest-derivative}.
\end{proof}

\section{Commutator estimate and weighted energy estimate}  \label{commutator}

In this section, we derive some estimates for the commutators between the weight function $W_{l}$ and the various linear, bilinear and trilinear operators whose un-weighted estimates are already given in Section \ref{linear}, \ref{bilinear} and \ref{trilinear}. These estimates together produce some weighted energy estimates that will be used in the next two sections.

Direct computation shows that the various commutators share a common term $\mathrm{D}(W_{l})$. We collect estimates
of various functionals involving  $\mathrm{D}(W_{l})$ in the following lemma.
In all the estimates for functionals involving $W_{l}$ in the article, the ``$\lesssim$'' could bring a constant $C_{l}$ depending on $l$ on the righthand side of ``$\lesssim$''. We do not specify this dependence for brevity.

\begin{lem}\label{full-integral-Wl-difference-g-h-f} Let $l \geq 0$. Let $\{a_{1},a_{2},a_{3}\} = \{0, \f12, 2\}$ and $s_{1},s_{2} \geq 0,s_{1}+s_{2}=\f{1}{2}$. The following estimates are valid.
\ben \label{Wl-difference-square-mu-star-g2-h2}
\int B\mathrm{D}^{2}(W_{l})\mu^{\frac{1}{4}}_{*} g^{2}_{*} h^{2}  \mathrm{d}V  &\lesssim& |\mu^{\frac{1}{64}}g|_{L^{2}}^{2} |h|_{L^{2}_{l+\gamma/2}}^{2} + |\mu^{\frac{1}{64}}g|_{H^{s_{1}}}^{2} |\mu^{\frac{1}{64}}h|_{H^{s_{2}}}^{2}.
\\ \label{Wl-difference-square-mu-star-prime-g2-h2}
\int B\mathrm{D}^{2}(W_{l})(\mu^{\frac{1}{4}})^{\prime}_{*} g^{2}_{*} h^{2}  \mathrm{d}V  &\lesssim& |g|_{L^{2}}^{2} |h|_{L^{2}_{l+\gamma/2}}^{2} + |\mu^{\frac{1}{64}}g|_{H^{s_{1}}}^{2} |\mu^{\frac{1}{64}}h|_{H^{s_{2}}}^{2}.
\\ \label{Wl-difference-order-1-mu-star-ghf}
|\int B\mathrm{D}(W_{l}) \mu^{\frac{1}{4}}_{*}g_{*} h f \mathrm{d}V|  &\lesssim& |\mu^{\frac{1}{64}}g|_{L^{2}} |h|_{L^{2}_{l+\gamma/2}} |f|_{L^{2}_{\gamma/2}} + |\mu^{\frac{1}{64}}g|_{H^{s_{1}}} |\mu^{\frac{1}{64}}h|_{H^{s_{2}}} |\mu^{\frac{1}{64}}f|_{L^{2}}.
\\ \label{Wl-difference-order-mu-star-g-hf-prime}
|\int B\mathrm{D}(W_{l}) \mu^{\frac{1}{4}}_{*}g_{*} (h f)^{\prime} \mathrm{d}V| &\lesssim& |\mu^{\frac{1}{64}}g|_{L^{2}} |h|_{L^{2}_{l+\gamma/2}} |f|_{L^{2}_{\gamma/2}}.
\\ \label{Wl-difference-order-mu-star-g-hfvarrho-prime}
|\int
B g_{*} h \varrho f \mathrm{D}(W_{l}) \mu^{\frac{1}{4}} \mu^{\frac{1}{4}}_{*} \mathrm{d}V| &\lesssim& |\mu^{\frac{1}{64}}g|_{H^{a_{1}}} |\mu^{\frac{1}{64}}h|_{H^{a_{2}}} |\mu^{\frac{1}{64}}\varrho|_{H^{a_{3}}} |\mu^{\frac{1}{64}}f|_{L^{2}}.
\\ \label{Wl-difference-order-mu-star-gvarrho-hf-prime}
|\int
B g_{*} h \varrho_{*} f \mathrm{D}(W_{l}) \mu^{\frac{1}{4}}_{*} \mathrm{d}V| &\lesssim& |\mu^{\frac{1}{64}}g|_{H^{a_{1}}} |h|_{H^{a_{2}}_{l+\gamma/2}} |\mu^{\frac{1}{64}}\varrho|_{H^{a_{3}}} |f|_{L^{2}_{\gamma/2}}.
\een
\end{lem}
The proof of Lemma \ref{full-integral-Wl-difference-g-h-f} is given in the Appendix \ref{appendix}. Based on the proof,
we give a remark to make Lemma \ref{full-integral-Wl-difference-g-h-f} more applicable.
\begin{rmk}\label{still-true-h-prime} In Lemma \ref{full-integral-Wl-difference-g-h-f}, if the exponent $\f{1}{4}$ is changed to some constant $a \geq \frac{1}{4}$, the results are still valid. If we replace $\mu$ with $N^{2}$ in \eqref{Wl-difference-order-1-mu-star-ghf}, \eqref{Wl-difference-order-mu-star-g-hf-prime}, \eqref{Wl-difference-order-mu-star-g-hfvarrho-prime} and \eqref{Wl-difference-order-mu-star-gvarrho-hf-prime}, the results are still valid.
Estimates \eqref{Wl-difference-square-mu-star-g2-h2} and \eqref{Wl-difference-square-mu-star-prime-g2-h2} are still valid if $h^{2}$ is replaced by $(h^{2})^{\prime}$.
\end{rmk}

\subsection{Commutator estimate between $W_{l}$ and $\mathcal{L}^{\rho}$}
In the following proposition, we give an estimate of the commutator between $W_{l}$ and $\mathcal{L}^{\rho}$.
\begin{prop}\label{commutator-linear}
Let $l\geq 0$. The following two estimates are valid.
\ben \label{commutator-Wl-f-f}
|\langle [W_{l}, \mathcal{L}^{\rho}]f, W_{l}f\rangle| \lesssim \rho |f|_{L^{2}_{l+\gamma/2}}^{2}.
\\ \label{commutator-Wl-g-f}
|\langle [W_{l}, \mathcal{L}^{\rho}]h, W_{l}f\rangle| \lesssim \rho |h|_{L^{2}_{l+\gamma/2}}|f|_{\mathcal{L}^{s}_{l+\gamma/2}}.
\een
\end{prop}
\begin{proof}
Recalling \eqref{linearized-quantum-Boltzmann-operator-UU}, \eqref{M--with-mathcal-M} and \eqref{L-epsilon-pm-tau-rest-part}, we have
\ben \label{L-into-two-sums}
\mathcal{L}^{\rho}f &=& \mathcal{L}^{\rho}_{m}f + \mathcal{L}^{\rho}_{r}f,
\\ \label{definition-of-lm}
(\mathcal{L}^{\rho}_{m}f)(v)  &\colonequals&    \rho \int  B N_{*} N^{\prime} N^{\prime}_{*}
\mathrm{D}(N^{-1}f)
\mathrm{d}\sigma \mathrm{d}v_{*}.
\een
By \eqref{L-into-two-sums}, we need to consider $\mathcal{L}^{\rho}_{m}$ and $\mathcal{L}^{\rho}_{r}$. For the operator $\mathcal{L}^{\rho}_{m}$, using \eqref{change-v-and-v-prime},
we derive
\beno
\langle [W_{l}, \mathcal{L}^{\rho}_{m}]f, W_{l}f\rangle = \rho \int B N_{*}
N^{\prime}_{*} f f^{\prime} W_{l}\mathrm{D}(W_{l}^{\prime})  \mathrm{d}V
= -\f{1}{2} \rho \int B N_{*}
N^{\prime}_{*} f f^{\prime} \mathrm{D}^{2}(W_{l})  \mathrm{d}V.
\eeno
Using $N_{*} N^{\prime}_{*} |f f^{\prime}| \leq \f{1}{2}(N_{*}^{2}f^{2}+(N^{2})^{\prime}_{*}(f^{2})^{\prime})$, \eqref{change-v-and-v-prime}, \eqref{M-rho-mu-N},
we get
\beno
|\langle [W_{l}, \mathcal{L}^{\rho}_{m}]f, W_{l}f\rangle| \leq  \f{1}{2} \rho \int B N_{*}^{2}  f^{2} \mathrm{D}^{2}(W_{l})  \mathrm{d}V
\lesssim   \rho \int B \mu_{*}
 f^{2} \mathrm{D}^{2}(W_{l})  \mathrm{d}V \lesssim \rho |f|_{L^{2}_{l+\gamma/2}}^{2}.
\eeno
where we use \eqref{Wl-difference-square-mu-star-g2-h2}(writing $\mu_{*}=\mu^{\frac{1}{4}}_{*}\mu^{\frac{3}{4}}_{*}$).
By Proposition \ref{ub-linearized-L2} and \eqref{deal-with-polynomial-weight}, we have
\ben \label{L-r-commutator-gh}
|\langle [W_{l}, \mathcal{L}^{\rho}_{r}]h, W_{l}f\rangle|
 \lesssim \rho |\mu^{\f{1}{128}}h|_{L^{2}} |\mu^{\f{1}{128}}f|_{L^{2}}.
\een
Patching together the previous two estimates, we get \eqref{commutator-Wl-f-f}.

We now set to prove \eqref{commutator-Wl-g-f}.
Recalling \eqref{L-into-two-sums}, we need to consider $\mathcal{L}^{\rho}_{m}$ and $\mathcal{L}^{\rho}_{r}$. For the operator $\mathcal{L}^{\rho}_{r}$, we already have \eqref{L-r-commutator-gh}. For the operator $\mathcal{L}^{\rho}_{m}$, recalling \eqref{definition-of-lm}, using \eqref{change-v-and-v-prime}, we deduce
\beno
\langle [W_{l}, \mathcal{L}^{\rho}_{m}]h, W_{l}f\rangle = -\rho \int B N_{*}
N^{\prime}_{*} f h^{\prime} W_{l}\mathrm{D}(W_{l})  \mathrm{d}V
= \rho \int B N_{*}
N^{\prime}_{*} h (W_{l}f)^{\prime}\mathrm{D}(W_{l})  \mathrm{d}V = \mathcal{I}_{1}+\mathcal{I}_{2},
\\
\mathcal{I}_{1} \colonequals     \rho \int B N^{2}_{*} h (W_{l}f)^{\prime}\mathrm{D}(W_{l})  \mathrm{d}V,
\quad
\mathcal{I}_{2} \colonequals     \rho \int B N_{*}
\mathrm{D}(N^{\prime}_{*})h (W_{l}f)^{\prime}\mathrm{D}(W_{l})  \mathrm{d}V.
\eeno
Note that $\mathcal{I}_{1} = -\rho \langle W_{l}Q_{c} (N^{2}, h) - Q_{c} (N^{2}, W_{l}h), W_{l}f\rangle$ by using  \eqref{commutator-inner-product-Q-epsilon-pm-0}. Then by Proposition \ref{Commutotator-Wl-and-Q-epsilon-pm-0}(which will be proved soon), we have
$
|\mathcal{I}_{1}| \lesssim \rho |h|_{L^{2}_{l+\gamma/2}}|f|_{\mathcal{L}^{s}_{l+\gamma/2}}.
$
By Cauchy-Schwartz inequality, \eqref{M-rho-mu-N}, \eqref{change-v-and-v-prime} and \eqref{change-v-and-v-star},
we have
\beno
|\mathcal{I}_{2}| \lesssim \rho \cs{\int B \mu_{*}
h^{2} \mathrm{D}^{2}(W_{l})  \mathrm{d}V}{\int B
\mathrm{D}^{2}(N) (W_{l}f)^{2}_{*}  \mathrm{d}V} \lesssim \rho |h|_{L^{2}_{l+\gamma/2}}|f|_{\mathcal{L}^{s}_{l+\gamma/2}}.
\eeno
where we use \eqref{Wl-difference-square-mu-star-g2-h2}, Theorem \ref{two-functional-ub-by-norm} and Remark \ref{mu-to-N-or-M-functional-N}. Patching together the estimates of $\mathcal{I}_{1}, \mathcal{I}_{2}$
and \eqref{L-r-commutator-gh},
we get \eqref{commutator-Wl-g-f}.
\end{proof}

Applying Theorem \ref{main-theorem} and Proposition \ref{commutator-linear}, we get
\begin{thm}
Let $l\geq 0$. If $0 \leq \rho \leq \rho_{0}$, there holds
\ben \label{wiehgted-coercivity-lb}
\langle \mathcal{L}^{\rho}f, W_{2l}f\rangle \geq  \frac{\lambda_{0}}{4} \rho |f|_{\mathcal{L}^{s}_{l+\gamma/2}}^{2} - C \rho |f|_{L^{2}_{\gamma/2}}^{2}, \een
where $\lambda_{0}$ is the constant appearing in Theorem \ref{main-theorem}. As a direct consequence, in the full $(x,v)$ space, there holds
\ben \label{coercivity-ub-and-lb-weight-xv}
(\mathcal{L}^{\rho}f, W_{2l}f)  \geq \frac{\lambda_{0}}{4} \rho \|f\|_{L^{2}_{x}\mathcal{L}^{s}_{l+\gamma/2}}^{2} - C \rho \|f\|_{L^{2}_{x}L^{2}_{\gamma/2}}^{2}.
\een
\end{thm}
\begin{proof} Recalling \eqref{linear-combination-of-basis} and \eqref{explicit-defintion-of-abc}, since $N \lesssim \mu^{\f12}$,
it is easy to see $|\mathbb{P}_{\rho}f|_{\mathcal{L}^{s}_{\gamma/2}} \lesssim |\mu^{\f18} f|_{L^{2}} \lesssim |f|_{L^{2}_{\gamma/2}}$. Then by the lower bound in Theorem \ref{main-theorem}, we have
\beno
\langle \mathcal{L}^{\rho}f , f \rangle \geq  \lambda_{0} \rho  (\f{1}{2}|f|_{\mathcal{L}^{s}_{\gamma/2}}^{2} - |\mathbb{P}_{\rho}f|_{\mathcal{L}^{s}_{\gamma/2}}^{2})
\geq \f{1}{2} \lambda_{0}\rho |f|_{\mathcal{L}^{s}_{\gamma/2}}^{2} - C |f|_{L^{2}_{\gamma/2}}^{2}.
\eeno
From which together with Proposition \ref{commutator-linear}, we have
\beno
\langle \mathcal{L}^{\rho}f, W_{2l}f\rangle = \langle \mathcal{L}^{\rho}W_{l}f, W_{l}f\rangle + \langle [W_{l}, \mathcal{L}^{\rho}]f, W_{l}f\rangle \geq  \f{1}{2} \lambda_{0} \rho |f|_{\mathcal{L}^{s}_{l+\gamma/2}}^{2} - C \rho |f|_{L^{2}_{l+\gamma/2}}^{2}.
\eeno
Using the interpolation $|f|_{L^{2}_{l+\gamma/2}}^{2} \leq \eta |f|_{L^{2}_{l+\gamma/2+s}}^{2} + C_{\eta} |f|_{L^{2}_{\gamma/2}}^{2} \leq \eta |f|_{\mathcal{L}^{s}_{l+\gamma/2}}^{2} + C_{\eta} |f|_{L^{2}_{\gamma/2}}^{2}$ and taking $\eta = \frac{\lambda_{0}}{4C}$, we finish the proof of \eqref{wiehgted-coercivity-lb}. Further taking integration over $x \in \mathbb{T}^{3}$, we get \eqref{coercivity-ub-and-lb-weight-xv}.
\end{proof}

In the $(x,v)$ space, we have the following lower bound estimates which yield the dissipation functional in later energy estimates.
\begin{thm} \label{L-energy-estimate} Recall $\partial^{\alpha} = \partial^{\alpha}_{x}, \partial^{\alpha}_{\beta} = \partial^{\alpha}_{x}\partial^{\beta}_{v}$. If $0 \leq \rho \leq \rho_{0}$, the following statements are valid. When there is only $x$ derivative, it holds that
\ben \label{no-v-derivative}
(\partial^{\alpha} \mathcal{L}^{\rho}f, W_{2l_{|\alpha|,0}} \partial^{\alpha}f) \geq
\frac{\lambda_{0}}{4}\rho\|\partial^{\alpha}f\|_{L^{2}_{x}\mathcal{L}^{s}_{l_{|\alpha|,0}+\gamma/2}}^{2} - C \rho
\|\partial^{\alpha}f\|_{L^{2}_{x}L^{2}_{\gamma/2}}^{2}.
\een
If $|\beta| \geq 1$, it holds that
\ben \label{with-v-derivative}
&& (\partial^{\alpha}_{\beta} \mathcal{L}^{\rho}f, W_{2l_{|\alpha|,|\beta|}} \partial^{\alpha}_{\beta}f)
\\ \nonumber &\geq&
\frac{\lambda_{0}}{8}\rho\|\partial^{\alpha}_{\beta}f\|_{L^{2}_{x}\mathcal{L}^{s}_{l_{|\alpha|,|\beta|}+\gamma/2}}^{2} -
C \rho \|\partial^{\alpha}_{\beta}f\|_{L^{2}_{x}L^{2}_{l_{|\alpha|,|\beta|}+\gamma/2}}^{2} - C \rho
\sum_{\beta_{1}<\beta} \|\partial^{\alpha}_{\beta_{1}}f\|_{L^{2}_{x}\mathcal{L}^{s}_{l_{|\alpha|,|\beta|}+\gamma/2}}^{2}.
\een
\end{thm}
\begin{proof} Since $\partial^{\alpha} \mathcal{L}^{\rho}f =  \mathcal{L}^{\rho}\partial^{\alpha}f$,
\eqref{no-v-derivative} is a direct result of \eqref{coercivity-ub-and-lb-weight-xv}.

Note that
$
(\partial^{\alpha}_{\beta} \mathcal{L}^{\rho}f, W_{2l_{|\alpha|,|\beta|}} \partial^{\alpha}_{\beta}f) =
\mathcal{I}_{1} + \mathcal{I}_{2}
$
where
\beno
\mathcal{I}_{1} \colonequals    (\mathcal{L}^{\rho} \partial^{\alpha}_{\beta}f, W_{2l_{|\alpha|,|\beta|}} \partial^{\alpha}_{\beta}f), \quad
\mathcal{I}_{2}\colonequals
([\partial_{\beta}, \mathcal{L}^{\rho}] \partial^{\alpha}f, W_{2l_{|\alpha|,|\beta|}} \partial^{\alpha}_{\beta}f).
\eeno
We use \eqref{coercivity-ub-and-lb-weight-xv} to get
\beno
\mathcal{I}_{1}  \geq  \frac{\lambda_{0}}{4} \rho   \|\partial^{\alpha}_{\beta}f\|_{L^{2}_{x}\mathcal{L}^{s}_{l_{|\alpha|,|\beta|}+\gamma/2}}^{2} - C \rho \|\partial^{\alpha}_{\beta}f\|_{L^{2}_{x}L^{2}_{l_{|\alpha|,|\beta|}+\gamma/2}}^{2}.
\eeno
By binomial formula and Remark \ref{sufficient-to-consider-123}, to estimate $([\partial_{\beta} ,\mathcal{L}^{\rho}] \partial^{\alpha}f, W_{2l_{|\alpha|,|\beta|}} \partial^{\alpha}_{\beta}f)$, it suffices to consider
\beno
(\mathcal{L}^{\rho} \partial^{\alpha}_{\beta_{1}}f, W_{2l_{|\alpha|,|\beta|}} \partial^{\alpha}_{\beta}f) =
(\mathcal{L}^{\rho} W_{l_{|\alpha|,|\beta|}} \partial^{\alpha}_{\beta_{1}} f, W_{l_{|\alpha|,|\beta|}} \partial^{\alpha}_{\beta}f) +
([W_{l_{|\alpha|,|\beta|}},\mathcal{L}^{\rho}] \partial^{\alpha}_{\beta_{1}}f, W_{l_{|\alpha|,|\beta|}} \partial^{\alpha}_{\beta}f),
\eeno
for all $\beta_{1} < \beta$.

Recalling \eqref{linear-combination-of-basis} and \eqref{explicit-defintion-of-abc},
it is easy to see $|\mathbb{P}_{\rho}f|_{\mathcal{L}^{s}_{\gamma/2}} \lesssim |\mu^{\f18} f|_{L^{2}} \lesssim |f|_{\mathcal{L}^{s}_{\gamma/2}}$. Then by the upper bound in Theorem \ref{main-theorem} and Remark \ref{up-no-small-on-rho}, for $\rho \leq\f{1}{2}(2 \pi)^{\frac{3}{2}}$,
 we have
\ben \label{rough-upper-bound-of-L}
\langle \mathcal{L}^{\rho}f , f \rangle \lesssim  \rho  (|f|_{\mathcal{L}^{s}_{\gamma/2}}^{2} + |\mathbb{P}_{\rho}f|_{\mathcal{L}^{s}_{\gamma/2}}^{2}) \lesssim  \rho  |f|_{\mathcal{L}^{s}_{\gamma/2}}^{2}.
\een
Recalling \eqref{self-joint-operator}, using Cauchy-Schwartz inequality and \eqref{rough-upper-bound-of-L}, for $\rho \leq\f{1}{2}(2 \pi)^{\frac{3}{2}}$,
we have
\ben \label{upper-bound-of-L}
|\langle \mathcal{L}^{\rho} g , h \rangle| \leq \cs{\langle \mathcal{L}^{\rho} g , g \rangle}{\langle \mathcal{L}^{\rho} h , h \rangle}
 \lesssim \rho |g|_{\mathcal{L}^{s}_{\gamma/2}} |h|_{\mathcal{L}^{s}_{\gamma/2}},
\een
which gives
\ben \label{ub-estimate-L-gh}
|(\mathcal{L}^{\rho} W_{l_{|\alpha|,|\beta|}} \partial^{\alpha}_{\beta_{1}} f, W_{l_{|\alpha|,|\beta|}} \partial^{\alpha}_{\beta}f)| \lesssim \rho\|\partial^{\alpha}_{\beta_{1}}f\|_{L^{2}_{x}\mathcal{L}^{s}_{l_{|\alpha|,|\beta|}+\gamma/2}}
\|\partial^{\alpha}_{\beta}f\|_{L^{2}_{x}\mathcal{L}^{s}_{l_{|\alpha|,|\beta|}+\gamma/2}}.
\een
By \eqref{commutator-Wl-g-f}, we have
\ben \label{commutator-estimate-L}
|([W_{l_{|\alpha|,|\beta|}},\mathcal{L}^{\rho}] \partial^{\alpha}_{\beta_{1}}f, W_{l_{|\alpha|,|\beta|}} \partial^{\alpha}_{\beta}f)| \lesssim \rho\|\partial^{\alpha}_{\beta_{1}}f\|_{L^{2}_{x}L^{2}_{l_{|\alpha|,|\beta|}+\gamma/2}}
\|\partial^{\alpha}_{\beta}f\|_{L^{2}_{x}\mathcal{L}^{s}_{l_{|\alpha|,|\beta|}+\gamma/2}}.
\een
Patching together \eqref{ub-estimate-L-gh} and \eqref{commutator-estimate-L}, using the basic inequality $2ab \leq  \eta a^{2} + \eta^{-1} b^{2}$, we have
\beno
|(\mathcal{L}^{\rho} \partial^{\alpha}_{\beta_{1}}f, W_{2l_{|\alpha|,|\beta|}} \partial^{\alpha}_{\beta}f)| \lesssim
\eta \rho \|\partial^{\alpha}_{\beta}f\|_{L^{2}_{x}\mathcal{L}^{s}_{l_{|\alpha|,|\beta|}+\gamma/2}}^{2} + \eta^{-1} \rho
\|\partial^{\alpha}_{\beta_{1}}f\|_{L^{2}_{x}\mathcal{L}^{s}_{l_{|\alpha|,|\beta|}+\gamma/2}}^{2}.
\eeno
Taking sum over  $\beta_{1} < \beta$ and taking $\eta$ suitably small, we finish the proof.
\end{proof}

\subsection{Estimate of the commutator $[W_{l},\Gamma_{2}^{\rho}(g, \cdot)]$} In this subsection,
we derive an estimate of the commutator between $W_{l}$ and the operator $\Gamma_{2}^{\rho}(g,\cdot)$. That is, we want to estimate
$
\langle W_{l}\Gamma_{2}^{\rho}(g,h) - \Gamma_{2}^{\rho}(g, W_{l}h), f\rangle.
$
Recalling \eqref{Gamma-main-remaining}, \eqref{definition-Gamma-main} and \eqref{Gamma-remaining-into-three-terms}, we have
\ben \label{G2-and-into-five-terms}
\Gamma_{2}^{\rho}(g,h) = \rho^{\f12}\big(Q_{c} (Ng, h) + I^{\rho}(g,h)\big) + \rho^{\f32}\big(\Gamma_{2,r,1}^{\rho}(g,h) + \Gamma_{2,r,2}^{\rho}(g,h) + \Gamma_{2,r,3}^{\rho}(g,h)\big).
\een
We will estimate the five terms on the right-hand side of \eqref{G2-and-into-five-terms} in the following five propositions. First, we give an estimate for the commutator $[W_{l}, Q_{c} (Ng, \cdot)]$ in the following proposition.
\begin{prop}\label{Commutotator-Wl-and-Q-epsilon-pm-0} Let $l, s_{1}, s_{2} \geq 0$  and  $s_{1} + s_{2} = \f{1}{2}$. It holds that
\beno
|\langle W_{l}Q_{c} (Ng, h) - Q_{c} (Ng, W_{l}h), f\rangle|  \lesssim (|\mu^{\frac{1}{64}}g|_{L^{2}} |h|_{L^{2}_{l+\gamma/2}} + |\mu^{\frac{1}{64}}g|_{H^{s_{1}}} |\mu^{\frac{1}{64}}h|_{H^{s_{2}}})|f|_{\mathcal{L}^{s}_{\gamma/2}}.
\eeno
\end{prop}
\begin{proof}
We observe that
\ben \label{commutator-inner-product-Q-epsilon-pm-0}
\langle W_{l}Q_{c} (Ng, h) - Q_{c} (Ng, W_{l}h), f\rangle  = \int B\mathrm{D}(W_{l}^{\prime})(Ng)_{*} h f^{\prime} \mathrm{d}V = \mathcal{I}_{1} + \mathcal{I}_{2},
\\  \nonumber
\mathcal{I}_{1}\colonequals     \int B\mathrm{D}(W_{l}^{\prime})(Ng)_{*} h \mathrm{D}(f^{\prime}) \mathrm{d}V, \quad
 \mathcal{I}_{2}\colonequals     \int B\mathrm{D}(W_{l}^{\prime})(Ng)_{*} h f \mathrm{d}V. \een
By Cauchy-Schwartz inequality and \eqref{M-rho-mu-N}, using \eqref{Wl-difference-square-mu-star-g2-h2} and Theorem \ref{two-functional-ub-by-norm}, we have
\beno  |\mathcal{I}_{1}| &\lesssim& \cs
{\int B\mathrm{D}^{2}(W_{l})\mu^{\f{1}{2}}_{*}g^{2}_{*} h^{2}  \mathrm{d}V}
{\int B \mu^{\f{1}{2}}_{*}\mathrm{D}^{2}(f) \mathrm{d}V}
\\  &\lesssim& (|\mu^{\frac{1}{64}}g|_{L^{2}} |h|_{L^{2}_{l+\gamma/2}} + |\mu^{\frac{1}{64}}g|_{H^{s_{1}}} |\mu^{\frac{1}{64}}h|_{H^{s_{2}}})|f|_{\mathcal{L}^{s}_{\gamma/2}}.
\eeno
By \eqref{Wl-difference-order-1-mu-star-ghf} and Remark \ref{still-true-h-prime}, we have
\beno
|\mathcal{I}_{2}| \lesssim |\mu^{\frac{1}{64}}g|_{L^{2}} |h|_{L^{2}_{l+\gamma/2}} |f|_{L^{2}_{\gamma/2}} + |\mu^{\frac{1}{64}}g|_{H^{s_{1}}} |\mu^{\frac{1}{64}}h|_{H^{s_{2}}} |\mu^{\frac{1}{64}}f|_{L^{2}}.
\eeno
Recalling \eqref{commutator-inner-product-Q-epsilon-pm-0}, patching together the estimates of $\mathcal{I}_{1}$ and $\mathcal{I}_{2}$, we finish the proof.
\end{proof}

The following proposition gives an estimate for the commutator $[W_{l},I^{\rho}(g,\cdot)]$.
\begin{prop}\label{Commutotator-Wl-and-I-epsilon-pm-tau}
	Let $l, s_{1}, s_{2} \geq 0$  and  $s_{1} + s_{2} = \f{1}{2}$. It holds that
		\beno
		|\langle W_{l}I^{\rho}(g,h) - I^{\rho}(g,W_{l}h), f\rangle| \lesssim
(|g|_{L^{2}}|h|_{L^{2}_{l+\gamma/2}}+|\mu^{\frac{1}{64}}g|_{H^{s_{1}}}|\mu^{\frac{1}{64}}h|_{H^{s_{2}}})
|f|_{\mathcal{L}^{s}_{\gamma/2}}.
		\eeno	
\end{prop}
\begin{proof}
Recalling \eqref{I-into-main-and-rest}, \eqref{definition-Im-epsilon}, \eqref{definition-Ir-epsilon}, we have
\ben \label{commutator-I-into-main-and-rest}
\langle [W_{l},I^{\rho}(g,\cdot)]h, f\rangle &=& \langle [W_{l},I^{\rho}_{m}(g,\cdot)]h, f\rangle+\langle [W_{l},I^{\rho}_{r}(g,\cdot)]h, f\rangle,
\\ \label{commutator-I-epsilon-pm-tau-m}
\langle [W_{l},I^{\rho}_{m}(g,\cdot)]h, f\rangle &=& -\int B
\mathrm{D}(\mu^{\f{1}{2}}_{*}) \mathrm{D}(W_{l}^{\prime}) (\frac{g}{1 - \rho \mu})_{*} (\frac{h}{1 - \rho \mu}) f^{\prime} \mathrm{d}V,
\\ \nonumber 
\langle [W_{l},I^{\rho}_{r}(g,\cdot)]h, f\rangle &=&  \rho \int B
\mathrm{D}(\mu^{\f{1}{2}}) \mu^{\f{1}{2}} \mu^{\f{1}{2}}_{*}\mathrm{D}(W_{l}^{\prime}) (\frac{g}{1 - \rho \mu})_{*} (\frac{h}{1 - \rho \mu}) f^{\prime} \mathrm{d}V.
\een
By Cauchy-Schwartz inequality, \eqref{change-v-and-v-star} and \eqref{change-v-and-v-prime},
using \eqref{Wl-difference-square-mu-star-g2-h2}, \eqref{Wl-difference-square-mu-star-prime-g2-h2} and Theorem \ref{two-functional-ub-by-norm}, we have
\ben \nonumber |\langle [W_{l},I^{\rho}_{m}(g,\cdot)]h, f\rangle| &\lesssim&
\cs{ \int  \mathrm{D}^{2}(W_{l}) \mathrm{A}(\mu^{\f{1}{2}}_{*}) g^{2}_{*} h^{2}
\mathrm{d}V}{
\int f^{2}_{*} \mathrm{D}^{2}(\mu^{\frac{1}{4}})
\mathrm{d}V}
\\\label{ub-commutator-main-part} &\lesssim& (|g|_{L^{2}} |h|_{L^{2}_{l+\gamma/2}}+|\mu^{\frac{1}{64}}g|_{H^{s_{1}}} |\mu^{\frac{1}{64}}h|_{H^{s_{2}}}) |f|_{\mathcal{L}^{s}_{\gamma/2}}.
\een
By Cauchy-Schwartz inequality, \eqref{change-v-and-v-star} and \eqref{change-v-and-v-prime}, using \eqref{Wl-difference-square-mu-star-g2-h2} and \eqref{estimate-of-Z-g2-h2-total-one-half-with-mu}, we have
\ben \nonumber
|\langle [W_{l},I^{\rho}_{r}(g,\cdot)]h, f\rangle| &\lesssim& \cs{\int B \mathrm{D}^{2}(W_{l}) \mu^{\f12} \mu^{\f{1}{2}}_{*} g_{*}^{2} h^{2} \mathrm{d}V}{\int B
\mathrm{D}^{2}(\mu^{\f{1}{2}}_{*}) \mu^{\f12} \mu^{\f{1}{2}}_{*} f^{2}_{*} \mathrm{d}V}
\\ \label{ub-commutator-rest-part} &\lesssim&
(|\mu^{\frac{1}{64}}g|_{L^{2}} |h|_{L^{2}_{l+\gamma/2}} + |\mu^{\frac{1}{64}}g|_{H^{s_{1}}} |\mu^{\frac{1}{64}}h|_{H^{s_{2}}})|\mu^{\frac{1}{4}}f|_{L^{2}}.
\een
Recalling \eqref{commutator-I-into-main-and-rest}, patching together \eqref{ub-commutator-main-part} and \eqref{ub-commutator-rest-part},
we finish the proof.
\end{proof}

The next proposition gives an estimate for the commutator $[W_{l},\Gamma_{2,r,1}^{\rho}(g,\cdot)]$.
\begin{prop}\label{Commutotator-Wl-and-Gamma-2-r-1-epsilon-pm-tau}
	Let $l, s_{1}, s_{2} \geq 0$  and  $s_{1} + s_{2} = \f{1}{2}$. It holds that
		\beno
		|\langle [W_{l},\Gamma_{2,r,1}^{\rho}(g,\cdot)]h, f\rangle| \lesssim
(|g|_{L^{2}}|h|_{L^{2}_{l+\gamma/2}}+|\mu^{\frac{1}{64}}g|_{H^{s_{1}}}|\mu^{\frac{1}{64}}h|_{H^{s_{2}}})|f|_{\mathcal{L}^{s}_{\gamma/2}}.
		\eeno	
\end{prop}
\begin{proof} Recalling \eqref{inner-product-remaining-1}, we have $\langle [W_{l},\Gamma_{2,r,1}^{\rho}(g,\cdot)]h, f\rangle = \mathcal{F}_{1} + \mathcal{F}_{2}$ where
\beno
\mathcal{F}_{1} \colonequals      \int B  (Ng)_{*}(Nh)(N^{-1}f)^{\prime} \mathrm{D}(W_{l}^{\prime}) M^{\prime}
\mathrm{d}V,
\quad
\mathcal{F}_{2} \colonequals   \int B  (Ng)_{*}(Nh)(N^{-1}f)^{\prime} \mathrm{D}(W_{l}^{\prime}) M^{\prime}_{*}
\mathrm{d}V.
\eeno
Recalling $M = \mu^{\f{1}{2}} N$ and \eqref{commutator-inner-product-Q-epsilon-pm-0}, we observe
\beno
\mathcal{F}_{1}= \int B  (Nh)(Ng)_{*}(\mu^{\f{1}{2}}f)^{\prime} \mathrm{D}(W_{l}^{\prime})
\mathrm{d}V = \langle W_{l}Q_{c} (Ng, Nh) - Q_{c} (Ng, W_{l}Nh), \mu^{\f{1}{2}}f\rangle.
\eeno
Then by Proposition \ref{Commutotator-Wl-and-Q-epsilon-pm-0}, we get
\ben \label{J-1-ub}
|\mathcal{F}_{1}|  \lesssim (|\mu^{\frac{1}{64}}g|_{L^{2}} |h|_{L^{2}_{l+\gamma/2}} + |\mu^{\frac{1}{64}}g|_{H^{s_{1}}} |\mu^{\frac{1}{64}}h|_{H^{s_{2}}})|f|_{\mathcal{L}^{s}_{\gamma/2}}.
\een
Note that $N^{-1} = \mu^{-\f{1}{2}} - \rho \mu^{\f{1}{2}}$ and $\mu\mu_{*}=\mu^{\prime}\mu_{*}^{\prime}$, we get
\beno
N N_{*} (N^{-1})^{\prime} M^{\prime}_{*} &=& \frac{\mu^{\f{1}{2}}}{1 - \rho \mu} (\frac{\mu^{\f{1}{2}}}{1 - \rho \mu})_{*}  (\mu^{-\f{1}{2}})^{\prime} - \rho (\mu^{\f{1}{2}})^{\prime}) M^{\prime}_{*}
\\&=&  \frac{1}{1 - \rho \mu} (\frac{1}{1 - \rho \mu})_{*}  (\mu^{\f{1}{2}}M)^{\prime}_{*} - \rho M M_{*} N^{\prime}_{*}
\\&=&  \frac{1}{1 - \rho \mu} (\frac{1}{1 - \rho \mu})_{*}  \mathrm{D}((\mu^{\f{1}{2}}M)^{\prime}_{*})
+ \frac{1}{1 - \rho \mu} (NM)_{*}
- \rho M M_{*} \mathrm{D}(N^{\prime}_{*}) - \rho M (NM)_{*},
\eeno
which gives $\mathcal{F}_{2} = \sum_{i=1}^{4} \mathcal{F}_{2,i}$ where
\ben \label{J-2-t1}
\mathcal{F}_{2,1} &\colonequals&   \int B
(\frac{g}{1 - \rho \mu})_{*} \frac{h}{1 - \rho \mu}  f^{\prime} \mathrm{D}((\mu^{\f{1}{2}}M)^{\prime}_{*})
\mathrm{D}(W_{l}^{\prime})
\mathrm{d}V,
\\
\label{J-2-t2}
\mathcal{F}_{2,2} &\colonequals&      \int B  (NMg)_{*} \frac{h}{1 - \rho \mu} f^{\prime} \mathrm{D}(W_{l}^{\prime})
\mathrm{d}V,
\\
\label{J-2-t3}
\mathcal{F}_{2,3} &\colonequals&  - \rho \int B  (Mg)_{*} M h f^{\prime} \mathrm{D}(N^{\prime}_{*})  \mathrm{D}(W_{l}^{\prime})
\mathrm{d}V,
\\ \label{J-2-t4}
\mathcal{F}_{2,4} &\colonequals&  - \rho  \int B   (NMg)_{*} M h f^{\prime} \mathrm{D}(W_{l}^{\prime})
\mathrm{d}V.
\een
Note that the two lines \eqref{J-2-t1} and \eqref{J-2-t3} resemble the line \eqref{commutator-I-epsilon-pm-tau-m}.
Then by \eqref{ub-commutator-main-part} we get
\ben \label{J-2-1-and-J-2-3}
|\mathcal{F}_{2,1}| + |\mathcal{F}_{2,3}| \lesssim (|g|_{L^{2}} |h|_{L^{2}_{l+\gamma/2}}+|\mu^{\frac{1}{64}}g|_{H^{s_{1}}} |\mu^{\frac{1}{64}}h|_{H^{s_{2}}}) |f|_{\mathcal{L}^{s}_{\gamma/2}}.
\een
Note that the two lines \eqref{J-2-t2} and \eqref{J-2-t4} resemble \eqref{commutator-inner-product-Q-epsilon-pm-0}. Then by Proposition \ref{Commutotator-Wl-and-Q-epsilon-pm-0}, we get
\ben \label{J-2-2-and-J-2-4}
|\mathcal{F}_{2,2}| + |\mathcal{F}_{2,4}| \lesssim (|\mu^{\frac{1}{64}}g|_{L^{2}} |h|_{L^{2}_{l+\gamma/2}} + |\mu^{\frac{1}{64}}g|_{H^{s_{1}}} |\mu^{\frac{1}{64}}h|_{H^{s_{2}}}) |f|_{\mathcal{L}^{s}_{\gamma/2}}.
\een
Recalling $\langle [W_{l},\Gamma_{2,r,1}^{\rho}(g,\cdot)]h, f\rangle = \mathcal{F}_{1} + \sum_{i=1}^{4} \mathcal{F}_{2,i}$,
patching together \eqref{J-1-ub}, \eqref{J-2-1-and-J-2-3} and \eqref{J-2-2-and-J-2-4}, we finish the proof.
\end{proof}

The next proposition gives an estimate for the commutator $[W_{l},\Gamma_{2,r,2}^{\rho}(g,\cdot)]$.
\begin{prop}\label{Commutotator-Wl-and-Gamma-2-r-2-epsilon-pm-tau}
	Let $l, s_{1}, s_{2} \geq 0$  and  $s_{1} + s_{2} = \f{1}{2}$. It holds that
		\beno
		|\langle [W_{l},\Gamma_{2,r,2}^{\rho}(g,\cdot)]h, f\rangle| \lesssim (|\mu^{\frac{1}{64}}g|_{H^{s_{1}}} |\mu^{\frac{1}{64}}h|_{H^{s_{2}}}+|g|_{L^{2}} |h|_{L^{2}_{l+\gamma/2}})|f|_{\mathcal{L}^{s}_{\gamma/2}}.
		\eeno	
\end{prop}
\begin{proof} Recalling \eqref{inner-product-remaining-3}, we have $\langle [\langle W_{l},\Gamma_{2,r,2}^{\rho}(g,\cdot)]h, f\rangle = \mathcal{G}_{1} + \mathcal{G}_{2}$ where
\beno  \mathcal{G}_{1} \colonequals      \int B  (Nh)(Ng)^{\prime}_{*}(N^{-1}f)^{\prime} \mathrm{D}(W_{l}^{\prime}) M_{*}
\mathrm{d}V,
\quad
\mathcal{G}_{2} \colonequals     - \int B  (Nh)(Ng)^{\prime}_{*}(N^{-1}f)^{\prime} \mathrm{D}(W_{l}^{\prime}) M^{\prime}
\mathrm{d}V.
\eeno
We first visit  $\mathcal{G}_{1}$.
By \eqref{change-v-and-v-prime} and the identity $N_{*} N^{\prime}  N^{-1}  M_{*}^{\prime} = M_{*} (1 - \rho \mu) N_{*}^{\prime} (\frac{1}{1 - \rho \mu})^{\prime}$, we get
\ben \nonumber
\mathcal{G}_{1}&=&  \int B  (Ng)_{*} (Nh)^{\prime}  (N^{-1}f)  M_{*}^{\prime} \mathrm{D}(W_{l})
\mathrm{d}V
\\ \label{G1-into-G11-12-13} &=&  \int B  (Mg)_{*} (\frac{h}{1 - \rho \mu})^{\prime} (1 - \rho \mu)f   N_{*}^{\prime} \mathrm{D}(W_{l}) \mathrm{d}V = \mathcal{G}_{1,1} + \mathcal{G}_{1,2} + \mathcal{G}_{1,3},
\\ \label{defintion-G11} \mathcal{G}_{1,1} &\colonequals    &  \int B  (Mg)_{*} (\frac{h}{1 - \rho \mu})^{\prime}
\mathrm{D}((1 - \rho \mu)f)N_{*}^{\prime} \mathrm{D}(W_{l}) \mathrm{d}V,
\\ \nonumber 
 \mathcal{G}_{1,2}&\colonequals    &  \int B  (Mg)_{*} (h f)^{\prime} \mathrm{D}(N_{*}^{\prime}) \mathrm{D}(W_{l}) \mathrm{d}V,
\\ \nonumber 
\mathcal{G}_{1,3}&\colonequals    &  \int B  (NMg)_{*} (h f)^{\prime} \mathrm{D}(W_{l}) \mathrm{d}V.
\een
By Cauchy-Schwartz inequality and \eqref{M-rho-mu-N}, using   \eqref{Wl-difference-square-mu-star-g2-h2}, Remark \ref{still-true-h-prime} and Theorem \ref{two-functional-ub-by-norm}, we get
\ben \nonumber
|\mathcal{G}_{1,1}| &\lesssim& \cs{\int B  (\mu g^{2})_{*} (h^{2})^{\prime}  \mathrm{D}^{2}(W_{l}) \mathrm{d}V}
{\int B  \mu_{*} \mathrm{D}^{2}((1 - \rho \mu)f) \mathrm{d}V}
\\ \label{result-G11} &\lesssim&
(|\mu^{\frac{1}{64}}g|_{L^{2}} |h|_{L^{2}_{l+\gamma/2}} + |\mu^{\frac{1}{64}}g|_{H^{s_{1}}} |\mu^{\frac{1}{64}}h|_{H^{s_{2}}})|f|_{\mathcal{L}^{s}_{\gamma/2}}.
\een
By Cauchy-Schwartz inequality, \eqref{M-rho-mu-N}, \eqref{change-v-and-v-star} and \eqref{change-v-and-v-prime}, we have
\ben \nonumber
|\mathcal{G}_{1,2}| &\lesssim& \cs{\int B  (\mu g)^{2}_{*} (h^{2})^{\prime}  \mathrm{D}^{2}(W_{l}) \mathrm{d}V}
{\int
B   f^{2}_{*} \mathrm{D}^{2}(N)
\mathrm{d}V}
\\ \label{result-G12} &\lesssim&
 (|\mu^{\frac{1}{64}}g|_{L^{2}} |h|_{L^{2}_{l+\gamma/2}} + |\mu^{\frac{1}{64}}g|_{H^{s_{1}}} |\mu^{\frac{1}{64}}h|_{H^{s_{2}}})|f|_{\mathcal{L}^{s}_{\gamma/2}}.
\een
where in the last line we use \eqref{Wl-difference-square-mu-star-g2-h2}, Remark \ref{still-true-h-prime}, Theorem \ref{two-functional-ub-by-norm} and Remark \ref{mu-to-N-or-M-functional-N}.
By \eqref{Wl-difference-order-mu-star-g-hf-prime} and Remark \ref{still-true-h-prime}, we have
\ben \label{result-G13}
|\mathcal{G}_{1,3}| \lesssim |\mu^{\frac{1}{64}}g|_{L^{2}} |h|_{L^{2}_{l+\gamma/2}} |f|_{L^{2}_{\gamma/2}}.
\een

We now go to see  $\mathcal{G}_{2}$. Recalling $M = \mu^{\f{1}{2}} N$ and using \eqref{change-v-and-v-prime}, we have
\ben
\mathcal{G}_{2}= - \int B  (Nh)(Ng)^{\prime}_{*}(\mu^{\f{1}{2}}f)^{\prime} \mathrm{D}(W_{l}^{\prime})
\mathrm{d}V
\label{G2-into-G21-G22} =  \int B  (Ng)_{*}(Nh)^{\prime}(\mu^{\f{1}{2}}f) \mathrm{D}(W_{l}^{\prime})
\mathrm{d}V = \mathcal{G}_{2,1} + \mathcal{G}_{2,2},
\\ \nonumber 
\mathcal{G}_{2,1} \colonequals      \int B  (Ng)_{*}(Nh)^{\prime} \mathrm{D}(\mu^{\f{1}{2}}f)  \mathrm{D}(W_{l}^{\prime}) \mathrm{d}V,
\quad 
\mathcal{G}_{2,2} \colonequals      \int B  (Ng)_{*}(\mu^{\f{1}{2}}Nhf)^{\prime}  \mathrm{D}(W_{l}^{\prime}) \mathrm{d}V.
\een
Note that the term $\mathcal{G}_{2,1}$ resembles $\mathcal{G}_{1,1}$ in \eqref{defintion-G11}. Similar to \eqref{result-G11}, we have
\ben \label{result-G21}
|\mathcal{G}_{2,1}|  \lesssim
(|\mu^{\frac{1}{64}}g|_{L^{2}} |h|_{L^{2}_{l+\gamma/2}} + |\mu^{\frac{1}{64}}g|_{H^{s_{1}}} |\mu^{\frac{1}{64}}h|_{H^{s_{2}}})|f|_{\mathcal{L}^{s}_{\gamma/2}}.
\een
By \eqref{Wl-difference-order-mu-star-g-hf-prime} and Remark \ref{still-true-h-prime}, we have
\ben \label{result-G22}
|\mathcal{G}_{2,2}| \lesssim |\mu^{\frac{1}{64}}g|_{L^{2}} |h|_{L^{2}_{l+\gamma/2}} |f|_{L^{2}_{\gamma/2}}.
\een

Recalling \eqref{G1-into-G11-12-13} and \eqref{G2-into-G21-G22} to see $\langle [\langle W_{l},\Gamma_{2,r,2}^{\rho}(g,\cdot)]h, f\rangle = \mathcal{G}_{1} + \mathcal{G}_{2} = \mathcal{G}_{1,1} + \mathcal{G}_{1,2} + \mathcal{G}_{1,3} + \mathcal{G}_{2,1} + \mathcal{G}_{2,2}$,
patching together \eqref{result-G11}, \eqref{result-G12}, \eqref{result-G13}, \eqref{result-G21} and \eqref{result-G22}, we finish the proof.
\end{proof}

The next proposition gives an estimate for the commutator $[W_{l},\Gamma_{2,r,3}^{\rho}(g,\cdot)]$.
\begin{prop}\label{upper-bound-of-Gamma-2-r-2-g-g-star-weighted} Let $l \geq 0$. It holds that
\beno
|\langle [W_{l},\Gamma_{2,r,3}^{\rho}(g,\cdot)]h, f\rangle| \lesssim  \min\{|g|_{H^{2}}  |h|_{\mathcal{L}^{s}_{\gamma/2}}, |g|_{L^{2}} ( |\mu^{\frac{1}{256}}h|_{H^{2}} + |h|_{H^{\frac{3}{2}}_{\gamma/2+s}})\} |f|_{\mathcal{L}^{s}_{\gamma/2}}.
\eeno
\end{prop}
\begin{proof} Recalling \eqref{Gamma-2-r-3-into-1-2}, we have $\langle [W_{l},\Gamma_{2,r,3}^{\rho}(g,\cdot)]h, f\rangle = \langle [W_{l},\Gamma_{2,r,3,1}^{\rho}(g,\cdot)]h, f\rangle + \langle [W_{l},\Gamma_{2,r,3,2}^{\rho}(g,\cdot)]h, f\rangle$.
Recalling \eqref{definition-Gamma-remaining-3-1}, it is easy to see $[W_{l},\Gamma_{2,r,3,1}^{\rho}(g,\cdot)]=0$.  For the term involving  $\Gamma_{2,r,3,2}^{\rho}$, it suffices to consider  $ \langle  \Gamma_{2,r,3,2}^{\rho}(g, W_{l_{1}}h), W_{l_{2}}f \rangle  $ for $l_{1}, l_{2} \geq 0$. By
\eqref{Gamma-2-r-3-2-into-2-terms}, \eqref{defintion-X1} and \eqref{defintion-X2}, we have $  \langle  \Gamma_{2,r,3,2}^{\rho}(g, W_{l_{1}}h), W_{l_{2}}f \rangle   = \mathcal{H}_{1}+\mathcal{H}_{2}$ where
\beno
\mathcal{H}_{1} &\colonequals&  \int
B (\frac{g}{1 - \rho \mu})^{\prime}_{*} (M W_{l_{1}}h)_{*} (1 - \rho \mu) W_{l_{2}}f \mathrm{D}(N^{\prime})
\mathrm{d}V,
\\  
\mathcal{H}_{2} &\colonequals&  \int
B (Ng)^{\prime}_{*}(N W_{l_{1}} h)_{*} W_{l_{2}}f\mathrm{D}((\mu^{\f{1}{2}})^{\prime})
\mathrm{d}V.
\eeno
Since the two quantities have a similar structure, it suffices to consider $\mathcal{H}_{2}$. Note that $\mathcal{H}_{2}$ only differs from $\mathcal{X}_{2}$ in \eqref{defintion-X2} by the two weight functions $(W_{l_{1}})_{*}$ and $W_{l_{2}}$. We can follow the derivation of the estimate \eqref{G-2-r-3-2-X-result} for $\mathcal{X}_{2}$.  Observing that in the derivation there is a factor $\mu^{\frac{1}{16}}_{*}\mu^{\frac{1}{16}}$. By using the fact $\mu^{\frac{1}{16}}_{*}\mu^{\frac{1}{16}} (W_{l_{1}})_{*} W_{l_{2}} \lesssim C(l_{1},l_{2}) \mu^{\frac{1}{32}}_{*}\mu^{\frac{1}{32}}$,
 we can get the same upper bound as that in \eqref{G-2-r-3-2-X-result}.
\end{proof}

By the above commutator estimates and the upper bound estimates in Section \ref{bilinear},
  weighted upper bounds of $\Gamma_{2,m}^{\rho}(\cdot, \cdot)$ and  $\Gamma_{2}^{\rho}(\cdot, \cdot)$ are given in the following theorem.
\begin{thm} \label{upper-bound-bilinear}
 Let $(s_{3},s_{4})=(2,0)$ or $(s_{3},s_{4})=(0,2)$. The following two estimates are valid.
\ben \label{Gamma-2-m-ub-weighted}
 |\langle \Gamma_{2,m}^{\rho}(g,h), W_{2l}f \rangle|
&\lesssim&  |g|_{H^{s_{3}}}  |h|_{\mathcal{L}^{s_{4},s}_{l+\gamma/2}} |f|_{\mathcal{L}^{s}_{l+\gamma/2}}.
\\ \label{Gamma-2-ub-weighted}
 |\langle \Gamma_{2}^{\rho}(g,h), W_{2l}f \rangle|
&\lesssim& \rho^{\f{1}{2}} |g|_{H^{s_{3}}}  |h|_{\mathcal{L}^{s_{4},s}_{l+\gamma/2}} |f|_{\mathcal{L}^{s}_{l+\gamma/2}}.
\een
\end{thm}
\begin{proof} Recall \eqref{definition-of-L-n} for the definition of $|\cdot|_{\mathcal{L}^{n,s}_{l}}$.
Recalling \eqref{definition-of-norm-L-epsilon-gamma}, we have
\ben \label{relation-with-weighted-Sobolev}
|h|_{H^{n}_{l+s}} \lesssim |h|_{\mathcal{L}^{n,s}_{l}}, \quad |h|_{H^{n+s}_{l}} \lesssim |h|_{\mathcal{L}^{n,s}_{l}}.
\een
Recalling \eqref{definition-Gamma-main},
by Proposition \ref{Commutotator-Wl-and-Q-epsilon-pm-0}, Proposition \ref{Commutotator-Wl-and-I-epsilon-pm-tau} and Theorem \ref{upper-bound-of-Gamma-2-m}, using \eqref{relation-with-weighted-Sobolev},
we get  \eqref{Gamma-2-m-ub-weighted}.

Recalling \eqref{G2-and-into-five-terms},
by Proposition \ref{Commutotator-Wl-and-Q-epsilon-pm-0}, Proposition \ref{Commutotator-Wl-and-I-epsilon-pm-tau}, Proposition \ref{Commutotator-Wl-and-Gamma-2-r-1-epsilon-pm-tau}, Proposition \ref{Commutotator-Wl-and-Gamma-2-r-2-epsilon-pm-tau},  Proposition \ref{upper-bound-of-Gamma-2-r-2-g-g-star-weighted} and Theorem \ref{Gamma-2-ub}, using \eqref{relation-with-weighted-Sobolev}, we get \eqref{Gamma-2-ub-weighted}.
\end{proof}

Recall \eqref{not-mix-x-v-norm-energy} and \eqref{not-mix-x-v-norm-dissipation} for the definition of $\|\cdot\|_{H^{n}_{x}H^{m}}$ and $\|\cdot\|_{H^{n}_{x}\mathcal{L}^{m,s}_{l}}$.

In the 3-dimensional space $\mathbb{T}^{3}$,
by imbedding $L^{\infty}_{x} \hookrightarrow H^{2}_{x}$ or $L^{p}_{x} \hookrightarrow H^{s}_{x}$ with $\frac{s}{3} = \f{1}{2} - \frac{1}{p}$, based on Theorem \ref{upper-bound-bilinear},
estimates of inner product in the full space $(x,v)$ are given in the following theorem.
\begin{thm} \label{upper-bound-bilinear-full-space}
 Let $ a,b \in \mathbb{N}, a+b=2, (s_{3},s_{4})=(2,0)$ or $(s_{3},s_{4})=(0,2)$. The following two estimates are valid.
\ben \label{Gamma-2-m-ub-weighted-full-space}
 |( \Gamma_{2,m}^{\rho}(g,h), W_{2l}f )|
&\lesssim&  \|g\|_{H^{a}_{x}H^{s_{3}}}  \|h\|_{H^{b}_{x}\mathcal{L}^{s_{4},s}_{l+\gamma/2}} \|f\|_{L^{2}_{x}\mathcal{L}^{s}_{l+\gamma/2}}.
\\ \label{Gamma-2-ub-weighted-full-space}
|( \Gamma_{2}^{\rho}(g,h), W_{2l}f )|
&\lesssim& \rho^{\f{1}{2}} \|g\|_{H^{a}_{x}H^{s_{3}}}  \|h\|_{H^{b}_{x}\mathcal{L}^{s_{4},s}_{l+\gamma/2}} \|f\|_{L^{2}_{x}\mathcal{L}^{s}_{l+\gamma/2}}.
\een
\end{thm}

Based on Theorem \ref{upper-bound-bilinear-full-space}, by making a suitable choice of parameters $a,b,s_{3},s_{4}$ to deal with different distribution of derivative order, we get
\begin{thm} \label{Gamma-2-energy-estimate-label}  Let $N \geq 5$. The following two estimates are valid.
\ben \label{Gamma-2-m-energy-estimate}
|\sum_{|\alpha|+|\beta| \leq N} ( \partial^{\alpha}_{\beta} \Gamma_{2,m}^{\rho}(g,h), W_{2l_{|\alpha|,|\beta|}} \partial^{\alpha}_{\beta}f )| \lesssim \mathcal{E}^{\f{1}{2}}_{N}(g)\mathcal{D}^{\f{1}{2}}_{N}(h)\mathcal{D}^{\f{1}{2}}_{N}(f).
\\
\label{Gamma-2-energy-estimate}
|\sum_{|\alpha|+|\beta| \leq N} ( \partial^{\alpha}_{\beta} \Gamma_{2}^{\rho}(g,h), W_{2l_{|\alpha|,|\beta|}} \partial^{\alpha}_{\beta}f )| \lesssim \rho^{\f{1}{2}}\mathcal{E}^{\f{1}{2}}_{N}(g)\mathcal{D}^{\f{1}{2}}_{N}(h)\mathcal{D}^{\f{1}{2}}_{N}(f)
.
\een
\end{thm}
\begin{proof} By binomial formula and Remark \ref{sufficient-to-consider}, we need to consider $(\Gamma_{2,m}^{\rho}(\partial^{\alpha_{1}}_{\beta_{1}}g,\partial^{\alpha_{2}}_{\beta_{2}}h), W_{2l_{|\alpha|,|\beta|}} \partial^{\alpha}_{\beta}f)$ for all combinations of $\alpha_{1}+\alpha_{2}=\alpha, \beta_{1}+\beta_{2} \leq \beta$ with $|\alpha|+|\beta| \leq N$. By  \eqref{Gamma-2-m-ub-weighted-full-space}, it suffices to prove that
the following inequality
\beno
\|\partial^{\alpha_{1}}_{\beta_{1}}g\|_{H^{a}_{x}H^{s_{3}}}  \|\partial^{\alpha_{2}}_{\beta_{2}}h\|_{H^{b}_{x}\mathcal{L}^{s_{4},s}_{l_{|\alpha|,|\beta|}+\gamma/2}} \|\partial^{\alpha}_{\beta}f\|_{L^{2}_{x}\mathcal{L}^{s}_{l_{|\alpha|,|\beta|}+\gamma/2}} \lesssim \mathcal{E}^{\f{1}{2}}_{N}(g)\mathcal{D}^{\f{1}{2}}_{N}(h)\mathcal{D}^{\f{1}{2}}_{N}(f).
\eeno
holds for some $a,b,s_{3},s_{4}$ verifying $a,b \in \mathbb{N}, a+b=2, (s_{3},s_{4})=(2,0)$ or $(s_{3},s_{4})=(0,2)$.

The following is divided into two cases: $|\alpha|+|\beta| \leq N-4$ and $|\alpha|+|\beta| = N-k$ for $k \in \{0,1,2,3\}$.

{\it {Case 1: $|\alpha|+|\beta| \leq N-4$.}} In this case, there holds $|\alpha_{1}| + |\beta_{1}| + 4 \leq N$ and
we take $a=2,b=0,s_{3}=2,s_{4}=0$ and use $l_{|\alpha|,|\beta|} \leq l_{|\alpha_{2}|,|\beta_{2}|}$ to get
\ben \label{G2-rho-up-to-3-order-directly}
 && \|\partial^{\alpha_{1}}_{\beta_{1}}g\|_{H^{2}_{x}H^{2}}
\|\partial^{\alpha_{2}}_{\beta_{2}}h\|_{L^{2}_{x}\mathcal{L}^{s}_{l_{|\alpha|,|\beta|}+\gamma/2}}
\|\partial^{\alpha}_{\beta}f\|_{L^{2}_{x}\mathcal{L}^{s}_{l_{|\alpha|,|\beta|}+\gamma/2}}
\\ \nonumber
&\lesssim& \|g\|_{H^{N}_{x,v}} \|\partial^{\alpha_{2}}_{\beta_{2}}h\|_{L^{2}_{x}\mathcal{L}^{s}_{l_{|\alpha_{2}|,|\beta_{2}|}+\gamma/2}} \|\partial^{\alpha}_{\beta}f\|_{L^{2}_{x}\mathcal{L}^{s}_{l_{|\alpha|,|\beta|}+\gamma/2}} \lesssim \mathcal{E}^{\f{1}{2}}_{N}(g)\mathcal{D}^{\f{1}{2}}_{N}(h)\mathcal{D}^{\f{1}{2}}_{N}(f).
\een
Here we recall \eqref{mix-x-v-norm-energy} for the definition of $\|\cdot\|_{H^{n}_{x,v}}$.

{\it {Case 2: $|\alpha|+|\beta| = N-k$ for $k \in \{0,1,2,3\}$.}}
We consider two subcases: $|\alpha_{1}| + |\beta_{1}| \leq N-4$ and $|\alpha_{1}| + |\beta_{1}| = N-j$ for $k \leq j \leq 3$.  In the first subcase $|\alpha_{1}| + |\beta_{1}| \leq N-4$, there holds $|\alpha_{1}| + |\beta_{1}| + 4\leq N$. As the same as \eqref{G2-rho-up-to-3-order-directly}, we get
\beno
 \|\partial^{\alpha_{1}}_{\beta_{1}}g\|_{H^{2}_{x}H^{2}}
\|\partial^{\alpha_{2}}_{\beta_{2}}h\|_{L^{2}_{x}\mathcal{L}^{s}_{l_{|\alpha|,|\beta|}+\gamma/2}}
\|\partial^{\alpha}_{\beta}f\|_{L^{2}_{x}\mathcal{L}^{s}_{l_{|\alpha|,|\beta|}+\gamma/2}} \lesssim \mathcal{E}^{\f{1}{2}}_{N}(g)\mathcal{D}^{\f{1}{2}}_{N}(h)\mathcal{D}^{\f{1}{2}}_{N}(f).
\eeno

Recall that the order sequence  $\{l_{|\alpha|, |\beta|}\}_{|\alpha|+|\beta| \leq N}$
satisfies \eqref{weight-order-condition-1} and \eqref{weight-order-condition-2}. As a result,
\ben
\label{weight-order-condition-on-total-regularity}
|\alpha_{1}|+|\beta_{1}| \geq |\alpha_{2}|+|\beta_{2}|+1 \quad  \Rightarrow \quad  l_{|\alpha_{1}|,|\beta_{1}|} \leq l_{|\alpha_{2}|,|\beta_{2}|}.
\een
Now we go to consider
the other subcase $|\alpha_{1}| + |\beta_{1}| = N-j$ for $k \leq j \leq 3$. Note that
$|\alpha_{1}| + |\beta_{1}| + j = N, |\alpha_{2}| + |\beta_{2}| \leq (N-k)-(N-j)=j-k$.
 By taking $a+s_{3}=j, b+s_{4}=4-j$, there holds $|\alpha_{1}| + |\beta_{1}| + a +s_{3} =N, |\alpha_{2}| + |\beta_{2}| + b +s_{4} \leq 4-k$. Since $N \geq 5$, then $|\alpha|+|\beta| \geq 5-k$ and thus (by \eqref{weight-order-condition-on-total-regularity})
$
 l_{|\alpha|,|\beta|} \leq l_{|\alpha_{2}|+b,|\beta_{2}|+s_{4}}
$
which gives
\beno 
&& \|\partial^{\alpha_{1}}_{\beta_{1}}g\|_{H^{a}_{x}H^{s_{3}}}  \|\partial^{\alpha_{2}}_{\beta_{2}}h\|_{H^{b}_{x}\mathcal{L}^{s_{4},s}_{l_{|\alpha|,|\beta|}+\gamma/2}} \|\partial^{\alpha}_{\beta}f\|_{L^{2}_{x}\mathcal{L}^{s}_{l_{|\alpha|,|\beta|}+\gamma/2}}
\\
&\lesssim& \|g\|_{H^{N}_{x,v}} \|\partial^{\alpha_{2}}_{\beta_{2}}h\|_{H^{b}_{x}\mathcal{L}^{s_{4},s}_{l_{|\alpha_{2}|+b,|\beta_{2}|+s_{4}}+\gamma/2}} \|\partial^{\alpha}_{\beta}f\|_{L^{2}_{x}\mathcal{L}^{s}_{l_{|\alpha|,|\beta|}+\gamma/2}} \lesssim \mathcal{E}^{\f{1}{2}}_{N}(g)\mathcal{D}^{\f{1}{2}}_{N}(h)\mathcal{D}^{\f{1}{2}}_{N}(f).
\eeno

Patching together all the above estimates, we obtain \eqref{Gamma-2-m-energy-estimate}. By \eqref{Gamma-2-ub-weighted-full-space} and the same derivation, we also have \eqref{Gamma-2-energy-estimate}.
\end{proof}

\subsection{Estimate of the commutator $[W_{l},\Gamma_{3}^{\rho}(g, \cdot,\varrho)]$} In this subsection, we derive a commutator estimate
between $W_{l}$ and the operator $\Gamma_{3}^{\rho}(g,\cdot,\varrho)$. More precisely, we will give some upper bound for the quantity
$ 
\langle W_{l} \Gamma_{3}^{\rho}(g,h,\varrho) - \Gamma_{3}^{\rho}(g, W_{l}h,\varrho), f\rangle.
$
Recalling \eqref{Gamma3-into-Gamma31-and-Gamma32},
it suffices to estimate
$ \langle  [W_{l},\Gamma_{3,1}^{\rho}(g, \cdot,\varrho)]h, f \rangle  $ and
$ \langle  [W_{l},\Gamma_{3,2}^{\rho}(g, \cdot,\varrho)]h, f \rangle  $.
We will estimate $ \langle  [W_{l},\Gamma_{3,1}^{\rho}(g, \cdot,\varrho)]h, f \rangle  $ in Proposition \ref{commu-Gamma-3-1-ub} and
$ \langle  [W_{l},\Gamma_{3,2}^{\rho}(g, \cdot,\varrho)]h, f \rangle  $ in Proposition \ref{commu-Gamma-3-2-ub}.

\begin{prop} \label{commu-Gamma-3-1-ub}  The following two estimates are valid.
\ben \label{rho-highest-derivative-commu-Gamma31}
| \langle  [W_{l},\Gamma_{3,1}^{\rho}(g, \cdot,\varrho)]h, f \rangle  | &\lesssim& |g|_{H^{2}} |h|_{H^{\f{1}{2}}_{l+\gamma/2}} |\varrho|_{\mathcal{L}^{s}_{\gamma/2}}|f|_{\mathcal{L}^{s}_{\gamma/2}}.
\\ \label{g-or-h-highest-derivative-commu-Gamma31}
| \langle  [W_{l},\Gamma_{3,1}^{\rho}(g, \cdot,\varrho)]h, f \rangle  | &\lesssim& |g|_{H^{s_{1}}}|h|_{H^{s_{2}}_{l+\gamma/2}} |\mu^{\frac{1}{64}}\varrho|_{H^{2}}|f|_{\mathcal{L}^{s}_{\gamma/2}}.
\een
\end{prop}
\begin{proof} We first prove  \eqref{rho-highest-derivative-commu-Gamma31}.
Recalling \eqref{cubic1-inner-product}, using $
N_{*} N (N)^{\prime}_{*} (N^{-1})^{\prime} = (\frac{1}{1 - \rho \mu})_{*}\frac{1}{1 - \rho \mu}M^{\prime}_{*} (1 - \rho \mu)^{\prime}$ and
 \eqref{change-v-and-v-star}, we have
\ben  \nonumber
\langle [W_{l},\Gamma_{3,1}^{\rho}(g,\cdot,\varrho)]h, f\rangle &=&  \int
B (Ng)_{*} (Nh) (N \varrho)^{\prime}_{*} (N^{-1} f)^{\prime}  \mathrm{D}(W_{l}^{\prime}) \mathrm{d}V
\\ \label{commutator-Gamma3-1-type1} &=&  \int B (\frac{g}{1 - \rho \mu})_{*}\frac{h}{1 - \rho \mu}(M \varrho)^{\prime}_{*} ((1 - \rho \mu)f)^{\prime}  \mathrm{D}(W_{l}^{\prime}) \mathrm{d}V
\\ \nonumber 
&=&   \int B (\frac{g}{1 - \rho \mu})^{\prime}_{*}(\frac{h}{1 - \rho \mu})^{\prime}(M \varrho)_{*} (1 - \rho \mu)f  \mathrm{D}(W_{l}) \mathrm{d}V.
\een
By rearrangement, we have $\langle [W_{l},\Gamma_{3,1}^{\rho}(g,\cdot,\varrho)]h, f\rangle = \mathcal{J}_{1} + \mathcal{J}_{2}+ \mathcal{J}_{3}$ where
\ben
\mathcal{J}_{1} &\colonequals&   \int B \mathrm{D}((\frac{g}{1 - \rho \mu})^{\prime}_{*})(\frac{h}{1 - \rho \mu})^{\prime}(M \varrho)_{*} (1 - \rho \mu)f  \mathrm{D}(W_{l}) \mathrm{d}V,
 \\ \nonumber 
\mathcal{J}_{2} &\colonequals&    \int B  (\frac{M \varrho g}{1 - \rho \mu})_{*} (\frac{h}{1 - \rho \mu})^{\prime}
\mathrm{D}((1 - \rho \mu)f) \mathrm{D}(W_{l}) \mathrm{d}V,
\\ \nonumber 
\mathcal{J}_{3} &\colonequals&   \int B  (\frac{M \varrho g}{1 - \rho \mu})_{*} (h f)^{\prime}  \mathrm{D}(W_{l}) \mathrm{d}V.
\een
Recalling \eqref{M-rho-mu-N}, we have $M \lesssim \mu$.
By Cauchy-Schwartz inequality, using \eqref{difference-g-over}, \eqref{mu-star-g-star-difference-h2}, \eqref{estimate-of-Z-g2-h2-total-one-half}, \eqref{Wl-difference-square-mu-star-g2-h2} and Remark \ref{still-true-h-prime},
 we have
\ben
|\mathcal{J}_{1}| \lesssim \cs{\int B \mathrm{D}^{2}((\frac{g}{1 - \rho \mu})_{*}) \mu_{*}f^{2}
\mathrm{d}V}{\int B (h^{2})^{\prime}(\mu \varrho^{2})_{*}  \mathrm{D}^{2}(W_{l})
\mathrm{d}V}
\label{result-K1} \lesssim |g|_{H^{\frac{3}{2}}}|h|_{H^{\f{1}{2}}_{l+\gamma/2}} |\mu^{\frac{1}{4}}\varrho|_{L^{2}} |f|_{L^{2}_{\gamma/2+s}}.
\een
By Cauchy-Schwartz inequality, use the imbedding $H^{2} \hookrightarrow L^{\infty}$ for $g$,
the estimate \eqref{Wl-difference-square-mu-star-g2-h2}, Remark \ref{still-true-h-prime}, Theorem \ref{two-functional-ub-by-norm} and \eqref{taking-out-with-a-weight},
we have
\ben
|\mathcal{J}_{2}| \lesssim \cs{\int B \mu_{*} g^{2}_{*} (h^{2})^{\prime} \varrho^{2}_{*}  \mathrm{D}^{2}(W_{l})
\mathrm{d}V}{\int B \mu_{*} \mathrm{D}^{2}((1 - \rho \mu)f)
\mathrm{d}V}
\label{result-K2} \lesssim |g|_{H^{2}}|h|_{H^{\f{1}{2}}_{l+\gamma/2}} |\mu^{\frac{1}{4}}\varrho|_{L^{2}} |f|_{\mathcal{L}^{s}_{\gamma/2}}.
\een
By \eqref{Wl-difference-order-mu-star-g-hf-prime} and Remark \ref{still-true-h-prime}, using \eqref{deal-with-product}, we have
\ben \label{result-K3}
|\mathcal{J}_{3}|  \lesssim |\mu^{\frac{1}{64}}g\varrho|_{L^{2}} |h|_{L^{2}_{l+\gamma/2}} |f|_{L^{2}_{\gamma/2}} \lesssim |g|_{H^{2}} |h|_{L^{2}_{l+\gamma/2}} |\mu^{\frac{1}{64}}\varrho|_{L^{2}} |f|_{L^{2}_{\gamma/2}}.
\een
Recalling $\langle [W_{l},\Gamma_{3,1}^{\rho}(g,\cdot,\varrho)]h, f\rangle = \mathcal{J}_{1} + \mathcal{J}_{2}+ \mathcal{J}_{3}$,
patching together \eqref{result-K1}, \eqref{result-K2} and \eqref{result-K3},  we arrive at \eqref{rho-highest-derivative-commu-Gamma31}.

We now prove  \eqref{g-or-h-highest-derivative-commu-Gamma31}.
Recalling the line \eqref{commutator-Gamma3-1-type1} and $M = N \mu^{\f12}$,
we have $\langle [W_{l},\Gamma_{3,1}^{\rho}(g,\cdot,\varrho)]h, f\rangle  = \mathcal{K}_{1} + \mathcal{K}_{2} + \mathcal{K}_{3} + \mathcal{K}_{4}$
where
\beno
\mathcal{K}_{1} &\colonequals    &  \int B (\frac{g}{1 - \rho \mu})_{*}\frac{h}{1 - \rho \mu} N^{\prime}_{*} \mathrm{D}((\mu^{\f{1}{2}}\varrho)_{*}) ((1 - \rho \mu)f)^{\prime}  \mathrm{D}(W_{l}) \mathrm{d}V,
 \\
\mathcal{K}_{2} &\colonequals    &  \int B (\frac{g}{1 - \rho \mu})_{*}\frac{h}{1 - \rho \mu} (\mu^{\f{1}{2}}\varrho)_{*} \mathrm{D}(N_{*}) ((1 - \rho \mu)f)^{\prime}  \mathrm{D}(W_{l}) \mathrm{d}V,
 \\
\mathcal{K}_{3} &\colonequals&   \int B (\frac{g M \varrho}{1 - \rho \mu})_{*}\frac{h}{1 - \rho \mu}
\mathrm{D}((1 - \rho \mu)f) \mathrm{D}(W_{l}) \mathrm{d}V,
\\ 
\mathcal{K}_{4} &\colonequals    & - \int B (\frac{g M \varrho}{1 - \rho \mu})_{*} hf \mathrm{D}(W_{l}) \mathrm{d}V.
\eeno
Recalling \eqref{M-rho-mu-N}, we have $N \lesssim \mu^{\f12}$.
By Cauchy-Schwartz inequality and \eqref{change-v-and-v-prime}, using \eqref{Wl-difference-square-mu-star-prime-g2-h2} and \eqref{mu-star-g-star-difference-h2}, we have
\ben
\mathcal{K}_{1} \lesssim \cs{\int B g^{2}_{*} h^{2}  (\mu^{\f{1}{2}})^{\prime}_{*}\mathrm{D}^{2}(W_{l})
\mathrm{d}V}{\int B  \mathrm{D}^{2}((\mu^{\f{1}{2}}\varrho)_{*})  \mu^{\f{1}{2}}_{*} f^{2}
\mathrm{d}V}
\label{result-M1} \lesssim |g|_{H^{s_{1}}}|h|_{H^{s_{2}}_{l+\gamma/2}} |\mu^{\f{1}{2}}\varrho|_{H^{\frac{3}{2}}} |f|_{L^{2}_{\gamma/2+s}}.
\een
By Cauchy-Schwartz inequality, using \eqref{Wl-difference-square-mu-star-g2-h2}, \eqref{estimate-of-Z-g2-h2-total-one-half} and Remark \ref{mu-to-N-or-M-X-Y}.  we have
\ben
\mathcal{K}_{2} \lesssim \cs{\int B g^{2}_{*} h^{2}   \mu^{\f{1}{2}}_{*} \mathrm{D}^{2}(W_{l})
\mathrm{d}V}{\int B  \mathrm{D}^{2}(N_{*}) \mu^{\f{1}{2}}_{*} \varrho_{*}^{2} (f^{2})^{\prime}
\mathrm{d}V}
\label{result-M2} \lesssim |g|_{H^{s_{1}}}|h|_{H^{s_{2}}_{l+\gamma/2}} |\mu^{\frac{1}{8}}\varrho|_{H^{\f{1}{2}}} |f|_{L^{2}_{\gamma/2+s}}.
\een
By Cauchy-Schwartz inequality, we have
\ben
\mathcal{K}_{3} \lesssim  \cs{\int B \mu_{*} g^{2}_{*} \varrho_{*}^{2}  h^{2}   \mathrm{D}^{2}(W_{l})
\mathrm{d}V}{\int B  \mu_{*} \mathrm{D}^{2}((1 - \rho \mu)f)
\mathrm{d}V}
\label{result-M3} \lesssim |g|_{H^{s_{1}}}|h|_{H^{s_{2}}_{l+\gamma/2}} |\mu^{\frac{1}{8}}\varrho|_{H^{2}} |f|_{\mathcal{L}^{s}_{\gamma/2}}.
\een
Here we use the imbedding $H^{2} \hookrightarrow L^{\infty}$ for $\mu^{\frac{1}{8}}\varrho$ and
the estimate \eqref{Wl-difference-square-mu-star-g2-h2} to deal with the former bracket. The latter bracket is estimated by using Theorem \ref{two-functional-ub-by-norm} and \eqref{taking-out-with-a-weight}. By \eqref{Wl-difference-order-mu-star-gvarrho-hf-prime} and \eqref{deal-with-product}, we get
\ben \label{result-M4}
|\mathcal{K}_{4}| \lesssim   |\mu^{\frac{1}{64}}g|_{H^{s_{1}}}|h|_{H^{s_{2}}_{l+\gamma/2}} |\mu^{\frac{1}{64}}\varrho|_{H^{2}}|f|_{L^{2}_{\gamma/2}}.
\een
Recalling $\langle [W_{l},\Gamma_{3,1}^{\rho}(g,\cdot,\varrho)]h, f\rangle  = \mathcal{K}_{1} + \mathcal{K}_{2} + \mathcal{K}_{3} + \mathcal{K}_{4}$,
patching together \eqref{result-M1}, \eqref{result-M2}, \eqref{result-M3} and \eqref{result-M4},
we arrive at \eqref{g-or-h-highest-derivative-commu-Gamma31}.
\end{proof}

\begin{prop} \label{commu-Gamma-3-2-ub}  The following two estimates are valid.
\beno
| \langle  [W_{l},\Gamma_{3,2}^{\rho}(g, \cdot,\varrho)]h, f \rangle  | &\lesssim& |\mu^{\frac{1}{256}}g|_{H^{2}} |\mu^{\frac{1}{256}}h|_{H^{2}} |\varrho|_{\mathcal{L}^{s}_{\gamma/2}}|f|_{\mathcal{L}^{s}_{\gamma/2}}.
\\
| \langle  [W_{l},\Gamma_{3,2}^{\rho}(g, \cdot,\varrho)]h, f \rangle  | &\lesssim& |\mu^{\frac{1}{256}}g|_{H^{s_{1}}}|\mu^{\frac{1}{256}}h|_{H^{s_{2}}} |\mu^{\frac{1}{256}}\varrho|_{H^{3}}|f|_{\mathcal{L}^{s}_{\gamma/2}}.
\eeno
\end{prop}
\begin{proof} Recalling \eqref{cubic2-inner-product}, we have
\ben \label{Gamma3-2-f-difference-commutator}
\langle [W_{l},\Gamma_{3,2}^{\rho}(g,\cdot,\varrho)]h, f\rangle =  \int
B (Ng)_{*} (Nh) (\varrho f)^{\prime}  \mathrm{D}(W_{l}^{\prime}) \mathrm{d}V = \mathcal{O}_{1} + \mathcal{O}_{2}+ \mathcal{O}_{3},
\\ \nonumber 
\mathcal{O}_{1} \colonequals      \int
B (Ng)_{*} (Nh) \varrho^{\prime} \mathrm{D}(f^{\prime}) \mathrm{D}(W_{l}^{\prime})  \mathrm{d}V,
\quad 
\mathcal{O}_{2} \colonequals      \int
B (Ng)_{*} (Nh) \mathrm{D}(\varrho^{\prime}) f \mathrm{D}(W_{l}^{\prime})  \mathrm{d}V,
\\ \nonumber 
\mathcal{O}_{3} \colonequals     \int
B (Ng)_{*} (Nh) \varrho f \mathrm{D}(W_{l}^{\prime})  \mathrm{d}V.
\een

Recalling \eqref{M-rho-mu-N}, we have $N_{*} N \lesssim \mu^{\f{1}{2}} \mu^{\f{1}{2}}_{*}$. We first consider $\mathcal{O}_{1}$.
By Cauchy-Schwartz inequality, we have
\beno
|\mathcal{O}_{1}| \lesssim \cs{\int
B   \mu^{\f{1}{2}} \mu^{\f{1}{2}}_{*} g^{2}_{*} h^{2} (\varrho^{2})^{\prime} \mathrm{D}^{2}(W_{l})
\mathrm{d}V}{\int B  \mu^{\f{1}{2}} \mu^{\f{1}{2}}_{*} \mathrm{D}^{2}(f)
 \mathrm{d}V}.
\eeno
The latter bracket is bounded by $|f|_{\mathcal{L}^{s}_{\gamma/2}}^{2}$ according to Theorem \ref{two-functional-ub-by-norm}.
On one hand,
using the imbedding $H^{2} \hookrightarrow L^{\infty}$ for $\mu^{\f14}h$, the estimate \eqref{Wl-difference-square-mu-star-g2-h2} and Remark \ref{still-true-h-prime}, we estimate the former bracket as
\beno
\int
B   \mu^{\f{1}{2}} \mu^{\f{1}{2}}_{*} g^{2}_{*} h^{2} (\varrho^{2})^{\prime} \mathrm{D}^{2}(W_{l})
\mathrm{d}V \lesssim |\mu^{\frac{1}{64}}g|_{H^{\f12}}^{2}|\mu^{\frac{1}{4}}h|_{H^{2}}^{2} |\varrho|_{L^{2}_{\gamma/2+s}}^{2}.
\eeno
On the other hand, since $\mu \mu_{*} = \mu^{\prime} \mu_{*}^{\prime}$,
using the imbedding $H^{2} \hookrightarrow L^{\infty}$ for $\mu^{\f18}\varrho$ and the estimate \eqref{Wl-difference-square-mu-star-g2-h2}, we get
\beno
\int
B   \mu^{\f{1}{2}} \mu^{\f{1}{2}}_{*} g^{2}_{*} h^{2} (\varrho^{2})^{\prime} \mathrm{D}^{2}(W_{l})
\mathrm{d}V \lesssim |\mu^{\frac{1}{64}}g|_{H^{s_{1}}}^{2}|\mu^{\frac{1}{64}}h|_{H^{s_{2}}}^{2} |\mu^{\frac{1}{8}}\varrho|_{H^{2}}^{2}.
\eeno
Therefore we have two estimates
\ben \label{result-J1}  |\mathcal{O}_{1}| \lesssim |\mu^{\frac{1}{64}}g|_{H^{2}}|\mu^{\frac{1}{4}}h|_{H^{2}} |h|_{\mathcal{L}^{s}_{\gamma/2}}|f|_{\mathcal{L}^{s}_{\gamma/2}}, \quad |\mathcal{O}_{1}| \lesssim |\mu^{\frac{1}{64}}g|_{H^{s_{1}}}|\mu^{\frac{1}{64}}h|_{H^{s_{2}}} |\mu^{\frac{1}{8}}\varrho|_{H^{2}}|f|_{\mathcal{L}^{s}_{\gamma/2}}.
\een

We next consider $\mathcal{O}_{2}$. By Cauchy-Schwartz inequality, on one hand, we have
\beno
|\mathcal{O}_{2}| \lesssim \cs{\int
B  \mu^{\f{1}{2}} \mu^{\f{1}{2}}_{*} g^{2}_{*} h^{2} f^{2} \mathrm{D}^{2}(W_{l})
\mathrm{d}V}{\int B  \mu^{\f{1}{2}} \mu^{\f{1}{2}}_{*} \mathrm{D}^{2}(\varrho)
 \mathrm{d}V}.
\eeno
The latter bracket is bounded by $|\varrho|_{\mathcal{L}^{s}_{\gamma/2}}^{2}$ according to Theorem \ref{two-functional-ub-by-norm}. Using the imbedding $H^{2} \hookrightarrow L^{\infty}$ for $\mu^{\f14}h$ and the estimate \eqref{Wl-difference-square-mu-star-g2-h2}, the former bracket is bounded by $|\mu^{\frac{1}{64}}g|_{H^{\f{1}{2}}}^{2}|\mu^{\frac{1}{4}}h|_{H^{2}}^{2} |f|_{L^{2}_{\gamma/2+s}}^{2}.$ Patching together the two estimates, we get
\ben \label{result-J2-for-rho}  |\mathcal{O}_{2}| \lesssim |\mu^{\frac{1}{64}}g|_{H^{\f{1}{2}}}|\mu^{\frac{1}{4}}h|_{H^{2}} |\varrho|_{\mathcal{L}^{s}_{\gamma/2}}|f|_{\mathcal{L}^{s}_{\gamma/2}}.
\een
By Cauchy-Schwartz inequality, on the other hand,  we have
\beno
|\mathcal{O}_{2}| \lesssim \cs{\int
B  \mu^{\f{1}{2}} \mu^{\f{1}{2}}_{*} g^{2}_{*} h^{2} \mathrm{D}^{2}(\varrho)
\mathrm{d}V}{\int B  \mu^{\f{1}{2}} \mu^{\f{1}{2}}_{*}  f^{2} \mathrm{D}^{2}(W_{l})
 \mathrm{d}V}.
\eeno
The latter bracket is bounded by $|\mu^{\frac{1}{16}}f|_{L^{2}}^{2}$ according to the estimate \eqref{Wl-difference-square-mu-star-g2-h2} and \eqref{deal-with-polynomial-weight}. Using \eqref{change-v-and-v-star} and \eqref{mu-star-g-star-difference-h2-f2-star-another}, the former bracket is bounded by $|\mu^{\frac{1}{256}}g|_{H^{s_{1}}}^{2}|\mu^{\frac{1}{256}}h|_{H^{s_{2}}}^{2} |\mu^{\frac{1}{256}}\varrho|_{H^{3}}^{2}.$ Patching together the two estimates, we get
\ben \label{result-J2-for-g-or-h}  |\mathcal{O}_{2}| \lesssim |\mu^{\frac{1}{256}}g|_{H^{s_{1}}}|\mu^{\frac{1}{256}}h|_{H^{s_{2}}} |\mu^{\frac{1}{256}}\varrho|_{H^{3}}|\mu^{\frac{1}{16}}f|_{L^{2}}.
\een

By \eqref{Wl-difference-order-mu-star-g-hfvarrho-prime} and Remark \ref{still-true-h-prime}, for $\{ a_{1}, a_{2}, a_{3}\} = \{ 2, \f{1}{2}, 0\}$,  we get
\ben \label{result-J3}
|\mathcal{O}_{3}|   \lesssim |\mu^{\frac{1}{64}}g|_{H^{a_{1}}}|\mu^{\frac{1}{64}}h|_{H^{a_{2}}} |\mu^{\frac{1}{64}}\varrho|_{H^{a_{3}}}|\mu^{\frac{1}{64}}f|_{L^{2}}.
\een

Recalling \eqref{Gamma3-2-f-difference-commutator},
patching together \eqref{result-J1}, \eqref{result-J2-for-rho}, \eqref{result-J2-for-g-or-h} and \eqref{result-J3},  we finish the proof.
\end{proof}

Patching together Proposition \ref{commu-Gamma-3-1-ub} and Proposition \ref{commu-Gamma-3-2-ub}, recalling \eqref{Gamma3-into-Gamma31-and-Gamma32},
we get the following proposition.
\begin{prop} \label{commu-Gamma-3-ub}  The following functional estimates are valid.
\beno 
| \langle  [W_{l},\Gamma_{3}^{\rho}(g, \cdot,\varrho)]h, f \rangle  | &\lesssim& \rho |g|_{H^{2}}|h|_{H^{2}_{l+\gamma/2}}|\varrho|_{\mathcal{L}^{s}_{\gamma/2}}|f|_{\mathcal{L}^{s}_{\gamma/2}}.
\\ 
| \langle  [W_{l},\Gamma_{3}^{\rho}(g, \cdot,\varrho)]h, f \rangle  | &\lesssim& \rho |g|_{H^{s_{1}}}|h|_{H^{s_{2}}_{l+\gamma/2}} |\mu^{\frac{1}{256}}\varrho|_{H^{3}}|f|_{\mathcal{L}^{s}_{\gamma/2}}.
\eeno
\end{prop}

Patching together Proposition \ref{commu-Gamma-3-ub} and Theorem \ref{Gamma-3-ub}, using \eqref{deal-with-polynomial-weight},
we arrive at the following theorem for the weighted upper bound of the trilinear term.
\begin{thm} \label{Gamma-3-ub-weighted}  The following functional estimates are valid.
\beno 
|\langle \Gamma_{3}^{\rho}(g,h,\varrho), W_{2l}f\rangle| &\lesssim& \rho|g|_{H^{2}}|h|_{H^{2}}|\varrho|_{\mathcal{L}^{s}_{\gamma/2}}|f|_{\mathcal{L}^{s}_{l+\gamma/2}} + \rho|g|_{H^{3}}|h|_{\mathcal{L}^{2,s}_{l+\gamma/2}}|\mu^{\frac{1}{256}}\varrho|_{L^{2}}|f|_{\mathcal{L}^{s}_{l+\gamma/2}}.
\\ 
|\langle \Gamma_{3}^{\rho}(g,h,\varrho), W_{2l}f\rangle| &\lesssim& \rho|g|_{H^{2}}|h|_{\mathcal{L}^{s}_{l+\gamma/2}}|\mu^{\frac{1}{256}}\varrho|_{H^{3}}|f|_{\mathcal{L}^{s}_{l+\gamma/2}}.
\\ 
|\langle \Gamma_{3}^{\rho}(g,h,\varrho), W_{2l}f\rangle| &\lesssim& \rho|g|_{L^{2}}|h|_{\mathcal{L}^{2,s}_{l+\gamma/2}}|\mu^{\frac{1}{256}}\varrho|_{H^{3}}|f|_{\mathcal{L}^{s}_{l+\gamma/2}}.
\eeno
\end{thm}

By Sobolev imbedding $L^{\infty}_{x} \hookrightarrow H^{2}_{x}$ in the 3-dimensional space $\mathbb{T}^{3}$, based on Theorem \ref{Gamma-3-ub-weighted},
we have the following result in the full space $(x,v)$.
\begin{thm} \label{Gamma-3-ub-weighted-x-v}  The following functional estimates are valid.
\ben \label{rho-highest-derivative-weighted-x-v}
|(\Gamma_{3}^{\rho}(g,h,\varrho), W_{2l}f)| &\lesssim& \rho \|g\|_{H^{2}_{x}H^{3}}(\|h\|_{H^{2}_{x}H^{2}}\|\varrho\|_{L^{2}_{x}\mathcal{L}^{s}_{\gamma/2}} + \|h\|_{H^{2}_{x}\mathcal{L}^{2,s}_{l+\gamma/2}}
\|\varrho\|_{L^{2}_{x}L^{2}})
\|f\|_{L^{2}_{x}\mathcal{L}^{s}_{l+\gamma/2}}.
\\ \label{h-highest-derivative-weighted-x-v}
|(\Gamma_{3}^{\rho}(g,h,\varrho), W_{2l}f)| &\lesssim& \rho \|g\|_{H^{2}_{x}H^{2}}\|h\|_{L^{2}_{x}\mathcal{L}^{s}_{l+\gamma/2}}\|\varrho\|_{H^{2}_{x}H^{3}}
\|f\|_{L^{2}_{x}\mathcal{L}^{s}_{l+\gamma/2}}.
\\ \label{g-highest-derivative-weighted-x-v}
|(\Gamma_{3}^{\rho}(g,h,\varrho), W_{2l}f)| &\lesssim& \rho \|g\|_{L^{2}_{x}L^{2}}\|h\|_{H^{2}_{x}\mathcal{L}^{2,s}_{l+\gamma/2}}\|\varrho\|_{H^{2}_{x}H^{3}}
\|f\|_{L^{2}_{x}\mathcal{L}^{s}_{l+\gamma/2}}.
\een
\end{thm}

Theorem \ref{Gamma-3-ub-weighted-x-v} allows us to derive the following weighted energy estimate.
\begin{thm} \label{Gamma-3-energy-estimate} Let $N \geq 9$, then
\beno
|\sum_{|\alpha|+|\beta| \leq N} ( \partial^{\alpha}_{\beta} \Gamma_{3}^{\rho}(g,h,\varrho), W_{2l_{|\alpha|,|\beta|}} \partial^{\alpha}_{\beta}f)| \lesssim \rho \mathcal{E}^{\f{1}{2}}_{N}(g)\big(\mathcal{E}^{\f{1}{2}}_{N}(h)\mathcal{D}^{\f{1}{2}}_{N}(\varrho)
+\mathcal{D}^{\f{1}{2}}_{N}(h)\mathcal{E}^{\f{1}{2}}_{N}(\varrho)\big)
\mathcal{D}^{\f{1}{2}}_{N}(f).
\eeno
\end{thm}
\begin{proof} By binomial formula and Remark \ref{sufficient-to-consider-123},
it suffices to establish the following estimate for all combinations of $\alpha_{1}+\alpha_{2}+\alpha_{3}=\alpha, \beta_{1}+\beta_{2}+\beta_{3} \leq \beta$ with $|\alpha|+|\beta| \leq N$,
\ben \label{Newton-Lebniz-typical-term}
|( \Gamma_{3}^{\rho}(\partial^{\alpha_{1}}_{\beta_{1}}g, \partial^{\alpha_{2}}_{\beta_{2}}h,\partial^{\alpha_{3}}_{\beta_{3}}\varrho), W_{2l_{|\alpha|,|\beta|}} \partial^{\alpha}_{\beta}f )|
\lesssim \rho \mathcal{E}^{\f{1}{2}}_{N}(g)\big(\mathcal{E}^{\f{1}{2}}_{N}(h)\mathcal{D}^{\f{1}{2}}_{N}(\varrho)
+\mathcal{D}^{\f{1}{2}}_{N}(h)\mathcal{E}^{\f{1}{2}}_{N}(\varrho)\big)
\mathcal{D}^{\f{1}{2}}_{N}(f).
\een
The following is divided into three cases: $|\alpha|+|\beta| \leq N-5; |\alpha|+|\beta| = N-4; |\alpha|+|\beta| = N-j, j=0,1,2,3$.

{\it {Case 1: $|\alpha|+|\beta| \leq N-5$.}} In this case, $|\alpha_{1}|+|\beta_{1}|+4 \leq N-1, |\alpha_{3}| + |\beta_{3}|+5 \leq N.$
We use \eqref{h-highest-derivative-weighted-x-v} and the condition $l_{|\alpha|,|\beta|} \leq l_{|\alpha_{2}|,|\beta_{2}|}$ to get
\ben \label{rho-up-to-3-order-directly}
&& |( \Gamma_{3}^{\rho}(\partial^{\alpha_{1}}_{\beta_{1}}g,\partial^{\alpha_{2}}_{\beta_{2}}h,\partial^{\alpha_{3}}_{\beta_{3}}\varrho),  W_{2l_{|\alpha|,|\beta|}} \partial^{\alpha}_{\beta}f )|
\\ \nonumber &\lesssim& \rho \|\partial^{\alpha_{1}}_{\beta_{1}}g\|_{H^{2}_{x}H^{2}}
\|\partial^{\alpha_{2}}_{\beta_{2}}h\|_{L^{2}_{x}\mathcal{L}^{s}_{l_{|\alpha|,|\beta|}+\gamma/2}}
\|\partial^{\alpha_{3}}_{\beta_{3}}\varrho\|_{H^{2}_{x}H^{3}}
\|\partial^{\alpha}_{\beta}f\|_{L^{2}_{x}\mathcal{L}^{s}_{l_{|\alpha|,|\beta|}+\gamma/2}}
\\ \nonumber &\lesssim& \rho \|g\|_{H^{N}_{x,v}} \|\partial^{\alpha_{2}}_{\beta_{2}}h\|_{L^{2}_{x}\mathcal{L}^{s}_{l_{|\alpha_{2}|,|\beta_{2}|}+\gamma/2}}\|\varrho\|_{H^{N}_{x,v}} \|\partial^{\alpha}_{\beta}f\|_{L^{2}_{x}\mathcal{L}^{s}_{l_{|\alpha|,|\beta|}+\gamma/2}}
\\ \nonumber &\lesssim& \rho \mathcal{E}^{\f{1}{2}}_{N}(g)\mathcal{D}^{\f{1}{2}}_{N}(h)\mathcal{E}^{\f{1}{2}}_{N}(\varrho)
\mathcal{D}^{\f{1}{2}}_{N}(f).
\een

{\it {Case 2: $|\alpha|+|\beta| = N-4$.}} In this case, $|\alpha_{1}|+|\beta_{1}|+4 \leq N$. We consider two subcases: $|\alpha_{3}| + |\beta_{3}| \leq N-5$ and $|\alpha_{3}| + |\beta_{3}| = N-4$.  In the first subcase,  $|\alpha_{3}| + |\beta_{3}| + 5\leq N$. As the same as \eqref{rho-up-to-3-order-directly},
 we use \eqref{h-highest-derivative-weighted-x-v} and the condition $l_{|\alpha|,|\beta|} \leq l_{|\alpha_{2}|,|\beta_{2}|}$ to get \eqref{Newton-Lebniz-typical-term}.

In the second subcase, $|\alpha_{3}| + |\beta_{3}| = |\alpha| + |\beta|$ gives $|\alpha_{1}| = |\alpha_{2}| =|\beta_{1}| = |\beta_{2}|= 0$. Since $N \geq 9$, then $|\alpha|+|\beta| \geq 5$
and so $l_{|\alpha|,|\beta|} \leq l_{2,2}$ by \eqref{weight-order-condition-on-total-regularity}. Then we use \eqref{rho-highest-derivative-weighted-x-v} to get
\beno 
&& |( \Gamma_{3}^{\rho}(\partial^{\alpha_{1}}_{\beta_{1}}g, \partial^{\alpha_{2}}_{\beta_{2}}h,\partial^{\alpha_{3}}_{\beta_{3}}\varrho),  W_{2l_{|\alpha|,|\beta|}} \partial^{\alpha}_{\beta}f )|
\\ 
&\lesssim& \rho
\|g\|_{H^{2}_{x}H^{3}}(\|h\|_{H^{2}_{x}H^{2}}\|\partial^{\alpha}_{\beta}\varrho\|_{L^{2}_{x}\mathcal{L}^{s}_{\gamma/2}} + \|h\|_{H^{2}_{x}\mathcal{L}^{2,s}_{l_{|\alpha|,|\beta|}+\gamma/2}}
\|\partial^{\alpha}_{\beta}\varrho\|_{L^{2}_{x}L^{2}})
\|\partial^{\alpha}_{\beta}f\|_{L^{2}_{x}\mathcal{L}^{s}_{l_{|\alpha|,|\beta|}+\gamma/2}},
\\ 
&\lesssim& \rho \|g\|_{H^{N}_{x,v}} (\|h\|_{H^{N}_{x,v}} \|\partial^{\alpha}_{\beta}\varrho\|_{L^{2}_{x}\mathcal{L}^{s}_{\gamma/2}} + \|h\|_{H^{2}_{x}\mathcal{L}^{2,s}_{l_{2,2}+\gamma/2}} \|\varrho\|_{H^{N}_{x,v}}) \|\partial^{\alpha}_{\beta}f\|_{L^{2}_{x}\mathcal{L}^{s}_{l_{|\alpha|,|\beta|}+\gamma/2}}
\\ 
&\lesssim& \rho \mathcal{E}^{\f{1}{2}}_{N}(g)\big(\mathcal{E}^{\f{1}{2}}_{N}(h)\mathcal{D}^{\f{1}{2}}_{N}(\varrho)
+\mathcal{D}^{\f{1}{2}}_{N}(h)\mathcal{E}^{\f{1}{2}}_{N}(\varrho)\big)
\mathcal{D}^{\f{1}{2}}_{N}(f) .
\eeno

{\it {Case 3: $|\alpha|+|\beta| = N-j$ for $j=0,1,2,3$.}}  We consider three subcases: $|\alpha_{3}| + |\beta_{3}| \leq N-5, |\alpha_{1}| + |\beta_{1}| \leq N-4$; $|\alpha_{3}| + |\beta_{3}| \leq N-5, |\alpha_{1}| + |\beta_{1}| \geq N-3$; $|\alpha_{3}| + |\beta_{3}| \geq N-4$.
In the first subcase, $|\alpha_{1}|+ |\beta_{1}|+ 4 \leq N, |\alpha_{3}| + |\beta_{3}| + 5\leq N$. As the same as \eqref{rho-up-to-3-order-directly}, we use \eqref{h-highest-derivative-weighted-x-v} and the condition $l_{|\alpha|,|\beta|} \leq l_{|\alpha_{2}|,|\beta_{2}|}$ to get \eqref{Newton-Lebniz-typical-term}.

In the second subcase, $|\alpha_{1}| + |\beta_{1}| \geq N-3$ implies $|\alpha_{2}| + |\alpha_{3}| +|\beta_{2}| + |\beta_{3}| \leq 3-j$ and so $|\alpha_{2}| +|\beta_{2}| + 4 \leq 7-j, |\alpha_{3}| +|\beta_{3}| + 5\leq 8-j \leq N$.
Since $N \geq 9$, then $|\alpha|+|\beta| = N-j \geq 9-j$ and so $l_{|\alpha|,|\beta|} \leq l_{|\alpha_{2}|+2,|\beta_{2}|+2}$ by \eqref{weight-order-condition-on-total-regularity}.
Therefore we use \eqref{g-highest-derivative-weighted-x-v} to get
\beno 
&&|(  \Gamma_{3}^{\rho}(\partial^{\alpha_{1}}_{\beta_{1}}g, \partial^{\alpha_{2}}_{\beta_{2}}h,\partial^{\alpha_{3}}_{\beta_{3}}\varrho),  W_{2l_{|\alpha|,|\beta|}} \partial^{\alpha}_{\beta}f )|
\\ 
&\lesssim& \rho \|\partial^{\alpha_{1}}_{\beta_{1}}g\|_{L^{2}_{x}L^{2}}
\|\partial^{\alpha_{2}}_{\beta_{2}}h\|_{H^{2}_{x}\mathcal{L}^{2,s}_{l_{|\alpha|,|\beta|}+\gamma/2}}
\|\partial^{\alpha_{3}}_{\beta_{3}}\varrho\|_{H^{2}_{x}H^{3}}
\|\partial^{\alpha}_{\beta}f\|_{L^{2}_{x}\mathcal{L}^{s}_{l_{|\alpha|,|\beta|}+\gamma/2}}
\\ 
&\lesssim& \rho \|g\|_{H^{N}_{x,v}}
\|\partial^{\alpha_{2}}_{\beta_{2}}h\|_{H^{2}_{x}\mathcal{L}^{2,s}_{l_{|\alpha_{2}|+2,|\beta_{2}|+2}+\gamma/2}}
\|\varrho\|_{H^{N}_{x,v}}\|\partial^{\alpha}_{\beta}f\|_{L^{2}_{x}\mathcal{L}^{s}_{l_{|\alpha|,|\beta|}+\gamma/2}}
\\ 
&\lesssim& \rho \mathcal{E}^{\f{1}{2}}_{N}(g)\mathcal{D}^{\f{1}{2}}_{N}(h)\mathcal{E}^{\f{1}{2}}_{N}(\varrho)
\mathcal{D}^{\f{1}{2}}_{N}(f).
\eeno

In the third subcase, note that $|\alpha_{3}| + |\beta_{3}| \geq N-4$ gives $|\alpha_{1}|+|\alpha_{2}|+|\beta_{1}|+|\beta_{2}| \leq 4-j$ and so $|\alpha_{1}| +|\beta_{1}| + 5 \leq 9-j \leq N, |\alpha_{2}| +|\beta_{2}| + 4 \leq 8-j$.
Since $N \geq 9$, then $|\alpha|+|\beta|  \geq 9-j$ and so $l_{|\alpha|,|\beta|} \leq l_{|\alpha_{2}|+2,|\beta_{2}|+2}$ by \eqref{weight-order-condition-on-total-regularity}.
Therefore we can use \eqref{rho-highest-derivative-weighted-x-v} to get
\beno 
&& |( \Gamma_{3}^{\rho}(\partial^{\alpha_{1}}_{\beta_{1}}g, \partial^{\alpha_{2}}_{\beta_{2}}h,\partial^{\alpha_{3}}_{\beta_{3}}\varrho),  W_{2l_{|\alpha|,|\beta|}} \partial^{\alpha}_{\beta}f )|
\\ 
&\lesssim& \rho
\|\partial^{\alpha_{1}}_{\beta_{1}}g|_{H^{2}_{x}H^{3}}(\|\partial^{\alpha_{2}}_{\beta_{2}}h\|_{H^{2}_{x}H^{2}}
\|\partial^{\alpha_{3}}_{\beta_{3}}\varrho\|_{L^{2}_{x}\mathcal{L}^{s}_{\gamma/2}} + \|\partial^{\alpha_{2}}_{\beta_{2}}h\|_{H^{2}_{x}\mathcal{L}^{2,s}_{l_{|\alpha|,|\beta|}+\gamma/2}}
\|\partial^{\alpha_{3}}_{\beta_{3}}\varrho\|_{L^{2}_{x}L^{2}})
\|\partial^{\alpha}_{\beta}f\|_{L^{2}_{x}\mathcal{L}^{s}_{l_{|\alpha|,|\beta|}+\gamma/2}}
\\ 
&\lesssim& \rho \|g\|_{H^{N}_{x,v}} (\|h\|_{H^{N}_{x,v}}\|\partial^{\alpha_{3}}_{\beta_{3}}\varrho\|_{L^{2}_{x}\mathcal{L}^{s}_{\gamma/2}} + \|\partial^{\alpha_{2}}_{\beta_{2}}h\|_{H^{2}_{x}\mathcal{L}^{2,s}_{l_{|\alpha_{2}|+2,|\beta_{2}|+2}+\gamma/2}}
\|\varrho\|_{H^{N}_{x,v}}) \|\partial^{\alpha}_{\beta}f\|_{L^{2}_{x}\mathcal{L}^{s}_{l_{|\alpha|,|\beta|}+\gamma/2}}
\\ 
&\lesssim& \rho \mathcal{E}^{\f{1}{2}}_{N}(g)\big(\mathcal{E}^{\f{1}{2}}_{N}(h)\mathcal{D}^{\f{1}{2}}_{N}(\varrho)
+\mathcal{D}^{\f{1}{2}}_{N}(h)\mathcal{E}^{\f{1}{2}}_{N}(\varrho)\big)
\mathcal{D}^{\f{1}{2}}_{N}(f) .
\eeno

Patching together the three cases, we finish the proof.
\end{proof}

\section{Local well-posedness}   \label{local}
In this section, we will prove local existence of \eqref{quantum-Boltzmann-CP} or \eqref{linearized-quantum-Boltzmann-eq}. To this end, we need
local well-posedness of \eqref{quantum-Boltzmann-UU-linear}(or \eqref{linearized-quantum-Boltzmann-eq-linear}).
We first apply Proposition \ref{for-positivity} to show that the solution to \eqref{quantum-Boltzmann-UU-linear} is non-negative.
\begin{prop}\label{positivity} Let $T>0$. Let $F_{0} \geq 0$.
If $G \geq  0$ and $C_{T}\colonequals  \sup_{0 \leq t \leq T}\|\mu^{-\f14}G(t, \cdot, \cdot)\|_{H^{2}_{x}H^{4}}< \infty$. Let $F \in L^{\infty}([0,T]; L^{2}_{x}L^{2})$ be a solution to \eqref{quantum-Boltzmann-UU-linear}, then $F \geq 0$.
\end{prop}
\begin{proof} Let $F_{-} = \min \{ 0 ,  F\}, F_{+} = \max \{ 0 ,  F\}$. Taking inner product with $F_{-}$, we have
\beno
\f{1}{2}\frac{\mathrm{d}}{\mathrm{d}t}\|F_{-} \|^{2}_{L^{2}_{x}L^{2}} = (\tilde{Q}(G, F), F_{-}) &=&  \int B
 G_{*} F (1 + G_{*}^{\prime} +  G^{\prime})\mathrm{D}(F_{-}^{\prime}) \mathrm{d}V \mathrm{d}x
 \\&\leq &  \int B
 G_{*} F_{-} (1 + G_{*}^{\prime} +  G^{\prime})\mathrm{D}(F_{-}^{\prime}) \mathrm{d}V \mathrm{d}x
 \\  &=& (\tilde{Q}(G, F_{-}), F_{-}),
\eeno
where we use $G_{*} (1 + G_{*}^{\prime} +  G^{\prime}) \geq 0,  F F_{-} = F_{-}^{2}, F F_{-}^{\prime} =
F_{+}F_{-}^{\prime} + F_{-} F_{-}^{\prime} \leq  F_{-} F_{-}^{\prime}$ in the inequality. Since $(\tilde{Q}(G,F_{-}), F_{-}) = \int  \langle \tilde{Q}(G,F_{-}), F_{-}\rangle \mathrm{d}x$, by Proposition \ref{for-positivity} and the imbedding $H^{2}_{x} \hookrightarrow L^{\infty}_{x}$, we have
\beno
(\tilde{Q}(G, F_{-}), F_{-}) \lesssim \|\mu^{-\f14}G\|_{H^{2}_{x}H^{4}} (1+\|\mu^{-\f14}G\|_{H^{2}_{x}H^{4}}) \|F_{-}\|_{L^{2}_{x}L^{2}}^{2},
\eeno
which yields
\beno
\f{1}{2}\frac{\mathrm{d}}{\mathrm{d}t}\|F_{-} \|^{2}_{L^{2}_{x}L^{2}} \lesssim C_{T}(1+C_{T})\|F_{-} \|^{2}_{L^{2}_{x}L^{2}}.
\eeno
Therefore the initial condition $\|F_{-}(0) \|_{L^{2}_{x}L^{2}} = 0$ yields $\|F_{-}(t) \|_{L^{2}_{x}L^{2}} = 0$ for any $0 \leq t \leq T$.
\end{proof}

We next prepare two lemmas.
By Proposition \ref{ub-linearized-L2}, Proposition \ref{ub-for-correction-term} and Remark \ref{sufficient-to-consider-123},
using Cauchy-Schwartz inequality, \eqref{deal-with-polynomial-weight} and \eqref{taking-out-with-a-weight},
we have the following lemma for energy estimates involving $\mathcal{L}^{\rho}_{r}$ and $\mathcal{C}^{\rho}$.
\begin{lem} \label{correction-term-energy-estimate} Let $N \geq 0$ and $|\alpha|+|\beta| \leq N$, then
\beno
|(\partial^{\alpha}_{\beta}\mathcal{L}^{\rho}_{r}h, W_{2l_{|\alpha|,|\beta|}}\partial^{\alpha}_{\beta}f)| \lesssim \rho \mathcal{E}^{\f12}_{N}(h)\mathcal{E}^{\f12}_{N}(f), \quad
|(\partial^{\alpha}_{\beta} \mathcal{C}^{\rho} h, W_{2l_{|\alpha|,|\beta|}} \partial^{\alpha}_{\beta}f)| \lesssim \rho^{2} \mathcal{D}^{\f{1}{2}}_{N}(h)
\mathcal{D}^{\f{1}{2}}_{N}(f).
\eeno
\end{lem}
The following lemma gives
 two standard results  for energy estimate involving $[v \cdot \nabla_{x}, \partial^{\alpha}_{\beta}]$. One is controlled by dissipation, while the other is directly controlled by energy.
\begin{lem} \label{transport-commutator} Let $\beta = (\beta^{1},\beta^{2},\beta^{3}), |\beta| \geq 1$. For any $\eta>0$,
it holds that
\ben \label{bounded-by-previous-level-dissipation}
|([v \cdot \nabla_{x}, \partial^{\alpha}_{\beta}]f, W_{2l_{|\alpha|,|\beta|}} \partial^{\alpha}_{\beta}f)|
\leq \eta \|\partial^{\alpha}_{\beta}f\|_{L^{2}_{x}\mathcal{L}^{s}_{l_{|\alpha|,|\beta|}+\gamma/2}}^{2} + \frac{1}{\eta } \sum_{j=1}^{3} |\beta^{j}|^{2} \|W_{l_{|\alpha|+1,|\beta|-1}} \partial^{\alpha+e^{j}}_{\beta-e^{j}}f\|_{L^{2}_{x}\mathcal{L}^{s}_{\gamma/2}}^{2}.
\een
Here $e^{1}=(1,0,0), e^{2}=(0,1,0), e^{3}=(0,0,1)$. In addition,
\ben \label{bounded-by-energy}
|([v \cdot \nabla_{x}, \partial^{\alpha}_{\beta}]f, W_{2l_{|\alpha|,|\beta|}} \partial^{\alpha}_{\beta}f)|
\lesssim \mathcal{E}_{N}(f).
\een
\end{lem}
\begin{proof} Note that
$
[v \cdot \nabla_{x}, \partial^{\alpha}_{\beta}]f = -\sum_{j=1}^{3} \beta^{j}  \partial^{\alpha+e^{j}}_{\beta-e^{j}}f.
$
Then by Cauchy-Schwartz inequality, we have
\beno
|([v \cdot \nabla_{x}, \partial^{\alpha}_{\beta}]f, W_{2l_{|\alpha|,|\beta|}} \partial^{\alpha}_{\beta}f)| \leq
\sum_{j=1}^{3} \beta^{j} \|W_{l_{|\alpha|,|\beta|}-(\gamma/2+s)}\partial^{\alpha+e^{j}}_{\beta-e^{j}}f\|_{L^{2}_{x}L^{2}}
\|W_{l_{|\alpha|,|\beta|}+\gamma/2+s} \partial^{\alpha}_{\beta}f\|_{L^{2}_{x}L^{2}}.
\eeno
Using the condition $l_{|\alpha|,|\beta|}-(\gamma/2+s) \leq l_{|\alpha|+1,|\beta|-1}+(\gamma/2+s)$ by \eqref{weight-order-condition-1}, the fact $|\cdot|_{L^{2}_{\gamma/2+s}} \leq |\cdot|_{\mathcal{L}^{s}_{\gamma/2}}$, and the basic inequality
$ab \leq \eta a^{2}  + \frac{b^{2}}{4 \eta}$, we finish the proof of \eqref{bounded-by-previous-level-dissipation}.

By \eqref{weight-order-condition-1}, $l_{|\alpha|,|\beta|} \leq l_{|\alpha|+1,|\beta|-1}$. Then by Cauchy-Schwartz inequality, we can derive \eqref{bounded-by-energy}.
\end{proof}

Now we are ready to prove well-posedness of the equation \eqref{linearized-quantum-Boltzmann-eq-linear}.

\begin{prop} \label{a-priori-estimate-of-linear-equation} Let $N \geq 9$ and $T>0$. There are universal constants $0<\rho_{1}, \delta_{0}<1$ such that for $0< \rho \leq \rho_{1}, 0< \delta \leq \delta_{0}$ the following statement is valid.
Suppose the initial data $f_{0}$ and the given function $g$ verifies
\ben \nonumber 
\mathcal{E}_{N}(f_{0}) < \infty, \quad \mathcal{M} + \mathcal{N} f_{0} \geq 0,
\\ \label{assumption-on-g}
\sup_{0 \leq t \leq T} \mathcal{E}_{N}(g(t)) + \rho \int_{0}^{T} \mathcal{D}_{N}(g(t)) \mathrm{d}t \leq \delta \rho, \quad \inf_{0 \leq t \leq T} (\mathcal{M} + \mathcal{N} g(t)) \geq 0,
\een
then \eqref{linearized-quantum-Boltzmann-eq-linear} has a unique solution $f^{\rho} \in L^{\infty}([0,T]; \mathcal{E}_{N})$ verifying $\mathcal{M} + \mathcal{N} f^{\rho}(t) \geq 0$ for any $0 \leq t \leq T$ and
\ben \label{solution-also-verifies-estimate}
\sup_{0 \leq t \leq T} \mathcal{E}_{N}(f^{\rho}(t)) + \rho \int_{0}^{T} \mathcal{D}_{N}(f^{\rho}(t)) \mathrm{d}t \leq C \exp(CT+C\delta\rho) (\mathcal{E}_{N}(f_{0}) + \delta \rho^{2} T + \delta \rho^{2}),
\een
for some universal constant $C$(independent of $\delta, \rho, T$).
\end{prop}
\begin{proof} Based on the operator estimates in Section \ref{linear}, \ref{bilinear} and \ref{trilinear}, we can use Hahn-Banach Theorem like in  \cite{alexandre2011global,morimoto2016global} to prove existence.

Positivity is a direct consequence of Proposition \ref{positivity}. Indeed, note that $F = \mathcal{M} + \mathcal{N} f^{\rho}$ is the solution to \eqref{quantum-Boltzmann-UU-linear} with the given function $G=\mathcal{M} + \mathcal{N} g$ and initial data $F_{0}=\mathcal{M} + \mathcal{N} f_{0}$. Recall that $\mathcal{M} = \frac{\rho \mu}{1 - \rho \mu}, \mathcal{N} = \frac{\rho^{\f12} \mu^{\f12}}{1 - \rho \mu}$.
Recall \eqref{definition-energy-and-dissipation}  and \eqref{mix-x-v-norm-energy} for the definition of $\mathcal{E}_{N}(\cdot)$ and $\|\cdot \|_{H^{N}_{x,v}}$, we have
\beno
\sup_{0 \leq t \leq T} \|g(t)\|_{H^{N}_{x,v}}^{2} \leq \sup_{0 \leq t \leq T} \mathcal{E}_{N}(g(t)) \leq \delta \rho,
\eeno
which gives
\beno
\|\mu^{-\f14}G(t)\|_{H^{2}_{x}H^{4}} = \|\frac{\rho \mu^{\f34}}{1 - \rho \mu} + \frac{\rho^{\f12} \mu^{\f14}}{1 - \rho \mu} g(t)\|_{H^{2}_{x}H^{4}} \lesssim \rho + \rho^{\f12} \|g(t)\|_{H^{N}_{x,v}} \lesssim \rho.
\eeno
The second inequality in \eqref{assumption-on-g} gives $G \geq 0$. By Proposition \ref{positivity}, we get $F(t) = \mathcal{M} + \mathcal{N} f^{\rho}(t) \geq 0$ for any
$0 \leq t \leq T$.

As for uniqueness, it is standard to take difference and do energy estimate. Indeed, one can revise the following proof to get uniqueness naturally.

Now it remains to prove the {\it a priori} estimate \eqref{solution-also-verifies-estimate}.
For simplicity, let $f = f^{\rho}$ be the solution to \eqref{linearized-quantum-Boltzmann-eq-linear}.
Fix $0 \leq k \leq N$.
Take two indexes $\alpha$ and $\beta$ such that $|\alpha|\leq N-k$ and $|\beta|= k$. Set $q=l_{|\alpha|,|\beta|}$. Applying $\partial^{\alpha}_{\beta}$ to both sides of \eqref{linearized-quantum-Boltzmann-eq-linear}, taking inner product with $W_{2l_{|\alpha|,|\beta|}}\partial^{\alpha}_{\beta} f$ over $(x,v)$, using periodic condition,
we get
\ben \label{mix-x-v-weight-linear}
\f{1}{2}\frac{\mathrm{d}}{\mathrm{d}t}\|W_{l_{|\alpha|,|\beta|}}\partial^{\alpha}_{\beta}f \|^{2}_{L^{2}_{x}L^{2}}  + (\partial^{\alpha}_{\beta}\mathcal{L}^{\rho}f,W_{2l_{|\alpha|,|\beta|}}\partial^{\alpha}_{\beta}f) =
(\mathcal{R}, W_{2l_{|\alpha|,|\beta|}}\partial^{\alpha}_{\beta}f),  \een
where for simplicity in this proof
\beno
\mathcal{R}\colonequals    [v \cdot \nabla_{x}, \partial^{\alpha}_{\beta}]f + \partial^{\alpha}_{\beta}\mathcal{L}^{\rho}_{r}f + \partial^{\alpha}_{\beta}\mathcal{C}^{\rho}f +
\rho^{\f{1}{2}}\partial^{\alpha}_{\beta}\Gamma_{2,m}^{\rho}(g, f + \rho^{\f{1}{2}}\mu^{\f{1}{2}} ) + \partial^{\alpha}_{\beta}\Gamma_{3}^{\rho}(g, f + \rho^{\f{1}{2}}\mu^{\f{1}{2}}, g).
\eeno
By \eqref{bounded-by-energy}, we have
$ 
|([v \cdot \nabla_{x}, \partial^{\alpha}_{\beta}]f, W_{2l_{|\alpha|,|\beta|}}\partial^{\alpha}_{\beta}f)| \lesssim \mathcal{E}_{N}(f).
$
By Lemma \ref{correction-term-energy-estimate}, we have
\beno 
|(\partial^{\alpha}_{\beta}\mathcal{L}^{\rho}_{r}f, W_{2l_{|\alpha|,|\beta|}}\partial^{\alpha}_{\beta}f)| \lesssim \rho\mathcal{E}_{N}(f), \quad
|(\partial^{\alpha}_{\beta}\mathcal{C}^{\rho}f, W_{2l_{|\alpha|,|\beta|}}\partial^{\alpha}_{\beta}f)| \lesssim \rho^{2}\mathcal{D}_{N}(f).
\eeno
By  \eqref{Gamma-2-m-energy-estimate} and the basic inequality
$ab \leq \eta a^{2}  + \frac{b^{2}}{4 \eta}$,
we have
\beno 
|(\rho^{\f{1}{2}}\partial^{\alpha}_{\beta}\Gamma_{2,m}^{\rho}(g, f + \rho^{\f{1}{2}}\mu^{\f{1}{2}} ), W_{2l_{|\alpha|,|\beta|}}\partial^{\alpha}_{\beta}f)| &\lesssim& \rho^{\f{1}{2}}\mathcal{E}^{\f{1}{2}}_{N}(g)\big(\mathcal{D}^{\f{1}{2}}_{N}(f)
+\rho^{\f{1}{2}}\big)\mathcal{D}^{\f{1}{2}}_{N}(f)
\\ &\lesssim& \rho^{\f{1}{2}}\mathcal{E}^{\f{1}{2}}_{N}(g) \mathcal{D}_{N}(f) +  \eta \rho \mathcal{D}_{N}(f) +  \eta^{-1} \rho \mathcal{E}_{N}(g).
\eeno
By Theorem \ref{Gamma-3-energy-estimate} and the basic inequality
$ab \leq \f{1}{2} a^{2}  +  \f{1}{2} b^{2}$,
we have
\beno 
&& |(\partial^{\alpha}_{\beta}\Gamma_{3}^{\rho}(g, f + \rho^{\f{1}{2}}\mu^{\f{1}{2}},g), W_{2l_{|\alpha|,|\beta|}}\partial^{\alpha}_{\beta}f)|
\\ &\lesssim& \rho \mathcal{E}^{\f{1}{2}}_{N}(g)\big(\mathcal{E}^{\f{1}{2}}_{N}(f)\mathcal{D}_{N}^{\f{1}{2}}(g)
+\mathcal{D}^{\f{1}{2}}_{N}(f)\mathcal{E}^{\f{1}{2}}_{N}(g)\big)
\mathcal{D}^{\f{1}{2}}_{N}(f)
\\ &&+ \rho^{\frac{3}{2}} \mathcal{E}^{\f{1}{2}}_{N}(g)\big(\mathcal{D}_{N}^{\f{1}{2}}(g)
+\mathcal{E}^{\f{1}{2}}_{N}(g)\big)
\mathcal{D}^{\f{1}{2}}_{N}(f)
\\ &\lesssim& \rho \mathcal{E}_{N}(g) \mathcal{D}_{N}(f) + \rho \mathcal{D}_{N}(g) \mathcal{E}_{N}(f) + \rho^{2} \mathcal{E}_{N}(g) + \rho^{2} \mathcal{D}_{N}(g).
\eeno
Patching together the above estimates, back to \eqref{mix-x-v-weight-linear},
we get
\beno
&& \f{1}{2}\frac{\mathrm{d}}{\mathrm{d}t}\|\partial^{\alpha}_{\beta}f \|^{2}_{L^{2}_{x}L^{2}_{l_{|\alpha|,|\beta|}}}  + (\partial^{\alpha}_{\beta}\mathcal{L}^{\rho}f,W_{2l_{|\alpha|,|\beta|}}\partial^{\alpha}_{\beta}f)
\\ &\lesssim& \big( \eta \rho
 + \rho^{2} + \rho^{\f{1}{2}}\mathcal{E}^{\f{1}{2}}_{N}(g) + \rho \mathcal{E}_{N}(g) \big)\mathcal{D}_{N}(f) +
\big(1+ \rho \mathcal{D}_{N}(g)\big)\mathcal{E}_{N}(f)  + \big( \eta^{-1} \rho + \rho^{2} \big) \mathcal{E}_{N}(g) + \rho^{2} \mathcal{D}_{N}(g).
\eeno
We now apply Theorem \ref{L-energy-estimate} to deal with the term involving $\mathcal{L}^{\rho}$.
From now on  in this proof, we assume $0 < \rho \leq \rho_{0}$.
When $|\beta|=0$, by \eqref{no-v-derivative}, since $\gamma \leq 0$, we have
\ben \label{no-v-derivative-linear}
(\partial^{\alpha} \mathcal{L}^{\rho}f, W_{2l_{|\alpha|,0}} \partial^{\alpha}f) \geq
\frac{\lambda_{0}}{4}\rho\|\partial^{\alpha}f\|_{L^{2}_{x}\mathcal{L}^{s}_{l_{|\alpha|,0}+\gamma/2}}^{2} - C \rho
\mathcal{E}_{N}(f).
\een
If $|\beta| \geq 1$, by \eqref{with-v-derivative}, since $\gamma \leq 0$, we have
\ben \label{with-v-derivative-linear}
(\partial^{\alpha}_{\beta} \mathcal{L}^{\rho}f, W_{2l_{|\alpha|,|\beta|}} \partial^{\alpha}_{\beta}f) \geq
\frac{\lambda_{0}}{8}\rho\|\partial^{\alpha}_{\beta}f\|_{L^{2}_{x}\mathcal{L}^{s}_{l_{|\alpha|,|\beta|}+\gamma/2}}^{2} -
C \rho \mathcal{E}_{N}(f) - C \rho
\sum_{\beta_{1}<\beta} \|\partial^{\alpha}_{\beta_{1}}f\|_{L^{2}_{x}\mathcal{L}^{s}_{l_{|\alpha|,|\beta|}+\gamma/2}}^{2}.
\een
Note that the last term in \eqref{with-v-derivative-linear} only involves the following term with index $\beta_{1}<\beta$. Since $l_{|\alpha|,|\beta|} \leq l_{|\alpha|,|\beta_{1}|}$, there holds
\beno
\|\partial^{\alpha}_{\beta_{1}}f\|_{L^{2}_{x}\mathcal{L}^{s}_{l_{|\alpha|,|\beta|}+\gamma/2}}^{2} \leq \|\partial^{\alpha}_{\beta_{1}}f\|_{L^{2}_{x}\mathcal{L}^{s}_{l_{|\alpha|,|\beta_{1}|}+\gamma/2}}^{2}.
\eeno
Therefore, there are universal constants $\{C_{k}\}_{0 \leq k \leq N}$ with $1 \leq C_{k} \leq C_{k-1}$ for $k=1, \cdots, N$ and a large universal constant $C$
such that
\beno
&& \frac{\mathrm{d}}{\mathrm{d}t} \sum_{k=0}^{N} C_{k}\sum_{|\alpha|\leq N-k} \sum_{|\beta|=k}\|\partial^{\alpha}_{\beta}f \|^{2}_{L^{2}_{x}L^{2}_{l_{|\alpha|,|\beta|}}}  + \frac{\lambda_{0}}{16}\rho \mathcal{D}_{N}(f)
\\ &\leq& C\big( \eta
 + \rho + \rho^{-\f{1}{2}}\mathcal{E}^{\f{1}{2}}_{N}(g) +  \mathcal{E}_{N}(g) \big) \rho\mathcal{D}_{N}(f) +
C \big(1+ \rho \mathcal{D}_{N}(g)\big)\mathcal{E}_{N}(f)  + C \big( \eta^{-1} \rho + \rho^{2} \big) \mathcal{E}_{N}(g) + C \rho^{2} \mathcal{D}_{N}(g).
\eeno
Let us take $\eta$ such that $C \eta = \frac{\lambda_{0}}{128}$. Suppose $\rho, \delta$ verify
\ben \label{condition-on-delta-rho-1}
 C \rho \leq \frac{\lambda_{0}}{128}, \quad C^{2} \delta  \leq (\frac{\lambda_{0}}{128})^{2}, \quad C \delta  \leq \frac{\lambda_{0}}{128},
\een
then by \eqref{assumption-on-g} there holds
$
C \rho^{-\f{1}{2}}\mathcal{E}^{\f{1}{2}}_{N}(g) \leq \frac{\lambda_{0}}{128},  C \mathcal{E}_{N}(g)  \leq \frac{\lambda_{0}}{128}.
$
By these smallness conditions on $\rho, \delta$ and the choice of $\eta$, we have
\beno
&&\frac{\mathrm{d}}{\mathrm{d}t} \sum_{k=0}^{N} C_{k}\sum_{|\alpha|\leq N-k} \sum_{|\beta|=k}\|\partial^{\alpha}_{\beta}f \|^{2}_{L^{2}_{x}L^{2}_{l_{|\alpha|,|\beta|}}}  + \frac{\lambda_{0}}{32}\rho \mathcal{D}_{N}(f)
\\&\leq& C\big(1+\rho \mathcal{D}_{N}(g)\big) \mathcal{E}_{N}(f) +  C\rho \mathcal{E}_{N}(g)   + C\rho^{2} \mathcal{D}_{N}(g).
\eeno
Note that
\ben \label{combined-energy}
\mathcal{E}_{N}(f) \leq \tilde{\mathcal{E}}_{N}(f)\colonequals    \sum_{k=0}^{N} C_{k}\sum_{|\alpha|\leq N-k} \sum_{|\beta|=k}\|\partial^{\alpha}_{\beta}f \|^{2}_{L^{2}_{x}L^{2}_{l_{|\alpha|,|\beta|}}} \leq \tilde{C} \mathcal{E}_{N}(f),
\een
where $\tilde{C}\colonequals  \max_{0 \leq k \leq N} \{C_{k}\}$. We arrive at
\beno
\frac{\mathrm{d}}{\mathrm{d}t} \tilde{\mathcal{E}}_{N}(f)  + \frac{\lambda_{0}}{32}\rho \mathcal{D}_{N}(f)
\leq C\big(1+\rho \mathcal{D}_{N}(g)\big)\tilde{\mathcal{E}}_{N}(f) +  C\rho \mathcal{E}_{N}(g)   + C\rho^{2} \mathcal{D}_{N}(g).
\eeno
By Gr\"{o}nwall's lemma and the assumption \eqref{assumption-on-g},  for any $0 \leq t \leq T$, we have
\beno
&& \tilde{\mathcal{E}}_{N}(f(t))  + \frac{\lambda_{0}}{32} \rho \int_{0}^{t} \mathcal{D}_{N}(f(\tau)) \mathrm{d}\tau
\\ &\leq&
\exp(C\int_{0}^{t} (1 + \rho \mathcal{D}_{N}(g)) \mathrm{d}\tau) (\tilde{\mathcal{E}}_{N}(f_{0}) + C\int_{0}^{t} \rho \mathcal{E}_{N}(g) \mathrm{d}\tau + C\int_{0}^{t} \rho^{2} \mathcal{D}_{N}(g) \mathrm{d}\tau)
\\ &\leq&
\exp(Ct + C\delta\rho) (\tilde{\mathcal{E}}_{N}(f_{0}) + C \delta \rho^{2} t  + C \delta \rho^{2}).
\eeno
Recalling \eqref{combined-energy}, since $\tilde{C}$ and $\lambda_{0}$ are universal constants,
 we get the desired result \eqref{solution-also-verifies-estimate} with a different constant $C$.

Let $C$ be the largest constant appearing in the above proof.
According to \eqref{condition-on-delta-rho-1}, we set
\beno \rho_{1} \colonequals   \min \{ \rho_{0}, \frac{\lambda_{0}}{128C} \}, \quad \delta_{0} \colonequals   (\frac{\lambda_{0}}{128C})^{2} \leq \frac{\lambda_{0}}{128C}. \eeno
Then if $0< \rho \leq \rho_{1}$ and $0< \delta \leq \delta_{0}$, the above proof is valid.
\end{proof}


By iteration on the equation \eqref{linearized-quantum-Boltzmann-eq-linear}, we derive local well-posedness of the Cauchy problem \eqref{linearized-quantum-Boltzmann-eq}.
The local well-posedness result can be concluded as follows:
\begin{thm}\label{local-well-posedness-LBE} Let $N \geq 9$. There are universal constants $\rho_{2}, \delta_{1}, T^{*}>0$ such that for $0< \rho \leq \rho_{2}, 0< \delta \leq \delta_{1}$, if
\beno \mathcal{E}_{N}(f_{0}) \leq \delta \rho, \quad \mathcal{M} + \mathcal{N} f_{0} \geq 0, \eeno
then the Cauchy problem \eqref{linearized-quantum-Boltzmann-eq} admits a unique solution $f^{\rho} \in
L^{\infty}([0,T^{*}]; \mathcal{E}_{N})$ verifying
\ben \label{uniform-estimate-of-solution}
\sup_{0 \leq t \leq T^{*}} \mathcal{E}_{N}(f^{\rho}(t)) + \rho \int_{0}^{T^{*}} \mathcal{D}_{N}(f^{\rho}(t)) \mathrm{d}t \leq C \delta \rho, \quad \mathcal{M} + \mathcal{N} f^{\rho}(t) \geq 0,
\een	
for some universal constant $C$.
\end{thm}
\begin{proof}
Let $f^{\rho, 0} \equiv 0$ and for $n \geq 1$, $f^{\rho, n}$ is the solution to
the following problem
\beno 
\partial _t f^{n} +  v \cdot \nabla_{x} f^{n} + \mathcal{L}^{\rho}f^{n} &=&  \mathcal{L}^{\rho}_{r}f^{n} + \mathcal{C}^{\rho}f^{n}+
\rho^{\f{1}{2}}\Gamma_{2,m}^{\rho}(f^{n-1},f^{n}+\rho^{\f{1}{2}}\mu^{\f{1}{2}})
\\&&+ \Gamma_{3}^{\rho}(f^{n-1},f^{n}+\rho^{\f{1}{2}}\mu^{\f{1}{2}},f^{n-1}),
\eeno
with the initial condition $f^{n}|_{t=0} = f_{0}.$

Let $\tilde{C}\colonequals     C \exp(C+C\delta_{0}\rho_{1}) (1 +  2 \rho_{1}) \geq 1, \tilde{\delta}_{1} \colonequals    \delta_{0}/\tilde{C}, \tilde{\rho}_{1}\colonequals    \rho_{1}/\tilde{C}$. Here and in the rest of this paragraph $C, \rho_{1}, \delta_{0}$ are the constants appearing in Proposition \ref{a-priori-estimate-of-linear-equation}. We now set to prove the following statement.
If $0< \rho \leq \tilde{\rho}_{1}, 0< \delta \leq \tilde{\delta}_{1}$ and $\mathcal{E}_{N}(f_{0}) \leq \delta \rho$, then the sequence $\{f^{\rho, n}\}_{n \geq 0}$ is well-defined on the time interval $[0, 1]$ and has uniform(in $n$) estimate
\ben \label{uniform-estimate-of-fn}
\sup_{0 \leq t \leq 1} \mathcal{E}_{N}(f^{\rho, n}(t)) + \rho \int_{0}^{1} \mathcal{D}_{N}(f^{\rho, n}(t)) \mathrm{d}t \leq \tilde{C} \delta \rho.
\een
Obviously  \eqref{uniform-estimate-of-fn}  is valid when $n=0$  since
$f^{\rho, 0} \equiv 0$. We now use mathematical induction over $n$. Let us assume that
the partial sequence $\{f^{\rho, n}\}_{0 \leq n \leq k}$ is well-defined on the time interval $[0, 1]$ and satisfies the uniform estimate \eqref{uniform-estimate-of-fn}. In particular, for $n=k$, it holds that
\beno 
\sup_{0 \leq t \leq 1} \mathcal{E}_{N}(f^{\rho, k}(t)) + \rho \int_{0}^{1} \mathcal{D}_{N}(f^{\rho, k}(t)) \mathrm{d}t \leq \tilde{C} \delta \rho.
\eeno
Since $0<\delta \leq \tilde{\delta}_{1}, 0< \rho \leq \tilde{\rho}_{1}$, by the definition of  $\tilde{C}, \tilde{\delta}_{1}, \tilde{\rho}_{1}$,  we have  $0< \tilde{C} \delta \leq \delta_{0}, 0 < \rho \leq \rho_{1}$. Now applying Proposition \ref{a-priori-estimate-of-linear-equation}(in which $\delta$ is replaced with $\tilde{C} \delta$), there is a unique solution $f^{\rho, k+1} \in L^{\infty}([0,1]; \mathcal{E}_{N})$ verifying
\beno 
&& \sup_{0 \leq t \leq 1} \mathcal{E}_{N}(f^{\rho, k+1}(t)) + \rho \int_{0}^{1} \mathcal{D}_{N}(f^{\rho, k+1}(t)) \mathrm{d}t
\\ 
&\leq& C \exp(C+C\tilde{C}\delta\rho) (\mathcal{E}_{N}(f_{0}) + 2\tilde{C} \delta \rho^{2})
\leq C \exp(C+C\delta_{0}\rho_{1}) (1 +  2 \rho_{1})\delta\rho = \tilde{C}\delta\rho,
\eeno
where we use $\tilde{C}\rho \leq \rho_{1}$ and $\delta \leq \tilde{\delta}_{1} \leq \delta_{0}$ in the last inequality. Now the statement about the sequence $\{f^{\rho, n}\}_{n \geq 0}$ is proved.

Moreover by Proposition \ref{a-priori-estimate-of-linear-equation}, we have for any $n \geq 0$ and $0 \leq t \leq 1$,
\ben \label{solution-positivity}
\mathcal{M} + \mathcal{N} f^{\rho, n}(t) \geq 0.
\een

We now prove that $\{f^{\rho, n}\}_{n \geq 0}$ is a Cauchy sequence in $L^{\infty}([0,T_{0}]; \mathcal{E}_{N})$ for some $0<T_{0} \leq 1$. Let $w^{n}\colonequals    f^{\rho, n+1}-f^{\rho, n}$ for $n \geq 0$. Then for $n \geq 1$, the function $w^{n}$ solves
\beno
\partial _t w^{n} +  v \cdot \nabla_{x} w^{n} + \mathcal{L}^{\rho} w^{n} =  \mathcal{L}^{\rho}_{r} w^{n} + \mathcal{C}^{\rho} w^{n} + \mathcal{Y}, \quad w^{n}(0,x,v) \equiv 0,
\eeno
where
\ben \label{definition-of-X-8-terms}
\mathcal{Y} &\colonequals    &\rho^{\f{1}{2}}\Gamma_{2,m}^{\rho}(f^{\rho, n},w^{n}) +
\rho^{\f{1}{2}}\Gamma_{2,m}^{\rho}(w^{n-1},f^{\rho, n}) + \rho \Gamma_{2,m}^{\rho}(w^{n-1}, \mu^{\f12})
\\ \nonumber &&+ \Gamma_{3}^{\rho}(f^{\rho, n}, w^{n}, f^{\rho, n}) + \Gamma_{3}^{\rho}(w^{n-1}, f^{\rho, n}, f^{\rho, n})
+ \Gamma_{3}^{\rho}(f^{\rho, n-1}, f^{\rho, n}, w^{n-1})
\\ \nonumber &&+ \rho^{\f{1}{2}} \Gamma_{3}^{\rho}(w^{n-1}, \mu^{\f12}, f^{\rho, n})
+ \rho^{\f{1}{2}} \Gamma_{3}^{\rho}(f^{\rho, n-1}, \mu^{\f12}, w^{n-1}).
\een
By basic energy estimate(similar to \eqref{mix-x-v-weight-linear}), we have
\ben \label{mix-x-v-weight-solution-diff}
\f{1}{2}\frac{\mathrm{d}}{\mathrm{d}t}\|\partial^{\alpha}_{\beta}w^{n} \|^{2}_{L^{2}_{x}L^{2}_{l_{|\alpha|,|\beta|}}}  + (\partial^{\alpha}_{\beta}\mathcal{L}^{\rho} w^{n}, W_{2l_{|\alpha|,|\beta|}}\partial^{\alpha}_{\beta} w^{n}) =
(\mathcal{R}, W_{2l_{|\alpha|,|\beta|}}\partial^{\alpha}_{\beta}w^{n}),  \een
where for simplicity in this proof
\beno
\mathcal{R}\colonequals    [v \cdot \nabla_{x}, \partial^{\alpha}_{\beta}]w^{n} + \partial^{\alpha}_{\beta}\mathcal{L}^{\rho}_{r} w^{n} + \partial^{\alpha}_{\beta}\mathcal{C}^{\rho} w^{n} + \partial^{\alpha}_{\beta} \mathcal{Y}.
\eeno

Let us first deal with
the inner product involving $\mathcal{R}$ on the right-hand side of \eqref{mix-x-v-weight-solution-diff}.
By \eqref{bounded-by-energy}, we have
\beno 
|([v \cdot \nabla_{x}, \partial^{\alpha}_{\beta}]w^{n}, W_{2l_{|\alpha|,|\beta|}}\partial^{\alpha}_{\beta}w^{n})| \lesssim \mathcal{E}_{N}(w^{n}).
\eeno
By Lemma \ref{correction-term-energy-estimate}, we have
\beno 
|(\partial^{\alpha}_{\beta}\mathcal{L}^{\rho}_{r}w^{n}, W_{2l_{|\alpha|,|\beta|}}\partial^{\alpha}_{\beta}w^{n})| \lesssim \rho\mathcal{E}_{N}(w^{n}), \quad
|(\partial^{\alpha}_{\beta}\mathcal{C}^{\rho}w^{n}, W_{2l_{|\alpha|,|\beta|}}\partial^{\alpha}_{\beta}w^{n})| \lesssim \rho^{2}\mathcal{D}_{N}(w^{n}).
\eeno
Recalling \eqref{definition-of-X-8-terms}, there are three terms involving $\Gamma_{2,m}^{\rho}$ and five terms involving $\Gamma_{3}^{\rho}$. For the three terms involving $\Gamma_{2,m}^{\rho}$,
by  \eqref{Gamma-2-m-energy-estimate}, we have
\beno
|\rho^{\f{1}{2}}(\partial^{\alpha}_{\beta} \Gamma_{2,m}^{\rho}(f^{\rho, n}, w^{n}) , W_{2l_{|\alpha|,|\beta|}} \partial^{\alpha}_{\beta} w^{n})| &\lesssim& \rho^{\f{1}{2}}
\mathcal{E}^{\f{1}{2}}_{N}(f^{\rho, n})\mathcal{D}_{N}(w^{n}),
\\ |\rho^{\f{1}{2}}(\partial^{\alpha}_{\beta} \Gamma_{2,m}^{\rho}(w^{n-1},f^{\rho, n}), W_{2l_{|\alpha|,|\beta|}} \partial^{\alpha}_{\beta} w^{n})| &\lesssim& \rho^{\f{1}{2}} \mathcal{E}^{\f{1}{2}}_{N}(w^{n-1})
\mathcal{D}^{\f{1}{2}}_{N}(f^{\rho, n})\mathcal{D}^{\f{1}{2}}_{N}(w^{n})
\\ &\lesssim& \eta \rho \mathcal{D}_{N}(w^{n}) + \eta^{-1} \mathcal{E}_{N}(w^{n-1}) \mathcal{D}_{N}(f^{\rho, n}),
\\ |\rho (\partial^{\alpha}_{\beta} \Gamma_{2,m}^{\rho} (w^{n-1}, \mu^{\f12}), W_{2l_{|\alpha|,|\beta|}} \partial^{\alpha}_{\beta} w^{n})| &\lesssim& \rho \mathcal{E}^{\f{1}{2}}_{N}(w^{n-1}) \mathcal{D}^{\f{1}{2}}_{N}(w^{n}) \lesssim  \eta \rho \mathcal{D}_{N}(w^{n}) + \eta^{-1} \rho \mathcal{E}_{N}(w^{n-1}).
\eeno
For the five terms involving $\Gamma_{3}^{\rho}$,
by Theorem \ref{Gamma-3-energy-estimate}, we have
\beno
&& |(\partial^{\alpha}_{\beta} \Gamma_{3}^{\rho}(f^{\rho, n}, w^{n}, f^{\rho, n}), W_{2l_{|\alpha|,|\beta|}} \partial^{\alpha}_{\beta}w^{n})|
\\ &\lesssim& \rho \mathcal{E}^{\f{1}{2}}_{N}(f^{\rho, n})\big(\mathcal{E}^{\f{1}{2}}_{N}(w^{n})\mathcal{D}^{\f{1}{2}}_{N}(f^{\rho, n})
+\mathcal{D}^{\f{1}{2}}_{N}(w^{n})\mathcal{E}^{\f{1}{2}}_{N}(f^{\rho, n})\big)
\mathcal{D}^{\f{1}{2}}_{N}(w^{n})
\\&\lesssim&    \rho \mathcal{E}_{N}(f^{\rho, n}) \mathcal{D}_{N}(w^{n}) + \rho  \mathcal{E}_{N}(w^{n}) \mathcal{D}_{N}(f^{\rho, n}),
\\
&& |(\partial^{\alpha}_{\beta} \Gamma_{3}^{\rho}(w^{n-1}, f^{\rho, n}, f^{\rho, n}), W_{2l_{|\alpha|,|\beta|}} \partial^{\alpha}_{\beta}w^{n})|
\\ &\lesssim& \rho \mathcal{E}^{\f{1}{2}}_{N}(w^{n-1})\mathcal{E}^{\f{1}{2}}_{N}(f^{\rho, n})\mathcal{D}^{\f{1}{2}}_{N}(f^{\rho, n})
\mathcal{D}^{\f{1}{2}}_{N}(w^{n})
\\&\lesssim&  \eta \rho \mathcal{D}_{N}(w^{n}) + \eta^{-1} \rho \mathcal{E}_{N}(f^{\rho, n}) \mathcal{E}_{N}(w^{n-1}) \mathcal{D}_{N}(f^{\rho, n}),
\\
&&|( \partial^{\alpha}_{\beta} \Gamma_{3}^{\rho}(f^{\rho, n-1}, f^{\rho, n}, w^{n-1}), W_{2l_{|\alpha|,|\beta|}} \partial^{\alpha}_{\beta}w^{n})|
\\&\lesssim& \rho \mathcal{E}^{\f{1}{2}}_{N}(f^{\rho, n-1})\big(\mathcal{E}^{\f{1}{2}}_{N}(f^{\rho, n})\mathcal{D}^{\f{1}{2}}_{N}(w^{n-1})
+\mathcal{D}^{\f{1}{2}}_{N}(f^{\rho, n})\mathcal{E}^{\f{1}{2}}_{N}(w^{n-1})\big)
\mathcal{D}^{\f{1}{2}}_{N}(w^{n})
\\&\lesssim&   \eta \rho \mathcal{D}_{N}(w^{n}) + \eta^{-1} \rho \mathcal{E}_{N}(f^{\rho, n-1}) \mathcal{E}_{N}(f^{\rho, n}) \mathcal{D}_{N}(w^{n-1}) + \eta^{-1} \rho \mathcal{E}_{N}(f^{\rho, n-1}) \mathcal{D}_{N}(f^{\rho, n}) \mathcal{E}_{N}(w^{n-1}),
\\
&&|( \partial^{\alpha}_{\beta} \Gamma_{3}^{\rho}(w^{n-1}, \rho^{\f{1}{2}} \mu^{\f12}, f^{\rho, n}), W_{2l_{|\alpha|,|\beta|}} \partial^{\alpha}_{\beta}w^{n})|
\\ &\lesssim& \rho^{\f32} \mathcal{E}^{\f{1}{2}}_{N}(w^{n-1})\big(\mathcal{E}^{\f{1}{2}}_{N}(f^{\rho, n}) +\mathcal{D}^{\f{1}{2}}_{N}(f^{\rho, n})\big)
\mathcal{D}^{\f{1}{2}}_{N}(w^{n})
\\&\lesssim& \eta \rho \mathcal{D}_{N}(w^{n}) + \eta^{-1} \rho^{2} \mathcal{E}_{N}(w^{n-1})( \mathcal{E}_{N}(f^{\rho, n}) +\mathcal{D}_{N}(f^{\rho, n})),
\\
&&|(\partial^{\alpha}_{\beta} \Gamma_{3}^{\rho}(f^{\rho, n-1}, \rho^{\f{1}{2}} \mu^{\f12}, w^{n-1}), W_{2l_{|\alpha|,|\beta|}} \partial^{\alpha}_{\beta}w^{n})|
\\&\lesssim& \rho^{\f32} \mathcal{E}^{\f{1}{2}}_{N}(f^{\rho, n-1})\big(\mathcal{D}^{\f{1}{2}}_{N}(w^{n-1})
+\mathcal{E}^{\f{1}{2}}_{N}(w^{n-1})\big)
\mathcal{D}^{\f{1}{2}}_{N}(w^{n})
\\&\lesssim& \eta \rho \mathcal{D}_{N}(w^{n}) + \eta^{-1} \rho^{2} \mathcal{E}_{N}(f^{\rho, n-1})( \mathcal{E}_{N}(w^{n-1}) +\mathcal{D}_{N}(w^{n-1})).
\eeno

Now let us deal with the inner product involving $\mathcal{L}^{\rho}$ on the left-hand side of \eqref{mix-x-v-weight-solution-diff}. Since $0< \rho \leq \tilde{\rho}_{1} \leq \rho_{1} \leq \rho_{0}$, we can apply Theorem
\ref{L-energy-estimate}. As the same as \eqref{no-v-derivative-linear} and \eqref{with-v-derivative-linear}, it holds that
\beno 
 (\partial^{\alpha} \mathcal{L}^{\rho}w^{n}, W_{2l_{|\alpha|,0}} \partial^{\alpha}w^{n}) \geq
\frac{\lambda_{0}}{4}\rho\|\partial^{\alpha}w^{n}\|_{L^{2}_{x}\mathcal{L}^{s}_{l_{|\alpha|,0}+\gamma/2}}^{2} - C \rho
\mathcal{E}_{N}(w^{n}),
\\ (\partial^{\alpha}_{\beta} \mathcal{L}^{\rho}w^{n}, W_{2l_{|\alpha|,|\beta|}} \partial^{\alpha}_{\beta}w^{n})
\geq
\frac{\lambda_{0}}{8}\rho\|\partial^{\alpha}_{\beta}w^{n}\|_{L^{2}_{x}\mathcal{L}^{s}_{l_{|\alpha|,|\beta|}+\gamma/2}}^{2} -
C \rho \mathcal{E}_{N}(w^{n}) - C \rho
\sum_{\beta_{1}<\beta} \|\partial^{\alpha}_{\beta_{1}}w^{n}\|_{L^{2}_{x}\mathcal{L}^{s}_{l_{|\alpha|,|\beta|}+\gamma/2}}^{2}.
\eeno
Therefore as in the proof of Proposition \ref{a-priori-estimate-of-linear-equation} making a combination of the energy inequality of order $|\beta|=0,1, \cdots, N$,
using the fact $\rho \leq 1$, we get
\beno
&& \frac{\mathrm{d}}{\mathrm{d}t} \sum_{k=0}^{N} C_{k}\sum_{|\alpha|\leq N-k} \sum_{|\beta|=k}\|\partial^{\alpha}_{\beta}w^{n} \|^{2}_{L^{2}_{x}L^{2}_{l_{|\alpha|,|\beta|}}}  + \frac{\lambda_{0}}{16} \rho \mathcal{D}_{N}(w^{n})
\\ &\leq& C \big( \eta
 + \rho + \rho^{-\f{1}{2}}\mathcal{E}^{\f{1}{2}}_{N}(f^{\rho, n}) +  \mathcal{E}_{N}(f^{\rho, n}) \big) \rho \mathcal{D}_{N}(w^{n}) + C
\big(1+\rho\mathcal{D}_{N}(f^{\rho, n})\big)\mathcal{E}_{N}(w^{n})
\\ && + C\big( \eta^{-1} \rho \mathcal{E}_{N}(f^{\rho, n-1}) \mathcal{E}_{N}(f^{\rho, n}) + \eta^{-1} \rho^{2} \mathcal{E}_{N}(f^{\rho, n-1}) \big) \mathcal{D}_{N}(w^{n-1})
\\&&+ C \eta^{-1} \big(\mathcal{D}_{N}(f^{\rho, n}) + \rho + \rho \mathcal{E}_{N}(f^{\rho, n}) \mathcal{D}_{N}(f^{\rho, n}) + \rho \mathcal{E}_{N}(f^{\rho, n-1}) \mathcal{D}_{N}(f^{\rho, n})
\\&& \quad \quad \quad \quad + \rho^{2} \mathcal{E}_{N}(f^{\rho, n}) + \rho^{2} \mathcal{E}_{N}(f^{\rho, n-1}) \big) \mathcal{E}_{N}(w^{n-1}),
\eeno
where $C$ is a universal constant. Let us take $\eta$ such that $C \eta = \frac{\lambda_{0}}{128}$. Suppose $\rho, \delta$ verify
\ben \label{condition-on-delta-rho}
 C \rho \leq \frac{\lambda_{0}}{128}, \quad C^{2} \tilde{C} \delta  \leq (\frac{\lambda_{0}}{128})^{2},
\een
then by \eqref{uniform-estimate-of-fn} there holds
$
C \rho^{-\f{1}{2}}\mathcal{E}^{\f{1}{2}}_{N}(f^{\rho, n}) \leq \frac{\lambda_{0}}{128},  C \mathcal{E}_{N}(f^{\rho, n})  \leq \frac{\lambda_{0}}{128}.
$
By these smallness conditions on $\rho, \delta$ and the choice of $\eta$, we have
\beno
&& \frac{\mathrm{d}}{\mathrm{d}t} \sum_{k=0}^{N} C_{k}\sum_{|\alpha|\leq N-k} \sum_{|\beta|=k}\|\partial^{\alpha}_{\beta}w^{n} \|^{2}_{L^{2}_{x}L^{2}_{l_{|\alpha|,|\beta|}}}  + \frac{\lambda_{0}}{32}\rho \mathcal{D}_{N}(w^{n})
\\ &\leq& C
(1+\rho\mathcal{D}_{N}(f^{\rho, n}))\mathcal{E}_{N}(w^{n}) + C (\rho+\mathcal{D}_{N}(f^{\rho, n}))\mathcal{E}_{N}(w^{n-1}) +
C \rho^{2} \mathcal{E}_{N}(f^{\rho, n-1}) \mathcal{D}_{N}(w^{n-1}),
\eeno
where we use $\sup_{0 \leq t \leq 1, n\geq 0}\mathcal{E}_{N}(f^{\rho, n}(t)) \lesssim \rho$ and $\rho \leq 1$. Recalling \eqref{combined-energy}, we have
\beno
 \frac{\mathrm{d}}{\mathrm{d}t} \tilde{\mathcal{E}}_{N}(w^{n})   + \frac{\lambda_{0}}{32}\rho \mathcal{D}_{N}(w^{n})
\leq C
(1+ \mathcal{D}_{N}(f^{\rho, n})) \big(\tilde{\mathcal{E}}_{N}(w^{n}) + \tilde{\mathcal{E}}_{N}(w^{n-1}) \big) +
C \rho^{2} \mathcal{E}_{N}(f^{\rho, n-1}) \mathcal{D}_{N}(w^{n-1}).
\eeno
For simplicity, let us define $y^{n}(t)\colonequals    \tilde{\mathcal{E}}_{N}(w^{n}(t)), x^{n}(t)\colonequals    \mathcal{D}_{N}(w^{n}(t)),  a^{n}(t)\colonequals    1+ \mathcal{D}_{N}(f^{\rho, n}(t))$, then the above inequality is
\beno
 \frac{\mathrm{d}}{\mathrm{d}t} y^{n}   + \frac{\lambda_{0}}{32}\rho x^{n}
\leq C
a^{n} \big(y^{n}+ y^{n-1} \big) +
C \rho^{2} \mathcal{E}_{N}(f^{\rho, n-1}) x^{n-1}.
\eeno
Let $g^{n}(t)\colonequals    \exp(-C\int_{0}^{t} a^{n}(\tau) \mathrm{d}\tau)$, then $\frac{\mathrm{d}}{\mathrm{d}t}g^{n}(t)\colonequals    -Ca^{n}(t) g^{n}(t)$ and thus
\beno
 \frac{\mathrm{d}}{\mathrm{d}t} (g^{n}y^{n})   + \frac{\lambda_{0}}{32}\rho g^{n}x^{n}
\leq C
a^{n} g^{n} y^{n-1}  +
C \rho^{2} \mathcal{E}_{N}(f^{\rho, n-1}) g^{n} x^{n-1}.
\eeno
When $n \geq 1$, recall that $y^{n}(0) = 0$ and thus
\beno
g^{n}(t)y^{n}(t)   + \frac{\lambda_{0}}{32}\rho \int_{0}^{t}g^{n}x^{n} \mathrm{d}\tau
\leq C \int_{0}^{t}
a^{n} g^{n} y^{n-1} \mathrm{d}\tau +
C \rho^{2} \int_{0}^{t} \mathcal{E}_{N}(f^{\rho, n-1}) g^{n} x^{n-1} \mathrm{d}\tau.
\eeno
Note that $\int_{0}^{1} a^{n}(\tau) \mathrm{d}\tau \leq 1 + \tilde{C}\delta \lesssim 1$ by \eqref{uniform-estimate-of-fn}, then $1 \lesssim g^{n}(t) \leq 1$ for any $0 \leq t \leq 1$. Then with a different constant $C$, we get for any $0 \leq t \leq 1$,
\beno
y^{n}(t)  + \frac{\lambda_{0}}{32} \rho \int_{0}^{t} x^{n} \mathrm{d}\tau
\leq C \int_{0}^{t}
a^{n} y^{n-1} \mathrm{d}\tau +
C \rho^{2} \int_{0}^{t} \mathcal{E}_{N}(f^{\rho, n-1}) x^{n-1} \mathrm{d}\tau.
\eeno
With a different $C$, we get for any $0 \leq t \leq 1$,
\beno
\sup_{0 \leq \tau \leq t}y^{n}(\tau)  + \frac{\lambda_{0}}{32} \rho \int_{0}^{t} x^{n} \mathrm{d}\tau
\leq C \int_{0}^{t}
a^{n}  \mathrm{d}\tau  \sup_{0 \leq \tau \leq t}y^{n-1}(\tau) +
C \rho^{2} \int_{0}^{t} \mathcal{E}_{N}(f^{\rho, n-1}) x^{n-1} \mathrm{d}\tau.
\eeno
If $\delta$ is small enough such that
\ben \label{condition-on-delta}
C \tilde{C} \delta \leq \frac{\lambda_{0}}{64},
\een
then by \eqref{uniform-estimate-of-fn} there holds
$
C \rho \sup_{0 \leq t \leq 1} \mathcal{E}_{N}(f^{\rho, n}(t)) \leq C \tilde{C} \delta \rho^{2} \leq \frac{\lambda_{0}}{64},
$
which gives
\beno
\sup_{0 \leq \tau \leq t}y^{n}(\tau)  +  \frac{\lambda_{0}}{32} \rho \int_{0}^{t} x^{n} \mathrm{d}\tau
\leq C \int_{0}^{t}
a^{n}  \mathrm{d}\tau  \sup_{0 \leq \tau \leq t}y^{n-1}(\tau)  +
\frac{\lambda_{0}}{64} \rho \int_{0}^{t}  x^{n-1} \mathrm{d}\tau.
\eeno
Suppose $t, \delta$ are small enough such that
\ben \label{condition-on-t-and-delta}
C t \leq \frac{1}{4}, \quad C \tilde{C} \delta \leq \frac{1}{4},
\een
then by \eqref{uniform-estimate-of-fn} there holds
\beno
C \int_{0}^{t}
a^{n}  \mathrm{d}\tau = C \int_{0}^{t}  \big(1+ \mathcal{D}_{N}(f^{\rho, n}(\tau))\big) \mathrm{d}\tau \leq C(t + \tilde{C} \delta) \leq \f{1}{2}.
\eeno
Let $T^{*} = \frac{1}{4C}$, we get
\beno
\sup_{0 \leq \tau \leq T^{*}} y^{n}(\tau)  + \frac{\lambda_{0}}{32} \rho \int_{0}^{T^{*}} x^{n} \mathrm{d}\tau
\leq \f{1}{2} \bigg( \sup_{0 \leq \tau \leq T^{*}}y^{n-1}(\tau)  +
\frac{\lambda_{0}}{32} \rho \int_{0}^{T^{*}}  x^{n-1} \mathrm{d}\tau \bigg).
\eeno
Note that the above inequality is valid for all $n \geq 1$. As a result, we get for all $n \geq 1$,
\beno
\sup_{0 \leq \tau \leq T^{*}} y^{n}(\tau)  + \frac{\lambda_{0}}{32} \rho \int_{0}^{T^{*}} x^{n} \mathrm{d}\tau
\leq \frac{1}{2^{n}} (\sup_{0 \leq \tau \leq T^{*}} y^{0}(\tau)  + \frac{\lambda_{0}}{32} \rho \int_{0}^{T^{*}} x^{0} \mathrm{d}\tau).
\eeno
Note that $w^{0} = f^{\rho, 1}-f^{\rho, 0} = f^{\rho, 1}$.
By \eqref{uniform-estimate-of-fn}, it holds that
\beno
y^{0} = \tilde{\mathcal{E}}_{N}(w^{0}) = \tilde{\mathcal{E}}_{N}(f^{\rho, 1}) \lesssim \mathcal{E}_{N}(f^{\rho, 1})  \lesssim \delta \rho,
\\ \rho \int_{0}^{T^{*}} x^{0} \mathrm{d}\tau  = \rho \int_{0}^{T^{*}} \mathcal{D}_{N}(w^{0}) \mathrm{d}\tau = \rho \int_{0}^{T^{*}} \mathcal{D}_{N}(f^{\rho, 1}) \mathrm{d}\tau
\lesssim \delta \rho.
\eeno
Recalling \eqref{combined-energy} and $y^{n} = \tilde{\mathcal{E}}_{N}(w^{n}) =
\tilde{\mathcal{E}}_{N}(f^{\rho, n+1}-f^{\rho, n})$, we arrive at
\beno
\sup_{0 \leq \tau \leq T^{*}} \mathcal{E}_{N}(f^{\rho, n+1}(\tau)-f^{\rho, n}(\tau))
\leq \sup_{0 \leq \tau \leq T^{*}} y^{n}(\tau)  + \frac{\lambda_{0}}{32} \rho \int_{0}^{T^{*}} x^{n} \mathrm{d}\tau
\leq \frac{C}{2^{n}} \delta \rho,
\eeno
and so $\{f^{\rho, n}\}_{n \geq 0}$ is a Cauchy sequence in $L^{\infty}([0,T^{*}]; \mathcal{E}_{N})$. The sequence $\{f^{\rho, n}\}_{n \geq 0}$ has a limit $ f^{\rho}  \in L^{\infty}([0,T^{*}]; \mathcal{E}_{N})$ which is the solution of the Cauchy problem \eqref{linearized-quantum-Boltzmann-eq} verifying \eqref{uniform-estimate-of-solution} thanks to \eqref{uniform-estimate-of-fn} and  \eqref{solution-positivity}.
As for uniqueness, see the proof of Theorem \ref{global-well-posedness}  at the end of Section \ref{global}.
Let $C$ be the largest constant in \eqref{condition-on-delta-rho}, \eqref{condition-on-delta} and \eqref{condition-on-t-and-delta}, we set
\beno
\tilde{\rho}_{2} \colonequals   \frac{\lambda_{0}}{128C}, \quad \tilde{\delta}_{2}\colonequals    \min\{ \frac{1}{C^{2} \tilde{C}} (\frac{\lambda_{0}}{128})^{2},   \frac{\lambda_{0}}{64C \tilde{C}},   \frac{1}{4 C \tilde{C}}\}.
\eeno
Finally, set $\rho_{2} \colonequals   \min \{\tilde{\rho}_{1}, \tilde{\rho}_{2}\}, \delta_{1} \colonequals   \min \{\tilde{\delta}_{1}, \tilde{\delta}_{2}\}$.
Then if $0< \rho \leq \rho_{2}$ and $0< \delta \leq \delta_{1}$,
all the conditions in the proof are fulfilled.
\end{proof}

\section{A priori estimate and global well-posedness}   \label{global}
This section is devoted to the proof to Theorem \ref{global-well-posedness}.
We first provide the {\it a priori} estimate for the equation \eqref{linearized-quantum-Boltzmann-eq} in Theorem \ref{a-priori-estimate-LBE}.
From which together with the local existence result in Theorem \ref{local-well-posedness-LBE}, the global well-posedness result
Theorem \ref{global-well-posedness} is constructively established by a standard continuity argument at the end of this section.

\subsection{A priori estimate of a general equation} This subsection is devoted to some {\it{a priori}} estimate for
 the following equation
\ben \label{lBE} \partial_{t}f + v\cdot \nabla_{x} f + \mathcal{L}^{\rho}f= g,  \quad t>0, x \in \mathbb{T}^{3}, v \in \mathbb{R}^{3},
\een
where $g$ is a given function.

Let $f$ be a solution to \eqref{lBE}.
Recalling the formula \eqref{linear-combination-of-basis} of projection operator $\mathbb{P}_{\rho}$, we denote
\ben \label{definition-f-1} (\mathbb{P}_{\rho} f)(t,x,v) = (a(t,x) + b(t,x) \cdot v + c(t,x)|v|^{2})N,\een
where $N=N_{\rho}$ and
\ben \label{a-b-c-for-t-x}
(a(t,x), b(t,x), c(t,x)) \colonequals     (a^{f(t,x,\cdot)}_{\rho}, b^{f(t,x,\cdot)}_{\rho}, c^{f(t,x,\cdot)}_{\rho}).
\een
Here we recall \eqref{explicit-defintion-of-abc} for the definition of $(a^{f(t,x,\cdot)}_{\rho}, b^{f(t,x,\cdot)}_{\rho}, c^{f(t,x,\cdot)}_{\rho})$ for fixed $t,x$.
Note that in \eqref{a-b-c-for-t-x} we omit $f, \rho$. However, we should
always keep in mind that $(a, b, c)$ are functions of $(t,x)$ originating from the solution $f$ to \eqref{lBE} for fixed $\rho$.

We set $f_{1} \colonequals     \mathbb{P}_{\rho} f$ and $f_{2} \colonequals     f - \mathbb{P}_{\rho} f$.
We first recall some basics of macro-micro decomposition. Plugging the macro-micro decomposition $f = f_{1} + f_{2}$ into \eqref{lBE}
and using the fact $ \mathcal{L}^{\rho} f_{1} =  0$, we get
\ben \label{macro-micro-LBE-2} \partial_{t}f_{1} + v\cdot \nabla_{x} f_{1}  = - \partial_{t}f_{2} - v\cdot \nabla_{x} f_{2} - \mathcal{L}^{\rho}f_{2} + g.\een

Recalling \eqref{definition-f-1},
the left-hand of \eqref{macro-micro-LBE-2} is
\ben \label{left-write-out}
\partial_{t}f_{1} + v\cdot \nabla_{x} f_{1} = (\partial_{t}a + \sum_{i=1}^{3}\partial_{t}b_{i} v_{i} + \partial_{t}c |v|^{2} )N +  (\sum_{i=1}^{3} \partial_{i}a v_{i} + \sum_{i \neq j} \partial_{i}b_{j} v_{i}v_{j} + \sum_{i=1}^{3} \partial_{i}c v_{i}  |v|^{2} )N.
\een
Here $\partial_{i} = \partial_{x_{i}}$ for $i=1,2,3, b=(b_{1},b_{2},b_{3})$ and $v=(v_{1},v_{2},v_{3})$.
We order the 13 functions of $v$ on the right-hand side of \eqref{left-write-out} as
\beno e_{1} = N, \quad e_{2} = v_{1}N, \quad e_{3} = v_{2}N, \quad e_{4} = v_{3}N, \quad e_{5} = v_{1}^{2}N, \quad e_{6} = v_{2}^{2}N, \quad e_{7} = v_{3}^{2}N,  \\ e_{8} =v_{1}v_{2}N, \quad e_{9} = v_{2}v_{3}N, \quad e_{10} = v_{3}v_{1}N, \quad  e_{11} = |v|^{2}v_{1}N, \quad e_{12} = |v|^{2}v_{2}N, \quad e_{13} = |v|^{2}v_{3}N.
\eeno
We emphasize that $e_{i}$ depends on $\rho$ through $N = N_{\rho}$. We also order the 13 functions of $(t, x)$ on the right-hand side of \eqref{left-write-out} as
\beno x_{1} \colonequals    \partial_{t}a, \quad x_{2} \colonequals     \partial_{t}b_{1}+ \partial_{1} a, \quad x_{3} \colonequals     \partial_{t}b_{2}+ \partial_{2} a, \quad x_{4} \colonequals     \partial_{t}b_{3}+ \partial_{3} a,
\\ x_{5} \colonequals     \partial_{t}c+ \partial_{1} b_{1}, \quad x_{6} \colonequals     \partial_{t}c+ \partial_{2} b_{2}, \quad x_{7} \colonequals     \partial_{t}c+ \partial_{3} b_{3},
\\ x_{8} \colonequals    \partial_{1}b_{2} + \partial_{2}b_{1}, \quad x_{9} \colonequals    \partial_{2}b_{3} + \partial_{3}b_{2}, \quad x_{10} \colonequals    \partial_{3}b_{1} + \partial_{1}b_{3}, \\ x_{11} \colonequals     \partial_{1}c, \quad x_{12} \colonequals     \partial_{2}c, \quad x_{13} \colonequals     \partial_{3}c.
\eeno
Use $\mathsf{T}$ to denote vector transpose. For simplicity, we define two column vectors
\beno
E \colonequals     (e_{1}, \cdots, e_{13})^{\mathsf{T}}, \quad X \colonequals     (x_{1}, \cdots, x_{13})^{\mathsf{T}}.
\eeno
With these two column vectors, \eqref{left-write-out}  becomes $\partial_{t}f_{1} + v\cdot \nabla_{x} f_{1} = E^{\mathsf{T}} X$ and thus
 \eqref{macro-micro-LBE-2} can be written as
\beno
 E^{\mathsf{T}} X = - \partial_{t}f_{2} - v\cdot \nabla_{x} f_{2} - \mathcal{L}^{\rho}f_{2} + g.
\eeno
Taking inner product with $E$
in the space $L^{2}(\mathbb{R}^{3})$, since $X$ depends on $(t,x)$ but not on $v$,
we get
\beno
\langle E, E^{\mathsf{T}} \rangle  X =  \langle E, - \partial_{t}f_{2} - v\cdot \nabla_{x} f_{2} - \mathcal{L}^{\rho}f_{2} + g \rangle.
\eeno
The $13 \times 13$ matrix $\langle E, E^{\mathsf{T}} \rangle = (\langle e_{i}, e_{j} \rangle)_{1\leq i \leq 13, 1\leq j \leq 13}$ is invertible for small $\rho$ and so
\beno 
  X =  (\langle E, E^{\mathsf{T}} \rangle)^{-1} \langle E, - \partial_{t}f_{2} - v\cdot \nabla_{x} f_{2} - \mathcal{L}^{\rho}f_{2} + g \rangle.
\eeno
For simplicity, let us define
\beno
Y = (Y^{(0)}, \{Y^{(1)}_{i}\}_{1\leq i \leq 3}, \{Y^{(2)}_{i}\}_{1\leq i \leq 3}, \{Y^{(2)}_{ij} \}_{1\leq i < j  \leq 3}, \{Y^{(3)}_{i}\}_{1\leq i \leq 3})^{\mathsf{T}} \colonequals     (\langle E, E^{\mathsf{T}} \rangle)^{-1} \langle E, f_{2} \rangle,
\\
Z  \colonequals     (\langle E, E^{\mathsf{T}} \rangle)^{-1} \langle E, - v\cdot \nabla_{x} f_{2} - \mathcal{L}^{\rho}f_{2} + g \rangle.
\eeno
Then we get the following system that consists of 13 macroscopic equations
\ben \label{macroscopic-system}
 X =  - \partial_{t} Y + Z.
\een

Observing  $|E| \lesssim \mu^{\f14}$ and using \eqref{upper-bound-of-L},
the functions $Y, Z$  can be controlled as:
\begin{lem} \label{estimate-on-fln-tilde} Let $0< \rho \leq \f{1}{2}(2 \pi)^{\f32}$. It holds that
\beno  \sum_{|\alpha|\leq N}|\partial^{\alpha}Y|^{2}_{L^{2}_{x}} \lesssim \|f_{2}\|^{2}_{H^{N}_{x}\mathcal{L}^{s}_{\gamma/2}}, \quad \sum_{|\alpha|\leq N-1}|\partial^{\alpha}Z|^{2}_{L^{2}_{x}} \lesssim \|f_{2}\|^{2}_{H^{N}_{x}\mathcal{L}^{s}_{\gamma/2}} + \sum_{|\alpha| \leq N-1}
|\langle E, \partial^{\alpha}g \rangle|^{2}_{L^{2}_{x}}. \eeno
\end{lem}


We now estimate the dynamics of $(a,b,c)$ in the following lemma.
\begin{lem}\label{estimate-for-ptabc} Let $0<\rho \leq \frac{1}{160}$. It holds that
\ben \label{simple-estimate-of-pt-a}
\sum_{|\alpha|\leq N-1}|\partial^{\alpha}\partial_{t}a|^{2}_{L^{2}_{x}} \lesssim \|f_{2}\|_{H^{N}_{x}\mathcal{L}^{s}_{\gamma/2}}^{2} +\sum_{|\alpha| \leq N-1}
|\langle E, \partial^{\alpha}g \rangle|^{2}_{L^{2}_{x}},
\\ \label{estimate-of-pt-b-c}
\sum_{|\alpha|\leq N-1}|\partial^{\alpha}\partial_{t} (b,c)|^{2}_{L^{2}_{x}} \lesssim |\nabla_{x}(a,b,c)|^{2}_{H^{N-1}_{x}}
+ \|f_{2}\|_{H^{N}_{x}\mathcal{L}^{s}_{\gamma/2}}^{2} +\sum_{|\alpha| \leq N-1}
|\langle E, \partial^{\alpha}g \rangle|^{2}_{L^{2}_{x}}.
\een
\end{lem}
\begin{proof} Recall \eqref{explicit-defintion-of-abc}.
In $L^{2}(\mathbb{R}^{3})$, taking inner products between equation \eqref{lBE} and the five functions $l_{1}N -  l_{2}N|v|^{2}, l_{3}N v_{j}, l_{4} N|v|^{2} - l_{2} N$ respectively, using $
\langle  \mathcal{L}^{\rho}f, N \rangle = \langle  \mathcal{L}^{\rho}f, N v_{j} \rangle = \langle  \mathcal{L}^{\rho}f, N|v|^{2} \rangle =0,
$
we get
\beno \partial_{t} a +  \langle v\cdot \nabla_{x} f, l_{1}N -  l_{2}N|v|^{2} \rangle = \langle g, l_{1}N -  l_{2}N|v|^{2}\rangle, \quad \partial_{t} b +  \langle v\cdot \nabla_{x} f, l_{3}N v  \rangle = \langle g,  l_{3}N v \rangle,
\\ \partial_{t} c + \langle v\cdot \nabla_{x} f, l_{4} N|v|^{2} - l_{2} N  \rangle = \langle g,  l_{4} N|v|^{2} - l_{2} N \rangle.
\eeno
Since $\langle v_{i}N, N\rangle = \langle v_{i} v_{j}N, v_{k}N\rangle = \langle v_{i}|v|^{2}N, |v|^{2}N\rangle = 0$ for $i, j, k \in \{1,2,3\}$.
Recalling the definition of $l_{1}, l_{2}$ in \eqref{defintion-of-l-i}, it is straightforward to see
\beno
\langle v\cdot \nabla_{x} f_{1}, l_{1}N -  l_{2}N|v|^{2} \rangle = 0, \quad  \langle v\cdot \nabla_{x} f_{1}, l_{3}N v  \rangle = \langle v \cdot \nabla_{x} (a N + c N |v|^{2}), l_{3}N v  \rangle,
\\ \langle v\cdot \nabla_{x} f_{1}, l_{4} N|v|^{2} - l_{2} N  \rangle = \langle v \cdot \nabla_{x} (b \cdot v N), l_{4} N|v|^{2} - l_{2} N  \rangle,
\eeno
which gives the local conservation laws
\ben \label{local-conservation-laws} \left\{ \begin{aligned}
 & \partial_{t} a = \langle g - v\cdot \nabla_{x} f_{2}, l_{1}N -  l_{2}N|v|^{2}\rangle,
\\&
\partial_{t} b +  l_{3} \langle Nv_{1}, N v_{1}  \rangle \nabla_{x} a +  l_{3} \langle N|v|^{2}v_{1}, N v_{1}  \rangle \nabla_{x} c
 = \langle g - v\cdot \nabla_{x} f_{2},  l_{3}N v \rangle,
\\&
\partial_{t} c + \langle N v_{1}^{2}, l_{4} N|v|^{2} - l_{2} N  \rangle \nabla_{x} \cdot b = \langle g - v\cdot \nabla_{x} f_{2},  l_{4} N|v|^{2} - l_{2} N \rangle.
\end{aligned} \right.
\een
By \eqref{boundedness-of-li}, $0 \leq l_{i} \lesssim 1$ for $i=1,2,3,4$. Recalling $N \lesssim \mu^{\f{1}{2}}$,
we get
 \eqref{simple-estimate-of-pt-a} and \eqref{estimate-of-pt-b-c} directly from \eqref{local-conservation-laws}.
\end{proof}
Let us recall the temporal  energy functional  $\mathcal{I}_{N}(f)$ in \cite{duan2008cauchy} as
\ben \label{interactive-INf} \mathcal{I}_{N}(f) &\colonequals&     \sum_{|\alpha|\leq N-1} (\mathcal{I}^{a}_{\alpha}(f) +
\sum_{j=1}^{3}\mathcal{I}^{b_{j}}_{\alpha}(f)+\mathcal{I}^{c}_{\alpha}(f)),
\\ \nonumber 
\mathcal{I}^{a}_{\alpha}(f) &\colonequals&    -\sum_{j=1}^{3} \langle \partial_{j} \pa^{\alpha}b_{j}, \pa^{\alpha} a\rangle_{x} - \sum_{j=1}^{3} \langle \partial_{j} \pa^{\alpha} Y^{(1)}_{j}, \pa^{\alpha} a\rangle_{x},
\\ \nonumber 
\mathcal{I}^{b_{j}}_{\alpha}(f) &\colonequals&    \langle \sum_{i\neq j, i=1}^{3} \partial_{j}\pa^{\alpha}Y^{(2)}_{i} - \sum_{i\neq j, i=1}^{3} \partial_{i} \pa^{\alpha}Y^{(2)}_{ij} - 2 \partial_{j}\pa^{\alpha}Y^{(2)}_{j}, \pa^{\alpha} b_{j}\rangle_{x},
\\ \nonumber 
\mathcal{I}^{c}_{\alpha}(f) &\colonequals&    - \sum_{j=1}^{3} \langle \partial_{j} \pa^{\alpha}Y^{(3)}_{j}, \pa^{\alpha} c\rangle_{x}.
\een

With Lemma \ref{estimate-on-fln-tilde} and Lemma \ref{estimate-for-ptabc}, based on the macroscopic system \eqref{macroscopic-system}, using integration by parts to balance derivative
it is standard to derive
\begin{lem}\label{estimate-for-highorder-abc}  Let $N \geq 1, T>0, 0<\rho \leq \frac{1}{160}$. Let $f \in L^{\infty}([0,T]; H^{N}_{x}L^{2})$ be a solution to \eqref{lBE},
then there exists a universal constant $C > 0$ such that
\ben \label{solution-property-part2} \frac{\mathrm{d}}{\mathrm{d}t}\mathcal{I}_{N}(f) + \f{1}{2}|\nabla_{x}(a,b,c)|^{2}_{H^{N-1}_{x}} \leq C(\|f_{2}\|_{H^{N}_{x}\mathcal{L}^{s}_{\gamma/2}}^{2}
+\sum_{|\alpha| \leq N-1}|\langle E, \partial^{\alpha}g \rangle|^{2}_{L^{2}_{x}}),\een
where $\mathcal{I}_{N}(f)$ is defined in \eqref{interactive-INf} and satisfying
\ben \label{temporal-bounded-by-norm}
|\mathcal{I}^{N}(f)| \leq C \|f\|^{2}_{H^{N}_{x}L^{2}}.
\een
\end{lem}
For a rigorous proof of Lemma \ref{estimate-for-highorder-abc}, one can refer to
\cite{duan2008cauchy}.
As a result of the Poincar\'{e} inequality and Lemma \ref{estimate-for-highorder-abc}, we have
\begin{lem}\label{estimate-for-highorder-abc-full} Let $N \geq 1, T>0, 0<\rho \leq \frac{1}{160}$. Let $f \in L^{\infty}([0,T]; H^{N}_{x}L^{2})$ be a solution to \eqref{lBE} verifying
\ben \label{a-priori-assumption}
\int (1,v,|v|^{2}) N_{\rho} f(t,x,v)\mathrm{d}x \mathrm{d}v = 0,
\een
then there exist two universal constant $C_{0}, C > 0$ such that
\ben \label{solution-property-part2-full}
\frac{\mathrm{d}}{\mathrm{d}t}\mathcal{I}_{N}(f) + C_{0}|(a,b,c)|^{2}_{H^{N}_{x}} \leq C(\|f_{2}\|_{H^{N}_{x}\mathcal{L}^{s}_{\gamma/2}}^{2}
+\sum_{|\alpha| \leq N-1}|\langle E, \partial^{\alpha}g \rangle|^{2}_{L^{2}_{x}}).\een
\end{lem}
\begin{proof} Recalling \eqref{a-b-c-for-t-x} and \eqref{explicit-defintion-of-abc}, by \eqref{a-priori-assumption}, we have
$
\int (a,b,c) \mathrm{d}x = 0.
$
Then by Poincar\'{e} inequality on the torus $\mathbb{T}^{3}$, there is a universal constant $C$ such that
$|(a,b,c)|_{L^{2}_{x}}  \leq C |\nabla(a,b,c)|_{L^{2}_{x}}$ and thus there is a universal constant $C_{0}>0$ such that
\ben \label{given-by-Poincare-inequality}
\f{1}{2}|\nabla_{x}(a,b,c)|^{2}_{H^{N-1}_{x}} \geq C_{0} |(a,b,c)|^{2}_{H^{N}_{x}}.
\een
Plugging \eqref{given-by-Poincare-inequality} into \eqref{solution-property-part2}, we get the desired result.
\end{proof}

We now give an interpolation result.
\begin{lem}\label{interpolation-to-symbol}
Let $l \geq 0, m \geq 1$. For any $\eta>0$, there exists a constant $C_{\eta}$ such that
\beno |f|_{H^{m}_{l}}^{2} \leq \eta \sum_{|\beta|= m} |\partial_{\beta}f|_{\mathcal{L}^{s}_{l}}^{2} + C_{\eta} |f|_{\mathcal{L}^{m-1,s}_{l}}^{2}. \eeno
\end{lem}
\begin{proof} By interpolation and
using  \eqref{relation-with-weighted-Sobolev}, we get
\beno
|f|_{H^{m}_{l}}^{2} \lesssim \eta |f|_{H^{m+s}_{l}}^{2} + C_{\eta} |f|_{H^{m-1+s}_{l}}^{2}  \lesssim \eta |f|_{\mathcal{L}^{m,s}_{l}}^{2} + C_{\eta} |f|_{\mathcal{L}^{m-1,s}_{l}}^{2}.
\eeno
Recalling \eqref{definition-of-L-n}, there holds $|f|_{\mathcal{L}^{m,s}_{l}}^{2} = \sum_{|\beta|= m} |\partial_{\beta}f|_{\mathcal{L}^{s}_{l}}^{2} + |f|_{\mathcal{L}^{m-1,s}_{l}}^{2}$ which ends the proof.
\end{proof}

We derive the following {\it a priori} estimate for equation \eqref{lBE}.

\begin{prop}\label{essential-estimate-of-micro-macro} Let $N \geq 1, T>0, 0 < \rho \leq \rho_{0}$. Let $f \in L^{\infty}([0,T]; \mathcal{E}_{N})$ be a solution to \eqref{lBE} verifying \eqref{a-priori-assumption},
then for any $0 \leq t \leq T$ there holds
\ben \label{essential-micro-macro-result-final} \frac{\mathrm{d}}{\mathrm{d}t}\Xi^{\rho}_{N}(f)+ \frac{\lambda_{0}}{16}\rho\mathcal{D}_{N}(f) &\lesssim& \sum_{|\alpha| \leq N}
|(\pa^{\alpha}g, \pa^{\alpha}f)|+\sum_{|\alpha|+|\beta| \leq N}
|(\partial^{\alpha}_{\beta}g, W_{2l_{|\alpha|,|\beta|}}\partial^{\alpha}_{\beta}f)|
\\&&+\sum_{|\alpha| \leq N-1}|\langle E, \partial^{\alpha}g \rangle|^{2}_{L^{2}_{x}},  \nonumber \een
where
\beno \Xi^{\rho}_{N}(f)\colonequals    K_{-2}\rho^{-2N+1}\mathcal{I}_{N}(f)
+ K_{-1}\rho^{-2N}\|f\|^{2}_{H^{N}_{x}L^{2}}+\sum_{j=0}^{N}K_{j}\rho^{2j-2N}\sum_{|\alpha|\leq N-j, |\beta|=j}
\| \partial^{\alpha}_{\beta} f\|^{2}_{L^{2}_{x}L^{2}_{l_{|\alpha|,|\beta|}}},\eeno
for some universal constants $K_{j} \geq 1, -2 \leq j \leq N$. Recalling \eqref{definition-energy-and-dissipation}, it holds that
\ben \label{equivalence-between-energy}
\mathcal{E}_{N}(f) \leq \Xi^{\rho}_{N}(f) \lesssim \rho^{-2N} \mathcal{E}_{N}(f).
\een
\end{prop}

\begin{proof}
Note that $\Xi^{\rho}_{N}(f)$ is a combination of several functionals. We already have $\mathcal{I}_{N}(f)$ from Lemma \ref{estimate-for-highorder-abc-full}. That is, the solution $f$ verifies \eqref{solution-property-part2-full}.
We add the other functionals step by step.

{\textit {Step 1: Pure $x$-derivative without weight $\|f\|^{2}_{H^{N}_{x}L^{2}}$.}}
Applying $\partial^{\alpha}$ to equation \eqref{lBE}, taking inner product with $\partial^{\alpha}f$, we have
\beno \f{1}{2} \frac{\mathrm{d}}{\mathrm{d}t} \|\partial^{\alpha}f\|^{2}_{L^{2}} + (\mathcal{L}^{\rho}\partial^{\alpha}f, \partial^{\alpha}f) =  (\partial^{\alpha}g, \partial^{\alpha}f).\eeno
Recalling $f_{2} = (\mathbb{I}-\mathbb{P}_{\rho})f$ and $(\partial^{\alpha}f)_{2} = \partial^{\alpha}f_{2}$.
By the lower bound in Theorem \ref{main-theorem}, taking sum over $|\alpha|\leq N$, we have
\ben \label{solution-property-part-g}\f{1}{2}\frac{\mathrm{d}}{\mathrm{d}t}\|f\|^{2}_{H^{N}_{x}L^{2}} + \lambda_{0} \rho \|f_{2}\|^{2}_{H^{N}_{x}\mathcal{L}^{s}_{\gamma/2}} \leq \sum_{|\alpha| \leq N}
|(\pa^{\alpha}g, \pa^{\alpha}f)|.\een
Then $\eqref{solution-property-part-g} \times 2C_{1} + \eqref{solution-property-part2-full} \times \rho$ gives
\beno 
&&\frac{\mathrm{d}}{\mathrm{d}t}(\rho \mathcal{I}_{N}(f) + C_{1}\|f\|^{2}_{H^{N}_{x}L^{2}})+ (C_{0}\rho|(a,b,c)|^{2}_{H^{N}_{x}}+2C_{1}\lambda_{0}\rho
\|f_{2}\|^{2}_{H^{N}_{x}\mathcal{L}^{s}_{\gamma/2}})
\\ 
&\leq& 2C_{1}\sum_{|\alpha| \leq N}
|(\pa^{\alpha}g, \pa^{\alpha}f)|+
C\rho(\|f_{2}\|_{H^{N}_{x}\mathcal{L}^{s}_{\gamma/2}}^{2}
+\sum_{|\alpha| \leq N-1}|\langle E, \partial^{\alpha}g \rangle|^{2}_{L^{2}_{x}}).
\eeno
Thanks to \eqref{temporal-bounded-by-norm}, we can choose $C_{1}$ large enough such that
\beno
\f{1}{2} C_{1}\|f\|^{2}_{H^{N}_{x}L^{2}} \leq \rho \mathcal{I}_{N}(f) + C_{1}\|f\|^{2}_{H^{N}_{x}L^{2}} \leq \frac{3}{2} C_{1}\|f\|^{2}_{H^{N}_{x}L^{2}}, \quad C_{1}\lambda_{0} \geq C, \quad C_{1}\lambda_{0} \geq C_{0},
\eeno
and then
\ben \label{essential-micro-macro-result-2} &&\frac{\mathrm{d}}{\mathrm{d}t}(\rho \mathcal{I}_{N}(f) + C_{1}\|f\|^{2}_{H^{N}_{x}L^{2}})+ C_{0}\rho(|(a,b,c)|^{2}_{H^{N}_{x}}+
\|f_{2}\|^{2}_{H^{N}_{x}\mathcal{L}^{s}_{\gamma/2}}) \\&\lesssim& \sum_{|\alpha| \leq N}
|(\pa^{\alpha}g, \pa^{\alpha}f)|+
\sum_{|\alpha| \leq N-1}|\langle E, \partial^{\alpha}g \rangle|^{2}_{L^{2}_{x}}. \nonumber \een

 {\textit {Step 2: Pure $x$-derivative with weight
 $\sum_{|\alpha|\leq N} \|\partial^{\alpha} f\|^{2}_{L^{2}_{x}L^{2}_{l_{|\alpha|,0}}}$.}}
Applying $\partial^{\alpha}$ to equation \eqref{lBE}, taking inner product with $W_{2l_{|\alpha|,0}}\partial^{\alpha}f$, taking sum over $|\alpha|\leq N$, we have
\beno&& \f{1}{2}\frac{\mathrm{d}}{\mathrm{d}t} \sum_{|\alpha|\leq N}\|\partial^{\alpha} f\|^{2}_{L^{2}_{x}L^{2}_{l_{|\alpha|,0}}}  + \sum_{|\alpha| \leq N} (\mathcal{L}^{\rho}\pa^{\alpha}f, W_{2l_{|\alpha|,0}}\pa^{\alpha}f) = \sum_{|\alpha| \leq N}(\pa^{\alpha}g, W_{2l_{|\alpha|,0}}\pa^{\alpha}f). \nonumber \eeno
By \eqref{no-v-derivative}, we have
\ben \label{H-Nx-L2l}
\frac{\mathrm{d}}{\mathrm{d}t}\sum_{|\alpha|\leq N}\|\partial^{\alpha} f\|^{2}_{L^{2}_{x}L^{2}_{l_{|\alpha|,0}}}  + \frac{\lambda_{0}}{2} \rho \sum_{|\alpha|\leq N}\|\partial^{\alpha} f\|^{2}_{L^{2}_{x}\mathcal{L}^{s}_{l_{|\alpha|,0}+\gamma/2}} \leq C \rho \|f\|^{2}_{H^{N}_{x}L^{2}_{\gamma/2}} +
2 \sum_{|\alpha| \leq N}|(\pa^{\alpha}g, W_{2l_{|\alpha|,0}}\pa^{\alpha}f)|.  \een
Note that
$
\|f\|^{2}_{H^{N}_{x}L^{2}_{\gamma/2}} \leq C (|(a,b,c)|^{2}_{H^{N}_{x}}+\|f_{2}\|^{2}_{H^{N}_{x}\mathcal{L}^{s}_{\gamma/2}})
$ for some universal constant $C$.
Taking $C_{2}$ large enough such that
$
C_{0} C_{2} \geq 2 C^{2}$ and $ \f{1}{2} C_{0} C_{2} \geq  \frac{\lambda_{0}}{16},
$
the combination  $\eqref{essential-micro-macro-result-2} \times C_{2} + \eqref{H-Nx-L2l}$ gives
\ben \label{essential-micro-macro-result-3} &&\frac{\mathrm{d}}{\mathrm{d}t} \bigg(C_{2}\rho\mathcal{I}_{N}(f)+C_{1}C_{2}\|f\|^{2}_{H^{N}_{x}L^{2}}+\sum_{|\alpha|\leq N}\|\partial^{\alpha} f\|^{2}_{L^{2}_{x}L^{2}_{l_{|\alpha|,0}}}\bigg)
\\ \nonumber  &&+  \frac{\lambda_{0}}{16}\rho \bigg( |(a,b,c)|^{2}_{H^{N}_{x}}
+\|f_{2}\|^{2}_{H^{N}_{x}\mathcal{L}^{s}_{\gamma/2}} +
\sum_{|\alpha|\leq N}\|\partial^{\alpha} f\|^{2}_{L^{2}_{x}\mathcal{L}^{s}_{l_{|\alpha|,0}+\gamma/2}} \bigg) \\&\lesssim& \sum_{|\alpha| \leq N}
|(\pa^{\alpha}g, \pa^{\alpha}f)|+\sum_{|\alpha| \leq N}|(\pa^{\alpha}g, W_{2l_{|\alpha|,0}}\pa^{\alpha}f)|+
\sum_{|\alpha| \leq N-1}|\langle E, \partial^{\alpha}g \rangle|^{2}_{L^{2}_{x}}. \nonumber \een

\noindent {\textit {Step 3: Weighted mixed derivatives.}}
We prove by mathematical induction that for any $0\leq i \leq N$, there exist some constants $K^{i}_{j} \geq 1, -2 \leq j \leq i$, such that
\ben \label{essential-micro-macro-result-4} \frac{\mathrm{d}}{\mathrm{d}t} \bigg(K^{i}_{-2}\rho^{-2i+1}\mathcal{I}_{N}(f)+ K^{i}_{-1}\rho^{-2i}\|f\|^{2}_{H^{N}_{x}L^{2}}+\sum_{j=0}^{i}K^{i}_{j}\rho^{2j-2i}\sum_{|\alpha|\leq N-j, |\beta|=j}\| \partial^{\alpha}_{\beta} f\|^{2}_{L^{2}_{x}L^{2}_{l_{|\alpha|,|\beta|}}} \bigg)\\+ \frac{\lambda_{0}}{16} \rho \bigg(|(a,b,c)|^{2}_{H^{N}_{x}}+\|f_{2}\|^{2}_{H^{N}_{x}\mathcal{L}^{s}_{\gamma/2}}
+ \sum_{j=0}^{i} \sum_{|\alpha|\leq N-j, |\beta|=j}
\|\partial^{\alpha}_{\beta} f\|^{2}_{L^{2}_{x}\mathcal{L}^{s}_{l_{|\alpha|,|\beta|}+\gamma/2}} \bigg) \nonumber\\ \lesssim \sum_{|\alpha| \leq N}
|(\pa^{\alpha}g, \pa^{\alpha}f)|+\sum_{j=0}^{i}\sum_{|\alpha|\leq N-j,|\beta|=j}
|(\partial^{\alpha}_{\beta}g, W_{2l_{|\alpha|,|\beta|}} \partial^{\alpha}_{\beta}f)|+ \sum_{|\alpha| \leq N-1}|\langle E, \partial^{\alpha}g \rangle|^{2}_{L^{2}_{x}}. \nonumber \een
Recalling \eqref{definition-energy-and-dissipation},
our final goal \eqref{essential-micro-macro-result-final} is a result of  \eqref{essential-micro-macro-result-4} by taking $i=N$.

Note that \eqref{essential-micro-macro-result-4} is true when $i=0$. Indeed, by \eqref{essential-micro-macro-result-3} we can take $K^{0}_{-2}=C_{2}, K^{0}_{-1}=C_{1}C_{2}, K^{0}_{0}=1$.

We prove \eqref{essential-micro-macro-result-4} by induction on $i$. Suppose \eqref{essential-micro-macro-result-4} is true when $i = k$ for some $0 \leq k \leq N-1$, we prove it is also valid when $i=k+1$.

Take two indexes $\alpha$ and $\beta$ such that $|\alpha|\leq N-(k+1)$ and $|\beta|= k+1 \geq 1$. Set $q=l_{|\alpha|,|\beta|}$. Applying $\partial^{\alpha}_{\beta}$ to both sides of \eqref{lBE}, we have
\beno 
\partial_{t}\partial^{\alpha}_{\beta}f + v\cdot \nabla_{x} \partial^{\alpha}_{\beta}f + \partial^{\alpha}_{\beta}\mathcal{L}^{\rho}f = [v \cdot \nabla_{x}, \partial^{\alpha}_{\beta}]f +\partial^{\alpha}_{\beta}g. \eeno
Taking inner product with $W_{2l_{|\alpha|,|\beta|}}\partial^{\alpha}_{\beta} f$, one has
\ben \label{mix-x-v-weight}
\f{1}{2}\frac{\mathrm{d}}{\mathrm{d}t}\|\partial^{\alpha}_{\beta}f \|^{2}_{L^{2}_{x}L^{2}_{l_{|\alpha|,|\beta|}}}  + (\partial^{\alpha}_{\beta}\mathcal{L}^{\rho}f,W_{2l_{|\alpha|,|\beta|}}\partial^{\alpha}_{\beta}f) =
 ([v \cdot \nabla_{x}, \partial^{\alpha}_{\beta}]f,W_{2l_{|\alpha|,|\beta|}}\partial^{\alpha}_{\beta}f) +
(\partial^{\alpha}_{\beta}g,W_{2l_{|\alpha|,|\beta|}}\partial^{\alpha}_{\beta}f).  \een
By \eqref{with-v-derivative}, since $l_{|\alpha|,|\beta|} \leq l_{|\alpha|,|\beta|-1} \leq l_{|\alpha|,|\beta_{1}|}$ if $\beta_{1}<\beta$, we get
\beno
&& (\partial^{\alpha}_{\beta} \mathcal{L}^{\rho}f, W_{2l_{|\alpha|,|\beta|}} \partial^{\alpha}_{\beta}f)
\\ &\geq&
\frac{\lambda_{0}}{8} \rho \|\partial^{\alpha}_{\beta}f\|_{L^{2}_{x}\mathcal{L}^{s}_{l_{|\alpha|,|\beta|}+\gamma/2}}^{2} -
C \rho \|\partial^{\alpha}_{\beta}f\|_{L^{2}_{x}L^{2}_{l_{|\alpha|,|\beta|}+\gamma/2}}^{2} -
C \rho \sum_{\beta_{1}<\beta} \|\partial^{\alpha}_{\beta_{1}}f\|_{L^{2}_{x}\mathcal{L}^{s}_{l_{|\alpha|,|\beta_{1}|}+\gamma/2}}^{2}.
\eeno
By \eqref{bounded-by-previous-level-dissipation}, we have
\beno
|([v \cdot \nabla_{x}, \partial^{\alpha}_{\beta}]f, W_{2l_{|\alpha|,|\beta|}} \partial^{\alpha}_{\beta}f)|
\leq \eta \rho \|\partial^{\alpha}_{\beta}f\|_{L^{2}_{x}\mathcal{L}^{s}_{l_{|\alpha|,|\beta|}+\gamma/2}}^{2} + \frac{1}{\eta \rho} \sum_{j=1}^{3} |\beta^{j}|^{2} \| \partial^{\alpha+e^{j}}_{\beta-e^{j}}f\|_{L^{2}_{x}\mathcal{L}^{s}_{l_{|\alpha|+1,|\beta|-1}+\gamma/2}}^{2}.
\eeno
Plugging the previous two inequalities into \eqref{mix-x-v-weight}, by taking $\eta = \frac{\lambda_{0}}{16}$,
taking sum over $|\alpha|\leq N-(k+1)$ and $|\beta|= k+1$,
 we get
\beno 
&& \f{1}{2}\frac{\mathrm{d}}{\mathrm{d}t}\sum_{|\alpha|\leq N-(k+1),|\beta|= k+1}\|\partial^{\alpha}_{\beta}f \|^{2}_{L^{2}_{x}L^{2}_{l_{|\alpha|,|\beta|}}}  + \frac{\lambda_{0}}{16} \rho \sum_{|\alpha|\leq N-(k+1),|\beta|= k+1} \|\partial^{\alpha}_{\beta}f\|_{L^{2}_{x}\mathcal{L}^{s}_{l_{|\alpha|,|\beta|}+\gamma/2}}^{2}
\\&\leq&
C \rho^{-1}
\sum_{j=0}^{k} \sum_{|\alpha|\leq N-j, |\beta|=j}
\|\partial^{\alpha}_{\beta} f\|^{2}_{L^{2}_{x}\mathcal{L}^{s}_{l_{|\alpha|,|\beta|}+\gamma/2}} + C \rho \sum_{|\alpha|\leq N-(k+1),|\beta|= k+1}\|\partial^{\alpha}_{\beta}f \|^{2}_{L^{2}_{x}L^{2}_{l_{|\alpha|,|\beta|}+\gamma/2}}
\\&&+ \sum_{|\alpha|\leq N-(k+1),|\beta|= k+1}(\partial^{\alpha}_{\beta}g,W_{2l_{|\alpha|,|\beta|}}\partial^{\alpha}_{\beta}f).
\eeno
By the interpolation Lemma \ref{interpolation-to-symbol}, we have
\beno \sum_{|\alpha|\leq N-(k+1),|\beta|= k+1}\|\partial^{\alpha}_{\beta}f \|^{2}_{L^{2}_{x}L^{2}_{l_{|\alpha|,|\beta|}+\gamma/2}} &\leq& \eta \sum_{|\alpha|\leq N-(k+1),|\beta|= k+1} \|\partial^{\alpha}_{\beta}f\|_{L^{2}_{x}\mathcal{L}^{s}_{l_{|\alpha|,|\beta|}+\gamma/2}}^{2}
\\&&+ C_{\eta}
\sum_{|\alpha|\leq N-(k+1),|\beta| \leq k} \|\partial^{\alpha}_{\beta}f\|_{L^{2}_{x}\mathcal{L}^{s}_{l_{|\alpha|,|\beta|}+\gamma/2}}^{2}.
\eeno
By taking $\eta$ small, we have
\ben \label{essential-micro-macro-result-pure-2}
&& \frac{\mathrm{d}}{\mathrm{d}t}\sum_{|\alpha|\leq N-(k+1),|\beta|= k+1}\|\partial^{\alpha}_{\beta}f \|^{2}_{L^{2}_{x}L^{2}_{l_{|\alpha|,|\beta|}}}  +  \frac{\lambda_{0}}{16} \rho \sum_{|\alpha|\leq N-(k+1),|\beta|= k+1} \|\partial^{\alpha}_{\beta}f\|_{L^{2}_{x}\mathcal{L}^{s}_{l_{|\alpha|,|\beta|}+\gamma/2}}^{2}
\\ \nonumber &\leq&
C\rho^{-1} \sum_{j=0}^{k} \sum_{|\alpha|\leq N-j, |\beta|=j}\|\partial^{\alpha}_{\beta} f\|^{2}_{L^{2}_{x}\mathcal{L}^{s}_{l_{|\alpha|,|\beta|}+\gamma/2}} + 2\sum_{|\alpha|\leq N-(k+1),|\beta|= k+1}(\partial^{\alpha}_{\beta}g,W_{2l_{|\alpha|,|\beta|}}\partial^{\alpha}_{\beta}f).  \een
By our induction assumption, \eqref{essential-micro-macro-result-4} is true when $i=k$. That is,
\ben \label{essential-micro-macro-result-k} \frac{\mathrm{d}}{\mathrm{d}t}(K^{k}_{-2}\rho^{-2k+1}\mathcal{I}_{N}(f)+ K^{k}_{-1}\rho^{-2k}\|f\|^{2}_{H^{N}_{x}L^{2}}+\sum_{j=0}^{k}K^{k}_{j}\rho^{2j-2k}\sum_{|\alpha|\leq N-j, |\beta|=j}\| \partial^{\alpha}_{\beta} f\|^{2}_{L^{2}_{x}L^{2}_{l_{|\alpha|,|\beta|}}}) \\+   \frac{\lambda_{0}}{16} \rho (|(a,b,c)|^{2}_{H^{N}_{x}}+\|f_{2}\|^{2}_{H^{N}_{x}\mathcal{L}^{s}_{\gamma/2}}
+ \sum_{j=0}^{k} \sum_{|\alpha|\leq N-j, |\beta|=j}\|\partial^{\alpha}_{\beta} f\|^{2}_{L^{2}_{x}\mathcal{L}^{s}_{l_{|\alpha|,|\beta|}+\gamma/2}}) \nonumber\\ \lesssim \sum_{|\alpha| \leq N}
|(\pa^{\alpha}g, \pa^{\alpha}f)|+\sum_{j=0}^{k}\sum_{|\alpha|\leq N-j,|\beta|=j}|(\partial^{\alpha}_{\beta}g, W_{2l_{|\alpha|,|\beta|}} \partial^{\alpha}_{\beta}f)|+ \sum_{|\alpha| \leq N-1}|\langle E, \partial^{\alpha}g \rangle|^{2}_{L^{2}_{x}}. \nonumber \een
Let $M_{k}$ be large enough such that
$
 \frac{\lambda_{0}}{16} M_{k} \geq 2 C$ and  $M_{k} \geq 2$,
then $\eqref{essential-micro-macro-result-k} \times M_{k}\rho^{-2} + \eqref{essential-micro-macro-result-pure-2}$ gives
\beno \frac{\mathrm{d}}{\mathrm{d}t}(K^{k+1}_{-2} \rho^{-2k-1}\mathcal{I}_{N}(f)+ K^{k+1}_{-1}\rho^{-2(k+1)}\|f\|^{2}_{H^{N}_{x}L^{2}}+\sum_{j=0}^{k+1}K^{k+1}_{j}\rho^{2j-2(k+1)}\sum_{|\alpha|\leq N-j, |\beta|=j}\| \partial^{\alpha}_{\beta} f\|^{2}_{L^{2}_{x}L^{2}_{l_{|\alpha|,|\beta|}}})\\+   \frac{\lambda_{0}}{16} \rho (|(a,b,c)|^{2}_{H^{N}_{x}}+\|f_{2}\|^{2}_{H^{N}_{x}\mathcal{L}^{s}_{\gamma/2}}
+ \sum_{j=0}^{k+1} \sum_{|\alpha|\leq N-j, |\beta|=j}\|\partial^{\alpha}_{\beta} f\|^{2}_{L^{2}_{x}\mathcal{L}^{s}_{l_{|\alpha|,|\beta|}+\gamma/2}})
\\ \lesssim \sum_{|\alpha| \leq N}
|(\pa^{\alpha}g, \pa^{\alpha}f)|+\sum_{j=0}^{k+1}\sum_{|\alpha|\leq N-j,|\beta|=j}|(\partial^{\alpha}_{\beta}g, W_{2l_{|\alpha|,|\beta|}} \partial^{\alpha}_{\beta}f)|+ \sum_{|\alpha| \leq N-1}|\langle E, \partial^{\alpha}g \rangle|^{2}_{L^{2}_{x}}.  \eeno
where $K^{k+1}_{j}=M_{k} K^{k}_{j}$ for $-2\leq j\leq k$,
$K^{k+1}_{k+1}=1$. Thus \eqref{essential-micro-macro-result-4} is proved when $i=k+1$ and so we finish the proof.
\end{proof}

\subsection{A priori estimate of the quantum Boltzmann equation.}
In this subsection, we derive the following { \it a priori} estimate for solutions to the Cauchy problem \eqref{linearized-quantum-Boltzmann-eq}.
\begin{thm}\label{a-priori-estimate-LBE} Let $N \geq 9$. Let $0< \rho \leq \rho_{0}$. Fix $T>0$.
There exists a constant  $\delta_{2}>0$ which is independent of $\rho$ and $T$, such that if
a solution $f^{\rho}$ to the Cauchy problem \eqref{linearized-quantum-Boltzmann-eq}  satisfies
\beno \sup_{0 \leq t \leq T} \mathcal{E}_{N}(f^{\rho}(t)) \leq \delta_{2} \rho,\eeno
then $f^{\rho}$ verifies
\ben \label{uniform-estimate-propagation} \sup_{0 \leq t \leq T} \mathcal{E}_{N}(f^{\rho}(t)) + \rho \int_{0}^{T}\mathcal{D}_{N}(f^{\rho}(t)) \mathrm{d}t \leq C \rho^{-2N} \mathcal{E}_{N}(f_{0}),\een	
for some universal constant $C$.
\end{thm}
\begin{proof} Observe that $f^{\rho}$ solves \eqref{lBE} with $g=\Gamma_{2}^{\rho}(f^{\rho},f^{\rho}) + \Gamma_{3}^{\rho}(f^{\rho},f^{\rho},f^{\rho})$. Note that \eqref{conversation-mass-momentum-energy} and \eqref{F0-determines-equilibrium} give
\ben \label{conversation-on-linearized-equation}
\int (1, v, |v|^{2}) \mathcal{N}_{\rho} f(t,x,v) \mathrm{d}x \mathrm{d}v = \int (1, v, |v|^{2}) \mathcal{N}_{\rho} f_{0}(x,v) \mathrm{d}x \mathrm{d}v = 0.
\een
Thanks to \eqref{conversation-on-linearized-equation} and recalling $\mathcal{N}_{\rho} = \rho^{\f12} N_{\rho}$, $f^{\rho}$  verifies \eqref{a-priori-assumption}.
Then we can apply  Proposition  \ref{essential-estimate-of-micro-macro} and the three terms on the righthand of \eqref{essential-micro-macro-result-final} are
\beno \mathcal{I}_{1} &\colonequals& \sum_{|\alpha| \leq N}
|(\pa^{\alpha}\Gamma_{2}^{\rho}(f^{\rho},f^{\rho}) + \pa^{\alpha} \Gamma_{3}^{\rho}(f^{\rho},f^{\rho},f^{\rho}), \pa^{\alpha}f^{\rho})|,\\
\mathcal{I}_{2} &\colonequals& \sum_{|\alpha|+|\beta| \leq N}
|(\partial^{\alpha}_{\beta}\Gamma_{2}^{\rho}(f^{\rho},f^{\rho}) + \partial^{\alpha}_{\beta} \Gamma_{3}^{\rho}(f^{\rho},f^{\rho},f^{\rho}), W_{2l_{|\alpha|,|\beta|}}\partial^{\alpha}_{\beta}f^{\rho})|,
\\\mathcal{I}_{3} &\colonequals& \sum_{|\alpha| \leq N-1}  |\langle  \pa^{\alpha}\Gamma_{2}^{\rho}(f^{\rho},f^{\rho}) + \pa^{\alpha} \Gamma_{3}^{\rho}(f^{\rho},f^{\rho},f^{\rho}), E \rangle|_{L^{2}_{x}}^{2} .\eeno
It is not necessary to estimate $\mathcal{I}_{1}$, since the upper bound of $\mathcal{I}_{2}$ controls $\mathcal{I}_{1}$ naturally.
By \eqref{Gamma-2-energy-estimate} and Theorem \ref{Gamma-3-energy-estimate}, we have
\beno \mathcal{I}_{2} \lesssim \rho^{\f{1}{2}}\mathcal{E}^{\f{1}{2}}_{N}(f^{\rho})\mathcal{D}_{N}(f^{\rho}) + \rho \mathcal{E}_{N}(f^{\rho})\mathcal{D}_{N}(f^{\rho}). \eeno
In view of the proof of \eqref{Gamma-2-energy-estimate} and Theorem \ref{Gamma-3-energy-estimate}, it is much easier to check
\beno \mathcal{I}_{3} \lesssim \rho \mathcal{E}_{N}(f^{\rho})\mathcal{D}_{N}(f^{\rho}) + \rho^{2} \mathcal{E}^{2}_{N}(f^{\rho})\mathcal{D}_{N}(f^{\rho}).\eeno
Therefore by \eqref{essential-micro-macro-result-final}, we get
\beno 
&& \frac{\mathrm{d}}{\mathrm{d}t}\Xi^{\rho}_{N}(f^{\rho}) +  \frac{\lambda_{0}}{16} \rho \mathcal{D}_{N}(f^{\rho}) \leq C \big( \rho^{\f{1}{2}}\mathcal{E}^{\f{1}{2}}_{N}(f^{\rho}) + \rho\mathcal{E}_{N}(f^{\rho}) + \rho^{2} \mathcal{E}^{2}_{N}(f^{\rho}) \big)\mathcal{D}_{N}(f^{\rho}).
\eeno
We take $\delta_{2}$ small enough such that
$
\rho_{0} \delta_{2}^{\f{1}{2}}  \leq 1, C \delta_{2}^{\f{1}{2}}  \leq \frac{\lambda_{0}}{96}.
$
If $\rho \leq \rho_{0}, \sup_{0 \leq t \leq T} \mathcal{E}_{N}(f^{\rho}(t)) \leq \delta_{2} \rho$, then
\beno
C(\rho^{\f{1}{2}}\mathcal{E}^{\f{1}{2}}_{N}(f^{\rho}) + \rho\mathcal{E}_{N}(f^{\rho}) + \rho^{2} \mathcal{E}^{2}_{N}(f^{\rho}) ) \leq C(\delta_{2}^{\f{1}{2}}\rho + \delta_{2}\rho^{2} + \delta_{2}^{2}\rho^{4} ) = C\delta_{2}^{\f{1}{2}} \rho(1 + \delta_{2}^{\f{1}{2}} \rho + \delta_{2}^{\frac{3}{2}}\rho^{3}) \leq \frac{\lambda_{0}}{32}\rho
\eeno
and thus
\ben \label{a-priori-estimate-1} \frac{\mathrm{d}}{\mathrm{d}t}\Xi^{\rho}_{N}(f^{\rho})+ \frac{\lambda_{0}}{32}\rho \mathcal{D}_{N}(f^{\rho}) \leq 0. \een
Integrating \eqref{a-priori-estimate-1} w.r.t. time, we finish the proof by recalling \eqref{equivalence-between-energy}.
\end{proof}

\subsection{Global well-posedness.} \label{global-proof}
With Theorem \ref{local-well-posedness-LBE}(local well-posedness) and Theorem \ref{a-priori-estimate-LBE}({\it a priori} estimate) in hand, we are ready to prove Theorem \ref{global-well-posedness}
for global well-posedness.

\begin{proof}[Proof of Theorem \ref{global-well-posedness}]  Recall the constants $\rho_{2}, \delta_{1}$ in Theorem \ref{local-well-posedness-LBE}. Recall the constant $\delta_{2}$ in Theorem \ref{a-priori-estimate-LBE}.
Take $\rho_{*} = \rho_{2}$. Note that $\rho_{2} \leq  \rho_{1} \leq \rho_{0}$ and so both Theorems are valid if $0< \rho \leq \rho_{*}$.
Denote still by $C$ the larger one of the two universal constants in \eqref{uniform-estimate-of-solution} and  \eqref{uniform-estimate-propagation}. Take
\beno 
\delta_{*} = \min \{\frac{\delta_{1}}{C}, \frac{\delta_{2}}{C^{2}}\}.
\eeno
Now we assume $\mathcal{E}_{N}(f_{0}) \leq \delta_{*} \rho^{2N+1}$ and set to establish global existence and the estimate \eqref{uniform-estimate-global}. Note that $\rho \leq 1, C \geq 1$.
First since $\mathcal{E}_{N}(f_{0}) \leq \delta_{*} \rho^{2N+1}$ and $\delta_{*} \rho^{2N} \leq \delta_{1}$, we can apply Theorem \ref{local-well-posedness-LBE}(taking $\delta = \delta_{*} \rho^{2N}$) to conclude that
the Cauchy problem \eqref{linearized-quantum-Boltzmann-eq} admits a unique solution $f^{\rho} \in L^{\infty}([0,T^{*}]; \mathcal{E}_{N})$ verifying
\beno 
\sup_{0 \leq t \leq T^{*}} \mathcal{E}_{N}(f^{\rho}(t)) \leq C \delta_{*} \rho^{2N+1} \leq \delta_{2}\rho.
\eeno	
Then by Theorem \ref{a-priori-estimate-LBE}(taking $T=T^{*}$), the solution verifies
\ben \label{priori-period-1} \sup_{0 \leq t \leq T^{*}} \mathcal{E}_{N}(f^{\rho}(t)) + \rho \int_{0}^{T^{*}}\mathcal{D}_{N}(f^{\rho}(t)) \mathrm{d}t \leq C \rho^{-2N} \mathcal{E}_{N}(f_{0}).\een
Now we go to establish the following result. For any $n \geq 1$, the Cauchy problem \eqref{linearized-quantum-Boltzmann-eq} admits a solution $f^{\rho} \in L^{\infty}([0, n T^{*}]; \mathcal{E}_{N})$ verifying
\ben \label{priori-period-n} \sup_{0 \leq t \leq n T^{*}} \mathcal{E}_{N}(f^{\rho}(t)) + \rho \int_{0}^{n T^{*}}\mathcal{D}_{N}(f^{\rho}(t)) \mathrm{d}t \leq C \rho^{-2N} \mathcal{E}_{N}(f_{0}).\een
We will prove \eqref{priori-period-n} by mathematical induction on $n$. First, $n=1$ is given by \eqref{priori-period-1}. Suppose \eqref{priori-period-n} is valid for $n=k \geq 1$, we now prove it is also valid for $n=k + 1$. By the assumption, the Cauchy problem \eqref{linearized-quantum-Boltzmann-eq} admits a solution $f^{\rho} \in L^{\infty}([0, k T^{*}]; \mathcal{E}_{N})$ verifying
\ben \label{priori-period-k} \sup_{0 \leq t \leq k T^{*}} \mathcal{E}_{N}(f^{\rho}(t)) + \rho \int_{0}^{k T^{*}}\mathcal{D}_{N}(f^{\rho}(t)) \mathrm{d}t \leq C \rho^{-2N} \mathcal{E}_{N}(f_{0})
\leq C \delta_{*} \rho.
\een
In particular, $\mathcal{E}_{N}(f^{\rho}(k T^{*})) \leq C \delta_{*} \rho$ and $C \delta_{*} \leq \delta_{1}$, we can apply Theorem \ref{local-well-posedness-LBE}(taking $\delta = C \delta_{*}$) to conclude that
the Cauchy problem \eqref{linearized-quantum-Boltzmann-eq} admits a unique solution $f^{\rho} \in L^{\infty}([kT^{*}, (k+1)T^{*}]; \mathcal{E}_{N})$ verifying
\ben \label{local-period-k-plus-1}
\sup_{kT^{*} \leq t \leq (k+1)T^{*}} \mathcal{E}_{N}(f^{\rho}(t)) \leq C^{2} \delta_{*} \rho.
\een	
By \eqref{priori-period-k} and \eqref{local-period-k-plus-1}, the solution verifying
\beno 
\sup_{0 \leq t \leq (k+1)T^{*}} \mathcal{E}_{N}(f^{\rho}(t)) \leq C^{2} \delta_{*} \rho \leq \delta_{2} \rho.
\eeno
Then by Theorem \ref{a-priori-estimate-LBE}(taking $T=(k+1)T^{*}$), the solution verifies
\beno 
\sup_{0 \leq t \leq (k+1)T^{*}} \mathcal{E}_{N}(f^{\rho}(t)) + \rho \int_{0}^{(k+1)T^{*}}\mathcal{D}_{N}(f^{\rho}(t)) \mathrm{d}t \leq C \rho^{-2N} \mathcal{E}_{N}(f_{0}).
\eeno
That is, \eqref{priori-period-n} is valid for $n=k+1$. Sending $n \rightarrow \infty$, we get a global solution $f^{\rho} \in L^{\infty}([0,\infty); \mathcal{E}_{N})$ satisfying \eqref{uniform-estimate-global}.

Note that positivity $\frac{\rho \mu}{1- \rho \mu} + \frac{\rho^{\f{1}{2}}\mu^{\f{1}{2}}}{1- \rho \mu}f^{\rho}(t) \geq 0$ follows from Theorem \ref{local-well-posedness-LBE}.

In the last, we prove uniqueness. Uniqueness can be easily proved by using the arguments in the proof of Theorem \ref{local-well-posedness-LBE}. Indeed, suppose $f,g$ are two solutions to the Cauchy problem \eqref{linearized-quantum-Boltzmann-eq} satisfying \eqref{uniform-estimate-global}. Let $h=f-g$, then $h$ solves
\beno
\partial _t h +  v \cdot \nabla_{x} h + \mathcal{L}^{\rho}h = \Gamma_{2}^{\rho}(f,h) + \Gamma_{2}^{\rho}(h,g) + \Gamma_{3}^{\rho}(f,f,h) + \Gamma_{3}^{\rho}(f,h,g) + \Gamma_{3}^{\rho}(h,g,g),
\eeno
with the initial condition $h(0,x,v) \equiv 0$.
By energy estimates like in the proof of Theorem \ref{local-well-posedness-LBE}, one will get
\beno
 \frac{\mathrm{d}}{\mathrm{d}t} \tilde{\mathcal{E}}_{N}(h)   + \frac{\lambda_{0}}{32}\rho \mathcal{D}_{N}(h)
\leq C
(1+ \mathcal{D}_{N}(f) + \mathcal{D}_{N}(g)) \tilde{\mathcal{E}}_{N}(h).
\eeno
Since $f,g$ are two solutions satisfying \eqref{uniform-estimate-global},
the quantity $1+ \mathcal{D}_{N}(f) + \mathcal{D}_{N}(g)$ is integrable over $[0,T]$ for any $T>0$.
Then by Gr\"{o}nwall's inequality
and using the initial condition $h(0,x,v) \equiv 0$,  one has $h(t) = 0$ in $\mathcal{E}_{N}$ for any $t \geq 0$ and so $f(t)=g(t)$.
Now the proof is complete.
\end{proof}

\section{Appendix} \label{appendix}

Recalling \eqref{definition-of-di} and \eqref{definition-of-e-pm-rho} for the definition of $\{ d^{\rho}_{i} \}_{1 \leq i \leq 5}$ and $\{ e^{\rho}_{i} \}_{1 \leq i \leq 5}$. Recalling \eqref{definition-of-di}, one has
\beno
\{ d^{0}_{i} \}_{1 \leq i \leq 5} = \{ \mu^{\f{1}{2}}, \mu^{\f{1}{2}} v_{1}, \mu^{\f{1}{2}} v_{2}, \mu^{\f{1}{2}} v_{3}, \mu^{\f{1}{2}}(|v|^{2} - 3) \}.
\eeno
For $\alpha>0$, we recall
\ben \label{integral-moment-0}
\int  e^{- \alpha^{2}|v|^{2}} \mathrm{d}v = \alpha^{-3} \pi^{\frac{3}{2}}, \quad
\int  e^{- \alpha^{2}|v|^{2}} |v|^{2} \mathrm{d}v = \alpha^{-5} \frac{3}{2} \pi^{\frac{3}{2}}, \quad
\int  e^{- \alpha^{2}|v|^{2}} |v|^{4} \mathrm{d}v = \alpha^{-7} \frac{15}{4} \pi^{\frac{3}{2}}.
\een
From which, for $i=2,3,4,$
 it is easy to see
\ben \label{d-0-1-l2}
|d^{0}_{1}|_{L^{2}}^{2} = \int  \mu \mathrm{d}v = 1, \quad
|d^{0}_{i}|_{L^{2}}^{2} = \frac{1}{3} \int  \mu |v|^{2} \mathrm{d}v = 1, \quad
|d^{0}_{5}|_{L^{2}}^{2} = \int  \mu(|v|^{2} - 3)^{2} \mathrm{d}v =  6.
\een
Recalling  \eqref{definition-of-e-pm-rho},
the orthonormal basis $\{ e^{0}_{i} \}_{1 \leq i \leq 5} $ for $\mathrm{Null}^{0}$ is
\beno
\{ e^{0}_{i} \}_{1 \leq i \leq 5} = \{ \mu^{\f{1}{2}}, \mu^{\f{1}{2}} v_{1}, \mu^{\f{1}{2}} v_{2}, \mu^{\f{1}{2}} v_{3}, \frac{1}{\sqrt{6}}\mu^{\f{1}{2}}(|v|^{2} - 3)\}.
\eeno

Intuitively, when $\rho$ is small, we expect the orthogonal basis $\{ d^{\rho}_{i} \}_{1 \leq i \leq 5}$ is close to $\{ d^{0}_{i} \}_{1 \leq i \leq 5}$, so is $\{ e^{\rho}_{i} \}_{1 \leq i \leq 5}$ to $\{ e^{0}_{i} \}_{1 \leq i \leq 5}$. The following lemma reflects this expectation mathematically.

\begin{lem} \label{basis-close}
Let $0 \leq \rho \leq \frac{1}{160}$, the following estimates are valid.
\ben \label{order-0-near-each-other}
\f{1}{2} \leq |d^{\rho}_{1}|_{L^{2}} \leq  2, \quad  ||d^{\rho}_{1}|_{L^{2}} - 1| \leq \frac{4}{3} \rho, \quad |e^{\rho}_{1} - e^{0}_{1}|_{L^{2}} \leq \frac{5}{3} \rho.
\\ \label{order-1-near-each-other}
\f{1}{2}  \leq |d^{\rho}_{i}|_{L^{2}} \leq  2, \quad  ||d^{\rho}_{i}|_{L^{2}} - 1| \leq \frac{4}{3} \rho, \quad |e^{\rho}_{i} - e^{0}_{i}|_{L^{2}} \leq  \frac{5}{3} \rho, \quad 2 \leq i \leq 4  .
\\ \label{order-2-near-each-other}
\sqrt{3} \leq |d^{\rho}_{5}|_{L^{2}} \leq  3, \quad ||d^{\rho}_{5}|_{L^{2}} - \sqrt{6}| \leq  160 \rho, \quad |e^{\rho}_{5} - e^{0}_{5}|_{L^{2}} \leq \frac{160}{\sqrt{6}} \rho.
\een
Recall \eqref{definition-of-m-i} and \eqref{defintion-of-l-i}. As a byproduct, we have
\ben \label{boundedness-of-li}
0 \leq l_{\rho,i} \lesssim 1, \quad  i=1,2,3,4.
\een
\end{lem}
\begin{proof} If $0 \leq \rho \leq \f{1}{2}(2 \pi)^{\frac{3}{2}}$, one has
\ben \label{tau-small-control}
\f{1}{2} \leq 1  - \rho \mu \leq 1 .\een
Recall \eqref{equilibrium-rho-not-vanish} and \eqref{definition-of-di}, we define
\beno
h(\rho) \colonequals    |d^{\rho}_{1}|_{L^{2}}^{2} =  \int  \frac{\mu}{(1  - \rho \mu)^{2}} \mathrm{d}v.
\eeno
Taking derivative w.r.t. $\rho$, using \eqref{tau-small-control} and $|\mu|_{L^{\infty}} \leq \frac{1}{8}$, recalling \eqref{d-0-1-l2},
we get
\ben \label{h-derivative-bound}
|h^{\prime}(\rho)| =  2 \int   \frac{\mu^{2}}{(1  - \rho \mu)^{3}} \mathrm{d}v \leq  2 \int  \mu \mathrm{d}v = 2.
\een
Note that $h(0)=1$ by \eqref{d-0-1-l2}. If $\rho \leq \frac{1}{4}$, by mean value theorem,
\ben \label{order-0-norm-squre-diff}
|h(\rho) - h(0)| \leq 2 \rho \leq \f{1}{2} \quad  \Rightarrow \quad \f{1}{2}  \leq h(\rho) \leq \frac{3}{2} \quad \Rightarrow \quad
\f{1}{2} \leq |d^{\rho}_{1}|_{L^{2}} \leq  2.
\een
From which we get
\ben \label{order-0-norm-difference}
| |d^{\rho}_{1}|_{L^{2}} - 1 | = \frac{||d^{\rho}_{1}|_{L^{2}}^{2} - 1|}{|d^{\rho}_{1}|_{L^{2}} + 1} = \frac{|h(\rho) - h(0)|}{|d^{\rho}_{1}|_{L^{2}} + 1} \leq \frac{4}{3}\rho.
\een
Recalling \eqref{definition-of-di} and \eqref{definition-of-e-pm-rho} for the definition of $d^{\rho}_{1}$ and $e^{\rho}_{1}$, we have
\beno
e^{\rho}_{1} - e^{0}_{1} = \frac{d^{\rho}_{1}}{|d^{\rho}_{1}|_{L^{2}}} - \mu^{\f{1}{2}} =  (1-|d^{\rho}_{1}|_{L^{2}})  \frac{d^{\rho}_{1}}{|d^{\rho}_{1}|_{L^{2}}} + \frac{\rho \mu^{\frac{3}{2}}}{1 - \rho \mu}.
\eeno
Using \eqref{tau-small-control}, since $|\mu|_{L^{\infty}} \leq \f18$,
it is easy to see $|\frac{ \mu^{\frac{3}{2}}}{1 - \rho \mu}|_{L^{2}} \leq |\frac{ \mu}{1 - \rho \mu}|_{L^{\infty}}|\mu^{\f{1}{2}} |_{L^{2}} \leq \f14 \leq \frac{1}{3}$. From which together with \eqref{order-0-norm-difference}, we get
\ben \label{order-0-norm-difference-basis}
|e^{\rho}_{1} - e^{0}_{1}|_{L^{2}} \leq  |1-|d^{\rho}_{1}|_{L^{2}}| + \rho |\frac{ \mu^{\frac{3}{2}}}{1 - \rho \mu}|_{L^{2}} \leq \frac{4}{3} \rho + \frac{1}{3} \rho = \frac{5}{3} \rho.
\een
Putting together \eqref{order-0-norm-squre-diff}, \eqref{order-0-norm-difference} and \eqref{order-0-norm-difference-basis}, we get \eqref{order-0-near-each-other}.

Fix $2 \leq i \leq 4$,
let $g(\rho)\colonequals    |d^{\rho}_{i}|_{L^{2}}^{2}$ for $\rho \geq 0$. Note that $g(\rho)$ is independent of $2 \leq i \leq 4$ and
$g(0)=1$ by \eqref{d-0-1-l2}.
Recalling the definition of $d^{\rho}_{i}$ in \eqref{definition-of-di}, there holds
\beno
g(\rho) =   \int  \frac{\mu v_{1}^{2}}{(1  - \rho \mu)^{2}} \mathrm{d}v = \frac{1}{3}  \int  \frac{\mu |v|^{2}}{(1  - \rho \mu)^{2}} \mathrm{d}v.
\eeno
Taking derivative w.r.t. $\rho$, using \eqref{tau-small-control} and $|\mu|_{L^{\infty}} \leq \f18$, we get
\ben \label{g-derivative-bound}
|g^{\prime}(\rho)| =  \frac{2}{3}\int   \frac{\mu^{2}|v|^{2}}{(1  - \rho \mu)^{3}} \mathrm{d}v \leq \frac{2}{3} \int  \mu|v|^{2} \mathrm{d}v  = 2.
\een
Recalling  $g(0) = 1$. When $\rho \leq 1/4$, by mean value theorem, we have
\ben \label{order-1-norm-squre-diff}
|g(\rho) - g(0)| \leq 2 \rho  \leq \f{1}{2} \quad  \Rightarrow \quad \f{1}{2}  \leq g(\rho) \leq \frac{3}{2} \quad \Rightarrow \quad
\f{1}{2} \leq |d^{\rho}_{i}|_{L^{2}} \leq  2.
\een
From which we get
\ben \label{order-1-norm-difference}
||d^{\rho}_{i}|_{L^{2}} - 1|  = \frac{||d^{\rho}_{i}|_{L^{2}}^{2} - 1|}{|d^{\rho}_{i}|_{L^{2}} + 1} = \frac{|g(\rho) - g(0)|}{|d^{\rho}_{i}|_{L^{2}} + 1} \leq \frac{4}{3}\rho.
\een
Recalling \eqref{definition-of-di} and \eqref{definition-of-e-pm-rho} for the definition of $e^{\rho}_{i}$ and $d^{\rho}_{i}$, for $2 \leq i \leq 4$,  we have
\beno
e^{\rho}_{i} - e^{0}_{i} = \frac{d^{\rho}_{i}}{|d^{\rho}_{i}|_{L^{2}}} - \mu^{\f{1}{2}} =  (1-|d^{\rho}_{i}|_{L^{2}})  \frac{d^{\rho}_{i}}{|d^{\rho}_{i}|_{L^{2}}}  + \frac{\rho \mu^{\frac{3}{2}} v_{i-1}}{1 - \rho \mu}.
\eeno
It is easy to check $|\frac{ \mu^{\frac{3}{2}}v_{i-1}}{1 - \rho \mu}|_{L^{2}} \leq |\frac{ \mu}{1 - \rho \mu}|_{L^{\infty}}|\mu^{\f{1}{2}}v_{i-1} |_{L^{2}} \leq \frac{1}{3}$. From which together with \eqref{order-1-norm-difference}, we have
\ben \label{order-1-norm-difference-basis}
|e^{\rho}_{i} - e^{0}_{i}|_{L^{2}} \leq  |1-|d^{\rho}_{i}|_{L^{2}}| + \rho |\frac{ \mu^{\frac{3}{2}} v_{i-1}}{1 - \rho \mu}|_{L^{2}} \leq \frac{4}{3} \rho + \frac{1}{3} \rho = \frac{5}{3}   \rho.
\een
Putting together \eqref{order-1-norm-squre-diff}, \eqref{order-1-norm-difference} and \eqref{order-1-norm-difference-basis}, we get \eqref{order-1-near-each-other}.

Let $\varphi(\rho) \colonequals     C^{\rho}_{5,1} = \langle N_{\rho}|v|^{2} , N_{\rho}\rangle|N_{\rho}|^{-2}_{L^{2}} = 3g(\rho) h^{-1}(\rho)$. Recalling \eqref{h-derivative-bound}, \eqref{order-0-norm-squre-diff}, \eqref{g-derivative-bound} and \eqref{order-1-norm-squre-diff}, we have
\ben \label{estimates-g-h}
\f{1}{2} \leq h(\rho) \leq \frac{3}{2}, \quad |h^{\prime}(\rho)| \leq 2, \quad \f{1}{2}  \leq g(\rho) \leq \frac{3}{2}, \quad |g^{\prime}(\rho)| \leq 2.\een
Taking derivative w.r.t. $\rho$ and using \eqref{estimates-g-h},  we get
\beno
|\varphi^{\prime}(\rho)| =|3 h^{-1}(\rho) g^{\prime}(\rho) - 3 g(\rho)  h^{-2}(\rho) h^{\prime}(\rho)| \leq 60.
\eeno
Note that $\varphi(0) = 3$ by recalling $g(0) = h(0) = 1$.
When $|\rho| \leq \frac{1}{40}$, by mean value theorem, we have
\ben \label{C-tau-5-1}
|\varphi(\rho) - \varphi(0)| \leq 60 \rho \leq \frac{3}{2}  \quad  \Rightarrow \quad
\frac{3}{2} \leq \varphi(\rho) = C^{\rho}_{5,1} \leq  \frac{9}{2}.
\een
Next, recalling \eqref{definition-of-di} and using $|d^{\rho}_{1}|^{2}_{L^{2}} C^{\rho}_{5,1} = \langle d^{\rho}_{1}|v|^{2} , d^{\rho}_{1}\rangle$, we have
\beno
f(\rho) \colonequals    |d^{\rho}_{5}|^{2}_{L^{2}} &=& \int (d^{\rho}_{1}|v|^{2} - C^{\rho}_{5,1} d^{\rho}_{1})^{2} \mathrm{d}v
\\ &=&  \int  \left( (d^{\rho}_{1})^{2} |v|^{4} - (d^{\rho}_{1})^{2} (C^{\rho}_{5,1})^{2} \right) \mathrm{d}v =  \int  \frac{\mu |v|^{4}}{(1  - \rho \mu)^{2}}  \mathrm{d}v -  \frac{9 g^{2}(\rho)}{h(\rho)}.
\eeno
Taking derivative w.r.t. $\rho$, using \eqref{tau-small-control} and $|\mu|_{L^{\infty}} \leq \f18$ to get
$\int   \frac{\mu^{2}|v|^{4}}{(1  - \rho \mu)^{3}} \mathrm{d}v \leq \int  \mu |v|^{4} \mathrm{d}v = 15,$
recalling \eqref{estimates-g-h},
 we get
\beno
|f^{\prime}(\rho)| = |2\int   \frac{\mu^{2}|v|^{4}}{(1  - \rho \mu)^{3}} \mathrm{d}v
- 18 g(\rho) g^{\prime}(\rho) h^{-1}(\rho) + 9 g^{2}(\rho) h^{-2}(\rho) h^{\prime}(\rho)| \leq 480.
\eeno
Note that $f(0) = 6$ by \eqref{d-0-1-l2}.
If $\rho \leq \frac{1}{160}$, by mean value theorem, we have
\ben \label{d-tau-5}
|f(\rho) - f(0)| \leq 480 \rho \leq 3 \quad  \Rightarrow \quad 3 \leq f(\rho) \leq 9
\quad  \Rightarrow \quad \sqrt{3} \leq |d^{\rho}_{5}| \leq 3.
\een
From which we get
\ben \label{order-2-norm-difference}
||d^{\rho}_{5}|_{L^{2}} - \sqrt{6}|  = \frac{||d^{\rho}_{5}|_{L^{2}}^{2} - 6|}{|d^{\rho}_{5}|_{L^{2}} + \sqrt{6}} = \frac{|f(\rho) - f(0)|}{|d^{\rho}_{5}|_{L^{2}} + \sqrt{6}} \leq 160\rho.
\een
By \eqref{d-tau-5} and \eqref{order-2-norm-difference}, the first two inequalities in \eqref{order-2-near-each-other} are proved.
Recalling \eqref{definition-of-di} and \eqref{definition-of-e-pm-rho} for the definition of $e^{\rho}_{5}$ and $d^{\rho}_{5}$,  we have
\beno
e^{\rho}_{5} - e^{0}_{5} &=& \frac{d^{\rho}_{5}}{|d^{\rho}_{5}|_{L^{2}}} - \frac{d^{0}_{5}}{|d^{0}_{5}|_{L^{2}}} =
\frac{|d^{0}_{5}|_{L^{2}} d^{\rho}_{5} - |d^{\rho}_{5}|_{L^{2}} d^{0}_{5} }{|d^{0}_{5}|_{L^{2}}|d^{\rho}_{5}|_{L^{2}}}
\\&=& \frac{(|d^{0}_{5}|_{L^{2}} - |d^{\rho}_{5}|_{L^{2}}) d^{\rho}_{5}  + |d^{\rho}_{5}|_{L^{2}}(d^{\rho}_{5}-d^{0}_{5})}{|d^{0}_{5}|_{L^{2}}|d^{\rho}_{5}|_{L^{2}}},
\eeno
and thus
\ben \label{e-L2-to-d-L2}
|e^{\rho}_{5} - e^{0}_{5}|_{L^{2}} \leq \frac{||d^{0}_{5}|_{L^{2}} - |d^{\rho}_{5}|_{L^{2}}|+|d^{\rho}_{5}-d^{0}_{5}|_{L^{2}}}{|d^{0}_{5}|_{L^{2}}} \leq \frac{2|d^{\rho}_{5}-d^{0}_{5}|_{L^{2}}}{|d^{0}_{5}|_{L^{2}}}.
\een
Now it remains to consider $|d^{\rho}_{5}-d^{0}_{5}|_{L^{2}}$. Recall
\beno
d^{\rho}_{5}  = \frac{\mu^{\f{1}{2}}|v|^{2}}{1 - \rho \mu} - C^{\rho}_{5,1} \frac{\mu^{\f{1}{2}}}{1 - \rho \mu}, \quad d^{0}_{5}  = \mu^{\f{1}{2}}|v|^{2} - 3 \mu^{\f{1}{2}}.
\eeno
Then
\beno
d^{\rho}_{5}-d^{0}_{5} =  \mu \rho \frac{\mu^{\f{1}{2}}|v|^{2}}{1 - \rho \mu} + (3-C^{\rho}_{5,1})\mu^{\f{1}{2}} -  \mu \rho C^{\rho}_{5,1} \frac{\mu^{\f{1}{2}}}{1 - \rho \mu}.
\eeno
Recalling \eqref{tau-small-control} and $|\mu|_{L^{\infty}} \leq \frac{1}{8}$, using $|\mu^{\f{1}{2}}|v|^{2}|_{L^{2}} = \sqrt{15}, |\mu^{\f{1}{2}}|_{L^{2}}=1$ and the results in \eqref{C-tau-5-1},
we have
\beno
|d^{\rho}_{5}-d^{0}_{5}|_{L^{2}} &\leq&   \frac{1}{4} \rho |\mu^{\f{1}{2}}|v|^{2}|_{L^{2}} + |3-C^{\rho}_{5,1}||\mu^{\f{1}{2}}|_{L^{2}} +  \frac{1}{4} \rho |C^{\rho}_{5,1}| |\mu^{\f{1}{2}}|_{L^{2}}
\\ &\leq&   \frac{1}{4} \rho \times \sqrt{15} + 60 \rho +  \frac{1}{4} \rho \times \f92 \leq 80 \rho.
\eeno
By recalling \eqref{e-L2-to-d-L2} and $|d^{0}_{5}|_{L^{2}} = \sqrt{6}$, we get the last inequality in \eqref{order-2-near-each-other}.

Recalling \eqref{definition-of-m-i} and \eqref{defintion-of-l-i}, using \eqref{order-0-near-each-other}, \eqref{order-1-near-each-other},  \eqref{order-2-near-each-other} and \eqref{C-tau-5-1},
we get \eqref{boundedness-of-li}.
\end{proof}

We can revise the proof of Lemma \ref{basis-close} to get
\begin{lem}\label{more-regularity-order-1} Let $m \geq 0, l \in \mathbb{R}, 0 \leq \rho \leq \frac{1}{160}$. For $1 \leq i \leq 5$, there holds
\beno
 |\mu^{-\frac{1}{4}}(e^{\rho}_{i} - e^{0}_{i})|_{L^{2}} \leq  \tilde{C} \rho, \quad
 |e^{\rho}_{i} - e^{0}_{i}|_{H^{m}_{l}} \leq  \tilde{C}_{m,l} \rho,
\eeno
for some universal constant $\tilde{C}$ and a constant $\tilde{C}_{m,l}$ depending only on $m,l$.
\end{lem}

We now apply Lemma \ref{more-regularity-order-1} to prove Lemma \ref{projection-close}.
\begin{proof}[Proof of Lemma \ref{projection-close}.] Observe
\beno
\mathbb{P}_{\rho}f - \mathbb{P}_{0}f  =  \sum_{i=1}^{5} \langle f, e^{\rho}_{i}\rangle e^{\rho}_{i} - \sum_{i=1}^{5} \langle f, e^{0}_{i}\rangle e^{0}_{i}  = \sum_{i=1}^{5} \langle f, e^{\rho}_{i}-e^{0}_{i} \rangle e^{\rho}_{i} + \sum_{i=1}^{5} \langle f, e^{0}_{i} \rangle (e^{\rho}_{i} - e^{0}_{i}).
\eeno
By Cauchy-Schwartz inequality, we have
\beno
|\langle f, e^{\rho}_{i} - e^{0}_{i} \rangle| \leq  |\mu^{\frac{1}{4}} f |_{L^{2}}| \mu^{-\frac{1}{4}}(e^{\rho}_{i} - e^{0}_{i})|_{L^{2}}, \quad |\langle f, e^{0}_{i} \rangle| \leq  |\mu^{\frac{1}{4}} f |_{L^{2}}| \mu^{-\frac{1}{4}} e^{0}_{i}|_{L^{2}},
\eeno
which gives
\beno
|\mathbb{P}_{\rho}f - \mathbb{P}_{0}f|_{H^{m}_{l}}  \leq  |\mu^{\frac{1}{4}} f |_{L^{2}}| \sum_{i=1}^{5} | \mu^{-\frac{1}{4}}(e^{\rho}_{i} - e^{0}_{i})|_{L^{2}} |e^{\rho}_{i}|_{H^{m}_{l}} + |\mu^{\frac{1}{4}} f |_{L^{2}}|\sum_{i=1}^{5} | \mu^{-\frac{1}{4}} e^{0}_{i}|_{L^{2}} |e^{\rho}_{i} - e^{0}_{i}|_{H^{m}_{l}}.
\eeno
With the constants $\tilde{C}$ and $\tilde{C}_{m,l}$ in Lemma \ref{more-regularity-order-1}, define $C_{m,l}$ in the following way,
\ben \label{constant-C-s-l}
C_{m,l} \colonequals    (\tilde{C}+\tilde{C}_{m,l}) (\sum_{i=1}^{5} |\mu^{-\frac{1}{4}} e^{0}_{i}|_{L^{2}} + \sup_{0 \leq \rho \leq 160^{-1}} \sum_{i=1}^{5} |e^{\rho}_{i}|_{H^{m}_{l}}). \een
Then by Lemma \ref{more-regularity-order-1},
we get the desired result.
\end{proof}

When $v$ and $v_{*}$ are close, we can exchange negative exponential weight freely. For example, if $|v-v_{*}| \leq 1$, then $|v|^{2} \geq \f{1}{2}|v_{*}|^{2} - 1$ and thus $\mu(v) \lesssim \mu^{\f12}(v_{*})$. More general, we have
\begin{lem} \label{mu-weight-transfer} Recall $v(\kappa) \colonequals   v + \kappa(v^{\prime} - v), v_{*}(\iota) \colonequals   v_{*} + \iota(v^{\prime}_{*} - v_{*})$.
For $\kappa_{1}, \iota_{1}, \kappa_{2}, \iota_{2} \in [0,1]$, if $|v-v_{*}| \leq 1$, then
\beno
\mu (v(\kappa_{1})) \lesssim \mu^{\f12} (v_{*}(\iota_{1})), \quad \mu (v_{*}(\iota_{1})) \lesssim \mu^{\f12} (v(\kappa_{1})), \quad \mu (v(\kappa_{1})) \lesssim \mu^{\f12} (v(\kappa_{2})), \quad \mu (v_{*}(\iota_{1})) \lesssim \mu^{\f12} (v_{*}(\iota_{2})).
\eeno
\end{lem}
Lemma \ref{mu-weight-transfer} is obvious since  $|v(\kappa_{1})-v_{*}(\iota_{1})|, |v(\kappa_{1})-v(\kappa_{2})|, |v_{*}(\iota_{1})-v_{*}(\iota_{2})| \leq |v-v_{*}| \leq 1$.

We estimate various integrals over $\mathbb{S}^{2}$ involving the difference $\mathrm{D}(\mu^{\f{1}{4}}_{*})$ in the following lemma.
\begin{lem} \label{sigma-integral-mu-difference} The following estimates are valid for any $v, v_{*} \in \mathbb{R}^{3}$.
\ben \label{mu-order-1-cancelation}
|\int B \mathrm{D}(\mu^{\f{1}{4}}_{*}) \mathrm{d}\sigma| &\lesssim& \mathrm{1}_{|v-v_{*}| \geq 1} \langle v \rangle^{\gamma+2s} + \mathrm{1}_{|v-v_{*}|\leq 1} |v-v_{*}|^{-2}  \mu^{\frac{1}{32}}\mu^{\frac{1}{32}}_{*}.
\\ \label{mu-order-1-cancelation-with-mu}
|\int B \mathrm{D}(\mu^{\f{1}{4}}_{*}) \mu^{\frac{1}{4}} \mathrm{d}\sigma| &\lesssim& \mathrm{1}_{|v-v_{*}| \geq 1} \mu^{\frac{1}{32}}\mu^{\frac{1}{32}}_{*} + \mathrm{1}_{|v-v_{*}|\leq 1} |v-v_{*}|^{-2}  \mu^{\frac{1}{32}}\mu^{\frac{1}{32}}_{*}.
\\ \label{mu-order-2-cancelation}
\int B \mathrm{D}^{2}(\mu^{\f{1}{4}}_{*}) \mathrm{d}\sigma &\lesssim& \mathrm{1}_{|v-v_{*}| \geq 1} \langle v \rangle^{\gamma+2s} + \mathrm{1}_{|v-v_{*}|\leq 1} |v-v_{*}|^{-1}  \mu^{\frac{1}{32}}\mu^{\frac{1}{32}}_{*}.
\\ \label{sigma-integral-mu-mu-star-difference}
\int B \mathrm{D}^{2}(\mu^{\f{1}{4}}_{*}) \mu^{\frac{1}{4}} \mathrm{d}\sigma &\lesssim&  \mathrm{1}_{|v-v_{*}| \geq 1} \mu^{\frac{1}{32}}\mu^{\frac{1}{32}}_{*} + \mathrm{1}_{|v-v_{*}|\leq 1} |v-v_{*}|^{-2}  \mu^{\frac{1}{32}}\mu^{\frac{1}{32}}_{*}.
\een
\end{lem}
\begin{proof} Note that when $|v-v_{*}|\leq 1$,
all the results contain a $\mu\mu_{*}$ type weight. This is given by Lemma \ref{mu-weight-transfer}.
We will give a detailed proof to  \eqref{mu-order-1-cancelation} and then only sketch the proof of the other three results.

Applying Taylor expansion, we have
\ben \label{Taylor-at-v-star}
-\mathrm{D}(\mu^{\f{1}{4}}_{*})=(\nabla \mu^{\frac{1}{4}})(v_{*}) \cdot\left(v^{\prime}_{*}-v_{*}\right)+ \int_{0}^{1}(1-\iota)(\nabla^{2} \mu^{\frac{1}{4}})(v_{*}(\iota)):\left(v^{\prime}_{*}-v_{*}\right) \otimes\left(v^{\prime}_{*}-v_{*}\right) \mathrm{d}\iota,
\een
where $v_{*}(\iota)\colonequals     v_{*} + \iota(v^{\prime}_{*}-v_{*})$.
We consider two cases: $|v-v_{*}| \leq 1$ and $|v-v_{*}| \geq 1$.

Case 1: $|v-v_{*}| \leq 1$. By \eqref{Taylor-at-v-star},
$
|\int
B \mathrm{D}(\mu^{\f{1}{4}}_{*}) \mathrm{d}\sigma| = |\mathcal{P}_{1} + \mathcal{P}_{2}|
$
where
\beno
\mathcal{P}_{1} \colonequals     \int
B (\nabla \mu^{\frac{1}{4}})(v_{*}) \cdot\left(v^{\prime}_{*}-v_{*}\right) \mathrm{d}\sigma,
\quad
\mathcal{P}_{2} \colonequals     \int
B (1-\iota)(\nabla^{2} \mu^{\frac{1}{4}})(v_{*}(\iota)):\left(v^{\prime}_{*}-v_{*}\right) \otimes\left(v^{\prime}_{*}-v_{*}\right) \mathrm{d}\iota \mathrm{d}\sigma .
\eeno
Noting that $v^{\prime}_{*}-v_{*} = v -v^{\prime}$, using \eqref{theta-squre-out}, the estimate \eqref{cancel-angular-singularity}, the fact
$|\nabla \mu^{\frac{1}{4}}| \lesssim \mu^{\frac{1}{8}}$ and Lemma \ref{mu-weight-transfer},
we get
\beno
|\mathcal{P}_{1}| = |\int
B \sin^{2} \frac{\theta}{2}(\nabla \mu^{\frac{1}{4}})(v_{*}) \cdot (v-v_{*}) \mathrm{d}\sigma|
\lesssim \mu^{\frac{1}{8}}(v_{*})|v-v_{*}|^{\gamma+1} \lesssim   |v-v_{*}|^{-2}  \mu^{\frac{1}{32}}\mu^{\frac{1}{32}}_{*}.
\eeno
By the fact $|\nabla^{2} \mu^{\frac{1}{4}}| \lesssim \mu^{\frac{1}{8}}$, Lemma \ref{mu-weight-transfer} and the estimate \eqref{cancel-angular-singularity}, we have
\beno
|\mathcal{P}_{2}| \leq \int
B \sin^{2} \frac{\theta}{2} |v-v_{*}|^{2} |(\nabla^{2} \mu^{\frac{1}{4}})(v_{*}(\iota))| \mathrm{d}\iota \mathrm{d}\sigma
\lesssim  |v-v_{*}|^{\gamma+2}  \mu^{\frac{1}{32}}\mu^{\frac{1}{32}}_{*} \lesssim |v-v_{*}|^{-1}  \mu^{\frac{1}{32}}\mu^{\frac{1}{32}}_{*}.
\eeno
Patching together the previous two estimates,  since $|v-v_{*}| \leq 1$,
we get
\ben \label{leq-1-case}
|\int
B \mathrm{D}(\mu^{\f{1}{4}}_{*}) \mathrm{d}\sigma| \leq |\mathcal{P}_{1}| + |\mathcal{P}_{2}|  \lesssim   |v-v_{*}|^{-2}  \mu^{\frac{1}{32}}\mu^{\frac{1}{32}}_{*}.
\een

Case 2: $|v-v_{*}| \geq 1$. We divide the integral into two parts according to $\sin\frac{\theta}{2} \leq |v-v_{*}|^{-1}$ and $\sin\frac{\theta}{2} \geq |v-v_{*}|^{-1}$ as $\int B \mathrm{D}(\mu^{\f{1}{4}}_{*}) \mathrm{d}\sigma = \mathcal{Q}_{\leq}+\mathcal{Q}_{\geq},$ where
\beno
\mathcal{Q}_{\leq}\colonequals    \int \mathrm{1}_{\sin\frac{\theta}{2} \leq |v-v_{*}|^{-1}}B \mathrm{D}(\mu^{\f{1}{4}}_{*}) \mathrm{d}\sigma, \quad
 \mathcal{Q}_{\geq}\colonequals     \int \mathrm{1}_{\sin\frac{\theta}{2} \geq |v-v_{*}|^{-1}}B \mathrm{D}(\mu^{\f{1}{4}}_{*}) \mathrm{d}\sigma.
\eeno
It is obvious $|\mathcal{Q}_{\geq}| \leq \mathcal{Q}_{\geq,1} + \mathcal{Q}_{\geq,2}$ where
\beno
\mathcal{Q}_{\geq,1}\colonequals    \int \mathrm{1}_{\sin\frac{\theta}{2} \geq |v-v_{*}|^{-1}}B (\mu^{\frac{1}{4}})^{\prime}_{*}  \mathrm{d}\sigma, \quad
\mathcal{Q}_{\geq,2}\colonequals    \int \mathrm{1}_{\sin\frac{\theta}{2} \geq |v-v_{*}|^{-1}}B  \mu_{*}^{\frac{1}{4}} \mathrm{d}\sigma.
\eeno
Note that
\ben \label{theta-large-critical}
\int \mathrm{1}_{\sin\frac{\theta}{2} \geq |v-v_{*}|^{-1}}  \sin^{-2-2s}\frac{\theta}{2} \mathrm{d}\sigma \lesssim |v-v_{*}|^{2s},
\een
which gives $
\mathcal{Q}_{\geq,2} \lesssim \mu_{*}^{\frac{1}{4}} |v-v_{*}|^{\gamma+2s} \lesssim \langle v \rangle^{\gamma+2s}.
$
The analysis of $\mathcal{Q}_{\geq,1}$ is a little bit technical. We will make use of the nice weight $(\mu^{\frac{1}{4}})^{\prime}_{*}$ on $v^{\prime}_{*}$. Since $|v-v^{\prime}_{*}| \sim |v-v_{*}|$ and $|v-v_{*}| \geq 1$, we have
\beno
|v-v_{*}|^{\gamma} (\mu^{\frac{1}{4}})^{\prime}_{*} \lesssim |v-v^{\prime}_{*}|^{\gamma+2s} (\mu^{\frac{1}{4}})^{\prime}_{*} |v-v_{*}|^{-2s} \lesssim \langle v \rangle^{\gamma+2s} |v-v_{*}|^{-2s}.
\eeno
Plugging which into $\mathcal{Q}_{\geq,1}$ and using \eqref{theta-large-critical}, we get
\beno
\mathcal{Q}_{\geq,1} \lesssim \int \mathrm{1}_{\sin\frac{\theta}{2} \geq |v-v_{*}|^{-1}}  \sin^{-2-2s}\frac{\theta}{2} \langle v \rangle^{\gamma+2s} |v-v_{*}|^{-2s}  \mathrm{d}\sigma  \lesssim \langle v \rangle^{\gamma+2s}.
\eeno
Note that $\mathcal{Q}_{\geq,1}$ and $\mathcal{Q}_{\geq,2}$ share the same upper bound $\langle v \rangle^{\gamma+2s}$ and so
$
|\mathcal{Q}_{\geq}| \lesssim \langle v \rangle^{\gamma+2s}.
$
Plugging \eqref{Taylor-at-v-star} into $\mathcal{Q}_{\leq}$, we get
$
|\mathcal{Q}_{\leq}| = |\mathcal{Q}_{\leq,1} + \mathcal{Q}_{\leq,2}|
$
where
\beno
\mathcal{Q}_{\leq,1} &\colonequals&     \int \mathrm{1}_{\sin\frac{\theta}{2} \leq |v-v_{*}|^{-1}}B (\nabla \mu^{\frac{1}{4}})(v_{*}) \cdot\left(v^{\prime}_{*}-v_{*}\right) \mathrm{d}\sigma,
\\
\mathcal{Q}_{\leq,2} &\colonequals&     \int \mathrm{1}_{\sin\frac{\theta}{2} \leq |v-v_{*}|^{-1}}B (1-\iota)(\nabla^{2} \mu^{\frac{1}{4}})(v_{*}(\iota)):\left(v^{\prime}_{*}-v_{*}\right) \otimes\left(v^{\prime}_{*}-v_{*}\right) \mathrm{d}\iota \mathrm{d}\sigma.
\eeno
Using \eqref{theta-squre-out} and $v^{\prime}_{*}-v_{*} = v -v^{\prime}$, we have
\beno
|\mathcal{Q}_{\leq,1}| = |\int \mathrm{1}_{\sin\frac{\theta}{2} \leq |v-v_{*}|^{-1}}B \sin^{2} \frac{\theta}{2}(\nabla \mu^{\frac{1}{4}})(v_{*}) \cdot (v-v_{*}) \mathrm{d}\sigma| \lesssim \mu^{\frac{1}{8}}(v_{*})|v-v_{*}|^{\gamma+2s-1} \lesssim \langle v \rangle^{\gamma+2s-1},
\eeno
where we use
\ben \label{theta-small-critical}
\int \mathrm{1}_{\sin\frac{\theta}{2} \leq |v-v_{*}|^{-1}} \sin^{-2s}\frac{\theta}{2}   \mathrm{d}\sigma \lesssim  |v-v_{*}|^{2s-2}.
\een
We will make use of the nice weight $(\nabla^{2} \mu^{\frac{1}{4}})(v_{*}(\iota))$ on $v_{*}(\iota)$ in $\mathcal{Q}_{\leq,2} $. Since $|v-v_{*}(\iota)| \sim |v-v_{*}|$ and $|v-v_{*}| \geq 1$, we have
\beno
|v-v_{*}|^{\gamma} |(\nabla^{2} \mu^{\frac{1}{4}})(v_{*}(\iota))| \lesssim |v-v_{*}(\iota)|^{\gamma+2s} \mu^{\frac{1}{8}}(v_{*}(\iota))|v-v_{*}|^{-2s} \lesssim \langle v \rangle^{\gamma+2s} |v-v_{*}|^{-2s}.
\eeno
From which together with \eqref{theta-small-critical}, we get
\beno
|\mathcal{Q}_{\leq,2}| \lesssim \langle v \rangle^{\gamma+2s} \int \mathrm{1}_{\sin\frac{\theta}{2} \leq |v-v_{*}|^{-1}}\sin^{-2s}\frac{\theta}{2}   |v-v_{*}|^{2-2s} \mathrm{d}\sigma \lesssim \langle v \rangle^{\gamma+2s}.
\eeno
Patching together the estimates of $|\mathcal{Q}_{\leq,1}|$ and $|\mathcal{Q}_{\leq,2}|$, we get
$
|\mathcal{Q}_{\leq}| \lesssim \langle v \rangle^{\gamma+2s}.
$
Therefore when $|v-v_{*}| \geq 1$, we get
\ben  \label{geq-1-case}
|\int B \mathrm{D}(\mu^{\f{1}{4}}_{*}) \mathrm{d}\sigma| \leq |\mathcal{Q}_{\leq}| + |\mathcal{Q}_{\geq}| \lesssim \langle v \rangle^{\gamma+2s}.
\een
Patching together  \eqref{leq-1-case} and \eqref{geq-1-case}, we finish the proof of \eqref{mu-order-1-cancelation}.

Thanks to the additional factor $\mu^{\f14}$ on the left-hand side of \eqref{mu-order-1-cancelation-with-mu},
we can use \eqref{mu-weight-result} to retain the $\mu \mu_{*}$ weight for the case $|v-v_{*}| \geq 1$ and get the desired result \eqref{mu-order-1-cancelation-with-mu}.

Now we explain how to prove \eqref{mu-order-2-cancelation}. Note that $\mathrm{D}^{2}(\mu^{\f{1}{4}}_{*})$ contains order 2 cancelation by first order Taylor expansion.
For the case $|v-v_{*}| \leq 1$, applying Taylor expansion to $\mathrm{D}(\mu^{\f{1}{4}}_{*})$ up to first order,
using the fact
$|\nabla \mu^{\frac{1}{4}}| \lesssim \mu^{\frac{1}{8}}$ and Lemma \ref{mu-weight-transfer}, we can get the result. For the case $|v-v_{*}| \geq 1$, like in the proof of \eqref{mu-order-1-cancelation}, we divide the integral into two parts according to $\sin\frac{\theta}{2} \leq |v-v_{*}|^{-1}$ and $\sin\frac{\theta}{2} \geq |v-v_{*}|^{-1}$. For the part with  $\sin\frac{\theta}{2} \geq |v-v_{*}|^{-1}$, we can use the same estimate for $\mathcal{Q}_{\geq}$. For the part with  $\sin\frac{\theta}{2} \leq |v-v_{*}|^{-1}$, applying Taylor expansion to $\mathrm{D}(\mu^{\f{1}{4}}_{*})$ up to first order,
we can use the same estimate for $\mathcal{Q}_{\leq,2}.$

At last, let us see how to prove \eqref{sigma-integral-mu-mu-star-difference}. Thanks to the additional factor $\mu^{\f14}$ on the left-hand side of \eqref{sigma-integral-mu-mu-star-difference},
we can use \eqref{mu-weight-result} to retain the $\mu \mu_{*}$ weight in the proof of \eqref{mu-order-2-cancelation}
for the case $|v-v_{*}| \geq 1$ and get the desired result \eqref{sigma-integral-mu-mu-star-difference}.
\end{proof}

The exponents of $\mu$ in Lemma \ref{sigma-integral-mu-difference} can be relaxed as long as they have a lower bound. Since $\frac{1}{16}$ is enough for our purpose, we give
\begin{lem} \label{sigma-integral-mu-difference-more} Let $\frac{1}{16} \leq a, b \leq 1$, then the following estimates are valid.
\ben \label{mu-order-1-cancelation-a-star}
|\int B \mathrm{D}(\mu^{a}_{*}) \mathrm{d}\sigma| &\lesssim& \mathrm{1}_{|v-v_{*}| \geq 1} \langle v \rangle^{\gamma+2s} + \mathrm{1}_{|v-v_{*}|\leq 1} |v-v_{*}|^{-2}  \mu^{\frac{1}{128}}\mu^{\frac{1}{128}}_{*}.
\\ \label{mu-order-1-cancelation-a-star-with-mu}
|\int B \mathrm{D}(\mu^{a}_{*}) \mu^{b} \mathrm{d}\sigma| &\lesssim&
\mathrm{1}_{|v-v_{*}| \geq 1} \mu^{\frac{1}{128}}\mu^{\frac{1}{128}}_{*} + \mathrm{1}_{|v-v_{*}|\leq 1} |v-v_{*}|^{-2}  \mu^{\frac{1}{128}}\mu^{\frac{1}{128}}_{*}.
\\ \label{mu-order-2-cancelation-a-star}
\int B \mathrm{D}^{2}(\mu^{a}_{*}) \mathrm{d}\sigma &\lesssim& \mathrm{1}_{|v-v_{*}| \geq 1} \langle v \rangle^{\gamma+2s} + \mathrm{1}_{|v-v_{*}|\leq 1} |v-v_{*}|^{-1}  \mu^{\frac{1}{128}}\mu^{\frac{1}{128}}_{*}.
\\ \label{mu-order-2-cancelation-a-star-with-mu}
\int B \mathrm{D}^{2}(\mu^{a}_{*}) \mu^{b} \mathrm{d}\sigma &\lesssim& \mathrm{1}_{|v-v_{*}| \geq 1} \mu^{\frac{1}{128}}\mu^{\frac{1}{128}}_{*} + \mathrm{1}_{|v-v_{*}|\leq 1} |v-v_{*}|^{-1}  \mu^{\frac{1}{128}}\mu^{\frac{1}{128}}_{*}.
\een
\end{lem}
\begin{proof} Since $\mu^{b} \lesssim \mu^{\frac{1}{16}}, \mu^{a} \lesssim \mu^{\frac{1}{16}}, |\nabla \mu^{a}| \lesssim \mu^{\frac{1}{32}}, |\nabla^{2} \mu^{a}| \lesssim \mu^{\frac{1}{32}}$, we can follow the proof of Lemma \ref{sigma-integral-mu-difference} to get the desired results.
\end{proof}

Based on the proofs of Lemma \ref{sigma-integral-mu-difference} and Lemma \ref{sigma-integral-mu-difference-more}, recalling \eqref{M-rho-mu-N}, we give the following remark.
\begin{rmk} \label{mu-to-N-or-M} If we replace $\mu$ with $N^{2}$ or $M$ in Lemma \ref{sigma-integral-mu-difference-more},
all the results are still valid. Let $P_{1}, P_{2}$ be two polynomials on $\mathbb{R}^{3}$. If we replace $\mu^{a}$ and $\mu^{b}$ with  $P_{1}\mu^{a}$ and $P_{2}\mu^{b}$ respectively in Lemma \ref{sigma-integral-mu-difference-more},
all the results are still valid. Moreover, If we replace $\mu^{a}$ and $\mu^{b}$ with  $P_{1}N^{2a}$(or $P_{1}M^{a}$) and $P_{2}N^{2b}$(or $P_{2}M^{b}$) respectively in Lemma \ref{sigma-integral-mu-difference-more},
all the results are still valid.
\end{rmk}

By Hardy-Sobolev-Littlewood inequality, we can derive
\begin{lem}\label{cancel-singularity-by-H-one-half} Let $s_{1},s_{2}\geq 0,s_{1}+s_{2}=\f{1}{2}$. There holds
\ben \label{minus-1-sigularity-g2-h2}
\int |v-v_{*}|^{-1} g_{*}^{2} h^{2} \mathrm{d}v_{*} \mathrm{d}v \lesssim |g|_{H^{s_{1}}}^{2}|h|_{H^{s_{2}}}^{2}.
\een
Let $a_{1},a_{2},a_{3}\geq 0,a_{1}+a_{2}+a_{3}=\f{1}{2}$. There holds
\ben \label{minus-2-sigularity-g-h-f}
\int |v-v_{*}|^{-2} |g_{*} h f| \mathrm{d}v_{*} \mathrm{d}v \lesssim |g|_{H^{a_{1}}} |h|_{H^{a_{2}}} |f|_{H^{a_{3}}}.
\een
\end{lem}

Based on Lemma \ref{sigma-integral-mu-difference-more} and Lemma \ref{cancel-singularity-by-H-one-half},
we now give a proof to Lemma \ref{functional-X-g-h}.
\begin{proof}[Proof of Lemma \ref{functional-X-g-h}] We first prove \eqref{estimate-of-X-g-h-f-total-one-half}.
By \eqref{mu-order-1-cancelation-a-star} and \eqref{minus-2-sigularity-g-h-f},
we have
\beno
&& |\int
 B \mathrm{D}(\mu^{a}_{*}) g_{*} h f \mathrm{d}V|
\\&\lesssim&  \int (\mathrm{1}_{|v-v_{*}| \geq 1} \langle v \rangle^{\gamma+2s} + \mathrm{1}_{|v-v_{*}|\leq 1} |v-v_{*}|^{-2}  \mu^{\frac{1}{128}}\mu^{\frac{1}{128}}_{*})  |g_{*} h f|  \mathrm{d}v_{*} \mathrm{d}v
\\ &\lesssim& |g|_{L^{1}} |h|_{L^{2}_{\gamma/2+s}} |f|_{L^{2}_{\gamma/2+s}} + |\mu^{\frac{1}{256}}g|_{H^{s_{1}}} |\mu^{\frac{1}{256}}h|_{H^{s_{2}}} |\mu^{\frac{1}{256}}f|_{H^{s_{3}}}.
\eeno
By the same argument as the above, using the imbedding $H^{2} \hookrightarrow L^{\infty}$(for $g$ or $W_{\gamma/2+s}h$ or $\varrho$), we can get \eqref{estimate-of-A-g-h-varrho-f-total-one-half-plus-2}.
By \eqref{mu-order-1-cancelation-a-star-with-mu} and \eqref{minus-2-sigularity-g-h-f}, we can similarly prove \eqref{estimate-of-X-g-h-f-total-one-half-with-mu} and \eqref{estimate-of-Y-g-h-f-total-one-half-with-mu}.
By \eqref{mu-order-1-cancelation-a-star-with-mu}, using the imbedding $H^{2} \hookrightarrow L^{\infty}$ and \eqref{minus-2-sigularity-g-h-f}, we can similarly prove \eqref{estimate-of-A-g-h-varrho-f-with-mu}.
By \eqref{mu-order-1-cancelation-a-star-with-mu}, using the imbedding $H^{2} \hookrightarrow L^{\infty}$ for $\mu^{\frac{1}{256}}\varrho$ and \eqref{minus-2-sigularity-g-h-f}, we can similarly get \eqref{estimate-of-B-g-h-varrho-f-with-mu}. By the same argument, we can also get \eqref{estimate-of-C-g-h-varrho-f-with-mu}. By \eqref{mu-order-1-cancelation-a-star-with-mu}, using the imbedding $H^{2} \hookrightarrow L^{\infty}$ for $\mu^{\frac{1}{256}}g$ and Lemma \ref{cancel-velocity-regularity}, we can get \eqref{estimate-of-D-g-h-f-total-one-half-with-mu}.

We now go to prove \eqref{estimate-of-Z-g2-h2-total-one-half}.
By \eqref{mu-order-2-cancelation-a-star} and \eqref{minus-1-sigularity-g2-h2}, using \eqref{deal-with-product} and \eqref{deal-with-polynomial-weight}, we have
\beno
\int
B  \mathrm{D}^{2}(\mu^{a}_{*}) g_{*}^{2} h^{2} \mathrm{d}V
&\lesssim&  \int (\mathrm{1}_{|v-v_{*}| \geq 1} \langle v \rangle^{\gamma+2s} + \mathrm{1}_{|v-v_{*}|\leq 1} |v-v_{*}|^{-1}  \mu^{\frac{1}{128}}\mu^{\frac{1}{128}}_{*})   g_{*}^{2} h^{2} \mathrm{d}v_{*} \mathrm{d}v
\\&\lesssim& |g|_{L^{2}}^{2}|h|_{L^{2}_{\gamma/2+s}}^{2} + |\mu^{\frac{1}{256}} g|_{H^{s_{1}}}^{2}|\mu^{\frac{1}{256}} h|_{H^{s_{2}}}^{2} \lesssim |g|_{H^{s_{1}}}^{2}|h|_{H^{s_{2}}_{\gamma/2+s}}^{2}.
\eeno
By \eqref{mu-order-2-cancelation-a-star-with-mu} and \eqref{minus-1-sigularity-g2-h2}, we can similarly prove \eqref{estimate-of-Z-g2-h2-total-one-half-with-mu}. Now the proof is complete.
\end{proof}

We set to prove Lemma \ref{optimal-weight-estimate} by adopting some arguments in the proof of Lemma \ref{sigma-integral-mu-difference}.
\begin{proof}[Proof of Lemma \ref{optimal-weight-estimate}] We will give a detailed proof to \eqref{mu-star-g-star-difference-h2} and sketch the proof of the other three results by pointing out the differences.

According to $|v-v_{*}| \leq 1$ and $|v-v_{*}| \geq 1$, we write $\int B \mu^{\f{1}{16}}_{*} \mathrm{D}^{2}(g_{*}) h^{2} \mathrm{d}V =
\mathcal{P} + \mathcal{Q}$ where
\beno
\mathcal{P}\colonequals     \int B \mu^{\f{1}{16}}_{*} \mathrm{1}_{|v-v_{*}| \leq 1}\mathrm{D}^{2}(g_{*}) h^{2}
\mathrm{d}V, \quad
\mathcal{Q}\colonequals     \int B \mu^{\f{1}{16}}_{*} \mathrm{1}_{|v-v_{*}| \geq 1}\mathrm{D}^{2}(g_{*}) h^{2} \mathrm{d}V.
\eeno
By Taylor expansion, we have
\ben \label{Taylor-at-v-star-g}
-\mathrm{D}(g_{*}) = \int_{0}^{1} \nabla g(v_{*}(\iota)) \cdot \left(v^{\prime}_{*}-v_{*}\right) \mathrm{d}\iota.
\een

{\it Estimate of $\mathcal{P}$.}
By \eqref{Taylor-at-v-star-g}, we get
\ben \label{term-1-taylor-g}
\mathcal{P} \lesssim
\int B |v-v_{*}|^{2}\sin^{2}\frac{\theta}{2} \mathrm{1}_{|v-v_{*}| \leq 1} \mu^{\f{1}{16}}_{*} |\nabla g(v_{*}(\iota))|^{2}
 h^{2} \mathrm{d}V \mathrm{d}\iota.
\een
Since $|v-v_{*}| \leq 1$, by Lemma \ref{mu-weight-transfer},
 we have $\mu_{*} \lesssim \mu^{\f{1}{2}}(v_{*}(\iota)), \mu_{*} \lesssim \mu^{\f{1}{2}}$.
For any fixed $\iota$, using the following change of variable
\ben \label{change-v-star-to-v-star-iota}
v_{*} \to v_{*}(\iota), \quad |\frac{\partial v_{*}}{\partial v_{*}(\iota)}| \lesssim 1, \quad \cos \theta(\iota) = \frac{v-v_{*}(\iota)}{|v-v_{*}(\iota)|} \cdot \sigma, \quad \frac{\theta}{2} \leq \theta(\iota) \leq \theta,
\een
the fact $\frac{\sqrt{2}}{2}|v-v_{*}| \leq |v-v_{*}(\iota)| \leq |v-v_{*}|$  and the estimate \eqref{cancel-angular-singularity},
we get
\beno
\mathcal{P} &\lesssim&
\int B |v-v_{*}|^{2}\sin^{2}\frac{\theta}{2} \mathrm{1}_{|v-v_{*}| \leq 1} \mu^{\frac{1}{64}}_{*}  \mu^{\frac{1}{64}} |(\nabla g)_{*}|^{2}
 h^{2} \mathrm{d}V
 \\ &\lesssim&
\int |v-v_{*}|^{-1} \mathrm{1}_{|v-v_{*}| \leq 1} \mu^{\frac{1}{64}}_{*}  \mu^{\frac{1}{64}} |(\nabla g)_{*}|^{2}
 h^{2} \mathrm{d}v \mathrm{d}v_{*}  \lesssim |\mu^{\frac{1}{128}} g|_{H^{s_{1}+1}}^{2} |\mu^{\frac{1}{128}} h|_{H^{s_{2}}}^{2},
\eeno
where we use \eqref{minus-1-sigularity-g2-h2} in the last inequality.

{\it Estimate of $\mathcal{Q}$.}  According to $\sin\frac{\theta}{2} \leq |v-v_{*}|^{-1}$ and $\sin\frac{\theta}{2} \geq |v-v_{*}|^{-1}$, we write $\mathcal{Q} = \mathcal{Q}_{\leq}+\mathcal{Q}_{\geq}$ where
\ben  \label{definition-mathcal-P-leq}
\mathcal{Q}_{\leq}\colonequals     \int B \mu^{\f{1}{16}}_{*} \mathrm{1}_{|v-v_{*}| \geq 1} \mathrm{1}_{\sin\frac{\theta}{2} \leq |v-v_{*}|^{-1}}\mathrm{D}^{2}(g_{*}) h^{2} \mathrm{d}V,
\\ \label{definition-mathcal-P-geq}
 \mathcal{Q}_{\geq}\colonequals     \int B \mu^{\f{1}{16}}_{*} \mathrm{1}_{|v-v_{*}| \geq 1} \mathrm{1}_{\sin\frac{\theta}{2} \geq |v-v_{*}|^{-1}}\mathrm{D}^{2}(g_{*}) h^{2} \mathrm{d}V.
\een
It is obvious $|\mathcal{Q}_{\geq}| \leq 2\mathcal{Q}_{\geq,1} + 2\mathcal{Q}_{\geq,2}$ where
\beno
\mathcal{Q}_{\geq,1}\colonequals    \int B \mu^{\f{1}{16}}_{*} \mathrm{1}_{|v-v_{*}| \geq 1} \mathrm{1}_{\sin\frac{\theta}{2} \geq |v-v_{*}|^{-1}} g_{*}^{2} h^{2} \mathrm{d}V,
\quad
\mathcal{Q}_{\geq,2}\colonequals    \int B \mu^{\f{1}{16}}_{*} \mathrm{1}_{|v-v_{*}| \geq 1} \mathrm{1}_{\sin\frac{\theta}{2} \geq |v-v_{*}|^{-1}} (g^{2})^{\prime}_{*}  h^{2} \mathrm{d}V.
\eeno
Recalling \eqref{theta-large-critical}, we have
\beno
\mathcal{Q}_{\geq,1} \lesssim \int \mu^{\f{1}{16}}_{*} |v-v_{*}|^{\gamma+2s} \mathrm{1}_{|v-v_{*}| \geq 1}  g_{*}^{2} h^{2} \mathrm{d}V \lesssim |\mu^{\f{1}{64}}g|_{L^{2}}^{2} |h|_{L^{2}_{\gamma+2s}}^{2},
\eeno
where we use
\ben \label{need-to-revise-1}
\mu^{\f{1}{16}}_{*} |v-v_{*}|^{\gamma+2s} \mathrm{1}_{|v-v_{*}| \geq 1} \lesssim \mu^{\f{1}{16}}_{*} \langle v-v_{*} \rangle^{\gamma+2s} \lesssim  \mu^{\f{1}{32}}_{*} \langle v \rangle^{\gamma+2s}.
\een
The analysis of $\mathcal{Q}_{\geq,2}$ is a little bit technical. We will make use of the change of variable $v_{*} \rightarrow v^{\prime}_{*}$ and the nice weight $\mu^{\f{1}{16}}_{*}$ on $v_{*}$. Since $|v-v^{\prime}_{*}| \sim |v-v_{*}|$ and $|v-v_{*}| \geq 1$, we have
\ben \label{need-to-revise-2}
|v-v_{*}|^{\gamma} \mu^{\f{1}{16}}_{*} \lesssim |v-v_{*}|^{\gamma+2s} \mu^{\f{1}{16}}_{*} |v-v^{\prime}_{*}|^{-2s} \lesssim \langle v \rangle^{\gamma+2s} |v-v^{\prime}_{*}|^{-2s}.
\een
Plugging which into $\mathcal{Q}_{\geq,2}$, by the change of variable $v_{*} \rightarrow v^{\prime}_{*}$(take $\iota=1$ in \eqref{change-v-star-to-v-star-iota} and let $\theta^{\prime} = \frac{\theta}{2}$) and using \eqref{theta-large-critical}, we get
\beno
\mathcal{Q}_{\geq,2} &\lesssim& \int \langle v \rangle^{\gamma+2s} |v-v^{\prime}_{*}|^{-2s} \sin^{-2-2s}\theta^{\prime} \mathrm{1}_{|v-v^{\prime}_{*}| \geq \sqrt{2}/2} \mathrm{1}_{\sin \theta^{\prime} \geq |v-v^{\prime}_{*}|^{-1}} (g^{2})^{\prime}_{*}  h^{2} \mathrm{d}V
\\ &\lesssim& \int \langle v \rangle^{\gamma+2s} |v-v_{*}|^{-2s} \sin^{-2-2s}\theta \mathrm{1}_{|v-v_{*}| \geq \sqrt{2}/2} \mathrm{1}_{\sin \theta \geq |v-v_{*}|^{-1}} g^{2}_{*}  h^{2} \mathrm{d}V
\\ &\lesssim& \int \langle v \rangle^{\gamma+2s}  g^{2}_{*}  h^{2} \mathrm{d}v \mathrm{d}v_{*} \lesssim |g|_{L^{2}}^{2} |h|_{L^{2}_{\gamma+2s}}^{2}.
\eeno
Patching together the estimates of $\mathcal{Q}_{\geq,1}$ and $\mathcal{Q}_{\geq,2}$, we arrive at
$
\mathcal{Q}_{\geq} \lesssim |g|_{L^{2}}^{2} |h|_{L^{2}_{\gamma+2s}}^{2}.
$
We now go to see $\mathcal{Q}_{\leq}$.
Plugging \eqref{Taylor-at-v-star-g} into $\mathcal{Q}_{\leq}$, we get
\ben \label{P-leq-after-Taylor}
\mathcal{Q}_{\leq} \lesssim
\int B |v-v_{*}|^{2}\sin^{2}\frac{\theta}{2} \mathrm{1}_{|v-v_{*}| \geq 1} \mathrm{1}_{\sin\frac{\theta}{2} \leq |v-v_{*}|^{-1}}  \mu^{\f{1}{16}}_{*} |\nabla g(v_{*}(\iota))|^{2}
 h^{2} \mathrm{d}V \mathrm{d}\iota.
\een
We will use the change of variable $v_{*} \rightarrow v_{*}(\iota)$ according to \eqref{change-v-star-to-v-star-iota}. Recall
$\frac{\sqrt{2}}{2}|v-v_{*}| \leq |v-v_{*}(\iota)| \leq |v-v_{*}|$ and $\frac{\theta}{2} \leq \theta(\iota) \leq \theta$, then
\ben \label{need-to-revise-3}
|v-v_{*}|^{\gamma+2} \mu^{\f{1}{16}}_{*} \lesssim |v-v_{*}|^{\gamma+2s} \mu^{\f{1}{16}}_{*} |v-v_{*}(\iota)|^{2-2s} \lesssim \langle v \rangle^{\gamma+2s} |v-v_{*}(\iota)|^{2-2s},
\\ \nonumber
|v-v_{*}| \geq 1 \Rightarrow |v-v_{*}(\iota)| \geq  \frac{\sqrt{2}}{2}, \quad \sin\frac{\theta}{2} \leq |v-v_{*}|^{-1} \Rightarrow \sin\frac{\theta(\iota)}{2} \leq |v-v_{*}(\iota)|^{-1},
\een
which yield
\beno
\mathcal{Q}_{\leq} &\lesssim&
\int  \langle v \rangle^{\gamma+2s} |v-v_{*}(\iota)|^{2-2s} \sin^{-2s}\frac{\theta(\iota)}{2} \mathrm{1}_{|v-v_{*}(\iota)| \geq \sqrt{2}/2} \mathrm{1}_{\sin\frac{\theta(\iota)}{2} \leq |v-v_{*}(\iota)|^{-1}}   |\nabla g(v_{*}(\iota))|^{2}
 h^{2} \mathrm{d}V \mathrm{d}\iota
 \\ &\lesssim&
\int  \langle v \rangle^{\gamma+2s} |v-v_{*}|^{2-2s} \sin^{-2s}\frac{\theta}{2} \mathrm{1}_{|v-v_{*}| \geq \sqrt{2}/2} \mathrm{1}_{\sin\frac{\theta}{2} \leq |v-v_{*}|^{-1}}   |\nabla g(v_{*})|^{2}
 h^{2} \mathrm{d}V
 \\ &\lesssim&
\int  \langle v \rangle^{\gamma+2s}   |\nabla g(v_{*})|^{2}
 h^{2} \mathrm{d}v \mathrm{d}v_{*} \lesssim  |g|_{H^{1}}^{2} |h|_{L^{2}_{\gamma+2s}}^{2},
\eeno
where we use the change of variable $v_{*} \rightarrow v_{*}(\iota)$ in the second line and
\eqref{theta-small-critical} in the last line.
Patching  together the estimates of $\mathcal{Q}_{\geq}$ and $\mathcal{Q}_{\leq}$, we get
$
|\mathcal{Q}| \lesssim |g|_{H^{1}}^{2} |h|_{L^{2}_{\gamma+2s}}^{2}.
$
Patching  together the estimates of $\mathcal{P}$ and $\mathcal{Q}$, we finish the proof of \eqref{mu-star-g-star-difference-h2}.

Now let us see how to get \eqref{mu-star-g-star-difference-h2-f2-star} by revising the proof of \eqref{mu-star-g-star-difference-h2}. Keep in mind the additional term $f^{2}_{*}$.
For $\mathcal{P}$ where $|v-v_{*}| \leq 1$, from \eqref{term-1-taylor-g} we just take out $|\nabla g|_{L^{\infty}} \lesssim |g|_{H^{3}}$ and thus we do not need the change of variable $v_{*} \to v_{*}^{\prime}(\iota)$. Finally we will
get $|\mathcal{P}| \lesssim |g|_{H^{3}}^{2} |\mu^{\frac{1}{128}} f|_{H^{s_{1}}}^{2} |\mu^{\frac{1}{128}} h|_{H^{s_{2}}}^{2}$.
Let us turn to $\mathcal{Q}$ where $|v-v_{*}| \geq 1$. Let us first see  $\mathcal{Q}_{\geq}$ in
\eqref{definition-mathcal-P-geq}. We simply use $\mathrm{D}^{2}(g_{*}) \lesssim |g|_{L^{\infty}}^{2} \lesssim |g|_{H^{2}}^{2}$ and then just use the same argument as that for $\mathcal{Q}_{\geq,1}$.
Finally we will
get $|\mathcal{Q}_{\geq}| \lesssim |g|_{H^{2}}^{2} |\mu^{\frac{1}{64}} f|_{L^{2}}^{2} |h|_{L^{2}_{\gamma+2s}}^{2}$.
Now it remains to see $\mathcal{Q}_{\leq}$ defined in
\eqref{definition-mathcal-P-leq}. From \eqref{P-leq-after-Taylor} we just take out $|\nabla g|_{L^{\infty}} \lesssim |g|_{H^{3}}$ and thus we do not need the change of variable $v_{*} \to v_{*}(\iota)$. The line \eqref{need-to-revise-3} is revised to
\beno 
|v-v_{*}|^{\gamma+2} \mu^{\f{1}{16}}_{*} \lesssim |v-v_{*}|^{\gamma+2s} \mu^{\f{1}{16}}_{*} |v-v_{*}|^{2-2s} \lesssim \langle v \rangle^{\gamma+2s} |v-v_{*}|^{2-2s} \mu_{*}^{\f{1}{32}},
\eeno
and we will
get $|\mathcal{Q}_{\leq}| \lesssim |g|_{H^{3}}^{2} |\mu^{\frac{1}{64}} f|_{L^{2}}^{2} |h|_{L^{2}_{\gamma+2s}}^{2}$.
Patching together the three parts, we get \eqref{mu-star-g-star-difference-h2-f2-star}.

Comparing to \eqref{mu-star-g-star-difference-h2-f2-star}, there is an additional factor $\mu^{\f{1}{16}}$ in \eqref{mu-star-g-star-difference-h2-f2-star-another} and so we can also get $\mu$ weight for $g$ and $h$ by using \eqref{mu-weight-result}. We omit the details here.

Now let us see how to get \eqref{key-estimate-only-star} by revising the proof of \eqref{mu-star-g-star-difference-h2}. Now we have the factor $\mu^{\f{1}{16}}$ instead of $\mu^{\f{1}{16}}_{*}$. In $\mathcal{P}$ where $|v-v_{*}| \leq 1$, by using Lemma \ref{mu-weight-transfer} to get
$\mu \lesssim \mu^{\f{1}{2}}(v_{*}(\iota))$, we can get exactly the same upper bound. As for $\mathcal{Q}$ where $|v-v_{*}| \geq 1$, we indicate the differences. Firstly, \eqref{need-to-revise-1} is revised to
\ben \label{need-to-revise-1-new}
\mu^{\f{1}{16}} |v-v_{*}|^{\gamma+2s} \mathrm{1}_{|v-v_{*}| \geq 1} \lesssim \mu^{\f{1}{16}} \langle v-v_{*} \rangle^{\gamma+2s} \lesssim  \mu^{\f{1}{32}} \langle v_{*} \rangle^{\gamma+2s}.
\een
Secondly, \eqref{need-to-revise-2} is revised to
\ben \label{need-to-revise-2-new}
|v-v_{*}|^{\gamma} \mu^{\f{1}{16}} \lesssim |v-v^{\prime}_{*}|^{\gamma+2s} \mu^{\f{1}{16}} |v-v^{\prime}_{*}|^{-2s} \lesssim \langle v^{\prime}_{*} \rangle^{\gamma+2s} |v-v^{\prime}_{*}|^{-2s} \mu^{\f{1}{32}}.
\een
Thirdly, the line \eqref{need-to-revise-3} is revised to
\ben \label{need-to-revise-3-new}
|v-v_{*}|^{\gamma+2} \mu^{\f{1}{16}} \lesssim |v-v_{*}(\iota)|^{\gamma+2s} \mu^{\f{1}{16}} |v-v_{*}(\iota)|^{2-2s} \lesssim \langle v_{*}(\iota) \rangle^{\gamma+2s} |v-v_{*}(\iota)|^{2-2s} \mu^{\f{1}{32}}.
\een
Note that in all the three estimates \eqref{need-to-revise-1-new}, \eqref{need-to-revise-2-new} and \eqref{need-to-revise-3-new}, we have polynomial weight $\langle \cdot \rangle^{\gamma+2s}$ for $g$ and negative exponential weight $\mu^{\f{1}{32}}$ for $h$. Therefore, we can get \eqref{key-estimate-only-star}.
\end{proof}

The following lemma gives upper bounds of various integrals over $\mathbb{S}^{2}$ involving the difference
$\mathrm{D}(W_{l})$.
\begin{lem} 
Let $l \geq 0$. The following estimates are valid.
\ben \label{Wl-difference-order-mu-star}
|\int B \mathrm{D}(W_{l}) \mu^{\frac{1}{4}}_{*} \mathrm{d}\sigma|  &\lesssim& (\mathrm{1}_{|v-v_{*}|\geq 1} \langle v \rangle^{l+\gamma} + \mathrm{1}_{|v-v_{*}|\leq 1} |v-v_{*}|^{-2}  \mu^{\frac{1}{32}}) \mu^{\frac{1}{32}}_{*}.
\\ \label{Wl-difference-square-mu-star}
\int B \mathrm{D}^{2}(W_{l}) \mu^{\frac{1}{4}}_{*} \mathrm{d}\sigma  &\lesssim& (\mathrm{1}_{|v-v_{*}|\geq 1} \langle v \rangle^{2l+\gamma} + \mathrm{1}_{|v-v_{*}|\leq 1} |v-v_{*}|^{-1}  \mu^{\frac{1}{32}}) \mu^{\frac{1}{32}}_{*}.
\\ \label{Wl-difference-square-mu-star-prime}
\int B \mathrm{D}^{2}(W_{l}) (\mu^{\frac{1}{4}})^{\prime}_{*} \mathrm{d}\sigma  &\lesssim& \mathrm{1}_{|v-v_{*}|\geq 1} \langle v \rangle^{2l+\gamma} + \mathrm{1}_{|v-v_{*}|\leq 1} |v-v_{*}|^{-1}  \mu^{\frac{1}{32}}\mu^{\frac{1}{32}}_{*}.
\een
\end{lem}
\begin{proof} Applying \eqref{Taylor-at-v} to $W_{l}$, we get
$
 |\int B \mathrm{D}(W_{l}) \mu^{\frac{1}{4}}_{*} \mathrm{d}\sigma| = |\mathcal{P}_{1} + \mathcal{P}_{2}|
$
where
\beno
\mathcal{P}_{1} \colonequals    \int
B \mu^{\frac{1}{4}}_{*} \nabla W_{l}(v) \cdot (v^{\prime}-v)   \mathrm{d}\sigma,
\quad
\mathcal{P}_{2} \colonequals     \int
B \mu^{\frac{1}{4}}_{*} (1-\kappa) \nabla^{2} W_{l} (v(\kappa)):
  (v^{\prime}-v)  \otimes (v^{\prime}-v)  \mathrm{d}\kappa \mathrm{d}\sigma.
\eeno
Using \eqref{theta-squre-out} and  \eqref{cancel-angular-singularity}, since $|\nabla W_{l}| \lesssim l W_{l-1}$
we get
\beno
|\mathcal{P}_{1}| = |\int
B \sin^{2} \frac{\theta}{2} \mu^{\frac{1}{4}}_{*}  \nabla W_{l}(v) \cdot (v-v_{*}) \mathrm{d}\sigma|
\lesssim \mu^{\frac{1}{4}}_{*} \langle v \rangle^{l-1} |v-v_{*}|^{\gamma+1}.
\eeno
If $|v-v_{*}| \geq 1$, then $|v-v_{*}|^{\gamma+1}  \sim \langle v-v_{*}\rangle^{\gamma+1}  \lesssim \langle v\rangle^{\gamma+1}\langle v_{*}\rangle^{|\gamma+1|}$. If $|v-v_{*}|\leq 1$, then $\mu^{\frac{1}{4}}_{*} \lesssim \mu^{\frac{1}{16}}\mu^{\frac{1}{16}}_{*}$ by Lemma \ref{mu-weight-transfer} and $|v-v_{*}|^{\gamma+1} \lesssim |v-v_{*}|^{-2}$ since $\gamma > -3$. Patching  together these two cases, we get
\ben \label{estimate-of-O1}
|\mathcal{P}_{1}| \lesssim (\mathrm{1}_{|v-v_{*}|\geq 1} \langle v \rangle^{l+\gamma} + \mathrm{1}_{|v-v_{*}|\leq 1} |v-v_{*}|^{-2}  \mu^{\frac{1}{32}})\mu^{\frac{1}{32}}_{*}.
\een

Using  $|\nabla^{2} W_{l}| \lesssim (l^{2}+1) W_{l-2}$ to get
\ben \label{O-2-with-v-kappa}
|\mathcal{P}_{2}| \lesssim   \int
\mu^{\frac{1}{4}}_{*} \langle v(\kappa) \rangle^{l-2} |v-v_{*}|^{\gamma+2}\sin^{-2s}\frac{\theta}{2} \mathrm{d}\kappa \mathrm{d}\sigma.
\een
If $|v-v_{*}| \geq 1$, then $|v-v_{*}|^{\gamma+2}  \lesssim \langle v\rangle^{\gamma+2}\langle v_{*}\rangle^{|\gamma+2|}$ and
$\langle v(\kappa) \rangle^{l-2}  \lesssim \langle v(\kappa)-v_{*}\rangle^{l-2} \langle v_{*}\rangle^{|l-2|}
 \lesssim \langle v-v_{*}\rangle^{l-2} \langle v_{*}\rangle^{|l-2|} \lesssim \langle v \rangle^{l-2} \langle v_{*}\rangle^{2|l-2|}$.
As a result, we have
\beno
\mu^{\frac{1}{4}}_{*} \langle v(\kappa) \rangle^{l-2} |v-v_{*}|^{\gamma+2} \lesssim \mu^{\frac{1}{32}}_{*} \langle v \rangle^{l+\gamma}.
\eeno
If $|v-v_{*}|\leq 1$, then $\mu^{\frac{1}{4}}_{*} \lesssim \mu^{\frac{1}{16}}\mu^{\frac{1}{16}}_{*}$ by Lemma \ref{mu-weight-transfer} and $|v-v_{*}|^{\gamma+2} \lesssim |v-v_{*}|^{-1}$ since $\gamma > -3$.  Using $|v(\kappa)| \lesssim |v|+|v_{*}|$,
we get $| (\nabla^{2} W_{l})(v(\kappa)) | \lesssim \langle v \rangle^{|l-2|}\langle v_{*} \rangle^{|l-2|}$. As a result, we have
\beno
\mu^{\frac{1}{4}}_{*} \langle v(\kappa) \rangle^{l-2} |v-v_{*}|^{\gamma+2} \lesssim \mu^{\frac{1}{32}}_{*} \mu^{\frac{1}{32}} |v-v_{*}|^{-1}.
\eeno
Patching these two cases together, using \eqref{cancel-angular-singularity}, we get
\ben \label{estimate-of-O2}
|\mathcal{P}_{2}| \lesssim (\mathrm{1}_{|v-v_{*}|\geq 1} \langle v \rangle^{l+\gamma} + \mathrm{1}_{|v-v_{*}|\leq 1} |v-v_{*}|^{-1}  \mu^{\frac{1}{32}})\mu^{\frac{1}{32}}_{*}.
\een
Patching together \eqref{estimate-of-O1} and \eqref{estimate-of-O2}, we arrive at \eqref{Wl-difference-order-mu-star}.

Now we sketch the proof of \eqref{Wl-difference-square-mu-star}. By Taylor expansion up to order 1, we have
\ben \label{Taylor-at-order-1}
-\mathrm{D}(W_{l})= \int_{0}^{1} \nabla W_{l}(v(\kappa)) \cdot  (v^{\prime}-v)  \mathrm{d}\kappa.
\een
Since $|\nabla W_{l}| \lesssim l W_{l-1}$, we have
\ben \label{Wl-square-with-v-kappa}
|\int B \mathrm{D}^{2}(W_{l}) \mu^{\frac{1}{4}}_{*} \mathrm{d}\sigma| \lesssim   \int
\mu^{\frac{1}{4}}_{*} \langle v(\kappa) \rangle^{2l-2} |v-v_{*}|^{\gamma+2}\sin^{-2s}\frac{\theta}{2} \mathrm{d}\kappa \mathrm{d}\sigma.
\een
Comparing \eqref{Wl-square-with-v-kappa} with
\eqref{O-2-with-v-kappa}, by the same arguments as that for $\mathcal{P}_{2}$, we can get \eqref{Wl-difference-square-mu-star}.

We set to prove \eqref{Wl-difference-square-mu-star-prime}.
Recalling \eqref{Taylor-at-order-1}.
Since  $|\nabla W_{l}| \lesssim l W_{l-1}$ and $|v(\kappa)| \lesssim |v|+|v^{\prime}_{*}|$,
we get $| \nabla W_{l}(v(\kappa)) | \lesssim \langle v \rangle^{l-1}\langle v^{\prime}_{*} \rangle^{l-1}$. Note that $|v^{\prime}-v| = |v-v_{*}|\sin\frac{\theta}{2} \sim |v-v^{\prime}_{*}|\sin\frac{\theta}{2}$ and thus
\beno
|\mathrm{D}(W_{l})| \lesssim \langle v \rangle^{l-1}\langle v^{\prime}_{*} \rangle^{l-1} |v-v^{\prime}_{*}|\sin\frac{\theta}{2},
\eeno
which gives
\beno
|v-v_{*}|^{\gamma} \mathrm{D}^{2}(W_{l}) (\mu^{\frac{1}{4}})^{\prime}_{*} \lesssim |v-v^{\prime}_{*}|^{\gamma+2} \langle v \rangle^{2l-2}\langle v^{\prime}_{*} \rangle^{2l-2} (\mu^{\frac{1}{4}})^{\prime}_{*} \sin^{2}\frac{\theta}{2} .
\eeno
If $|v-v_{*}| \geq 1$, then $|v-v^{\prime}_{*}|^{\gamma+2}  \sim \langle v-v^{\prime}_{*}\rangle^{\gamma+2}  \lesssim \langle v\rangle^{\gamma+2}\langle v^{\prime}_{*}\rangle^{|\gamma+2|}$. If $|v-v_{*}|\leq 1$, then $(\mu^{\frac{1}{4}})^{\prime}_{*} \lesssim \mu^{\frac{1}{16}}\mu^{\frac{1}{16}}_{*}$ and $|v-v^{\prime}_{*}|^{\gamma+2} \lesssim |v-v_{*}|^{-1}$. Patching these two cases together, we get
\beno
|v-v_{*}|^{\gamma} \mathrm{D}^{2}(W_{l}) (\mu^{\frac{1}{4}})^{\prime}_{*} \lesssim (\mathrm{1}_{|v-v_{*}|\geq 1} \langle v \rangle^{2l+\gamma} + \mathrm{1}_{|v-v_{*}|\leq 1} |v-v_{*}|^{-1}  \mu^{\frac{1}{32}}\mu^{\frac{1}{32}}_{*})\sin^{2}\frac{\theta}{2}.
\eeno
From which together with \eqref{cancel-angular-singularity},
we get \eqref{Wl-difference-square-mu-star-prime} and finish the proof.
\end{proof}

\begin{proof}[Proof of Lemma \ref{full-integral-Wl-difference-g-h-f}.] By \eqref{Wl-difference-square-mu-star} and \eqref{minus-1-sigularity-g2-h2}, we directly obtain \eqref{Wl-difference-square-mu-star-g2-h2}.
By \eqref{Wl-difference-square-mu-star-prime} and \eqref{minus-1-sigularity-g2-h2}, we directly obtain \eqref{Wl-difference-square-mu-star-prime-g2-h2}.
By \eqref{Wl-difference-order-mu-star} and \eqref{minus-2-sigularity-g-h-f}, we directly obtain \eqref{Wl-difference-order-1-mu-star-ghf}.

By \eqref{Wl-difference-order-mu-star}, the imbedding $H^{2} \hookrightarrow L^{\infty}$
and \eqref{minus-2-sigularity-g-h-f}, we get \eqref{Wl-difference-order-mu-star-g-hfvarrho-prime} and \eqref{Wl-difference-order-mu-star-gvarrho-hf-prime}.

It remains to prove \eqref{Wl-difference-order-mu-star-g-hf-prime}.
By \eqref{type-2-cancel} and recalling $|\nabla^{2} W_{l}| \lesssim (l^{2}+1) W_{l-2}$, using $|v^{\prime}(\kappa) -v_{*}| \sim |v -v_{*}|,  \langle v(\kappa) \rangle^{l-2} \lesssim \langle v^{\prime} \rangle^{l-2} \langle v_{*} \rangle^{2|l-2|}$,
recalling \eqref{change-v-to-v-prime} for the change of variable $v \rightarrow v^{\prime}$, the fact $\int  \sin^{-2s}\theta
\mathrm{d}\sigma \lesssim \frac{1}{1-s}$,
 we have
\beno
&&|\int B \mathrm{D}(W_{l}) \mu^{\frac{1}{4}}_{*} g_{*} (h f)^{\prime} \mathrm{d}V|
\\&\lesssim& \int  \langle v(\kappa) \rangle^{l-2} |v^{\prime} -v_{*}|^{\gamma+2} \mathrm{1}_{0 \leq \theta \leq \pi/2} \sin^{-2s}\frac{\theta}{2}
  |\mu^{\frac{1}{4}}_{*} g_{*} (h f)^{\prime}| \mathrm{d}V
\\&\lesssim& \int \langle v \rangle^{l-2} \langle v_{*} \rangle^{2|l-2|} |v -v_{*}|^{\gamma+2} \mathrm{1}_{0 \leq \theta \leq \pi/4} \sin^{-2s}\theta
 | \mu^{\frac{1}{4}}_{*}g_{*} h f| \mathrm{d}V
\\&\lesssim& \int (\mathrm{1}_{|v-v_{*}|\geq 1} \langle v \rangle^{l+\gamma} + \mathrm{1}_{|v-v_{*}|\leq 1} |v-v_{*}|^{-1}  \mu^{\frac{1}{32}})\mu^{\frac{1}{32}}_{*} |g_{*} h f| \mathrm{d}v_{*} \mathrm{d}v
\\&\lesssim& |\mu^{\frac{1}{64}}g|_{L^{2}} |h|_{L^{2}_{l+\gamma/2}} |f|_{L^{2}_{\gamma/2}} + |\mu^{\frac{1}{64}}g|_{L^{2}} |\mu^{\frac{1}{64}}h|_{L^{2}} |\mu^{\frac{1}{64}}f|_{L^{2}} \lesssim |\mu^{\frac{1}{64}}g|_{L^{2}} |h|_{L^{2}_{l+\gamma/2}} |f|_{L^{2}_{\gamma/2}},
\eeno
where in the last line we use
\beno
\int \mathrm{1}_{|v-v_{*}|\leq 1} |v-v_{*}|^{-1} \mu^{\frac{1}{32}}_{*} |g_{*}| \mathrm{d}v_{*} \lesssim |\mu^{\frac{1}{64}}g|_{L^{2}}.
\eeno
The proof of the lemma is complete now.
\end{proof}

 {\bf Acknowledgments.} Yu-Long Zhou was supported by National Key R\&D Program of China under the grant 2021YFA1002100 and NSF of China under the grant 12001552.
The author is indebted to Prof. Ling-Bing He for his continuous encouragement and supervision.
Great gratitude goes to Prof. Xuguang Lu for his insightful comments, especially in the scope and focus of the article.

\bibliographystyle{siam}
\bibliography{quantum-B-IPL}

\begin{thebibliography}{10}

\bibitem{alexandre2000some}
{\sc R.~Alexandre}, {\em {On some related non homogeneous 3D Boltzmann models
  in the non cutoff case}}, Journal of Mathematics of Kyoto University, 40
  (2000), pp.~493--524.

\bibitem{alexandre2000entropy}
{\sc R.~Alexandre, L.~Desvillettes, C.~Villani, and B.~Wennberg}, {\em Entropy
  dissipation and long-range interactions}, Archive for Rational Mechanics and
  Analysis, 152 (2000), pp.~327--355.

\bibitem{alexandre2011global}
{\sc R.~Alexandre, Y.~Morimoto, S.~Ukai, C.-J. Xu, and T.~Yang}, {\em {Global
  existence and full regularity of the Boltzmann equation without angular
  cutoff}}, Communications in Mathematical Physics, 304 (2011), pp.~513--581.

\bibitem{alexandre2012boltzmann}
\leavevmode\vrule height 2pt depth -1.6pt width 23pt, {\em {The Boltzmann
  equation without angular cutoff in the whole space: I, Global existence for
  soft potential}}, Journal of Functional Analysis, 262 (2012), pp.~915--1010.

\bibitem{bae2021relativistic}
{\sc G.-C. Bae, J.~W. Jang, and S.-B. Yun}, {\em {The relativistic quantum
  Boltzmann equation near equilibrium}}, Archive for Rational Mechanics and
  Analysis, 240 (2021), pp.~1593--1644.

\bibitem{benedetto2004some}
{\sc D.~Benedetto, F.~Castella, R.~Esposito, and M.~Pulvirenti}, {\em {Some
  considerations on the derivation of the nonlinear quantum Boltzmann
  equation}}, Journal of Statistical Physics, 116 (2004), pp.~381--410.

\bibitem{benedetto2006some}
\leavevmode\vrule height 2pt depth -1.6pt width 23pt, {\em {Some considerations
  on the derivation of the nonlinear quantum Boltzmann equation II: the low
  density regime}}, Journal of Statistical Physics, 124 (2006), pp.~951--996.

\bibitem{benedetto2007short}
\leavevmode\vrule height 2pt depth -1.6pt width 23pt, {\em {A short review on
  the derivation of the nonlinear quantum Boltzmann equations}}, Communications
  in Mathematical Sciences, 5 (2007), pp.~55--71.

\bibitem{benedetto2008n}
\leavevmode\vrule height 2pt depth -1.6pt width 23pt, {\em {From the N-body
  Schr{\"o}dinger equation to the quantum Boltzmann equation: a term-by-term
  convergence result in the weak coupling regime}}, Communications in
  Mathematical Physics, 277 (2008), pp.~1--44.

\bibitem{benedetto2005weak}
{\sc D.~Benedetto, M.~Pulvirenti, F.~Castella, and R.~Esposito}, {\em {On the
  weak-coupling limit for bosons and fermions}}, Mathematical Models and
  Methods in Applied Sciences, 15 (2005), pp.~1811--1843.

\bibitem{briant2016cauchy}
{\sc M.~Briant and A.~Einav}, {\em {On the Cauchy problem for the homogeneous
  Boltzmann--Nordheim equation for bosons: local existence, uniqueness and
  creation of moments}}, Journal of Statistical Physics, 163 (2016),
  pp.~1108--1156.

\bibitem{cai2019spatially}
{\sc S.~Cai and X.~Lu}, {\em {The spatially homogeneous Boltzmann equation for
  Bose--Einstein particles: Rate of strong convergence to equilibrium}},
  Journal of Statistical Physics, 175 (2019), pp.~289--350.

\bibitem{chapman1990mathematical}
{\sc S.~Chapman and T.~G. Cowling}, {\em {The mathematical theory of
  non-uniform gases: an account of the kinetic theory of viscosity, thermal
  conduction and diffusion in gases}}, Cambridge university press, 1990.

\bibitem{dolbeault1994kinetic}
{\sc J.~Dolbeault}, {\em {Kinetic models and quantum effects: A modified
  Boltzmann equation for Fermi-Dirac Particles}}, Archive for Rational
  Mechanics and Analysis, 127 (1994), pp.~101--131.

\bibitem{duan2008cauchy}
{\sc R.~Duan}, {\em {On the Cauchy problem for the Boltzmann equation in the
  whole space: Global existence and uniform stability in
  $L^{2}_{\xi}(H^{N}_{x})$}}, Journal of Differential Equations, 244 (2008),
  pp.~3204--3234.

\bibitem{duan2021global}
{\sc R.~Duan, S.~Liu, S.~Sakamoto, and R.~M. Strain}, {\em {Global mild
  solutions of the Landau and non-cutoff Boltzmann equations}}, Communications
  on Pure and Applied Mathematics, 74 (2021), pp.~932--1020.

\bibitem{erdHos2004quantum}
{\sc L.~Erd{\H{o}}s, M.~Salmhofer, and H.-T. Yau}, {\em {On the quantum
  Boltzmann equation}}, Journal of Statistical Physics, 116 (2004),
  pp.~367--380.

\bibitem{escobedo2003homogeneous}
{\sc M.~Escobedo, S.~Mischler, and M.~A. Valle}, {\em {Homogeneous Boltzmann
  equation in quantum relativistic kinetic theory}}, Electronic Journal of
  Differential Equations, 4 (2003), pp.~1--85.

\bibitem{escobedo2007fundamental}
{\sc M.~Escobedo, S.~Mischler, and J.~J.~L. Vel{\'a}zquez}, {\em {On the
  fundamental solution of a linearized Uehling--Uhlenbeck equation}}, Archive
  for Rational Mechanics and Analysis, 186 (2007), pp.~309--349.

\bibitem{escobedo2008singular}
\leavevmode\vrule height 2pt depth -1.6pt width 23pt, {\em {Singular solutions
  for the Uehling--Uhlenbeck equation}}, Proceedings of the Royal Society of
  Edinburgh Section A: Mathematics, 138 (2008), pp.~67--107.

\bibitem{escobedo2014blow}
{\sc M.~Escobedo and J.~Vel{\'a}zquez}, {\em {On the blow up and condensation
  of supercritical solutions of the Nordheim equation for bosons}},
  Communications in Mathematical Physics, 330 (2014), pp.~331--365.

\bibitem{escobedo2015finite}
\leavevmode\vrule height 2pt depth -1.6pt width 23pt, {\em {Finite time blow-up
  and condensation for the bosonic Nordheim equation}}, Inventiones
  Mathematicae, 200 (2015), pp.~761--847.

\bibitem{gressman2011global}
{\sc P.~Gressman and R.~M. Strain}, {\em {Global classical solutions of the
  Boltzmann equation without angular cut-off}}, Journal of the American
  Mathematical Society, 24 (2011), pp.~771--847.

\bibitem{guo2003classical}
{\sc Y.~Guo}, {\em {Classical solutions to the Boltzmann equation for molecules
  with an angular cutoff}}, Archive for Rational Mechanics and Analysis, 169
  (2003), pp.~305--353.

\bibitem{guo2012vlasov}
\leavevmode\vrule height 2pt depth -1.6pt width 23pt, {\em {The
  Vlasov-Poisson-Landau system in a periodic box}}, Journal of the American
  Mathematical Society, 25 (2012), pp.~759--812.

\bibitem{he2022asymptotic}
{\sc L.-B. He and Y.-L. Zhou}, {\em Asymptotic analysis of the linearized
  boltzmann collision operator from angular cutoff to non-cutoff}, Annales de
  l'Institut Henri Poincar{\'e} C,  (2022).

\bibitem{jiang2021incompressible}
{\sc N.~Jiang, L.~Xiong, and K.~Zhou}, {\em {The incompressible
  Navier-Stokes-Fourier limit from Boltzmann-Fermi-Dirac equation}}, arXiv
  preprint arXiv:2102.02656,  (2021).

\bibitem{li2019global}
{\sc W.~Li and X.~Lu}, {\em {Global existence of solutions of the Boltzmann
  equation for Bose--Einstein particles with anisotropic initial data}},
  Journal of Functional Analysis, 276 (2019), pp.~231--283.

\bibitem{lions1994compactness}
{\sc P.~L. Lions}, {\em {Compactness in Boltzmann's equation via Fourier
  integral operators and applications. III}}, Journal of Mathematics of Kyoto
  University, 34 (1994), pp.~539--584.

\bibitem{lu2000modified}
{\sc X.~Lu}, {\em {A modified Boltzmann equation for Bose--Einstein particles:
  isotropic solutions and long-time behavior}}, Journal of Statistical Physics,
  98 (2000), pp.~1335--1394.

\bibitem{lu2001spatially}
\leavevmode\vrule height 2pt depth -1.6pt width 23pt, {\em {On spatially
  homogeneous solutions of a modified Boltzmann equation for Fermi--Dirac
  particles}}, Journal of Statistical Physics, 105 (2001), pp.~353--388.

\bibitem{lu2004isotropic}
\leavevmode\vrule height 2pt depth -1.6pt width 23pt, {\em {On isotropic
  distributional solutions to the Boltzmann equation for Bose-Einstein
  particles}}, Journal of Statistical Physics, 116 (2004), pp.~1597--1649.

\bibitem{lu2005boltzmann}
\leavevmode\vrule height 2pt depth -1.6pt width 23pt, {\em {The Boltzmann
  equation for Bose--Einstein particles: velocity concentration and convergence
  to equilibrium}}, Journal of Statistical Physics, 119 (2005), pp.~1027--1067.

\bibitem{lu2006boltzmann}
\leavevmode\vrule height 2pt depth -1.6pt width 23pt, {\em {On the Boltzmann
  equation for Fermi--Dirac particles with very soft potentials: Averaging
  compactness of weak solutions}}, Journal of Statistical Physics, 124 (2006),
  pp.~517--547.

\bibitem{lu2008boltzmann}
\leavevmode\vrule height 2pt depth -1.6pt width 23pt, {\em {On the Boltzmann
  equation for Fermi--Dirac particles with very soft potentials: Global
  existence of weak solutions}}, Journal of Differential equations, 245 (2008),
  pp.~1705--1761.

\bibitem{lu2013boltzmann}
\leavevmode\vrule height 2pt depth -1.6pt width 23pt, {\em {The Boltzmann
  equation for Bose-Einstein particles: condensation in finite time}}, Journal
  of Statistical Physics, 150 (2013), pp.~1138--1176.

\bibitem{lu2014boltzmann}
\leavevmode\vrule height 2pt depth -1.6pt width 23pt, {\em {The Boltzmann
  equation for Bose--Einstein particles: regularity and condensation}}, Journal
  of Statistical Physics, 156 (2014), pp.~493--545.

\bibitem{lu2016long}
\leavevmode\vrule height 2pt depth -1.6pt width 23pt, {\em {Long time
  convergence of the Bose--Einstein condensation}}, Journal of Statistical
  Physics, 162 (2016), pp.~652--670.

\bibitem{lu2018long}
\leavevmode\vrule height 2pt depth -1.6pt width 23pt, {\em {Long time strong
  convergence to Bose-Einstein distribution for low temperature}}, Kinetic and
  Related Models, 11 (2018), pp.~715--734.

\bibitem{lu2003stability}
{\sc X.~Lu and B.~Wennberg}, {\em {On stability and strong convergence for the
  spatially homogeneous Boltzmann equation for Fermi-Dirac particles}}, Archive
  for Rational Mechanics and Analysis, 168 (2003), pp.~1--34.

\bibitem{lukkarinen2009not}
{\sc J.~Lukkarinen and H.~Spohn}, {\em {Not to normal order--notes on the
  kinetic limit for weakly interacting quantum fluids}}, Journal of Statistical
  Physics, 134 (2009), pp.~1133--1172.

\bibitem{morimoto2016global}
{\sc Y.~Morimoto and S.~Sakamoto}, {\em {Global solutions in the critical Besov
  space for the non-cutoff Boltzmann equation}}, Journal of Differential
  Equations, 261 (2016), pp.~4073--4134.

\bibitem{nordhiem1928kinetic}
{\sc L.~Nordhiem}, {\em {On the kinetic method in the new statistics and
  application in the electron theory of conductivity}}, Proceedings of the
  Royal Society of London. Series A, Containing Papers of a Mathematical and
  Physical Character, 119 (1928), pp.~689--698.

\bibitem{ouyang2021quantum}
{\sc Z.~Ouyang and L.~Wu}, {\em {On the quantum Boltzmann equation near
  Maxwellian and vacuum}}, arXiv preprint arXiv:2102.00657,  (2021).

\bibitem{spohn2010kinetics}
{\sc H.~Spohn}, {\em {Kinetics of the Bose--Einstein condensation}}, Physica D:
  Nonlinear Phenomena, 239 (2010), pp.~627--634.

\bibitem{strain2006almost}
{\sc R.~M. Strain and Y.~Guo}, {\em {Almost exponential decay near
  Maxwellian}}, Communications in Partial Differential Equations, 31 (2006),
  pp.~417--429.

\bibitem{uehling1933transport}
{\sc E.~A. Uehling and G.~Uhlenbeck}, {\em {Transport phenomena in
  Einstein-Bose and Fermi-Dirac gases. i}}, Physical Review, 43 (1933), p.~552.

\bibitem{ukai1974existence}
{\sc S.~Ukai}, {\em {On the existence of global solutions of mixed problem for
  non-linear Boltzmann equation}}, Proceedings of the Japan Academy, 50 (1974),
  pp.~179--184.

\end{thebibliography}

\end{document}